\documentclass[12pt,reqno]{amsart}
\usepackage{etex,mathtools}
\usepackage{mathdots,longtable}
\usepackage{amscd,amssymb,amsmath,multicol,tikz,float,mathrsfs}
\usepackage[left=1in,right=1in]{geometry}
\begin{document}

\restylefloat{table}
\newtheorem{thm}[equation]{Theorem}
\numberwithin{equation}{section}
\newtheorem{cor}[equation]{Corollary}
\newtheorem{expl}[equation]{Example}
\newtheorem{rmk}[equation]{Remark}
\newtheorem{conv}[equation]{Convention}
\newtheorem{claim}[equation]{Claim}
\newtheorem{lem}[equation]{Lemma}
\newtheorem{sublem}[equation]{Sublemma}
\newtheorem{conj}[equation]{Conjecture}
\newtheorem{defin}[equation]{Definition}
\newtheorem{diag}[equation]{Diagram}
\newtheorem{prop}[equation]{Proposition}
\newtheorem{notation}[equation]{Notation}
\newtheorem{tab}[equation]{Table}
\newtheorem{fig}[equation]{Figure}
\newcounter{bean}
\renewcommand{\theequation}{\thesection.\arabic{equation}}

\raggedbottom \voffset=-.7truein \hoffset=0truein \vsize=8truein
\hsize=6truein \textheight=8truein \textwidth=6truein
\baselineskip=18truept
\def\mapleft#1{\smash{\mathop{\longleftarrow}\limits^{#1}}}
\def\mapright#1{\ \smash{\mathop{\longrightarrow}\limits^{#1}}\ }
\def\ml#1{\,\smash{\mathop{\leftarrow}\limits^{#1}}\,}
\def\mapup#1{\Big\uparrow\rlap{$\vcenter {\hbox {$#1$}}$}}
\def\mapdown#1{\Big\downarrow\rlap{$\vcenter {\hbox {$\ssize{#1}$}}$}}
\def\mapne#1{\nearrow\rlap{$\vcenter {\hbox {$#1$}}$}}
\def\mapse#1{\searrow\rlap{$\vcenter {\hbox {$\ssize{#1}$}}$}}
\def\mapr#1{\smash{\mathop{\rightarrow}\limits^{#1}}}
\def\Mt{\widetilde{\M}}
\def\ss{\smallskip}
\def\s{\sigma}
\def\Bt{\widetilde{\mathcal{B}}}
\def\l{\lambda}
\def\Ah{\widehat{A}}
\def\Bh{\widehat{B}}
\def\vp{v_1^{-1}\pi}
\def\at{{\widetilde\alpha}}
\def\At{\widetilde{\mathcal{A}}}
\def\as{\mathscr{A}}
\def\Ast{\widetilde{\as}}
\def\Mct{\widetilde{\mathcal{M}}}
\def\sm{\wedge}
\def\la{\langle}
\def\ra{\rangle}
\def\lar{\leftarrow}
\def\ev{\text{ev}}
\def\od{\text{od}}
\def\on{\operatorname}
\def\ol#1{\overline{#1}{}}
\def\spin{\on{Spin}}
\def\cat{\on{cat}}
\def\Lbar{\overline{\Lambda}}
\def\qed{\quad\rule{8pt}{8pt}\bigskip}
\def\ssize{\scriptstyle}
\def\a{\alpha}
\def\bz{{\Bbb Z}}
\def\Rhat{\hat{R}}
\def\im{\on{im}}
\def\ct{\widetilde{C}}
\def\ext{\on{Ext}}
\def\sq{\on{Sq}}
\def\eps{\epsilon}
\def\ar#1{\stackrel {#1}{\rightarrow}}
\def\br{{\bold R}}
\def\bC{{\bold C}}
\def\bA{{\bold A}}
\def\bB{{\bold B}}
\def\bD{{\bold D}}
\def\bC{{\bold C}}
\def\bh{{\bold H}}
\def\bQ{{\bold Q}}
\def\bP{{\bold P}}
\def\bx{{\bold x}}
\def\bo{{\bold{bo}}}
\def\dh{\widehat{d}}
\def\A{\mathcal{A}}
\def\B{\mathcal{B}}
\def\si{\sigma}
\def\Vbar{{\overline V}}
\def\dbar{{\overline d}}
\def\wbar{{\overline w}}
\def\Sum{\sum}
\def\tfrac{\textstyle\frac}

\def\tb{\textstyle\binom}
\def\Si{\Sigma}
\def\w{\wedge}
\def\equ{\begin{equation}}
\def\b{\beta}
\def\G{\Gamma}
\def\L{\Lambda}
\def\g{\gamma}
\def\d{\delta}
\def\k{\kappa}
\def\psit{\widetilde{\Psi}}
\def\tht{\widetilde{\Theta}}
\def\psiu{{\underline{\Psi}}}
\def\thu{{\underline{\Theta}}}
\def\aee{A_{\text{ee}}}
\def\aeo{A_{\text{eo}}}
\def\aoo{A_{\text{oo}}}
\def\aoe{A_{\text{oe}}}
\def\vbar{{\overline v}}
\def\endeq{\end{equation}}
\def\sn{S^{2n+1}}
\def\zp{\bold Z_p}
\def\cR{{\mathcal R}}
\def\P{{\mathcal P}}
\def\cQ{{\mathcal Q}}
\def\cj{{\cal J}}
\def\zt{{\bold Z}_2}
\def\bs{{\bold s}}
\def\bof{{\bold f}}
\def\bq{{\bold Q}}
\def\be{{\bold e}}
\def\Hom{\on{Hom}}
\def\ker{\on{ker}}
\def\kot{\widetilde{KO}}
\def\coker{\on{coker}}
\def\da{\downarrow}
\def\colim{\operatornamewithlimits{colim}}
\def\zphat{\bz_2^\wedge}
\def\io{\iota}
\def\om{\omega}
\def\Prod{\prod}
\def\e{{\cal E}}
\def\zlt{\Z_{(2)}}
\def\exp{\on{exp}}
\def\abar{{\overline a}}
\def\xbar{{\overline x}}
\def\ybar{{\overline y}}
\def\zbar{{\overline z}}
\def\mbar{{\overline m}}
\def\nbar{{\overline n}}
\def\sbar{{\overline s}}
\def\kbar{{\overline k}}
\def\bbar{{\overline b}}
\def\et{{\widetilde E}}
\def\ni{\noindent}
\def\tsum{\textstyle \sum}
\def\coef{\on{coef}}
\def\den{\on{den}}
\def\lcm{\on{l.c.m.}}
\def\Ext{\operatorname{Ext}}
\def\iso{\approx}
\def\lra{\longrightarrow}
\def\vi{v_1^{-1}}
\def\ot{\otimes}
\def\psibar{{\overline\psi}}
\def\thbar{{\overline\theta}}
\def\Mh{{\widehat M}}
\def\exc{\on{exc}}
\def\ms{\medskip}
\def\ehat{{\hat e}}
\def\etao{{\eta_{\text{od}}}}
\def\etae{{\eta_{\text{ev}}}}
\def\dirlim{\operatornamewithlimits{dirlim}}
\def\gt{\widetilde{L}}
\def\lt{\widetilde{\lambda}}
\def\st{\widetilde{s}}
\def\ft{\widetilde{f}}
\def\sgd{\on{sgd}}
\def\lfl{\lfloor}
\def\rfl{\rfloor}
\def\ord{\on{ord}}
\def\gd{{\on{gd}}}
\def\rk{{{\on{rk}}_2}}
\def\nbar{{\overline{n}}}
\def\MC{\on{MC}}
\def\lg{{\on{lg}}}
\def\cH{\mathcal{H}}
\def\cS{\mathcal{S}}
\def\cP{\mathcal{P}}
\def\N{{\Bbb N}}
\def\Z{{\Bbb Z}}
\def\Q{{\Bbb Q}}
\def\R{{\Bbb R}}
\def\C{{\Bbb C}}
\def\Lb{\overline\Lambda}
\def\mo{\on{mod}}
\def\xt{\times}
\def\notimm{\not\subseteq}
\def\Remark{\noindent{\it  Remark}}
\def\kut{\widetilde{KU}}
\def\Eb{\overline E}
\def\*#1{\mathbf{#1}}
\def\0{$\*0$}
\def\1{$\*1$}
\def\22{$(\*2,\*2)$}
\def\33{$(\*3,\*3)$}
\def\ss{\smallskip}
\def\ssum{\sum\limits}
\def\dsum{\displaystyle\sum}
\def\la{\langle}
\def\ra{\rangle}
\def\on{\operatorname}
\def\proj{\on{proj}}
\def\od{\text{od}}
\def\ev{\text{ev}}
\def\o{\on{o}}
\def\U{\on{U}}
\def\lg{\on{lg}}
\def\a{\alpha}
\def\bz{{\Bbb Z}}
\def\ccM{{\Bbb M}}
\def\E{\mathcal{E}}
\def\eps{\varepsilon}
\def\bc{{\bold C}}
\def\bN{{\bold N}}
\def\bB{{\bold B}}
\def\bW{{\bold W}}
\def\nut{\widetilde{\nu}}
\def\tfrac{\textstyle\frac}
\def\b{\beta}
\def\G{\Gamma}
\def\g{\gamma}
\def\zt{{\Bbb Z}_2}
\def\zth{{\bold Z}_2^\wedge}
\def\bs{{\bold s}}
\def\bx{{\bold x}}
\def\bof{{\bold f}}
\def\bq{{\bold Q}}
\def\be{{\bold e}}
\def\lline{\rule{.6in}{.6pt}}
\def\xb{{\overline x}}
\def\xbar{{\overline x}}
\def\ybar{{\overline y}}
\def\zbar{{\overline z}}
\def\ebar{{\overline e}}
\def\nbar{{\overline n}}
\def\ubar{{\overline u}}
\def\bbar{{\overline b}}
\def\et{{\widetilde e}}
\def\M{\mathcal{M}}
\def\lf{\lfloor}
\def\rf{\rfloor}
\def\ni{\noindent}
\def\ms{\medskip}
\def\Dhat{{\widehat D}}
\def\what{{\widehat w}}
\def\Yhat{{\widehat Y}}
\def\abar{{\overline{a}}}
\def\minp{\min\nolimits'}
\def\sb{{$\ssize\bullet$}}
\def\mul{\on{mul}}
\def\N{{\Bbb N}}
\def\Z{{\Bbb Z}}
\def\Q{{\Bbb Q}}
\def\R{{\Bbb R}}
\def\C{{\Bbb C}}
\def\Xb{\overline{X}}
\def\eb{\overline{e}}
\def\notint{\cancel\cap}
\def\cS{\mathcal S}
\def\cR{\mathcal R}
\def\el{\ell}
\def\TC{\on{TC}}
\def\GC{\on{GC}}
\def\wgt{\on{wgt}}
\def\Ht{\widetilde{H}}
\def\wbar{\overline w}
\def\dstyle{\displaystyle}
\def\Sq{\on{sq}}
\def\Om{\Omega}
\def\ds{\dstyle}
\def\tz{tikzpicture}
\def\zcl{\on{zcl}}
\def\bd{\bold{d}}
\def\cM{\mathcal{M}}
\def\io{\iota}
\def\Vb#1{{\overline{V_{#1}}}}
\def\Ebar{\overline{E}}
\def\lb{\,\begin{picture}(-1,1)(1,-1)\circle*{3.5}\end{picture}\ }
\def\rlb{\,\begin{picture}(-1,1)(1,-1) \circle*{4.5}\end{picture}\ }
\def\lbb{\,\begin{picture}(-1,1)(1,-1)\circle*{8}\end{picture}\ }
\def\zp{\Z_p}
\def\lbr{\,\begin{picture}(-1,1)(1,-1)[dashed]\circle*{3.5}\end{picture}\ }
\def\llb{\,\begin{picture}(-1,1)(1,-1)\circle*{2.6}\end{picture}\ }
\def\blb{\,\begin{picture}(-1,1)(1,-1) \circle*{5.8}\end{picture}\ }

\title
{The connective $KO$ theory of the Eilenberg-MacLane space $K(\Z/2,2)$}
\author{Donald M. Davis}
\address{Department of Mathematics, Lehigh University\\Bethlehem, PA 18015, USA}
\email{dmd1@lehigh.edu}
\date{February 19, 2025}

\begin{abstract} We compute $ko_*(K(\Z/2,2))$ and $ko^*(K(\Z/2,2))$, the connective $KO$-homology and -cohomology groups of the mod 2 Eilenberg MacLane space $K(\Z/2,2)$, using the Adams spectral sequence. The work relies heavily on work done several years earlier for the (complex) $ku$ groups by the author and W.S.Wilson. We illustrate an interesting duality relation between the $ko$-homology and -cohomology groups. We deduce a new result about Stiefel-Whitney classes in Spin manifolds.\end{abstract}
\maketitle 
\tableofcontents
\section{Introduction and overview}\label{intro}
Let $\zt=\Z/2$ and let $K_2$ denote the Eilenberg-MacLane space $K(\zt,2)$. In this paper, we use the Adams spectral sequence (ASS) to compute the connective $KO$-homology 
and -cohomology groups $ko_*(K_2)$ and $ko^*(K_2)$. The groups $ko_*(K_2)$ were initially studied long ago in \cite{W} and more recently in \cite{DWSW} because of their close relationship with Stiefel-Whitney classes of Spin-manifolds, but only fragmentary results were obtained. A consequence of our work here is the following new result, which is discussed and proved in Section \ref{SWsec}.

\begin{thm}\label{SWthm} There exists an $n$-dimensional Spin manifold with the dual Stiefel-Whitney class $\overline{w}_{n-2}\ne0$ if and only if $n$ is a $2$-power $\ge8$.\end{thm}

Our work here draws heavily from the computation in \cite{DW} of the complex analogue $ku_*(K_2)$ and $ku^*(K_2)$, which in turn relied on the connective Morava $K$-theory groups $k(1)_*(K_2)$ and $k(1)^*(K_2)$ determined in \cite{DRW}. Our primary focus in \cite{DW} was the $ku$-cohomology groups because of their product structure, but here our primary focus will be on the $ko$-homology groups, for historical reasons and because of the more familiar form of its ASS. We begin this introduction with a slight reformulation of results of \cite{DW} regarding $ku_*(K_2)$. This will enable us to describe the overall structure of $ko_*(K_2)$, and also to see the significant increase in complication of the $ko$ result as compared with $ku$.

The $ku$-cohomology groups $ku^*(K_2)$ were a combination of suspensions of three basic types of summands, $A_k$, $B_k$, and $S_{k,\ell}$, $1\le k<\ell$. It was observed in \cite[Section 7]{DW} that if $S_{k,\ell}$ is combined with two specific copies of $B_k$, we obtain something, which we now call $B_{k,\ell}$, that appears both in $ku^*(K_2)$ and, after switching to homology grading, in $ku_*(K_2)$. We denote by $\B_{k,\ell}$ the $ko$-homology analogue of this combination. The $ku$-homology analogue of $A_k$, which we denote by $\as_k$, was pictured when $k=5$ in \cite[Figure 4]{DW}, which we repeat here as Figure \ref{kuA5}. Short vertical lines indicate multiplication by 2, short diagonal lines are multiplication by $v\in ku_2$, and the long dashed lines, sometimes slightly curved, are exotic extensions ($\cdot2$).

\bigskip
 \begin{minipage}{6in}
\begin{fig}\label{kuA5}

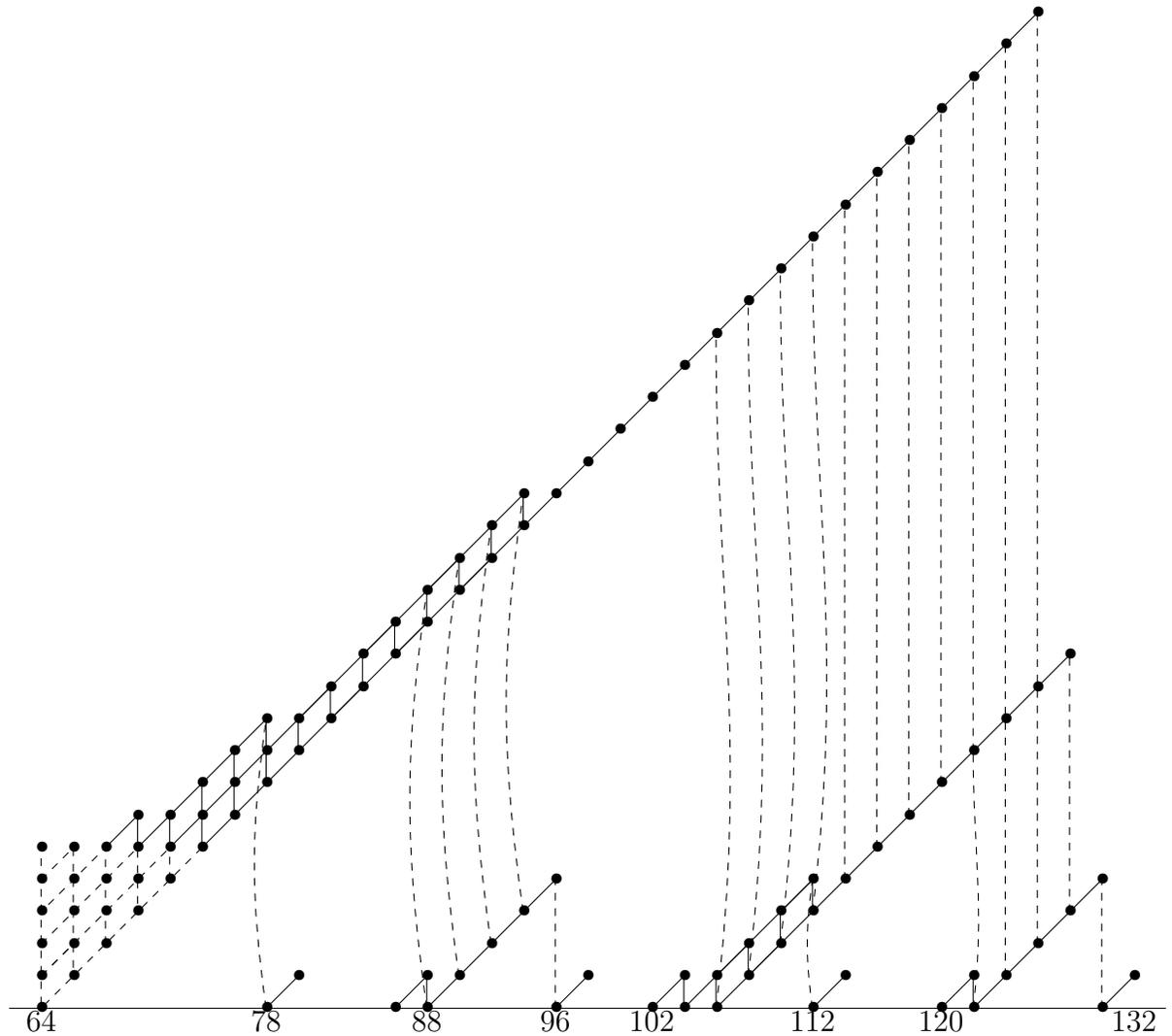
{\bf $\as_5$, the $ku$-homology analogue of $A_5$, from \cite{DW}}

\begin{center}

\begin{\tz}[scale=.44]
{\draw (0,5) -- (1,6) -- (1,5) -- (5,9) -- (5,7) -- (29,31);
\draw (-3,0) -- (33,0);
\draw (3,7) -- (3,5);
\draw (4,6) -- (4,8);
\draw [dashed] (3,5) -- (-2,0) -- (-2,5);
\draw [dashed] (-2,1) -- (2,5);
\draw [dashed] (-2,1) -- (-1,2) -- (-1,1);
\draw [dashed] (0,2) -- (0,5) -- (-2,3);
\draw [dashed] (-2,2) -- (1,5) -- (1,3);
\draw [dashed] (-2,4) -- (-1,5) -- (-1,2);
\draw [dashed] (2,4) -- (2,5);
\draw [dashed] (5,0) to[out=98, in=262] (5,9);
\draw [dashed] (10,0) to[out=98, in=262] (10,13);
\draw [dashed] (11,1) to[out=98, in=262] (11,14);
\draw [dashed] (12,2) to[out=98, in=262] (12,15);
\draw [dashed] (13,3) to[out=98, in=262] (13,16);
\draw [dashed] (22,0) to[out=98, in=262] (22,4);
\draw [dashed] (22,3) to[out=83, in=270] (22,24);
\draw [dashed] (21,2) to[out=83, in=270] (21,23);
\draw [dashed] (20,1) to[out=83, in=270] (20,22);
\draw [dashed] (19,0) to[out=83, in=270] (19,21);
\draw [dashed] (27,0) to[out=83, in=270] (27,8);
\draw [dashed] (14,0) -- (14,4);
\draw [dashed] (23,4) -- (23,25);
\draw [dashed] (24,5) -- (24,26);
\draw [dashed] (25,6) -- (25,27);
\draw [dashed] (26,7) -- (26,28);
\draw [dashed] (27,8) -- (27,29);
\draw [dashed] (28,1) -- (28,30);
\draw [dashed] (29,2) -- (29,31);
\draw [dashed] (30,3) -- (30,11);
\draw [dashed] (31,0) -- (31,4);
\node at (0,5) {\lb};
\node at (1,5) {\lb};
\node at (1,6) {\lb};
\node at (2,6) {\lb};
\node at (3,7) {\lb};
\node at (4,8) {\lb};
\node at (5,9) {\lb};
\node at (2,5) {\lb};
\node at (3,6) {\lb};
\node at (4,7) {\lb};
\node at (5,8) {\lb};
\node at (6,9) {\lb};
\node at (7,10) {\lb};
\node at (8,11) {\lb};
\node at (9,12) {\lb};
\node at (10,13) {\lb};
\node at (11,14) {\lb};
\node at (12,15) {\lb};
\node at (13,16) {\lb};
\node at (3,5) {\lb};
\node at (4,6) {\lb};
\node at (5,7) {\lb};
\node at (6,8) {\lb};
\node at (7,9) {\lb};
\node at (8,10) {\lb};
\node at (9,11) {\lb};
\node at (10,12) {\lb};
\node at (11,13) {\lb};
\node at (12,14) {\lb};
\node at (13,15) {\lb};
\node at (14,16) {\lb};
\node at (15,17) {\lb};
\node at (16,18) {\lb};
\node at (17,19) {\lb};
\node at (18,20) {\lb};
\node at (19,21) {\lb};
\node at (20,22) {\lb};
\node at (21,23) {\lb};
\node at (22,24) {\lb};
\node at (23,25) {\lb};
\node at (24,26) {\lb};
\node at (25,27) {\lb};
\node at (26,28) {\lb};
\node at (27,29) {\lb};
\node at (28,30) {\lb};
\node at (29,31) {\lb};
\node at (5,0) {\lb};
\node at (6,1) {\lb};
\node at (9,0) {\lb};
\node at (10,1) {\lb};
\node at (10,0) {\lb};
\node at (11,1) {\lb};
\node at (12,2) {\lb};
\node at (13,3) {\lb};
\node at (14,4) {\lb};
\node at (14,0) {\lb};
\node at (15,1) {\lb};
\node at (17,0) {\lb};
\node at (18,1) {\lb};
\node at (18,0) {\lb};
\node at (19,1) {\lb};
\node at (20,2) {\lb};
\node at (21,3) {\lb};
\node at (22,4) {\lb};
\node at (19,0) {\lb};
\node at (20,1) {\lb};
\node at (21,2) {\lb};
\node at (22,3) {\lb};
\node at (23,4) {\lb};
\node at (24,5) {\lb};
\node at (25,6) {\lb};
\node at (26,7) {\lb};
\node at (27,8) {\lb};
\node at (28,9) {\lb};
\node at (29,10) {\lb};
\node at (30,11) {\lb};
\node at (22,0) {\lb};
\node at (23,1) {\lb};
\node at (26,0) {\lb};
\node at (27,1) {\lb};
\node at (27,0) {\lb};
\node at (28,1) {\lb};
\node at (29,2) {\lb};
\node at (30,3) {\lb};
\node at (31,4) {\lb};
\node at (31,0) {\lb};
\node at (32,1) {\lb};

\draw (17,0) -- (18,1) -- (18,0) -- (22,4);
\draw (19,0) -- (30,11);
\draw (19,0) -- (19,1) -- (20,2) -- (20,1) -- (21,2) -- (21,3) -- (22,4) -- (22,3);
\node at (-2,0) {\lb};
\node at (-2,1) {\lb};
\node at (-2,2) {\lb};

\node at (-2,3) {\lb};
\node at (-2,4) {\lb};
\node at (-2,5) {\lb};
\node at (-1,1) {\lb};
\node at (-1,2) {\lb};
\node at (-1,3) {\lb};
\node at (-1,4) {\lb};
\node at (-1,5) {\lb};
\node at (0,2) {\lb};
\node at (0,3) {\lb};
\node at (0,4) {\lb};
\node at (1,3) {\lb};
\node at (1,4) {\lb};
\node at (2,4) {\lb};

\draw (2,6) -- (2,5) -- (13,16) -- (13,15);
\draw (6,8) -- (6,9) -- (7,10) -- (7,9) -- (8,10) -- (8,11) -- (9,12) -- (9,11) -- (10,12) -- (10,13) -- (11,14) -- (11,13) -- (12,14) -- (12,15);
\draw (3,5) -- (5,7);
\draw (9,0) -- (10,1) -- (10,0) -- (14,4);
\draw (5,0) -- (6,1);
\draw (14,0) -- (15,1);
\draw (22,0) -- (23,1);
\draw (26,0) -- (27,1) -- (27,0) -- (31,4);
\draw (31,0) -- (32,1);
\node at (-2,-.4) {$64$};
\node at (5,-.4) {$78$};
\node at (10,-.4) {$88$};
\node at (14,-.4) {$96$};
\node at (17,-.4) {$102$};
\node at (22,-.4) {$112$};
\node at (26,-.4) {$120$};
\node at (32,-.4) {$132$};
}

\end{\tz}
\end{center}
\end{fig}
\end{minipage}

\bigskip
The iterative structure of the $ku$ groups $A_k$ and $B_k$ can be seen by noting that the $ku$-homology analogue of $B_4$ is a desuspension of the portion of Figure \ref{kuA5} in grading $\ge102$, not including the upper edge, and the $A_4$ analogue can be obtained from that by placing a triangle similar to the one in the lower left corner of the chart, but one level smaller, beneath the classes in grading 102 to 106. The $ko$-homology analogue of $A_k$, which we denote by $\A_k$, is much more complicated. We picture $\A_k$ for $1\le k\le 6$ in Section \ref{Akexpls}, and in Section \ref{Akthmsec} state several results, \ref{Akthm}, \ref{E'def}, \ref{third}, \ref{E1'def}, and \ref{extnthm}, which give a complete determination of all $\A_k$.

In the $ku$ situation, each $B_{k,\ell}$ could be multiplied by any number of classes $z_j$, which just amounted to suspending by $\Si^{2^{j+2}+2}$ for each. In the $ko$ analogue, multiplying $\B_{k,\ell}$ by $z_j$ suspends by $2^{j+2}$ but also changes its form. We denote by $z^i\B_{k,\ell}$ the modified form of $\B_{k,\ell}$ after $i$ such changes of form, but without the various $\Si^{2^{j+2}}$'s. So $z^i$ can be thought of as a product of $i$ different $z_j$'s, each desuspended by $2^{j+2}$. It turns out that $z^4\B_{k,\ell}=\Si^8\B_{k,\ell}$. We will describe and illustrate $z^i\B_{k,\ell}$ in Section \ref{Bklsec}.

We can now state the theorem which, along with the aforementioned descriptions of $\A_k$ and $z^i
\B_{k,\ell}$, completely determines $ko_*(K_2)$. This will be proved in Section \ref{kosummandsec}.
\begin{thm}\label{main} There is an isomorphism of $ko_*$-modules
$$ko_*(K_2)\approx\bigoplus_{k\ge1} \bigoplus_{i\ge0}\Si^{2^{k+2}i}\A_k\oplus\ \bigoplus_{1\le k<\ell}\bigoplus_{i,j\ge0}\Si^{2^{k+2}i+2^{\ell+3}j}z^{\a(j)}\B_{k,\ell}$$
plus a trivial $ko_*$-module.\end{thm}
\noindent Here and elsewhere, $\a(j)$ denotes the number of 1's in the binary expansion of $j$. The trivial $ko_*$-module  could be calculated, but is not of interest. It is discussed at the end of Section \ref{kosummandsec}.

In Section \ref{cohsec}, we determine the $ko$-cohomology groups, $ko^*(K_2)$, and discuss the following interesting duality theorem.
\begin{thm} \label{dualitythm} There is an isomorphism of $ko_*$-modules
$$ko_*(K_2)\approx(ko^{*+6}(K_2))^\vee,$$
where $M^\vee=\Hom(M,\Z/2^\infty)$, the Pontryagin dual, localized at 2.
\end{thm}

In the first six sections, we describe our results in several ways. In Section \ref{outline}, we outline the proof, which occupies the subsequent seven sections. At the end of the paper, we discuss the application to Stiefel-Whitney classes of Spin manifolds, and the adaptation to
$ko$-cohomology groups of $K_2$.
\section{Examples of $\A_k$}\label{Akexpls}
In this section, we display charts of $\A_k$ for $1\le k\le 6$. We will make much use of these charts in later sections, as models for the general statements and proofs. These are ASS-type charts, but we frequently elevate filtrations\footnote{``Filtration'' refers to the vertical position in the chart.} of classes in order to better display the $ko_*$-module structure. In the transition from $E_2$ to $E_\infty$, there have been many differentials, which are not displayed here. The short diagonal lines represent multiplication by $\eta\in ko_1$.
In Figure \ref{A3chart}, we depict $\A_1$, $\A_2$, and $\A_3$ on the same set of axes. The classes in grading 8--12 are $\A_2$, and $\A_1$ and $\A_3$ are on the left and right, respectively.

\bigskip
\begin{minipage}{6in}
\begin{fig}\label{A3chart}

{\bf $\A_1$, $\A_2$, and $\A_3$}

\begin{center}

\begin{\tz}[scale=.5]

\draw (18,0) -- (49,0);
\draw (34,3) -- (34,0) -- (36,2);
\node at (37,3) {\lb};
\draw (38,3) -- (38,4) -- (37,3);
\draw (42,4) -- (44,6) -- (44,5) -- (45,6);
\node at (46,0) {\lb};
\node at (48,1) {\lb};
\node at (48,-.6) {$30$};
\node at (38,-.6) {$20$};
\node at (34,-.6) {$16$};
\node at (42,-.6) {$24$};
\node at (34,0) {\lb};
\node at (34,1) {\lb};
\node at (34,2) {\lb};
\node at (34,3) {\lb};
\node at (35,1) {\lb};
\node at (36,2) {\lb};
\node at (38,3) {\lb};
\node at (38,4) {\lb};
\node at (42,4) {\lb};
\node at (43,5) {\lb};
\node at (44,6) {\lb};
\node at (44,5) {\lb};
\node at (45,6) {\lb};
\node at (20,0) {\lb};
\node at (22,1) {\lb};
\node at (26,0) {\lb};
\node at (26,1) {\lb};
\node at (26,2) {\lb};
\node at (27,1) {\lb};
\node at (28,2) {\lb};
\node at (29,2) {\lb};
\node at (30,3) {\lb};
\draw (26,2) -- (26,0) -- (28,2);
\draw (29,2) -- (30,3);
\node at (20,-.6) {$2$};
\node at (22,-.6) {$4$};
\node at (26,-.6) {$8$};
\node at (30,-.6) {$12$};
\end{\tz}

\end{center}
\end{fig}
\end{minipage}
\bigskip

In Figure \ref{A4} and subsequently, the dashed lines are exotic extensions. One reason for our choice of filtrations is the structure of edges, which will be introduced in Section \ref{Akthmsec}.

\bigskip
\begin{minipage}{6in}
\begin{fig}\label{A4}

{\bf $\A_4$}

\begin{center}

\begin{\tz}[scale=.4]
\draw (31.5,0) -- (66.5,0);
\draw (32,4) -- (32,0) -- (34,2);
\draw (36,3) -- (36,5) -- (35,4);
\draw (40,5) -- (40,4) -- (42,6);
\draw (42,6.3) -- (43,7);
\draw (48,8) -- (50,10) -- (50,9) -- (52,11);
\draw (56,12) -- (58,14) -- (58,13);
\draw (58,1) -- (60,3);
\draw [dashed] (58,1) -- (58,13);
\draw [dashed] (62,0) -- (62,4);
\node at (32,0) {\lb};
\node at (32,1) {\lb};
\node at (32,2) {\lb};
\node at (32,3) {\lb};
\node at (32,4) {\lb};
\node at (33,1) {\lb};
\node at (34,2) {\lb};
\node at (36,3) {\lb};
\node at (36,4) {\lb};
\node at (36,5) {\lb};
\node at (35,4) {\lb};
\node at (40,4) {\lb};
\node at (40,5) {\lb};
\node at (41,5) {\lb};
\node at (42,6) {\lb};
\node at (42,6.3) {\lb};
\node at (43,7) {\lb};
\node at (44,0) {\lb};
\node at (46,1) {\lb};
\node at (44,7) {\lb};
\node at (48,8) {\lb};
\node at (49,9) {\lb};
\node at (50,10) {\lb};
\node at (50,9) {\lb};
\node at (51,10) {\lb};
\node at (52,11) {\lb};
\node at (56,12) {\lb};
\node at (57,13) {\lb};
\node at (58,14) {\lb};
\node at (58,13) {\lb};
\node at (52,0) {\lb};
\node at (54,0) {\lb};
\node at (56,1) {\lb};
\node at (58,1) {\lb};
\node at (59,2) {\lb};
\node at (60,3) {\lb};
\node at (62,0) {\lb};
\node at (62,4) {\lb};
\node at (66,1) {\lb};
\node at (32,-.6) {$32$};
\node at (36,-.6) {$36$};
\node at (40,-.6) {$40$};
\node at (44,-.6) {$44$};
\node at (48,-.6) {$48$};
\node at (52,-.6) {$52$};
\node at (56,-.6) {$56$};
\node at (60,-.6) {$60$};
\node at (64,-.6) {$64$};
\node at (66,-.6) {$66$};

\end{\tz}
\end{center}
\end{fig}
\end{minipage}
\bigskip

For $k\ge2$, it will be convenient for us to work with $\Si^{-2^{k+1}}\A_k$, which we denote by $\At_k$. The configuration of six classes which occurs in grading 48 to 52 in Figure \ref{A4} is often called a ``lightning flash,'' terminology which we will use frequently.

\begin{minipage}{6in}
\begin{fig}\label{A5}

{\bf $\At_5$}

\begin{center}

\begin{\tz}[scale=.44]
\draw (-.5,0) -- (34,0);
\node at (0,-.6) {$0$};
\node at (4,-.6) {$4$};
\node at (8,-.6) {$8$};
\node at (12,-.6) {$12$};
\node at (16,-.6) {$16$};
\node at (20,-.6) {$20$};
\node at (24,-.6) {$24$};
\node at (28,-.6) {$28$};
\node at (32,-.6) {$32$};
\draw (0,5) -- (0,0) -- (2,2);
\draw (3,5) -- (4,6) -- (4,3);
\draw (8,6) -- (8,4) -- (10,6);
\draw (10,7) -- (11,8);
\draw (12,8) -- (12,7);
\draw (16,9) -- (16,8) -- (18,10);
\draw (18,10.3) -- (20,12) -- (20,11);
\draw (24,13) -- (24,12) -- (26,14);
\node at (0,0) {\lb};
\node at (0,1) {\lb};
\node at (0,2) {\lb};
\node at (0,3) {\lb};
\node at (0,4) {\lb};
\node at (0,5) {\lb};
\node at (1,1) {\lb};
\node at (2,2) {\lb};
\node at (4,3) {\lb};
\node at (4,4) {\lb};
\node at (4,5) {\lb};
\node at (4,6) {\lb};
\node at (3,5) {\lb};
\node at (8,4) {\lb};
\node at (8,5) {\lb};
\node at (8,6) {\lb};
\node at (9,5) {\lb};
\node at (10,6) {\lb};
\node at (10,7) {\lb};
\node at (11,8) {\lb};
\node at (12,8) {\lb};
\node at (12,7) {\lb};
\node at (12,0) {\lb};
\node at (14,1) {\lb};
\node at (20,0) {\lb};
\node at (22,0) {\lb};
\node at (24,1) {\lb};
\node at (26,1) {\lb};
\node at (27,2) {\lb};
\node at (28,3) {\lb};
\node at (30,0) {\lb};
\node at (30,4) {\lb};
\node at (16,8) {\lb};
\node at (16,9) {\lb};
\node at (17,9) {\lb};
\node at (18,10) {\lb};
\node at (18,10.3) {\lb};
\node at (19,11.15) {\lb};
\node at (20,11) {\lb};
\node at (20,12) {\lb};
\node at (24,12) {\lb};
\node at (24,13) {\lb};
\node at (25,13) {\lb};
\node at (26,14) {\lb};
\node at (26,13.7) {\lb};
\draw (26,1) -- (28,3);
\draw [dashed] (26,1) -- (26,13.7);
\draw [dashed] (30,0) -- (30,4);
\node at (28,15) {\lb};
\draw (-1,17) -- (35.5,17);
\node at (0,16.4) {$32$};
\node at (4,16.4) {$36$};
\node at (8,16.4) {$40$};
\node at (12,16.4) {$44$};
\node at (16,16.4) {$48$};
\node at (20,16.4) {$52$};
\node at (24,16.4) {$56$};
\node at (28,16.4) {$60$};
\node at (32,16.4) {$64$};
\node at (34,1) {\lb};
\node at (4,17) {\lb};
\node at (8,17) {\lb};
\node at (8,18) {\lb};
\node at (10,18) {\lb};
\node at (11,19) {\lb};
\node at (12,20) {\lb};
\node at (12,19) {\lb};
\node at (14,17) {\lb};
\node at (14,20) {\lb};
\node at (16,21) {\lb};
\node at (18,18) {\lb};
\node at (18,21) {\lb};
\node at (19,22) {\lb};
\node at (20,23) {\lb};
\node at (22,24) {\lb};
\node at (26,25) {\lb};
\node at (30,28) {\lb};
\node at (24,17) {\lb};
\node at (26,18) {\lb};
\node at (27,18) {\lb};
\node at (28,19) {\lb};
\node at (30,20) {\lb};
\node at (34,21) {\lb};
\node at (34,17) {\lb};
\node at (35,18) {\lb};
\draw (8,17) -- (8,18);
\draw (10,18) -- (12,20) -- (12,19);
\draw (18,21) -- (20,23);
\draw (27,18) -- (28,19);
\draw (34,17) -- (35,18);
\node at (16,25) {\rm upper edge should be 12 higher};
\node at (12,20) {\lb};
\node at (0,21) {\lb};
\node at (1,22) {\lb};
\node at (2,23) {\lb};
\node at (2,22) {\lb};
\node at (3,23) {\lb};
\node at (4,24) {\lb};
\node at (8,25) {\lb};
\node at (9,26) {\lb};
\node at (10,27) {\lb};
\node at (10,26) {\lb};
\node at (16,29) {\lb};
\node at (17,30) {\lb};
\node at (18,31) {\lb};
\node at (18,30) {\lb};
\node at (24,33) {\lb};
\node at (25,34) {\lb};
\node at (26,35) {\lb};
\node at (26,34) {\lb};
\draw (0,21) -- (2,23) -- (2,22) -- (4,24);
\draw (8,25) -- (10,27) -- (10,26);
\draw (16,29) -- (18,31) -- (18,30);
\draw (24,33) -- (26,35) -- (26,34);
\draw [dashed] (18,21) -- (18,30);
\draw [dashed] (26,25) -- (26,34);
\draw [dashed] (34,17) -- (34,21);
\draw [dashed] (30,20) -- (30,28);
\node at (20,11) {\lb};
\end{\tz}
\end{center}
\end{fig}
\end{minipage}

\bigskip

\bigskip
\begin{minipage}{6in}
\begin{fig}\label{A6}

{\bf $\At_6$}

\begin{center}

\begin{\tz}[scale=.24]
\draw (68,0) -- (132,0);
\draw (72,0) -- (72,1);
\draw (74,1) -- (76,3) -- (76,1);
\node at (68,0) {\llb};
\node at (72,0) {\llb};
\node at (72,1) {\llb};
\node at (74,1) {\llb};
\node at (75,2) {\llb};
\node at (76,3) {\llb};
\node at (76,2) {\llb};
\node at (76,1) {\llb};
\node at (78,0) {\llb};
\node at (78,2) {\llb};
\node at (80,3) {\llb};
\draw (80,3) -- (80,4);
\node at (80,4) {\llb};
\node at (82,1) {\llb};
\node at (82,4) {\llb};
\node at (83,5) {\llb};
\node at (84,6) {\llb};
\node at (84,5) {\llb};
\node at (86,6) {\llb};
\node at (68,-.8) {$68$};
\node at (72,-.8) {$72$};
\node at (76,-.8) {$76$};
\node at (80,-.8) {$80$};
\node at (84,-.8) {$84$};
\node at (88,-.8) {$88$};
\node at (92,-.8) {$92$};
\node at (96,-.8) {$96$};
\node at (100,-.8) {$100$};
\node at (104,-.8) {$104$};
\node at (108,-.8) {$108$};
\node at (112,-.8) {$112$};
\node at (116,-.8) {$116$};
\node at (120,-.8) {$120$};
\node at (124,-.8) {$124$};
\node at (128,-.8) {$128$};
\node at (132,-.8) {$132$};
\draw (82,4) -- (84,6) -- (84,5);
\draw (91,1) -- (92,2);
\draw (98,0) -- (99,1);
\draw (98,11) -- (100,13);
\node at (88,7) {\llb};
\node at (90,8) {\llb};
\node at (92,9) {\llb};
\node at (94,10) {\llb};
\node at (96,11) {\llb};
\node at (98,11) {\llb};
\node at (99,12) {\llb};
\node at (100,13) {\llb};
\node at (102,14) {\llb};
\node at (106,15) {\llb};
\node at (110,18) {\llb};
\node at (114,19) {\llb};
\node at (118,22) {\llb};
\node at (122,23) {\llb};
\node at (126,26) {\llb};
\node at (72,8) {\llb};
\node at (73,9) {\llb};
\node at (74,10) {\llb};
\node at (74,9) {\llb};
\node at (80,12) {\llb};
\node at (81,13) {\llb};
\node at (82,14) {\llb};
\node at (82,13) {\llb};
\node at (96,20) {\llb};
\node at (97,21) {\llb};
\node at (98,22) {\llb};
\node at (98,21) {\llb};
\node at (88,16) {\llb};
\node at (89,17) {\llb};
\node at (90,18) {\llb};
\node at (90,17) {\llb};
\draw (72,8) -- (74,10) -- (74,9);
\draw (80,12) -- (82,14) -- (82,13);
\draw (88,16) -- (90,18) -- (90,17);
\draw (96,20) -- (98,22) -- (98,21);
\draw (104,24) -- (106,26) -- (106,25);
\draw (112,28) -- (114,30) -- (114,29);
\draw (120,32) -- (122,34) -- (122,33);
\node at (84,24) {\rm upper edge should be 28 higher};
\node at (104,24) {\llb};
\node at (105,25) {\llb};
\node at (106,26) {\llb};
\node at (106,25) {\llb};
\node at (112,28) {\llb};
\node at (113,29) {\llb};
\node at (114,30) {\llb};
\node at (114,29) {\llb};
\node at (120,32) {\llb};
\node at (121,33) {\llb};
\node at (122,34) {\llb};
\node at (122,33) {\llb};
\node at (88,0) {\llb};
\node at (90,1) {\llb};
\node at (91,1) {\llb};
\node at (92,2) {\llb};
\node at (94,3) {\llb};
\node at (98,0) {\llb};
\node at (99,1) {\llb};
\node at (98,4) {\llb};
\node at (104,0) {\llb};
\node at (106,1) {\llb};
\node at (107,1) {\llb};
\node at (108,2) {\llb};
\node at (108,1) {\llb};
\node at (110,2) {\llb};
\node at (112,3) {\llb};
\node at (114,0) {\llb};
\node at (115,1) {\llb};
\node at (114,4) {\llb};
\node at (115,4) {\llb};
\node at (116,5) {\llb};
\node at (118,6) {\llb};
\node at (122,7) {\llb};
\node at (126,10) {\llb};
\node at (130,11) {\llb};
\node at (122,0) {\llb};
\node at (123,1) {\llb};
\node at (124,1) {\llb};
\node at (126,2) {\llb};
\node at (130,3) {\llb};
\node at (131,0) {\llb};
\node at (132,1) {\llb};
\draw (107,1) -- (108,2) -- (108,1);
\draw (114,0) -- (115,1);
\draw (115,4) -- (116,5);
\draw (122,0) -- (123,1);
\draw (131,0) -- (132,1);
\draw [dashed] (98,0) -- (98,4);
\draw [dashed] (114,0) -- (114,4);
\draw [dashed] (98,11) -- (98,22);
\draw [dashed] (106,15) -- (106,26);
\draw [dashed] (114,19) -- (114,30);
\draw [dashed] (118,22) -- (118,6);
\draw [dashed] (122,0) -- (122,34);
\draw [dashed] (110,2) -- (110,18);
\draw [dashed] (126,2) -- (126,26);
\draw [dashed] (130,3) -- (130,11);
\draw (68,-38) -- (135,-38);
\node at (68,-38.8) {$0$};
\node at (72,-38.8) {$4$};
\node at (76,-38.8) {$8$};
\node at (80,-38.8) {$12$};
\node at (84,-38.8) {$16$};
\node at (88,-38.8) {$20$};
\node at (92,-38.8) {$24$};
\node at (96,-38.8) {$28$};
\node at (100,-38.8) {$32$};
\node at (104,-38.8) {$36$};
\node at (108,-38.8) {$40$};
\node at (112,-38.8) {$44$};
\node at (116,-38.8) {$48$};
\node at (120,-38.8) {$52$};
\node at (124,-38.8) {$56$};
\node at (128,-38.8) {$60$};
\node at (132,-38.8) {$64$};
\draw (70,-36)  -- (68,-38) -- (68,-32);
\draw (72,-35) -- (72,-31) -- (71,-32);
\draw (76,-31) -- (76,-34) -- (78,-32);
\draw (78,-30) -- (79,-29);
\draw (80,-31) -- (80,-29);
\draw (84,-28) -- (84,-30) -- (86,-28);
\draw (86,-27) -- (88,-25) -- (88,-27);
\draw (92,-24) -- (92,-26) -- (94,-24);
\draw (96,-23) -- (96,-22);
\draw (100,-21) -- (100,-22) -- (102,-20);
\draw (102.2,-20.2) -- (104,-18) -- (104,-19);
\draw (108,-17) -- (108,-18) -- (110,-16);
\draw (116,-13) -- (116,-14) -- (118,-12);
\draw (124,-9) -- (124,-10) -- (126,-8);
\draw (132,-6) -- (134,-4) -- (134,-5) -- (136,-3);
\node at (68,-38) {\llb};
\node at (68,-37) {\llb};
\node at (68,-36) {\llb};
\node at (68,-35) {\llb};
\node at (68,-34) {\llb};
\node at (68,-33) {\llb};
\node at (68,-32) {\llb};
\node at (69,-37) {\llb};
\node at (70,-36) {\llb};
\node at (72,-35) {\llb};
\node at (72,-34) {\llb};
\node at (72,-33) {\llb};
\node at (72,-32) {\llb};
\node at (72,-31) {\llb};
\node at (71,-32) {\llb};
\node at (76,-34) {\llb};
\node at (76,-33) {\llb};
\node at (76,-32) {\llb};
\node at (76,-31) {\llb};
\node at (77,-33) {\llb};
\node at (78,-32) {\llb};
\node at (78,-30) {\llb};
\node at (79,-29) {\llb};
\node at (80,-29) {\llb};
\node at (80,-30) {\llb};
\node at (80,-31) {\llb};
\node at (80,-38) {\llb};
\node at (82,-37) {\llb};
\node at (84,-30) {\llb};
\node at (84,-29) {\llb};
\node at (84,-28) {\llb};
\node at (85,-29) {\llb};
\node at (86,-28) {\llb};
\node at (86,-27) {\llb};
\node at (87,-26) {\llb};
\node at (88,-25) {\llb};
\node at (88,-26) {\llb};
\node at (88,-27) {\llb};
\node at (92,-26) {\llb};
\node at (92,-25) {\llb};
\node at (92,-24) {\llb};
\node at (93,-25) {\llb};
\node at (94,-24) {\llb};
\node at (94.3,-23) {\llb};
\node at (88,-38) {\llb};
\node at (90,-38) {\llb};
\node at (92,-37) {\llb};
\node at (94.3,-37) {\llb};
\node at (95.15,-36) {\llb};
\node at (96,-35) {\llb};
\node at (98,-38) {\llb};
\node at (98,-34) {\llb};
\node at (96,-23) {\llb};
\node at (96,-22) {\llb};
\node at (102,-37) {\llb};
\node at (104,-38) {\llb};
\draw (94.3,-37) -- (96,-35);
\draw [dashed] (94.3,-37) -- (94.3,-23);
\draw [dashed] (98,-38) -- (98,-34);
\node at (100,-22) {\llb};
\node at (100,-21) {\llb};
\node at (101,-21) {\llb};
\node at (102,-20) {\llb};
\node at (102.2,-20.2) {\llb};
\node at (103.1,-19.1) {\llb};
\node at (104,-18) {\llb};
\node at (108,-18) {\llb};
\node at (108,-17) {\llb};
\node at (109,-17) {\llb};
\node at (110,-16) {\llb};
\node at (110.2,-16.2) {\llb};
\node at (112,-15) {\llb};
\node at (116,-14) {\llb};
\node at (116,-13) {\llb};
\node at (117,-13) {\llb};
\node at (118,-12) {\llb};
\node at (118.2,-12.2) {\llb};
\node at (108,-38) {\llb};
\node at (108,-37) {\llb};
\node at (110,-37) {\llb};
\node at (111,-36) {\llb};
\node at (112,-35) {\llb};
\node at (112,-36) {\llb};
\node at (114,-38) {\llb};
\node at (114,-35) {\llb};
\node at (116,-34) {\llb};
\node at (118,-37) {\llb};
\node at (118.2,-34) {\llb};
\node at (119.1,-33) {\llb};
\node at (120,-32) {\llb};
\node at (120,-11) {\llb};
\node at (122,-31) {\llb};
\node at (126,-30) {\llb};
\node at (130,-27) {\llb};
\draw (108,-38) -- (108,-37);
\draw (110,-37) -- (112,-35) -- (112,-36);
\draw (118.2,-34) -- (120,-32);
\draw [dashed] (118.2,-34) -- (118.2,-12.2);
\node at (124,-10) {\llb};
\node at (124,-9) {\llb};
\node at (125,-9) {\llb};
\node at (126,-8) {\llb};
\node at (126.2,-8.2) {\llb};
\node at (128,-7) {\llb};
\node at (132,-6) {\llb};
\node at (133,-5) {\llb};
\node at (134,-4) {\llb};
\node at (134,-5) {\llb};
\node at (135,-4) {\llb};
\node at (136,-3) {\llb};
\node at (124,-38) {\llb};
\node at (126,-37) {\llb};
\node at (127,-37) {\llb};
\node at (128,-36) {\llb};
\node at (130,-35) {\llb};
\node at (134,-34) {\llb};
\node at (134,-38) {\llb};
\node at (135,-37) {\llb};
\node at (104,-19) {\llb};
\draw (127,-37) -- (128,-36);
\draw (134,-38) -- (135,-37);
\draw [dashed] (134,-38) -- (134,-34);
\draw [dashed] (130,-35) -- (130,-27);
\draw [dashed] (126.2,-30) -- (126.2,-8.2);

\end{\tz}
\end{center}
\end{fig}
\end{minipage}
\bigskip

\section{One description of $\A_k$}\label{Akthmsec}
Each $\A_k$ is structured in terms of edges and subedges. In this section, we describe the subedge structure, and give an explicit description of all edges. A more formulaic description of $\A_k$ is given in Section \ref{expsec}. The proof will be given in Section \ref{pfsec}, the culmination of several preliminary sections.

An edge $\E_{e,\ell}$ will always be suspended by some multiple of 8.
It occurs for the first time in $\At_\ell$, but, if $e>1$, it will occur again as a subedge of $\At_k$ for all $k>\ell$, often more than once for the same $\At_k$. If $e=4a+b>1$ with $b=(0,1,2,3)$, the bottom class of $\E_{e,\ell}$ is in grading $8a+(2,3,4,8)$ for any $\ell$, unless $e\equiv3$ mod 4 and $\ell-e=1$.

\begin{thm}\label{Akthm} $\At_k$ is built up recursively, beginning with $\E_{1,k}$, which begins in grading $0$. Then, for $1\le e\le k-2$, each occurrence of $\Si^D\E_{e,\ell}$ contains subedges $\Si^{D+2^{d+1}}\E_{e+1,e+d}$ for $2\le d\le \ell-e$.\end{thm}

For example, $\At_5$ has $\E_{1,5}$, then $\Si^8\E_{2,3}$, $\Si^{16}\E_{2,4}$, and $\Si^{32}\E_{2,5}$.
Then under $\Si^{16}\E_{2,4}$ is $\Si^{16+8}\E_{3,4}$, under $\Si^{32}\E_{2,5}$ are $\Si^{32+8}\E_{3,4}$ and $\Si^{32+16}\E_{3,5}$, and under $\Si^{48}\E_{3,5}$ is $\Si^{48+8}\E_{4,5}$.
Recalling that $\E_{2,\ell}$ begins in grading 4, $\E_{3,\ell}$ in grading 8, and $\E_{4,\ell}$ in grading 10, the structure of the subedges of $\At_5$ should be clear in Figure \ref{A5}.
Theorem \ref{Akthm} will be proved at the end of Section \ref{derivationsec}.

We are interested in the $ko_*$-module structure of  $ko_*(K_2)$, especially the action of $\eta$ and $v_1^4$. By $v_1^4$, we mean Adams or Bott periodicity, of bigrading $(8,4)$. Most of the time, we overlook the action of the generator of $ko_4$.
We include in Figure \ref{ko*} the well-known chart of $ko_*$, periodic with period $(8,4)$. The elements $\eta$ and $v_1^4$ are indicated with large dots.

\bigskip
\begin{minipage}{6in}
\begin{fig}\label{ko*}

{\bf $ko_*$}

\begin{center}

\begin{\tz}[scale=.4]
\draw (-.5,0) -- (11,0);
\draw [->] (0,0) -- (0,7);
\draw (0,0) -- (2,2);
\draw [->] (4,3) -- (4,7);
\draw [->] (8,4) -- (8,7);
\draw (8,4) -- (10,6);
\node at (11,7.7) {$\iddots$};
\node at (0,-.6) {$0$};
\node at (4,-.6) {$4$};
\node at (8,-.6) {$8$};
\node at (0,0) {\llb};
\node at (0,1) {\llb};
\node at (0,2) {\llb};
\node at (0,3) {\llb};
\node at (0,4) {\llb};
\node at (0,5) {\llb};
\node at (0,6) {\llb};
\node at (2,2) {\llb};
\node at (4,3) {\llb};
\node at (4,4) {\llb};
\node at (4,5) {\llb};
\node at (4,6) {\llb};
\node at (8,5) {\llb};
\node at (8,6) {\llb};
\node at (9,5) {\llb};
\node at (10,6) {\llb};
\node at (1,1) {\blb};
\node at (8,4) {\blb};

\end{\tz}
\end{center}
\end{fig}
\end{minipage}
\bigskip

We now work toward a description of $\E_{e,k}$ for $e\ge2$. We introduce some charts $M_k^i$ that will appear repeatedly throughout the paper. The derivation and significance of these charts will be discussed in Section \ref{kosummandsec}, when we begin our proof. By ``charts,'' we mean ASS diagrams, often involving filtration shifts of some elements.

We begin with charts $M^0_k$ for $k\ge4$, which are building blocks for our calculations.
  They are similar to familiar charts of $ko_*(RP^{2n})$ (e.g., \cite{D}). In fact, there are isomorphisms $M^0_{4\ell+4}\approx ko_*(RP^{8\ell+2})$ and $M^0_{4\ell+5}\approx ko_*(P^{8\ell+4})$. These charts were derived in \cite{DW2}; their derivation will be explained in Section \ref{kosummandsec}. For all $k$, all classes in $M^0_k$ are $v_1^4$-periodic.  All the charts $M^0_k$ have the same upper edge. The lower edge drops by 1 for each increase in $k$, with classes of negative filtration removed. The chart $M^0_4$ is a sequence of lightning flashes starting starting in position $(8i+1,4i)$ for $i\ge0$. In Figure \ref{Mcharts} we show the beginning of the charts for $5\le k\le7$.

\bigskip
\begin{minipage}{6in}
\begin{fig}\label{Mcharts}

{\bf $M^0_k$}

\begin{center}

\begin{\tz}[scale=.4]
\node at (8,4) {$\iddots$};
\node at (22,6.5) {$\iddots$};
\node at (37,5.5) {$\iddots$};

\draw (-.5,0) -- (6.5,0);
\draw (8.5,0) -- (21.5,0);
\draw (23.5,0) -- (36.5,0);
\node at (0,0) {\lb};
\draw (0,0) -- (2,2) -- (2,0) -- (4,2);
\node at (1,1) {\lb};
\node at (2,2) {\lb};
\node at (2,0) {\lb};
\node at (2,1) {\lb};
\node at (3,1) {\lb};
\node at (4,2) {\lb};
\node at (6,3) {\lb};
\node at (9,0) {\lb};
\node at (10,1) {\lb};
\node at (11,2) {\lb};
\node at (11,0) {\lb};\node at (11,1) {\lb};
\node at (12,0) {\lb};
\draw (9,0) -- (11,2) -- (11,0);
\node at (0,-.8) {$1$};
\node at (6,-.8) {$7$};
\node at (3,-2) {$k=5$};
\node at (13,1) {\lb};
\node at (15,2) {\lb};
\node at (15,3) {\lb};
\node at (19,3) {\lb};
\node at (19,4) {\lb};
\node at (19,5) {\lb};
\node at (19,6) {\lb};
\node at (20,4) {\lb};
\draw (15,2) -- (15,3);
\draw (17,4) -- (19,6) -- (19,3) -- (21,5);
\node at (9,-.8) {$1$};
\node at (15,-.8) {$7$};
\node at (19,-.8) {$11$};
\node at (15,-2) {$k=6$};
\node at (17,4) {\lb};
\node at (18,5) {\lb};
\node at (21,5) {\lb};
\node at (24,0) {\lb};
\node at (25,1) {\lb};
\node at (26,2) {\lb};
\node at (26,1) {\lb};
\node at (26,0) {\lb};
\node at (28,0) {\lb};
\node at (30,1) {\lb};
\draw (24,0) -- (26,2) -- (26,0);
\draw (30,1) -- (30,3);
\draw (32,4) -- (34,6) -- (34,2) -- (36,4);
\node at (24,-.8) {$1$};
\node at (30,-.8) {$7$};
\node at (34,-.8) {$11$};
\node at (10,1) {\lb};
\node at (30,2) {\lb};
\node at (30,3) {\lb};
\node at (32,4) {\lb};
\node at (33,5) {\lb};
\node at (34,6) {\lb};
\node at (34,5) {\lb};
\node at (34,4) {\lb};
\node at (34,3) {\lb};
\node at (34,2) {\lb};
\node at (35,3) {\lb};
\node at (36,4) {\lb};
\draw (12,0) -- (13,1);
\node at (30,-2) {$k=7$};

\end{\tz}
\end{center}
\end{fig}
\end{minipage}
\bigskip

Explicitly, $M^0_k$ has, for $i\ge0$,
\begin{itemize}
    \item $0$ in grading 0 and 6 mod 8,
    \item $\zt$ in grading $8i+1$ and $8i+2$ of filtration $4i$ and $4i+1$, respectively.
    \item $\zt$ in grading $8i+4$ and $8i+5$ of filtration $4i-k+6$ and $4i-k+7$, respectively, if the filtration is $\ge0$, else 0,
    \item $\Z/2^{k-4}$ in grading $8i+7$ with generator of filtration $4i-k+8$ if $4i-k+8\ge0$, else $\Z/2^{4i+4}$ with generator of filtration 0, and
    \item $\Z/2^{k-2}$ in grading $8i+3$ with generator of filtration $4i-k+5$ if $4i-k+5\ge0$, else $\Z/2^{4i+3}$ with generator of filtration 0.
\end{itemize}

Next we define $M_k^i$ for $k\ge4$ and $i\ge0$ to be $M^0_k$ with classes of filtration less than $i$ removed, and filtrations of other classes  decreased by $i$. In Figure \ref{M6}, we depict $\M_6^i$ for $1\le i\le3$. All elements are acted on freely by $v_1^4$. It follows that $M_k^{i+4}=\Si^8M_k^i$.

\bigskip
\begin{minipage}{6in}
\begin{fig}\label{M6}

{\bf $M_6^i$}

\begin{center}

\begin{\tz}[scale=.4]
\node at (14,6) {$\iddots$};
\node at (29,5) {$\iddots$};
\node at (42,5.5) {$\iddots$};
\draw (1,0) -- (13.5,0);
\draw (17.5,0) -- (28,0);
\draw (32.5,0) -- (41,0);
\draw (2,0) -- (3,1) -- (3,0);
\draw (7,1) -- (7,2);
\draw (9,3) -- (11,5) -- (11,2) -- (13,4);
\draw (22,0) -- (22,1);
\draw (24,2) -- (26,4) -- (26,1) -- (28,3);
\draw (35,1) -- (37,3) -- (37,0) -- (39,2);
\draw (41,3) -- (41,4);
\node at (2,0) {\lb};
\node at (3,1) {\lb};
\node at (3,0) {\lb};
\node at (5,0) {\lb};
\node at (7,1) {\lb};
\node at (7,2) {\lb};
\node at (9,3) {\lb};
\node at (10,4) {\lb};
\node at (11,5) {\lb};
\node at (11,4) {\lb};
\node at (11,3) {\lb};
\node at (11,2) {\lb};
\node at (12,3) {\lb};
\node at (13,4) {\lb};
\node at (18,0) {\lb};
\node at (22,0) {\lb};
\node at (22,1) {\lb};
\node at (24,2) {\lb};
\node at (25,3) {\lb};
\node at (26,4) {\lb};
\node at (26,3) {\lb};
\node at (26,2) {\lb};
\node at (26,1) {\lb};
\node at (27,2) {\lb};
\node at (28,3) {\lb};
\node at (33,0) {\lb};
\node at (35,1) {\lb};
\node at (36,2) {\lb};
\node at (37,3) {\lb};
\node at (37,2) {\lb};
\node at (37,1) {\lb};
\node at (37,0) {\lb};
\node at (38,1) {\lb};
\node at (39,2) {\lb};
\node at (41,3) {\lb};
\node at (41,4) {\lb};
\node at (2,-.7) {$2$};
\node at (7,-.7) {$7$};
\node at (11,-.7) {$11$};
\node at (18,-.7) {$3$};
\node at (22,-.7) {$7$};
\node at (26,-.7) {$11$};
\node at (33,-.7) {$7$};
\node at (37,-.7) {$11$};
\node at (41,-.7) {$15$};
\node at (8,-1.7) {$M_6^1$};
\node at (23,-1.7) {$M_6^2$};
\node at (38,-1.7) {$M_6^3$};

\end{\tz}
\end{center}
\end{fig}
\end{minipage}
\bigskip

A variant of $M_4^s$ that will be useful is given in the following definition.
\begin{defin}\label{Mhdef} For $0\le s\le3$ and $s\ge0$, the chart $\Mh_4^s$ is formed from $M_4^0$ by removing classes of grading $\le s$. If $s\ge4$, $\Mh_4^s=\Si^8\Mh_4^{s-4}$.\end{defin}

As $\Mh_4^s$ will always be combined into other charts, its filtration as an individual entity is irrelevant and undefined.
Note that $\Mh_4^s=M_4^s$ if $s\equiv0,1$ mod 4, while if $s\equiv2,3$ mod 4, $\Mh_4^s$ is formed from $M_4^s$ by adjoining one class. In Figure \ref{Mhatfig}, we picture $\Mh_4^s$ for $s\equiv2,3$ mod 4.

\bigskip
\begin{minipage}{6in}
\begin{fig}\label{Mhatfig}

{\bf $\Mh_4^s$ for $s\equiv2,3$ mod 4}

\begin{center}

\begin{\tz}[scale=.5]
\node at (0,-.8) {$8i+3$};
\node at (20,.2) {$8i+4$};
\draw (0,1) -- (0,0) -- (2,2);
\draw (6,3) -- (8,5) -- (8,4) -- (10,6);
\draw (20,1) -- (21,2);
\draw (25,3) -- (27,5) -- (27,4) -- (29,6);
\node at (0,0) {\lb};
\node at (0,1) {\lb};
\node at (1,1) {\lb};
\node at (2,2) {\lb};
\node at (6,3) {\lb};
\node at (7,4) {\lb};
\node at (8,5) {\lb};
\node at (8,4) {\lb};
\node at (9,5) {\lb};
\node at (10,6) {\lb};
\node at (20,1) {\lb};
\node at (21,2) {\lb};
\node at (25,3) {\lb};
\node at (26,4) {\lb};
\node at (27,5) {\lb};
\node at (27,4) {\lb};
\node at (28,5) {\lb};
\node at (29,6) {\lb};
\node at (3,3.5) {$\Mh_4^{4i+2}$};
\node at (23,3.5) {$\Mh_4^{4i+3}$};
\node at (11,7) {$\iddots$};
\node at (30,7) {$\iddots$};

\end{\tz}
\end{center}
\end{fig}
\end{minipage}
\bigskip

The following definition will be useful.
\begin{defin} \label{stable} A chart is {\em stably} $\Si^iM_k$ if it agrees with $\Si^iM_k^0$ in sufficiently large  grading, without regard for filtration.\end{defin}

Now we define what we call pre-edges $\E'_{e,\ell}$. This formulation will be derived in Section \ref{derivationsec}.

\begin{defin}\label{E'def}For $2\le e< \ell$, $\E'_{e,\ell}$ is formed from the following sequence of charts.
$$\Si M_{\ell-e+3}^e\lar\Si^8M_4^{e-1}\lar\cdots\lar \Si^{2^i}M_4^{e-1}\lar\cdots\lar\Si^{2^{\ell-e+1}}M_4^{e-1}\lar \Si^{2^{\ell-e+2}+2}\Mh_4^{e-2}$$
Working from left to right, each $\Si^{2^i}M_4^{e-1}$  is placed so that there is a $d_1$ differential from its generators in grading $1$, $3$, $4$, and $5$ mod $8$ to the upper edge of the chart resulting from the preceding steps, and there are $\eta$ extensions on its top classes in grading $3$ mod $8$. The chart resulting after incorporating $\Si^{2^i}M_4^{e-1}$ is stably $\Si M_{\ell-e+5-i}$. The $\Si^{2^{\ell-e+2}+2}\Mh_4^{e-2}$ is placed so that all its classes support $d_1$ differentials.
 The pre-edge $\E'_{e,\ell}$ is the chart remaining after all these $d_1$'s. It vanishes in the range of $\Si^{2^{\ell-e+2}+2}\Mh_4^{e-2}$. \end{defin}

We illustrate in Figure \ref{E'expl} the forming of $\E'_{4,7}$, which is derived from
$$\Si M_6^4\lar \Si^8M_4^3\lar \Si^{16}M_4^3\lar \Si^{34}\Mh_4^2.$$
At each step, the chart being adjoined is indicated in red, while the black part is the result of the preceding step. Then $\E'_{4,7}$ consists of the classes remaining in the lower part of Figure \ref{E'expl} after the indicated differentials.

\bigskip
\begin{minipage}{6in}
\begin{fig}\label{E'expl}

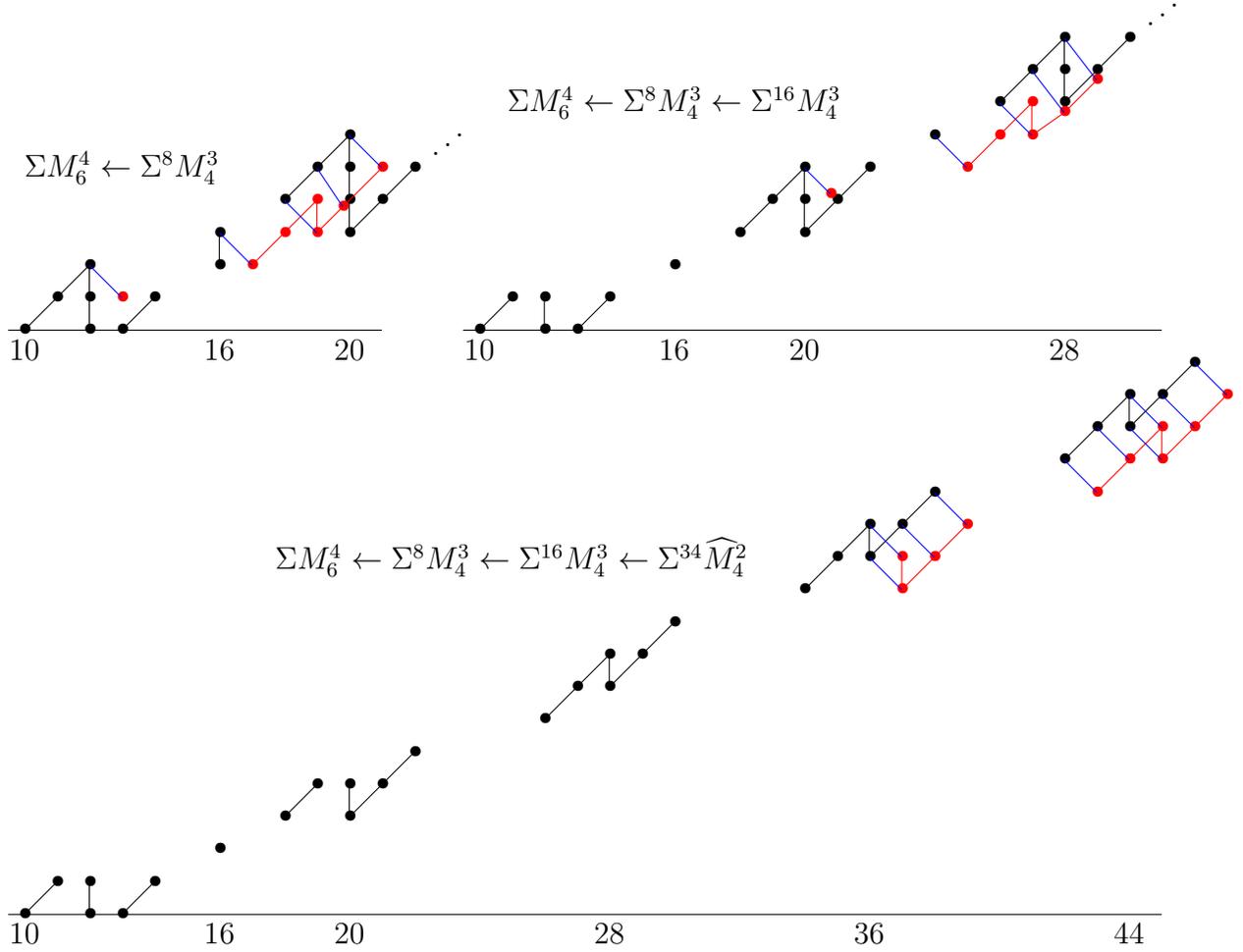
{\bf Forming $\E'_{4,7}$}

\begin{center}

\begin{\tz}[scale=.44]
\draw (10,0) -- (11,1);
\draw (12,0) -- (12,1);
\draw (13,0) -- (14,1);
\draw (18,3) -- (19,4);
\draw (22,5) -- (20,3) -- (20,4);
\draw (26,6) -- (28,8) -- (28,7) -- (30,9);
\draw (34,10) -- (36,12) -- (36,11) -- (38,13);
\draw (42,14) -- (44,16) -- (44,15) -- (46,17);
\draw (10,18) -- (12,20) -- (12,18);
\draw (13,18) -- (14,19);
\draw (16,20) -- (16,21);
\draw (18,22) -- (20,24) -- (20,21) -- (22,23);
\draw (24,18) -- (25,19);
\draw (26,18) -- (26,19);
\draw (27,18) -- (28,19);
\draw (32,21) -- (34,23) -- (34,21) -- (36,23);
\draw (40,25) -- (42,27) -- (42,25) -- (44,27);
\node at (10,0) {\lb};
\node at (11,1) {\lb};
\node at (12,0) {\lb};
\node at (12,1) {\lb};
\node at (13,0) {\lb};
\node at (14,1) {\lb};
\node at (16,2) {\lb};
\node at (18,3) {\lb};
\node at (19,4) {\lb};
\node at (20,3) {\lb};
\node at (20,4) {\lb};
\node at (21,4) {\lb};
\node at (22,5) {\lb};
\node at (26,6) {\lb};
\node at (27,7) {\lb};
\node at (28,8) {\lb};
\node at (28,7) {\lb};
\node at (29,8) {\lb};
\node at (30,9) {\lb};
\node at (34,10) {\lb};
\node at (35,11) {\lb};
\node at (36,12) {\lb};
\node at (36,11) {\lb};
\node at (37,12) {\lb};
\node at (38,13) {\lb};
\node at (42,14) {\lb};
\node at (43,15) {\lb};
\node at (44,16) {\lb};
\node at (44,15) {\lb};
\node at (45,16) {\lb};
\node at (46,17) {\lb};
\node at (10,18) {\lb};
\node at (11,19) {\lb};
\node at (12,20) {\lb};
\node at (12,19) {\lb};
\node at (12,18) {\lb};
\node at (13,18) {\lb};
\node at (14,19) {\lb};
\node at (16,20) {\lb};
\node at (16,21) {\lb};
\node at (18,22) {\lb};
\node at (19,23) {\lb};
\node at (20,24) {\lb};
\node at (20,23) {\lb};
\node at (20,22) {\lb};
\node at (20,21) {\lb};
\node at (21,22) {\lb};
\node at (22,23) {\lb};
\node at (24,18) {\lb};
\node at (25,19) {\lb};
\node at (26,18) {\lb};
\node at (26,19) {\lb};
\node at (27,18) {\lb};
\node at (28,19) {\lb};
\node at (30,20) {\lb};
\node at (32,21) {\lb};
\node at (33,22) {\lb};
\node at (34,23) {\lb};
\node at (34,22) {\lb};
\node at (34,21) {\lb};
\node at (35,22) {\lb};
\node at (36,23) {\lb};
\node at (38,24) {\lb};
\node at (40,25) {\lb};
\node at (41,26) {\lb};
\node at (42,27) {\lb};
\node at (42,26) {\lb};
\node at (42,25) {\lb};
\node at (43,26) {\lb};
\node at (44,27) {\lb};
\node at (23,24) {$\iddots$};
\node at (45,28) {$\iddots$};
\draw (9.5,0) -- (45,0);

\draw (9.5,18) -- (21,18);
\draw (23.5,18) -- (45,18);
\node at (10,-.6) {$10$};
\node at (16,-.6) {$16$};
\node at (20,-.6) {$20$};
\node at (28,-.6) {$28$};
\node at (36,-.6) {$36$};
\node at (44,-.6) {$44$};
\node at (10,17.4) {$10$};
\node at (16,17.4) {$16$};
\node at (20,17.4) {$20$};
\node at (24,17.4) {$10$};
\node at (30,17.4) {$16$};
\node at (34,17.4) {$20$};
\node at (42,17.4) {$28$};
\node at (13,23) {$\Si M_6^4\lar \Si^8M_4^3$};
\node at (30,25) {$\Si M_6^4\lar \Si^8M_4^3\lar \Si^{16}M_4^3$};
\node at (25,11)  {$\Si M_6^4\lar \Si^8M_4^3\lar \Si^{16}M_4^3\lar \Si^{34}\Mh_4^2$};
\draw [red] (17,20) -- (19,22) -- (19,21) -- (21,23);
\draw [red] (39,23) -- (41,25) -- (41,24) --  (42,24.7) -- (43,25.7);
\draw [red] (37,11) -- (37,10) -- (39,12);
\draw [red] (43,13) -- (45,15) -- (45,14) -- (47,16);
\node [red] at (37,11) {\lb};
\node [red] at (37,10) {\lb};
\node [red] at (38,11) {\lb};
\node [red] at (39,12) {\lb};
\node [red] at (43,13) {\lb};
\node [red] at (44,14) {\lb};
\node [red] at (45,15) {\lb};
\node [red] at (45,14) {\lb};
\node [red] at (46,15) {\lb};
\node [red] at (47,16) {\lb};
\node [red] at (13,19) {\lb};
\node [red] at (17,20) {\lb};
\node [red] at (18,21) {\lb};
\node [red] at (19,22) {\lb};
\node [red] at (19,21) {\lb};
\node [red] at (19.8,21.8) {\lb};
\node [red] at (21,23) {\lb};
\node [red] at (34.8,22.2) {\lb};
\node [red] at (39,23) {\lb};
\node [red] at (40,24) {\lb};
\node [red] at (41,25) {\lb};
\node [red] at (41,24) {\lb};
\node [red] at (42,24.7) {\lb};
\node [red] at (43,25.7) {\lb};
\draw [blue] (13,19) -- (12,20);
\draw [blue] (17,20) -- (16,21);
\draw [blue] (19,21) -- (18,22);
\draw [blue] (19.8,21.8) -- (19,23);
\draw [blue] (21,23) -- (20,24);
\draw [blue] (34.8,22.2) -- (34,23);
\draw [blue] (39,23) -- (38,24);
\draw [blue] (41,24) -- (40,25);
\draw [blue] (42,24.7) -- (41,26);
\draw [blue] (43,25.7) -- (42,27);
\draw [blue] (37,10) -- (36,11);
\draw [blue] (37,11) -- (36,12);
\draw [blue] (38,11) -- (37,12);
\draw [blue] (39,12) -- (38,13);
\draw [blue] (43,13) -- (42,14);
\draw [blue] (44,14) -- (43,15);
\draw [blue] (45,15) -- (44,16);
\draw [blue] (45,14) -- (44,15);
\draw [blue] (46,15) -- (45,16);
\draw [blue] (47,16) -- (46,17);

\end{\tz}
\end{center}
\end{fig}
\end{minipage}
\bigskip

Since $M_k^{e+4}=\Si^8 M_k^e$ and $\Mh_4^{e+4}=\Si^8\Mh_4^e$ for $e\ge0$, it follows that $\E'_{e+4,\ell+4}=\Si^8\E'_{e,\ell}$ for $e\ge2$. So it suffices to study $\E'_{e,\ell}$ for $e\le5$.

Recall from Theorem \ref{Akthm} that beneath any edge $\Si^D\E_{e,\ell}$ in $\At_k$ there occur subedges $\Si^{D+2^{d+1}}\E_{e+1,e+d}$ for $e+2\le e+d\le\ell$. The same is true for $\E'_{e,\ell}$.
The following theorem will be proved in Section \ref{pfsec}.
\begin{thm} For $e\ge1$, there are differentials from $\Si^{2^{d+1}}\E'_{e+1,e+d}$ into $\E'_{e,\ell}$ from all classes in grading $4$ or  $5$ mod $8$, except filtration-$0$ classes $x$ in grading $4$ mod $8$ with $\eta x=0$.  For $e\ge2$, $\E_{e,\ell}$ is formed from $\E'_{e,\ell}$ after removing all classes either supporting or hit by differentials. This is true, after appropriate suspension, of all occurrences of $\E_{e,\ell}$    as edges or subedges.\label{third}
\end{thm}

In Figures \ref{four}, \ref{2,6}, \ref{3,7}, \ref{4,8}, and \ref{5,9}, we display $\E'_{e,e+d}$ for $2\le e\le 5$ and $1\le d\le4$. We circle classes supporting differentials, and use a larger dot for classes hit by differentials.
Then $\E_{e,e+d}$ consists of uncircled small dots. 

\bigskip
\begin{minipage}{6in}
\begin{fig}\label{four}

{\bf $\E'_{e,\ell}$ and differentials}

\begin{center}

\begin{\tz}[scale=.37]
\draw (0,0) -- (1.5,0);
\draw (4.5,0) -- (14.5,0);
\draw (17.5,0) -- (43,0);
\draw (0,0) -- (1,1) -- (1,0);
\draw (5,0) -- (6,1) -- (6,0);
\draw (12,2) -- (14,4) -- (14,3);
\draw (18,0) -- (19,1) -- (19,0);
\draw (25,2) -- (27,4) -- (27,2) -- (29,4);
\draw (33,5) -- (35,7) -- (35,6) -- (37,8);
\draw (41,9) -- (43,11) -- (43,10);
\draw [dashed] (3,-.6) -- (3,5);
\draw [dashed] (16,-.6) -- (16,7);
\draw (0,12) -- (1.5,12);
\draw (4.5,12) -- (14,12);
\draw (17.5,12) -- (43,12);
\draw (34,18) -- (36,20) -- (36,19) -- (38,21);
\draw (42,22) -- (43,23);
\draw (0,24) -- (4,24);
\draw (6.5,24) -- (17,24);
\draw (20,24) -- (46,24);
\draw (38,31) -- (40,33) -- (40,32) -- (42,34);
\draw (0,37) -- (2,37);
\draw (4.5,37) -- (15,37);
\draw (18,37) -- (44,37);
\draw (40,45) -- (42,47) -- (42,46) -- (44,48);
\node at (0,-.6) {$11$};
\node at (5,-.6) {$11$};
\node at (14,-.6) {$20$};
\node at (8,-.6) {$14$};
\node at (18,-.6) {$11$};
\node at (23,-.6) {$16$};
\node at (27,-.6) {$20$};
\node at (35,-.6) {$28$};
\node at (43,-.6) {$36$};
\node at (0,11.4) {$10$};
\node at (5,11.4) {$10$};
\node at (7,11.4) {$12$};
\node at (9,11.4) {$14$};
\node at (14,11.4) {$19$};
\node at (18,11.4) {$10$};
\node at (20,11.4) {$12$};
\node at (24,11.4) {$16$};
\node at (28,11.4) {$20$};
\node at (36,11.4) {$28$};
\node at (43,12.6) {$35$};
\draw [dashed] (3,11.4) -- (3,16);
\draw [dashed] (16,11.4) -- (16,18);
\node at (0,23.4) {$6$};
\node at (4,23.4) {$10$};
\node at (7,23.4) {$8$};
\node at (11,23.4) {$12$};
\node at (17,23.4) {$18$};
\node at (20,23.4) {$8$};
\node at (24,23.4) {$12$};
\node at (28,23.4) {$16$};
\node at (32,23.4) {$20$};
\node at (40,23.4) {$28$};
\node at (46,23.4) {$34$};
\draw [dashed] (5.5,23.4) -- (5.5,26);
\draw [dashed] (18.5,23.4) -- (18.5,30);
\node at (0,36.4) {$4$};
\node at (2,36.4) {$6$};
\node at (5,36.4) {$4$};
\node at (9,36.4) {$8$};
\node at (13,36.4) {$12$};
\node at (18,36.4) {$4$};
\node at (22,36.4) {$8$};
\node at (26,36.4) {$12$};
\node at (34,36.4) {$20$};
\node at (42,36.4) {$28$};
\draw [dashed] (3.5,36.4) -- (3.5,40);
\draw [dashed] (16.5,36.4) -- (16.5,43);
\draw (7,0) -- (8,1);
\node at (0,0) {\llb};
\node at (1,1) {\llb};
\node at (5,0) {\llb};
\node at (6,1) {\llb};
\node at (6,0) {\llb};
\node at (8,1) {\llb};
\node at (12,2) {\llb};
\node at (13,3) {\llb};
\node at (18,0) {\llb};
\node at (19,1) {\llb};
\node at (19,0) {\llb};
\node at (21,0) {\llb};
\node at (23,1) {\llb};
\node at (25,2) {\llb};
\node at (26,3) {\llb};
\node at (27,3) {\llb};
\node at (29,4) {\llb};
\node at (33,5) {\llb};
\node at (34,6) {\llb};
\node at (35,7) {\llb};
\node at (37,8) {\llb};
\node at (41,9) {\llb};
\node at (14,4) {\blb};
\node at (27,4) {\blb};
\node at (42,10) {\blb};
\node at (43,11) {\blb};
\draw  (1,0) circle (.3);
\draw  (7,0) circle (.3);
\draw  (14,3) circle (.3);
\draw  (27,2) circle (.3);
\draw  (28,3) circle (.3);
\draw  (35,6) circle (.3);
\draw  (36,7) circle (.3);
\draw  (43,10) circle (.3);
\draw (0,12) -- (1,13);
\draw (5,12) -- (6,13);
\draw (7,13) -- (7,12) -- (9,14);
\draw (13,15) -- (14,16);
\draw (18,12) -- (19,13);
\draw (20,12) -- (20,13);
\draw (21,12) -- (22,13);
\draw (26,15) -- (27,16);
\draw (28,16) -- (28,15) -- (30,17);
\node at (0,12) {\llb};
\node at (1,13) {\llb};
\node at (5,12) {\llb};
\node at (6,13) {\llb};
\node at (7,13) {\llb};
\node at (9,14) {\llb};
\node at (13,15) {\llb};
\node at (18,12) {\llb};
\node at (19,13) {\llb};
\node at (20,12) {\llb};
\node at (20,13) {\llb};
\node at (22,13) {\llb};
\node at (24,14) {\llb};
\node at (26,15) {\llb};
\node at (28,16) {\llb};
\node at (30,17) {\llb};
\node at (34,18) {\llb};
\node at (35,19) {\llb};
\node at (38,21) {\llb};
\node at (42,22) {\llb};
\draw (7,12) circle (.3);
\draw (8,13) circle (.3);
\draw (21,12) circle (.3);
\draw (28,15) circle (.3);
\draw (29,16) circle (.3);
\draw (36,19) circle (.3);
\draw (37,20) circle (.3);
\node at (14,16) {\blb};
\node at (27,16) {\blb};
\node at (36,20) {\blb};

\node at (43,23) {\blb};
\node at (1,3) {$5,6$};
\node at (8,4) {$5,7$};
\node at (27,6) {$5,8$};
\node at (.5,15) {$4,5$};
\node at (9,16) {$4,6$};
\node at (28,19) {$4,7$};
\draw (10,25) -- (11,26) -- (11,25) -- (13,27);
\draw (23,25) -- (24,26) -- (24,24) -- (26,26);
\draw (31,28) -- (32,29) -- (32,28) -- (34,30);
\node at (0,24) {\llb};
\node at (4,25) {\llb};
\node at (7,24) {\llb};
\node at (9,25) {\llb};
\node at (10,25) {\llb};
\node at (11,26) {\llb};
\node at (13,27) {\llb};
\node at (17,28) {\llb};
\node at (20,24) {\llb};
\node at (22,25) {\llb};
\node at (23,25) {\llb};
\node at (24,26) {\llb};
\node at (24,25) {\llb};
\node at (26,26) {\llb};
\node at (28,27) {\llb};
\node at (30,28) {\llb};
\node at (31,28) {\llb};
\node at (32,29) {\llb};
\node at (34,30) {\llb};
\node at (38,31) {\llb};
\node at (42,34) {\llb};
\node at (46,35) {\llb};
\node at (39,32) {\blb};
\node at (40,33) {\blb};
\draw (11,25) circle (.3);
\draw (12,26) circle (.3);
\draw (24,24) circle (.3);
\draw (25,25) circle (.3);
\draw (32,28) circle (.3);
\draw (33,29) circle (.3);
\draw (40,32) circle (.3);
\draw (41,33) circle (.3);
\node at (2,27) {$3,4$};
\node at (10,28) {$3,5$};
\node at (28,30) {$3,6$};
\draw (1,37) -- (2,38);
\draw (11,38) -- (13,40) -- (13,39) -- (15,41);
\draw (22,37) -- (22,38);
\draw (24,38) -- (26,40) -- (26,38) -- (28,40);
\draw (32,41) -- (34,43) -- (34,42) -- (36,44);
\node at (0,37) {\llb};
\node at (1,37) {\llb};
\node at (2,38) {\llb};
\node at (5,37) {\llb};
\node at (7,37) {\llb};
\node at (9,38) {\llb};
\node at (11,38) {\llb};
\node at (12,39) {\llb};
\node at (13,40) {\llb};
\node at (15,41) {\llb};
\node at (18,37) {\llb};
\node at (22,37) {\llb};
\node at (22,38) {\llb};
\node at (24,38) {\llb};
\node at (25,39) {\llb};
\node at (26,40) {\llb};
\node at (28,40) {\llb};
\node at (30,41) {\llb};
\node at (32,41) {\llb};
\node at (33,42) {\llb};
\node at (34,43) {\llb};
\node at (36,44) {\llb};
\node at (40,45) {\llb};
\node at (44,48) {\llb};
\draw (1,37) circle (.3);
\draw (13,39) circle (.3);
\draw (14,40) circle (.3);
\node at (26,39) {\llb};
\draw (26,38) circle (.3);
\draw (27,39) circle (.3);
\draw (34,42) circle (.3);
\draw (35,43) circle (.3);
\draw (42,46) circle (.3);
\draw (43,47) circle (.3);
\node at (41,46) {\blb};
\node at (42,47) {\blb};
\node at (1,40) {$2,3$};
\node at (8,40) {$2,4$};
\node at (29,43) {$2,5$};

\end{\tz}
\end{center}
\end{fig}
\end{minipage}
\bigskip

 One added feature in Figures \ref{2,6}, \ref{3,7}, \ref{4,8}, and \ref{5,9} is that for the classes hit by differentials, we include below them the name of the subedge that supported the differential. The classes supporting those differentials can be seen in Figure \ref{four}, where they appear with circles. Note that these are complete figures; they are finite charts.

\bigskip
\begin{minipage}{6in}
\begin{fig}\label{2,6}

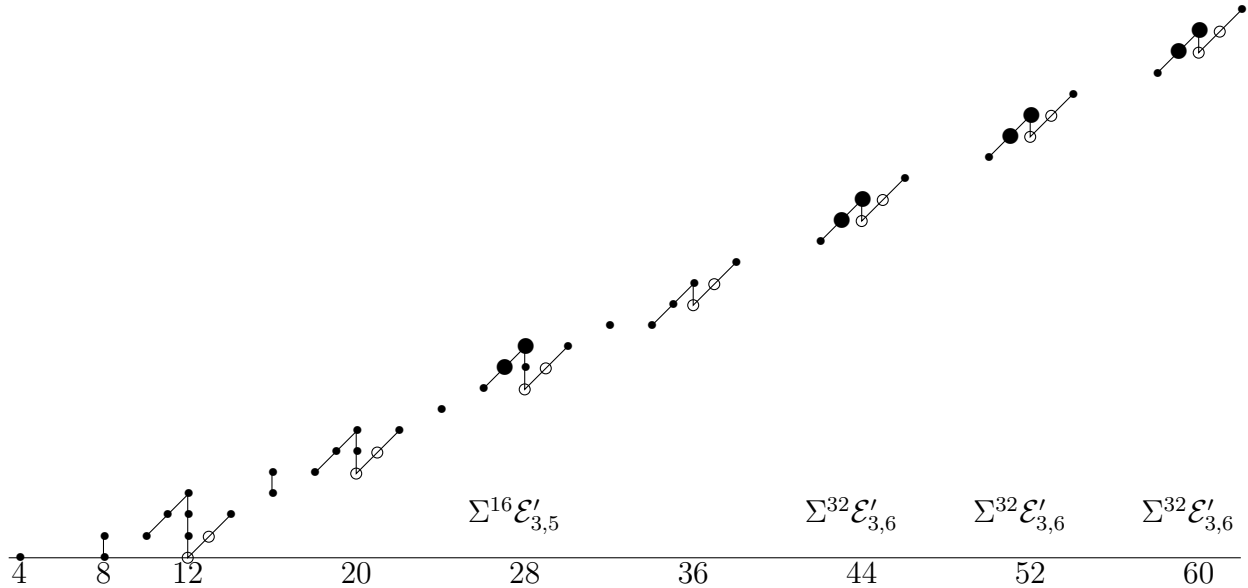
{\bf $\E'_{2,6}$ and differentials}

\begin{center}

\begin{\tz}[scale=.28]
\draw (3.5,0) -- (62,0);
\draw (8,0) -- (8,1);
\draw (10,1) -- (12,3) -- (12,0) -- (14,2);
\draw (16,4) -- (16,3);
\draw (18,4) -- (20,6) -- (20,4) -- (22,6);
\draw (26,8) -- (28,10) -- (28,8) -- (30,10);
\draw (34,11) -- (36,13) -- (36,12) -- (38,14);
\draw (42,15) -- (44,17) -- (44,16) -- (46,18);
\draw (50,19) -- (52,21) -- (52,20) -- (54,22);
\draw (58,23) -- (60,25) -- (60,24) -- (62,26);
\node at (4,0) {\llb};
\node at (8,0) {\llb};
\node at (8,1) {\llb};
\node at (10,1) {\llb};
\node at (11,2) {\llb};
\node at (12,3) {\llb};
\node at (12,2) {\llb};
\node at (12,1) {\llb};
\node at (14,2) {\llb};
\node at (16,3) {\llb};
\node at (16,4) {\llb};
\node at (18,4) {\llb};
\node at (19,5) {\llb};
\node at (20,6) {\llb};
\node at (20,5) {\llb};
\node at (22,6) {\llb};
\node at (26,8) {\llb};
\node at (28,9) {\llb};
\node at (30,10) {\llb};
\node at (32,11) {\llb};
\node at (34,11) {\llb};
\node at (35,12) {\llb};
\node at (36,13) {\llb};
\node at (38,14) {\llb};
\node at (42,15) {\llb};
\node at (50,19) {\llb};
\node at (54,22) {\llb};
\node at (46,18) {\llb};
\node at (58,23) {\llb};
\node at (62,26) {\llb};
\node at (4,-.7) {$4$};
\node at (8,-.7) {$8$};
\node at (12,-.7) {$12$};
\node at (20,-.7) {$20$};
\node at (28,-.7) {$28$};
\node at (36,-.7) {$36$};
\node at (44,-.7) {$44$};
\node at (52,-.7) {$52$};
\node at (60,-.7) {$60$};
\draw (12,0) circle (.26);
\draw (13,1) circle (.26);
\draw (20,4) circle (.26);
\draw (21,5) circle (.26);
\draw (28,8) circle (.26);
\draw (29,9) circle (.26);
\draw (36,12) circle (.26);
\draw (37,13) circle (.26);
\draw (44,16) circle (.26);
\draw (45,17) circle (.26);
\draw (52,20) circle (.26);
\draw (53,21) circle (.26);
\draw (60,24) circle (.26);
\draw (61,25) circle (.26);
\node at (27,9) {\blb};
\node at (28,10) {\blb};
\node at (43,16) {\blb};
\node at (44,17) {\blb};
\node at (51,20) {\blb};
\node at (52,21) {\blb};
\node at (59,24) {\blb};
\node at (60,25) {\blb};
\node at (27.5,2) {$\Si^{16}\E'_{3,5}$};
\node at (43.5,2) {$\Si^{32}\E'_{3,6}$};
\node at (51.5,2)  {$\Si^{32}\E'_{3,6}$};
\node at (59.5,2)  {$\Si^{32}\E'_{3,6}$};
\node at (24,7) {\llb};

\end{\tz}
\end{center}
\end{fig}
\end{minipage}
\bigskip

It is instructive to compare Figure \ref{2,6} with Figure \ref{A6}, in which $\Si^{64}\E_{2,6}$ appears prominently.

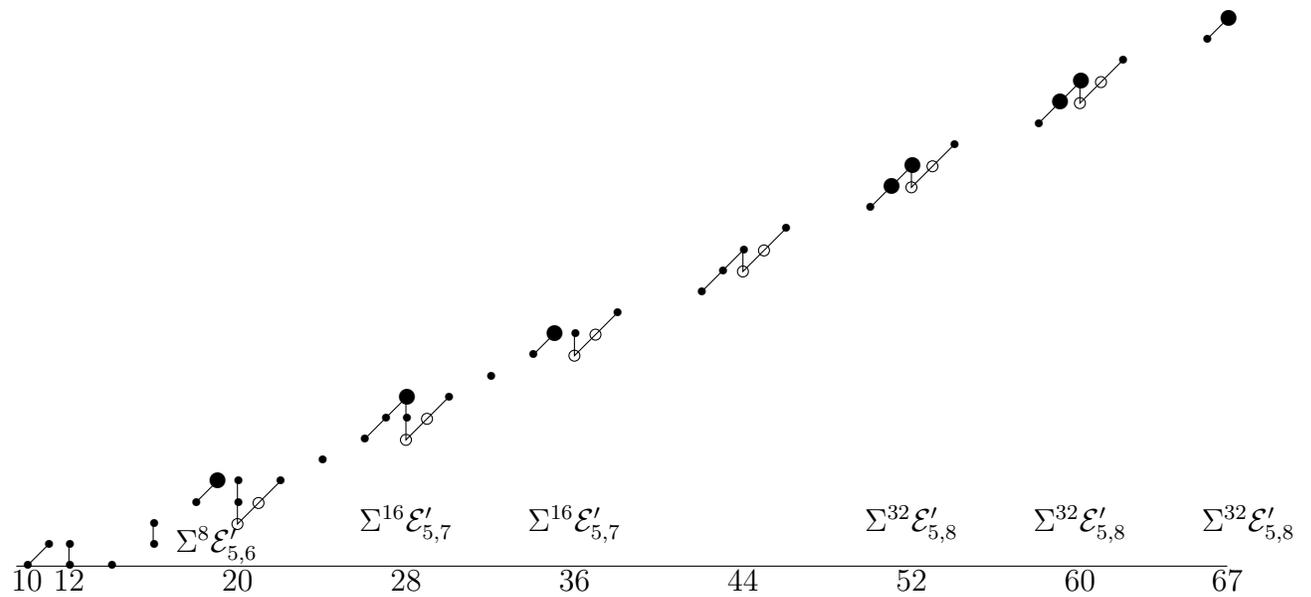
\begin{minipage}{6in}
\begin{fig}\label{3,7}

{\bf $\E'_{3,7}$ and differentials}

\begin{center}

\begin{\tz}[scale=.28]
\draw (7.5,0) -- (66,0);
\node at (8,-.7) {$8$};
\node at (12,-.7) {$12$};
\node at (20,-.7) {$20$};
\node at (28,-.7) {$28$};
\node at (36,-.7) {$36$};
\node at (44,-.7) {$44$};
\node at (52,-.7) {$52$};
\node at (60,-.7) {$60$};
\node at (66,-.7) {$66$};
\draw (11,1) -- (12,2) -- (12,0);
\draw (13,0) -- (14,1);
\draw (16,2) -- (16,3);
\draw (19,4) -- (20,5) -- (20,3) -- (22,5);
\draw (26,7) -- (28,9) -- (28,7) -- (30,9);
\draw (35,11) -- (36,12) -- (36,11) -- (38,13);
\draw (42,14) -- (44,16) -- (44,15) -- (46,17);
\draw (50,18) -- (52,20) -- (52,19) -- (54,21);
\draw (58,22) -- (60,24) -- (60,23) -- (62,25);
\node at (8,0) {\llb};
\node at (10,1) {\llb};
\node at (11,1) {\llb};
\node at (12,2) {\llb};
\node at (12,1) {\llb};
\node at (12,0) {\llb};
\node at (14,1) {\llb};
\node at (16,2) {\llb};
\node at (16,3) {\llb};
\node at (18,4) {\llb};
\node at (19,4) {\llb};
\node at (20,5) {\llb};
\node at (20,4) {\llb};
\node at (22,5) {\llb};
\node at (24,6) {\llb};
\node at (26,7) {\llb};
\node at (28,8) {\llb};
\node at (30,9) {\llb};
\node at (32,10) {\llb};
\node at (34,11) {\llb};
\node at (35,11) {\llb};
\node at (36,12) {\llb};
\node at (38,13) {\llb};
\node at (42,14) {\llb};
\node at (43,15) {\llb};
\node at (46,17) {\llb};
\node at (50,18) {\llb};
\node at (54,21) {\llb};
\node at (58,22) {\llb};
\node at (62,25) {\llb};
\node at (66,26) {\llb};
\draw (13,0) circle (.26);
\draw (20,3) circle (.26);
\draw (21,4) circle (.26);
\draw (28,7) circle (.26);
\draw (29,8) circle (.26);
\draw (36,11) circle (.26);
\draw (37,12) circle (.26);
\draw (44,15) circle (.26);
\draw (45,16) circle (.26);
\draw (52,19) circle (.26);
\draw (53,20) circle (.26);
\draw (60,23) circle (.26);
\draw (61,24) circle (.26);
\node at (27,8) {\blb};
\node at (28,9) {\blb};
\node at (44,16) {\blb};
\node at (51,19) {\blb};
\node at (52,20) {\blb};
\node at (59,23) {\blb};
\node at (60,24) {\blb};
\node at (27.5,2) {$\Si^{16}\E'_{4,6}$};
\node at (44,2) {$\Si^{32}\E'_{4,7}$};
\node at (51.5,2)  {$\Si^{32}\E'_{4,7}$};
\node at (59.5,2)  {$\Si^{32}\E'_{4,7}$};

\end{\tz}
\end{center}
\end{fig}
\end{minipage}

\bigskip

\begin{minipage}{6in}
\begin{fig}\label{4,8}

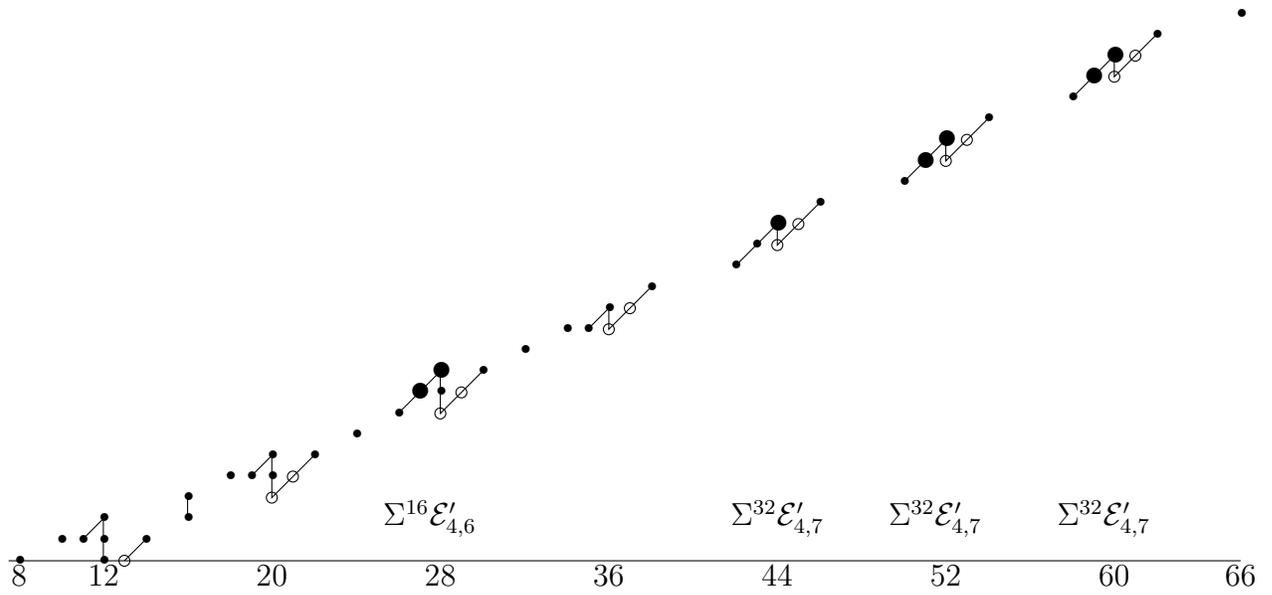
{\bf $\E'_{4,8}$ and differentials}

\begin{center}

\begin{\tz}[scale=.28]
\draw (9.5,0) -- (67,0);
\draw (10,0) -- (11,1);
\draw (12,0) -- (12,1);
\draw (16,1) -- (16,2);
\draw (18,3) -- (19,4);
\draw (20,4) -- (20,2) -- (22,4);
\draw (26,6) -- (28,8) -- (28,6) -- (30,8);
\draw (34,10) -- (35,11);
\draw (36,11) -- (36,10) -- (38,12);
\draw (42,13) -- (44,15) -- (44,14) -- (46,16);
\draw (50,17) -- (52,19) -- (52,18) -- (54,20);
\draw (58,21) -- (60,23) -- (60,22) -- (62,24);
\draw (66,25) -- (67,26);
\node at (10,-.7) {$10$};
\node at (12,-.7) {$12$};
\node at (20,-.7) {$20$};
\node at (28,-.7) {$28$};
\node at (36,-.7) {$36$};
\node at (44,-.7) {$44$};
\node at (52,-.7) {$52$};
\node at (60,-.7) {$60$};
\node at (67,-.7) {$67$};
\node at (10,0) {\llb};
\node at (11,1) {\llb};
\node at (12,0) {\llb};
\node at (12,1) {\llb};
\node at (14,0) {\llb};
\node at (16,1) {\llb};
\node at (16,2) {\llb};
\node at (18,3) {\llb};
\node at (20,3) {\llb};
\node at (20,4) {\llb};
\node at (22,4) {\llb};
\node at (24,5) {\llb};
\node at (26,6) {\llb};
\node at (27,7) {\llb};
\node at (28,7) {\llb};
\node at (30,8) {\llb};
\node at (32,9) {\llb};
\node at (34,10) {\llb};
\node at (36,11) {\llb};
\node at (38,12) {\llb};
\node at (42,13) {\llb};
\node at (43,14) {\llb};
\node at (44,15) {\llb};
\node at (46,16) {\llb};
\node at (50,17) {\llb};
\node at (54,20) {\llb};
\node at (58,21) {\llb};
\node at (62,24) {\llb};
\node at (66,25) {\llb};
\draw (20,2) circle (.26);
\draw (21,3) circle (.26);
\draw (28,6) circle (.26);
\draw (29,7) circle (.26);
\draw (36,10) circle (.26);
\draw (37,11) circle (.26);
\draw (44,14) circle (.26);
\draw (45,15) circle (.26);
\draw (52,18) circle (.26);
\draw (53,19) circle (.26);
\draw (60,22) circle (.26);
\draw (61,23) circle (.26);
\node at (19,4) {\blb};
\node at (28,8) {\blb};
\node at (35,11) {\blb};
\node at (51,18) {\blb};
\node at (52,19) {\blb};
\node at (59,22) {\blb};
\node at (60,23) {\blb};
\node at (67,26) {\blb};
\node at (19,1) {$\Si^8\E'_{5,6}$};
\node at (28,2) {$\Si^{16}\E'_{5,7}$};
\node at (36,2) {$\Si^{16}\E'_{5,7}$};
\node at (52,2) {$\Si^{32}\E'_{5,8}$};
\node at (60,2) {$\Si^{32}\E'_{5,8}$};
\node at (68,2) {$\Si^{32}\E'_{5,8}$};
\end{\tz}
\end{center}
\end{fig}
\end{minipage}

\bigskip

\begin{minipage}{6in}
\begin{fig}\label{5,9}

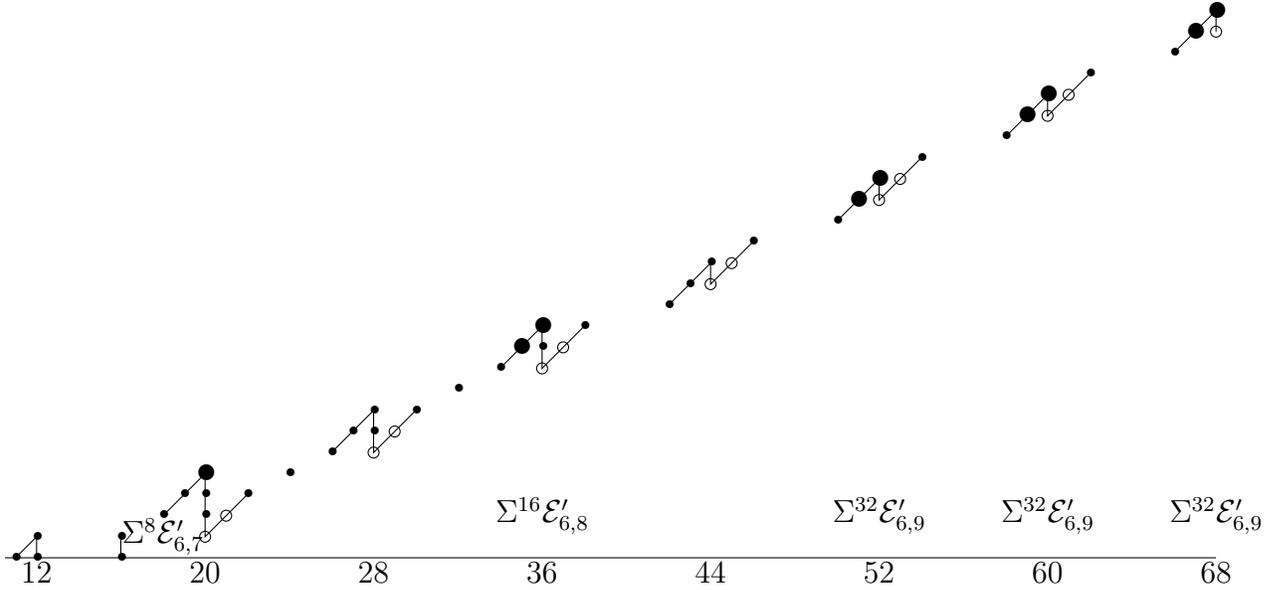
{\bf $\E'_{5,9}$ and differentials}

\begin{center}

\begin{\tz}[scale=.28]
\draw (10.5,0) -- (68,0);
\draw (11,0) -- (12,1) -- (12,0);
\draw (16,0) -- (16,1);
\draw (18,2) -- (20,4) -- (20,1) -- (22,3);
\draw (26,5) -- (28,7) -- (28,5) -- (30,7);
\draw (34,9) -- (36,11) -- (36,9) -- (38,11);
\draw (42,12) -- (44,14) -- (44,13) -- (46,15);
\draw (50,16) -- (52,18) -- (52,17) -- (54,19);
\draw (58,20) -- (60,22) -- (60,21) -- (62,23);
\draw (66,24) -- (68,26) -- (68,25);
\node at (12,-.7) {$12$};
\node at (20,-.7) {$20$};
\node at (28,-.7) {$28$};
\node at (36,-.7) {$36$};
\node at (44,-.7) {$44$};
\node at (52,-.7) {$52$};
\node at (60,-.7) {$60$};
\node at (68,-.7) {$68$};
\node at (11,0) {\llb};
\node at (12,1) {\llb};
\node at (12,0) {\llb};
\node at (16,0) {\llb};
\node at (16,1) {\llb};
\node at (18,2) {\llb};
\node at (19,3) {\llb};
\node at (20,2) {\llb};
\node at (20,3) {\llb};
\node at (24,4) {\llb};
\node at (26,5) {\llb};
\node at (27,6) {\llb};
\node at (28,7) {\llb};
\node at (28,6) {\llb};
\node at (30,7) {\llb};
\node at (32,8) {\llb};
\node at (34,9) {\llb};
\node at (36,10) {\llb};
\node at (38,11) {\llb};
\node at (42,12) {\llb};
\node at (43,13) {\llb};
\node at (44,14) {\llb};
\node at (46,15) {\llb};
\node at (50,16) {\llb};
\node at (54,19) {\llb};
\node at (58,20) {\llb};
\node at (62,23) {\llb};
\node at (66,24) {\llb};
\node at (20,4) {\blb};
\node at (35,10) {\blb};
\node at (36,11) {\blb};
\node at (51,17) {\blb};
\node at (52,18) {\blb};
\node at (59,21) {\blb};
\node at (60,22) {\blb};
\node at (67,25) {\blb};
\node at (68,26) {\blb};
\draw (20,1) circle (.26);
\draw (21,2) circle (.26);
\draw (28,5) circle (.26);
\draw (29,6) circle (.26);
\draw (36,9) circle (.26);
\draw (37,10) circle (.26);
\draw (44,13) circle (.26);
\draw (45,14) circle (.26);
\draw (52,17) circle (.26);
\draw (53,18) circle (.26);
\draw (60,21) circle (.26);
\draw (61,22) circle (.26);
\draw (68,25) circle (.26);
\node at (18,1) {$\Si^8\E'_{6,7}$};
\node at (36,2) {$\Si^{16}\E'_{6,8}$};
\node at (52,2) {$\Si^{32}\E'_{6,9}$};
\node at (60,2) {$\Si^{32}\E'_{6,9}$};
\node at (68,2) {$\Si^{32}\E'_{6,9}$};
\node at (22,3) {\llb};

\end{\tz}
\end{center}
\end{fig}
\end{minipage}
\bigskip

The sources of the differentials are very regular. There are differentials from $\Si^{2^{d+2}}\E'_{e,e+d}$ into $\E'_{e-1,\ell}$ for each $\ell\ge e+d$ on all  classes in grading 4 and 5 mod 8 in $\Si^{2^{d+2}+1}M^e_{d+3}$ (except for a class $x$ in filtration 0 and grading 4 mod 8 satisfying $\eta x=0$) until the $\E'_{e,e+d}$ is ended by the differentials from $\Si^{2^{d+2}+2}\Mh_4^{e-2}$.

Careful study of the above charts $\E'_{e,e+d}$ for $d=1$, 2, 3, and 4 can give great insight into the form of $\E'_{e,e+d}$ for arbitrary $d$. As $d$ increases, the upper edge of the chart stays fixed while the lower edge drops by 1 each time. Only $e=2$, 3, 4, and 5 need be considered, since $\E_{e+4,k+4}=\Si^8\E_{e,k}$.

The slightly differing forms of the upper edge of $\E'_{e,e+d}$ depending on the mod 4 value of $e$ are caused by the differing ways that $M_4^{e-1}$ begins. In Figure \ref{E35}, we show the formation of $\E'_{3,5}$. Comparison with Figure \ref{E'expl} is instructive.

\bigskip

\begin{minipage}{6in}
\begin{fig}\label{E35}

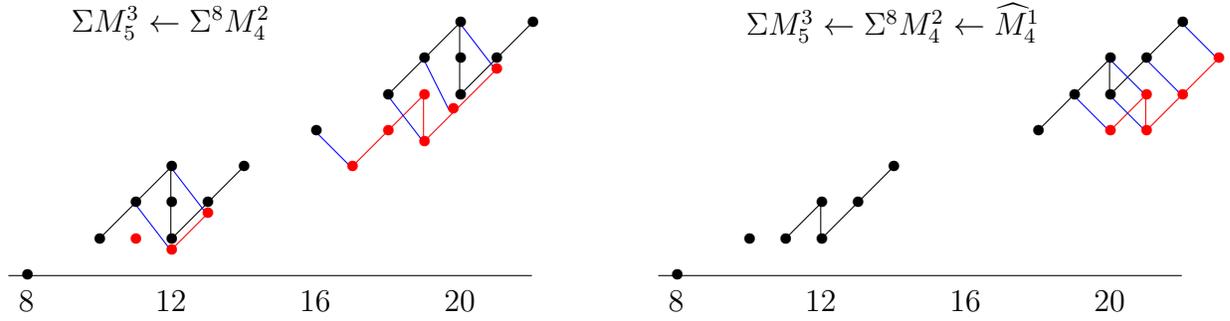
{\bf Formation of $\E'_{3,5}$}

\begin{center}

\begin{\tz}[scale=.48]
\draw (7.5,0) -- (22,0);
\draw (25.5,0) -- (40,0);
\draw (10,1) -- (12,3) -- (12,1) -- (14,3);
\draw (18,5) -- (20,7) -- (20,5) -- (22,7);
\draw (29,1) -- (30,2) -- (30,1) -- (32,3);
\draw (36,4) -- (38,6) -- (38,5) -- (40,7);
\draw [red] (12,.7) -- (13,1.7);
\draw [red] (17,3) -- (19,5) -- (19,3.7) -- (21,5.7);
\draw [red] (38,4) -- (39,5) -- (39,4) -- (41,6);
\draw [blue] (12,.7) -- (11,2);
\draw [blue] (13,1.7) -- (12,3);
\draw [blue] (17,3) -- (16,4);
\draw [blue] (19,3.7) -- (18,5);
\draw [blue] (19.7,4.6) -- (19,6);
\draw [blue] (21,5.7) -- (20,7);
\draw [blue] (38,4) -- (37,5);
\draw [blue] (39,5) -- (38,6);
\draw [blue] (39,4) -- (38,5);
\draw [blue] (40,5) -- (39,6);
\draw [blue] (41,6) -- (40,7);
\node at (8,0) {\lb};
\node at (10,1) {\lb};
\node at (11,2) {\lb};
\node at (12,3) {\lb};
\node at (12,2) {\lb};
\node at (12,1) {\lb};
\node at (13,2) {\lb};
\node at (14,3) {\lb};
\node at (16,4) {\lb};
\node at (18,5) {\lb};
\node at (19,6) {\lb};
\node at (20,7) {\lb};
\node at (20,6) {\lb};
\node at (20,5) {\lb};
\node at (21,6) {\lb};
\node at (22,7) {\lb};
\node at (26,0) {\lb};
\node at (28,1) {\lb};
\node at (29,1) {\lb};
\node at (30,2) {\lb};
\node at (30,1) {\lb};
\node at (31,2) {\lb};
\node at (32,3) {\lb};
\node at (36,4) {\lb};
\node at (37,5) {\lb};
\node at (38,6) {\lb};
\node at (38,5) {\lb};
\node at (39,6) {\lb};
\node at (40,7) {\lb};
\node [red] at (11,1) {\lb};
\node [red] at (12,.7) {\lb};
\node [red] at (13,1.7) {\lb};
\node [red] at (17,3) {\lb};
\node [red] at (18,4) {\lb};
\node [red] at (19,5) {\lb};
\node [red] at (19,3.7) {\lb};
\node [red] at (19.8,4.6) {\lb};
\node [red] at (21,5.7) {\lb};
\node [red] at (38,4) {\lb};
\node [red] at (39,5) {\lb};
\node [red] at (39,4) {\lb};
\node [red] at (40,5) {\lb};
\node [red] at (41,6) {\lb};
\node at (12,7) {$\Si M_5^3\lar\Si^8M_4^2$};
\node at (32,7) {$\Si M_5^3\lar\Si^8M_4^2\lar \Mh_4^1$};
\node at (8,-.7) {$8$};
\node at (12,-.7) {$12$};
\node at (16,-.7) {$16$};
\node at (20,-.7) {$20$};
\node at (26,-.7) {$8$};
\node at (30,-.7) {$12$};
\node at (34,-.7) {$16$};
\node at (38,-.7) {$20$};

\end{\tz}
\end{center}
\end{fig}
\end{minipage}
\bigskip

It is interesting to see how the ending of each $\E'_{e,e+d}$ occurs. Prior to the $\Si^{2^{d+2}+2}\Mh_4^{e-2}$ at the end of the sequence in Definition \ref{E'def}, the sequence will have stabilized to a sequence of lightning flashes with initial classes in grading 2 mod 8. In Figure \ref{ends}, we show how the lightning flash beginning in grading $2^{d+2}+2$ is hit by $\Si^{2^{d+2}+2}\Mh_4^{e-2}$ for $2\le e\le5$. Increasing $e$ by $4$ would push this behavior out by $8$ gradings. In Figure \ref{ends}, the remaining classes are indicated with bigger dots. Compare with the endings in Figures \ref{four}, \ref{2,6}, \ref{3,7}, \ref{4,8}, and \ref{5,9}.

\bigskip

\begin{minipage}{6in}
\begin{fig}\label{ends}

{\bf Termination of $\E'_{e,e+d}$}

\begin{center}

\begin{\tz}[scale=.45]
\draw (0,1) -- (2,3) -- (2,2) -- (4,4);
\draw (8,1) -- (10,3) -- (10,2) -- (12,4);
\draw (16,1) -- (18,3) -- (18,2) -- (20,4);
\draw (24,1) -- (26,3) -- (26,2) -- (28,4);
\draw [red] (1,0) -- (3,2) -- (3,1) -- (5,3);
\draw [red] (10,1) -- (11,2) -- (11,1) -- (13,3);
\draw [red] (19,2) -- (19,1) -- (21,3);
\draw [red] (28,2) -- (29,3);
\draw [blue] (1,0) -- (0,1);
\draw [blue] (2,1) -- (1,2);
\draw [blue] (3,2) -- (2,3);
\draw [blue] (3,1) -- (2,2);
\draw [blue] (4,2) -- (3,3);
\draw [blue] (5,3) -- (4,4);
\draw [blue] (10,1) -- (9,2);
\draw [blue] (11,2) -- (10,3);
\draw [blue] (11,1) -- (10,2);
\draw [blue] (12,2) -- (11,3);
\draw [blue] (13,3) -- (12,4);
\draw [blue] (19,1) -- (18,2);
\draw [blue] (19,2) -- (18,3);
\draw [blue] (20,2) -- (19,3);
\draw [blue] (21,3) -- (20,4);
\draw [blue] (28,2) -- (27,3);
\draw [blue] (29,3) -- (28,4);
\node at (2.5,-1) {$e=2$};
\node at (10.5,-1) {$e=3$};
\node at (18.5,-1) {$e=4$};
\node at (26.5,-1) {$e=5$};
\node at (0,1) {\lb};
\node at (1,2) {\lb};
\node at (2,3) {\lb};
\node at (2,2) {\lb};
\node at (3,3) {\lb};
\node at (4,4) {\lb};
\node at (9,2) {\lb};
\node at (10,3) {\lb};
\node at (10,2) {\lb};
\node at (11,3) {\lb};
\node at (12,4) {\lb};
\node at (18,2) {\lb};
\node at (18,3) {\lb};
\node at (19,3) {\lb};
\node at (20,4) {\lb};
\node at (27,3) {\lb};
\node at (28,4) {\lb};
\node [red] at (1,0) {\lb};
\node [red] at (2,1) {\lb};
\node [red] at (3,2) {\lb};
\node [red] at (3,1) {\lb};
\node [red] at (4,2) {\lb};
\node [red] at (5,3) {\lb};
\node [red] at (10,1) {\lb};
\node [red] at (11,2) {\lb};
\node [red] at (11,1) {\lb};
\node [red] at (12,2) {\lb};
\node [red] at (13,3) {\lb};
\node [red] at (19,1) {\lb};
\node [red] at (19,2) {\lb};
\node [red] at (20,2) {\lb};
\node [red] at (21,3) {\lb};
\node [red] at (29,3) {\lb};
\node [red] at (28,2) {\lb};
\node at (8,1) {\blb};
\node at (16,1) {\blb};
\node at (17,2) {\blb};
\node at (24,1) {\blb};
\node at (25,2) {\blb};
\node at (26,3) {\blb};
\node at (26,2) {\blb};

\end{\tz}
\end{center}
\end{fig}
\end{minipage}
\bigskip

\begin{defin}\label{Phidef} If $C$ is a chart, then $\Phi^iC$ is the chart obtained from $C$ by increasing filtration of all elements by $i$.\end{defin}

The upper edge $\E_{1,k}$ has a different form. The following definition will be illustrated in Figure \ref{E14} and justified at the end of Section \ref{derivationsec}. 

\begin{defin}\label{E1'def} For $k\ge2$, $V_k$ is a chart with, for $i\ge0$, classes $g_{8i}$ in position $(8i,4i)$, and $g_{8i+4}$ in position $(8i+4,4i+3)$, of order $2^{k-1}$ except that the order of $g_0$ is $2^{k+1}$. There are also elements $\eta g_{8i}$ and $\eta^2 g_{8i}$, $x_3$ such that $\eta x_3=2^{k-2}g_4$, and, for $i>0$, elements $x_{8i+2}$ and $\eta x_{8i+2}$ such that $\eta^2x_{8i+2}=g_{8i+4}$.  If $k=2$ and $i>0$, then $2x_{8i+2}=\eta^2g_{8i}$. In grading $\ge8$, $V_k$ agrees with $\Phi^4ko_*(M(2^{k-1}))$, where $M(n)$ is the mod $n$ Moore spectrum.

The chart $\E'_{1,k}$ is formed from the sequence
\begin{equation}\label{E1seq}V_k\lar \Si^8M_4^0\lar\cdots\lar\Si^{2^i}M_4^0\lar\cdots\lar \Si^{2^{k+1}}M_4^0.\end{equation}
Working from left to right, each $\Si^{2^i}M_4^0$ is placed so that there are $d_1$ differentials from its generators in grading $1$, $3$, $4$, and $5$ mod $8$ into the upper edge of the chart resulting from the preceding steps. For $2\le i\le k-1$, $\E'_{1,k}$ agrees with $\Phi^{2^{i-1}}ko_*(M(2^{k-i}))$ in grading $2^{i+1}$ through $2^{i+2}-1$.

The edge $\E_{1,k}$ is obtained from $\E'_{1,k}$ by removing the classes in grading $3$ or $4$ mod $8$ which are  hit by differentials from $\Si^{2^i}\E'_{2,i}$ for $3\le i\le k$.
\end{defin}

In Figure \ref{E14}, we show the formation of $\E'_{1,4}$ and $\E_{1,4}$. At each step, the chart being adjoined is shown in red. Then $\E'_{1,4}$ consists of all classes in the bottom chart of grading $\le28$. The classes with the big dots are hit by differentials, and so do not appear in $\E_{1,4}$. The chart hitting them is indicated below them. The classes supporting those differentials can be seen in Figure \ref{four}. The result of Figure \ref{E14} can be seen in Figure \ref{A4}.

\bigskip

\begin{minipage}{6in}
\begin{fig}\label{E14}

{\bf Formation of $\E'_{1,4}$ and $\E_{1,4}$}

\begin{center}

\begin{\tz}[scale=.44]
\draw (0,0) -- (36,0);
\draw (0,4) -- (0,0) -- (2,2);
\draw (3,4) -- (4,5) -- (4,3);
\draw (8,5) -- (8,4) -- (10,6);
\draw (10,6.3) -- (12,8) -- (12,7);
\draw (16,8) -- (18,10) -- (18,9) -- (20,11);
\draw (24,12) -- (26,14) -- (26,13) -- (28,15);
\draw (32,16) -- (34,18) -- (34,17) -- (36,19);
\draw [red] (33,15) -- (35,17) -- (35,16) -- (37,18);
\node at (0,-.7) {$0$};
\node at (4,-.7) {$4$};
\node at (8,-.7) {$8$};
\node at (12,-.7) {$12$};
\node at (16,-.7) {$16$};
\node at (20,-.7) {$20$};
\node at (24,-.7) {$24$};
\node at (28,-.7) {$28$};
\node at (32,-.7) {$32$};
\draw (0,21) -- (13,21);
\draw (16,21) -- (37,21);
\draw (0,25) -- (0,21) -- (2,23);
\draw (3,25) -- (4,26) -- (4,24);
\draw (8,27) -- (8,25) -- (10,27);
\draw (10,28) -- (12,30) -- (12,28);
\draw [red] (9,26.3) -- (11,28) -- (11,27) -- (13,28.6);
\draw [blue] (9,26.3) -- (8,27);
\draw [blue] (13,28.6) -- (12,30);
\draw [blue] (11,27) -- (10,28);
\draw [blue] (11.8,27.64) -- (11,29);
\draw (16,25) -- (16,21) -- (18,23);
\draw (19,25) -- (20,26) -- (20,24);
\draw (24,26) -- (24,25) -- (26,27);
\draw (26,27.3) -- (28,29) -- (28,28);
\draw (32,30) -- (32,29) -- (34,31);
\draw (34.1,31.4) -- (36,33) -- (36,32);
\draw [red] (33,29) -- (35,31) -- (35,30) -- (37,32);
\draw [blue] (33,29) -- (32,30);
\draw [blue] (35,30) -- (34.1,31.4);
\draw [blue] (35.05,32.2) -- (36,31);
\draw [blue] (37,32) -- (36,33);
\draw [blue] (33,15) -- (32,16);
\draw [blue] (34,16) -- (33,17);
\draw [blue] (35,16) -- (34,17);
\draw [blue] (35,17) -- (34,18);
\draw [blue] (36,17) -- (35,18);
\draw [blue] (37,18) -- (36,19);
\node at (4,29) {$V_4\lar\Si^8M_4^0$};
\node at (26,31) {$V_4\lar\Si^8M_4^0\lar\Si^{16}M_4^0$};
\node at (12,8) {\blb};
\node at (12,7) {\llb};
\node at (11,7.15) {\llb};
\node at (0,20.3) {$0$};
\node at (4,20.3) {$4$};
\node at (8,20.3) {$8$};
\node at (12,20.3) {$12$};
\node at (16,20.3) {$0$};
\node at (20,20.3) {$4$};
\node at (24,20.3) {$8$};
\node at (28,20.3) {$12$};
\node at (32,20.3) {$16$};
\node at (36,20.3) {$20$};
\node at (14,30) {$\iddots$};
\node at (38,33.5) {$\iddots$};
\node at (0,0) {\llb};
\node at (0,1) {\llb};
\node at (0,2) {\llb};
\node at (0,3) {\llb};
\node at (0,4) {\llb};
\node at (1,1) {\llb};
\node at (2,2) {\llb};
\node at (4,3) {\llb};
\node at (4,4) {\llb};
\node at (4,5) {\llb};
\node at (3,4) {\llb};
\node at (8,4) {\llb};
\node at (8,5) {\llb};
\node at (9,5) {\llb};
\node at (10,6) {\llb};
\node at (10,6.3) {\llb};
\node at (12,7) {\llb};
\node at (12,8) {\blb};
\node at (11,7.15) {\llb};
\node at (16,8) {\llb};
\node at (17,9) {\llb};
\node at (18,10) {\llb};
\node at (18,9) {\llb};
\node at (19,10) {\llb};
\node at (20,11) {\llb};
\node at (24,12) {\llb};
\node at (25,13) {\llb};
\node at (26,14) {\llb};
\node at (26,13) {\llb};
\node at (27,14) {\blb};
\node at (28,15) {\blb};
\node at (32,16) {\llb};
\node at (33,17) {\llb};
\node at (34,18) {\llb};
\node at (34,17) {\llb};
\node at (35,18) {\llb};
\node at (36,19) {\llb};
\node at (0,21) {\llb};
\node at (0,22) {\llb};
\node at (0,23) {\llb};
\node at (0,24) {\llb};
\node at (0,25) {\llb};
\node at (1,22) {\llb};
\node at (2,23) {\llb};
\node at (4,24) {\llb};
\node at (4,25) {\llb};
\node at (4,26) {\llb};
\node at (3,25) {\llb};
\node at (8,25) {\llb};
\node at (8,26) {\llb};
\node at (8,27) {\llb};
\node at (9,26) {\llb};
\node at (10,27) {\llb};
\node at (10,28) {\llb};
\node at (11,29) {\llb};
\node at (12,30) {\llb};
\node at (12,29) {\llb};
\node at (12,28) {\llb};
\node at (16,21) {\llb};
\node at (16,22) {\llb};
\node at (16,23) {\llb};
\node at (16,24) {\llb};
\node at (16,25) {\llb};
\node at (17,22) {\llb};
\node at (18,23) {\llb};
\node at (20,24) {\llb};
\node at (20,25) {\llb};
\node at (20,26) {\llb};
\node at (19,25) {\llb};
\node at (24,25) {\llb};
\node at (24,26) {\llb};
\node at (25,26) {\llb};
\node at (26,27) {\llb};
\node at (26,27.3) {\llb};
\node at (27,28.15) {\llb};
\node at (28,29) {\llb};
\node at (28,28) {\llb};
\node at (32,29) {\llb};
\node at (32,30) {\llb};
\node at (33,30) {\llb};
\node at (34,31) {\llb};
\node at (34,31.3) {\llb};
\node at (35,32.15) {\llb};
\node at (36,33) {\llb};
\node at (36,32) {\llb};
\node at (12,2) {$\Si^8\E'_{2,3}$};
\node at (27,2) {$\Si^{16}\E'_{2,4}$};
\node [red] at (33,15) {\llb};
\node [red] at (34,16) {\llb};
\node [red] at (35,17) {\llb};
\node [red] at (35,16) {\llb};
\node [red] at (36,17) {\llb};
\node [red] at (37,18) {\llb};
\node [red] at (33,29) {\llb};
\node [red] at (34,30) {\llb};
\node [red] at (35,31) {\llb};
\node [red] at (35,30) {\llb};
\node [red] at (36,31) {\llb};
\node [red] at (37,32) {\llb};
\node [red] at (9,26.3) {\llb};
\node [red] at (10,27.15) {\llb};
\node [red] at (11,28) {\llb};
\node [red] at (11,27) {\llb};
\node [red] at (11.8,27.64) {\llb};
\node [red] at (13,28.6) {\llb};

\node at (16,13) {$V_4\lar\Si^8M_4^0\lar\Si^{16}M_4^0\lar\Si^{32}M_4^0$};
\end{\tz}
\end{center}
\end{fig}
\end{minipage}
\bigskip

The following theorem about the exotic extensions will be proved in Section \ref{pfsec}.
\begin{thm}\label{extnthm} The only exotic extensions in $\At_k$ are as follows:
\begin{enumerate}
\item Into $\E_{1,k}$ from $\Si^{2^\ell}\E_{2,\ell}$ for $4\le\ell\le k$ in grading $8i+2$ for $3\cdot2^{\ell-4}\le i\le2^{\ell-2}-1$. The target is the class which is not in $\im(\eta)$.
\item For $e\ge2$ and $e+1<\ell\le k$, into $\Si^D\E_{e,\ell}$ from $\Si^{2^{\ell+1-e}+D}\E_{e+1,\ell}$ in grading $6$ mod $8$ throughout the entire extent of $\Si^{2^{\ell+1-e}+D}\E_{e+1,\ell}$.
\item For $e\ge2$, $\ell\le k$, and $e+1< m\le\ell$, into $\Si^D\E_{e,\ell}$ from $\Si^{2^{m+1-e}+D}\E_{e+1,m}$ in grading $2$ mod $8$ throughout the second half of  $\E_{e+1,m}$ (the range of the lightning flashes in $\E'_{e+1,m}$, including the final element in $2$ mod $8$ if $e\equiv0$  or $3$ mod $4$, but not if $e\equiv2$ mod $4$).
\end{enumerate}
\end{thm}

Since $\E_{e,\ell}$ has order $\le 2$ in grading 2 mod 4, we don't have to specify which element is involved in the extension in parts (2) and (3) of the theorem. In Figure \ref{A6}, the extensions into the upper edge are easily checked to agree with part (1). In (\ref{extnlist}), we list the other extensions in $\At_6$, with those of type (2) in the left column and those of type (3) in the right column. We list the grading followed by the edges involved, with dashes indicating the extension.

\def\ddd{\text{\,---\,}}
\begin{align} 
\nonumber30:\,\Si^{16}(\Si^8\E_{3,4}\ddd\E_{2,4})&\qquad66:\,\Si^{48}(\Si^8\E_{4,5}\ddd\E_{3,5})\\
\nonumber62:\,\Si^{32}(\Si^{16}\E_{3,5}\ddd\E_{2,5})&\qquad98:\,\Si^{80}(\Si^8\E_{4,5}\ddd\E_{3,5})\\
\nonumber110:\,\Si^{64}(\Si^{32}\E_{3,6}\ddd\E_{2,6})&\qquad114:\,\Si^{96}(\Si^8\E_{4,5}\ddd\E_{3,6})\\
\nonumber118:\,\Si^{64}(\Si^{32}\E_{3,6}\ddd\E_{2,6})&\qquad122:\,\Si^{64}(\Si^{32}\E_{3,6}\ddd\E_{2,6})\\
\nonumber126:\,\Si^{96}(\Si^{16}\E_{4,6}\ddd\E_{3,6})&\qquad122:\,\Si^{96}(\Si^{16}\E_{4,6}\ddd\E_{3,6})\\
\label{extnlist}126:\,\Si^{64}(\Si^{32}\E_{3,6}\ddd\E_{2,6})&\qquad130:\,\Si^{96}(\Si^{16}\E_{4,6}\ddd\E_{3,6})
\end{align}

\section{$z^i\B_{k,\ell}$}\label{Bklsec}
In this section, we describe the summands $z^i\B_{k,\ell}$, which appeared in Theorem \ref{main}. The proof will be in Section \ref{Bklpfsec}.

The description of $ku^*(K_2)$ in \cite{DW} was in terms of summands $A_k$, $B_k$, and $S_{k,\ell}$. It was stated in \cite[Section 7]{DW} that $S_{k,\ell}$ and two specific copies of $B_k$ 
combine together nicely, in the sense that differentials in the ASS that form them involve just the three\footnote{This will be proved in Theorem \ref{closed}.}, and that the contribution to $ku^*(K_2)$ of the three is isomorphic to the corresponding contribution to $ku_*(K_2)$ after dualizing the gradings, implying self-duality of $B_{k,\ell}$. The example there yielded the chart for $B_{3,4}\subset ku_*(K_2)$ in Figure \ref{kuB34}.

\bigskip
\begin{minipage}{6in}
\begin{fig}\label{kuB34}

{\bf $B_{3,4}\subset ku_*(K_2)$}

\begin{center}

\begin{\tz}[scale=.5]
\draw (-.5,0) -- (28.5,0);
\draw (0,0) -- (2,1) -- (2,0) -- (10,4);
\draw [dashed] (10,4) -- (10,0);
\draw (10,0) -- (12,1);
\draw (11,0) -- (17,3);
\draw (16,0) -- (18,1) -- (18,0) -- (26,4);
\draw [dashed] (26,4) -- (26,0);
\draw (26,0) -- (28,1);
\node at (0,0) {\lb};
\node at (2,1) {\lb};
\node at (2,0) {\lb};
\node at (4,1) {\lb};
\node at (6,2) {\lb};
\node at (8,3) {\lb};
\node at (10,4) {\lb};
\node at (10,0) {\lb};
\node at (12,1) {\lb};
\node at (11,0) {\lb};
\node at (13,1) {\lb};
\node at (15,2) {\lb};
\node at (17,3) {\lb};
\node at (16,0) {\lb};
\node at (18,1) {\lb};
\node at (18,0) {\lb};
\node at (20,1) {\lb};
\node at (22,2) {\lb};
\node at (24,3) {\lb};
\node at (26,4) {\lb};
\node at (26,0) {\lb};
\node at (28,1) {\lb};
\node at (0,-.7) {$70$};
\node at (10,-.7) {$80$};
\node at (11,-.7) {$81$};
\node at (16,-.7) {$86$};
\node at (26,-.7) {$96$};

\end{\tz}
\end{center}
\end{fig}
\end{minipage}
\bigskip

In Section \ref{Bklpfsec}, we explain how the $ko$ analogues of $B_{k,\ell}$ are defined and explain how the ASS of this part is computed. We denote by $\B_{k,\ell}$ the resulting chart for this summand of $ko_*(K_2)$. The result for $\B_{3,4}$ is shown in Figure \ref{B34}.

\bigskip
\begin{minipage}{6in}
\begin{fig}\label{B34}

{\bf $\B_{3,4}$}

\begin{center}

\begin{\tz}[scale=.4]
\draw (3.5,0) -- (34.5,0);
\draw (10,1) -- (12,3) -- (12,2) -- (14,4);
\draw [dashed] (14,0) -- (14,4);
\draw (17,4) -- (19,6) -- (19,5) -- (21,7);
\draw (25,8) -- (26,9);
\draw [dashed] (26,1) -- (26,9);
\draw  (26,1) -- (28,3);
\node at (4,0) {\lb};
\node at (6,0) {\lb};
\node at (8,1) {\lb};
\node at (10,1) {\lb};
\node at (11,2) {\lb};
\node at (12,3) {\lb};
\node at (12,2) {\lb};
\node at (13,3) {\lb};
\node at (14,4) {\lb};
\node at (14,0) {\lb};
\node at (18,1) {\lb};
\node at (20,0) {\lb};
\node at (17,4) {\lb};
\node at (18,5) {\lb};
\node at (19,6) {\lb};
\node at (19,5) {\lb};
\node at (20,6) {\lb};
\node at (21,7) {\lb};
\node at (24,0) {\lb};
\node at (24,1) {\lb};
\node at (30,0) {\lb};
\node at (34,1) {\lb};
\node at (32,4) {\lb};
\node at (26,1) {\lb};
\node at (27,2) {\lb};
\node at (28,3) {\lb};
\node at (26,9) {\lb};
\node at (25,8) {\lb};
\draw [red] [dotted] (24,0) -- (25,8);
\node at (4,-.6) {$68$};
\node at (8,-.6) {$72$};
\node at (12,-.6) {$76$};
\node at (16,-.6) {$80$};
\node at (20,-.6) {$84$};
\node at (24,-.6) {$88$};
\node at (28,-.6) {$92$};
\node at (32,-.6) {$96$};
\draw (24,0) -- (24,1);
\end{\tz}
\end{center}
\end{fig}
\end{minipage}
\bigskip

The lightning flash in the middle is the analogue of the $S_{3,4}$ part, which was the short $v$-tower in grading 81 to 87 in Figure \ref{kuB34}. We choose to raise the filtration of this for two reasons. (a) When we look at the differentials in the ASS it is convenient to increase some filtrations to make the extensions clearer and the pictures nicer. (b) The classes in 84 and 85 are $v_1^4$ times the classes in 76 and 77, and we like to have our chart depict this.

In the $ku$ version, the copies of $B_k$ on either side of the $S_{k,\ell}$ are isomorphic, as is clear in Figure \ref{kuB34}. This is not the case in the $ko$ version, as can be observed in Figure \ref{B34}. There is an exotic $\eta$ extension from 88 to 89, accompanying the exotic $\cdot 2$ in 90. One might prefer to lower the filtrations of the high classes in 89 and 90, but they are $v_1^4$ times the classes in 81 and 82. Also, when we show how these charts are obtained, it will be clear that the high filtrations are warranted.

Let $\E'_{e,\infty}=\ds\dirlim_{\ell\to\infty}\E'_{e,\ell}$ and $M^e_\infty=\ds\dirlim_{\ell\to\infty}M^e_\ell$. Then $\E'_{e,\infty}$ can also be defined as the chart obtained from the sequence of $d_1$ differentials, situated as in Definition \ref{E'def},
$$\Si M^e_\infty\lar\Si^8M_4^{e-1}\lar\Si^{16}M_4^{e-1}\lar\Si^{32}M_4^{e-1}\lar\cdots.$$

\begin{defin}\label{hdef} For $k\ge2$, define functions $h_k$ by
$$h_k(8a+b)=4a-k+\begin{cases}0&b=2\\
1&b=3\\ 2&b=4\\ 3&5\le b\le9.\end{cases}$$
Let $C_{i,k}$ denote the subchart of $\E'_{2+i,\infty}$ consisting of elements in position $(x,y)$ satisfying $y\ge h_{k+i}(x)$.
\end{defin}

In Figure \ref{C04}, we depict a portion of $\E'_{2,\infty}$ with all of $C_{0,4}$ indicated by big dots. The dashed line is the graph of $h_4$.

\bigskip
\begin{minipage}{6in}
\begin{fig}\label{C04}

{\bf $C_{0,4}$}

\begin{center}

\begin{\tz}[scale=.4]
\draw (35,0) -- (72.5,0);
\node at (72,-.7) {$72$};
\draw [ultra thick] (40,0) -- (40,1);
\draw [ultra thick] (58,8) -- (60,10);
\draw [ultra thick] (42,1) -- (44,3) -- (44,1.5);
\draw (44,1.5) -- (44,0);
\draw  (48,0) -- (48,2.5);
\draw [ultra thick] (48,2.5) -- (48,4);
\draw [ultra thick] (50,4) -- (52,6) -- (52,5.5);
\draw (52,5.5) -- (52,0);
\draw  (56,0) -- (56,6.5);
\draw [ultra thick] (56,6.5) -- (56,7);
\draw   (60,0) -- (60,10) -- (58,8);
\node at (59,9) {\blb};
\node at (60,10) {\blb};
\draw  (64,0) -- (64,10.5);
\draw [ultra thick] (64,10.5) -- (64,11);
\draw  (66,11) -- (68,13) -- (68,0);
\node at (66,11) {\llb};
\node at (67,12) {\llb};
\draw  (72,0) -- (72,14);
\draw [dashed] (37,-1.2) -- (41,-1.2) -- (45,3) -- (49,3) -- (53,7) -- (57,7) -- (61,11) -- (65,11) -- (69,15) -- (73,15);
\node at (36,0) {\blb};
\node at (40,0) {\blb};
\node at (40,1) {\blb};
\node at (42,1) {\blb};
\node at (43,2) {\blb};
\node at (44,3) {\blb};
\node at (44,2) {\blb};
\node at (48,3) {\blb};
\node at (48,4) {\blb};
\node at (50,4) {\blb};
\node at (51,5) {\blb};
\node at (52,6) {\blb};
\node at (56,7) {\blb};
\node at (58,8) {\blb};
\node at (64,11) {\blb};
\node  at (44,0) {\llb};
\node  at (44,1) {\llb};
\node  at (48,0) {\llb};
\node  at (48,1) {\llb};
\node  at (48,2) {\llb};
\node  at (52,0) {\llb};
\node  at (52,1) {\llb};
\node  at (52,2) {\llb};
\node  at (52,3) {\llb};
\node  at (52,4) {\llb};
\node  at (52,5) {\llb};
\node at (56,0) {\llb};
\node at (56,1) {\llb};
\node at (56,2) {\llb};
\node at (60,0) {\llb};
\node at (60,1) {\llb};
\node at (60,2) {\llb};
\node at (60,3) {\llb};
\node at (64,0) {\llb};
\node at (64,1) {\llb};
\node at (64,2) {\llb};
\node at (64,3) {\llb};
\node at (64,4) {\llb};
\node at (64,5) {\llb};
\node at (64,6) {\llb};
\node at (68,0) {\llb};
\node at (68,1) {\llb};
\node at (68,2) {\llb};
\node at (68,3) {\llb};
\node at (68,4) {\llb};
\node at (68,5) {\llb};
\node at (68,6) {\llb};
\node at (68,7) {\llb};
\node at (72,0) {\llb};
\node at (72,1) {\llb};
\node at (72,2) {\llb};
\node at (72,3) {\llb};
\node at (72,4) {\llb};
\node at (72,5) {\llb};
\node at (72,6) {\llb};
\node at (72,7) {\llb};
\node at (72,8) {\llb};
\node at (72,9) {\llb};
\node at (72,10) {\llb};
\node at (72,11) {\llb};
\node at (72,12) {\llb};
\node at (72,13) {\llb};
\node at (72,14) {\llb};

\node  at (56,3) {\llb};
\node  at (56,4) {\llb};
\node  at (56,5) {\llb};
\node  at (56,6) {\llb};
\node  at (60,4) {\llb};
\node  at (60,5) {\llb};
\node  at (60,6) {\llb};
\node at (60,7) {\llb};
\node at (60,8) {\llb};
\node  at (60,9) {\llb};

\node  at (64,7) {\llb};
\node at (64,8) {\llb};
\node  at (64,9) {\llb};
\node  at (64,10) {\llb};

\node  at (68,13) {\llb};
\node  at (68,12) {\llb};
\node  at (68,11) {\llb};
\node  at (68,10) {\llb};
\node  at (68,9) {\llb};
\node  at (68,8) {\llb};

\node at (36,-.7) {$4$};
\node at (40,-.7) {$8$};
\node at (44,-.7) {$12$};
\node at (48,-.7) {$16$};
\node at (52,-.7) {$20$};
\node at (56,-.7) {$24$};
\node at (60,-.7) {$28$};
\node at (64,-.7) {$32$};
\node at (68,-.7) {$36$};

\end{\tz}
\end{center}
\end{fig}
\end{minipage}
\bigskip

Theorem \ref{main} involved summands $z^i\B_{k,\ell}$. These charts will be derived in Section \ref{Bklpfsec}. Of course, $z^0\B_{k,\ell}=\B_{k,\ell}$. Just as we did with letting $\At_k=\Si^{-2^{k+1}}\A_k$, it will be convenient to let $\Bt_{k,\ell}=\Si^{-2^{\ell+2}}\B_{k,\ell}$.

\begin{thm} \label{ziBklthm} For $i\ge0$ and $1\le k<\ell$, the chart $z^i\Bt_{k,\ell}$ consists of the following four parts:
\begin{enumerate}
    \item The portion of $\Si^{2^{k+1}}M^i_{\ell-k+3}$ in filtration $0$ through $k$, with filtrations increased by $2^k-k-1$. If $M^i_{\ell-k+3}$ has a filtration-$k$ element in grading $1$ mod $8$, that element is omitted, since it appears in part ($4$).
    
    \item A modification of $\E_{\ell-k+i+1,\ell+i}$ and all edges under it, as enumerated in Theorem \ref{Akthm}, including  extensions among these edges. The modification is that the elements of $\E'_{\ell-k+i+1,\ell+i}$ in grading $4$ and $5$ mod $8$, which supported differentials in forming $\E_{\ell-k+i+1,\ell+i}$, do not support differentials in this case.
    \item A chart formed from $\Si^{2^{k+1}}C_{i,k}$ together with all the edges strictly under $\Si^{2^{k+1}}\E_{2+i,k+1+i}$, incorporating all the differentials and extensions among these lower edges, and from them into $\Si^{2^{k+1}}\E'_{2+i,k+i+1}$. The target elements in $\Si^{2^{k+1}}\E'_{2+i,k+i+1}$
    are part of $\Si^{2^{k+1}}C_{i,k}$.
    \item For $\bigl[\frac{k+3+i}4\bigr]\le j\le 2^{k-2}+\bigl[\frac{i+3}4\bigr]-1$, elements $x_j$ in $(2^{k+1}+1+8j,2^k-k-1+4j-i)$ and $\eta x_j$. If $i=4t+1$ for some $t$, and $j=2^{k-2}+t$, then $\eta x_j$ is not present. If $j\ge 2^{k-3}+\bigl[\frac{i+3}4\bigr]$, there is an exotic extension from the  element in $(2^{k+1}+8j+2,4j-k-i)$ to $\eta x_j$.
\end{enumerate}
\end{thm}

\begin{rmk} \label{rem}{\rm If $\ell-k+i\equiv0$ mod 4, there is an $\eta$-extension from the last element of $\E_{\ell-k+i+1,\ell+i}$ in $z^i\Bt_{k,\ell}$ to a $\zt$ in the lowest filtration of $\Si^{2^{k+1}}M^i_{\ell-k+3}$. The first element is like the circled element at the end of the $e=5$ part of Figure \ref{four} (which does not support a differential in $z^i\Bt_{k,\ell}$) and the second element is like the element in grading 5 in the $M_6^1$ chart in Figure \ref{M6}.}\end{rmk}

\begin{rmk}\label{rem2}{\rm There are exotic $\eta$ extensions in $z^i\Bt_{k,\ell}$ wherever the exotic $\cdot2$ extensions occur in part (4) of the theorem. If $2\a=\eta\b$ and $v_1\g=\a$, then $\eta\g=\b$. In Figure \ref{M11}, the classes $\a$ are the two classes supporting the exotic extensions, and if $\a$ is in position $(x,y)$, then $\g$ is in position $(x-2,y-1)$.}\end{rmk}

The middle lightning flash in Figure \ref{B34} furnishes one example of part (1) of the theorem, which  corresponds to the $S_{k,\ell}$ portion in $ku^*(K_2)$. In Figure \ref{Sexpls}, we provide three more examples, without indicating the  filtration.

\bigskip
\begin{minipage}{6in}
\begin{fig}\label{Sexpls}

{\bf Examples of Theorem \ref{ziBklthm}(1)}

\begin{center}

\begin{\tz}[scale=.3]
\draw (1,0) -- (3,2) -- (3,0) -- (5,2);
\draw (9,4) -- (10,5);
\draw (11,5) -- (11,4) -- (12,5);
\draw (15,0) -- (17,2) -- (17,0);
\draw (21,2) -- (21,3);
\draw (25,4) -- (25,3) -- (26,4);
\draw (28,0) -- (30,2) -- (30,0);
\draw (34,0) -- (34,3);
\draw (36,4) -- (37,5);
\draw (38,5) -- (38,1) -- (40,3);
\draw (42,4) -- (42,5);
\node at (1,0) {\lb};
\node at (2,1) {\lb};
\node at (3,2) {\lb};
\node at (3,0) {\lb};
\node at (4,1) {\lb};
\node at (5,2) {\lb};
\node at (7,3) {\lb};
\node at (9,4) {\lb};
\node at (10,5) {\lb};
\node at (11,5) {\lb};
\node at (11,4) {\lb};
\node at (12,5) {\lb};
\node at (15,0) {\lb};
\node at (16,1) {\lb};
\node at (17,2) {\lb};
\node at (17,1) {\lb};
\node at (17,0) {\lb};
\node at (21,2) {\lb};
\node at (21,3) {\lb};
\node at (25,4) {\lb};
\node at (25,3) {\lb};
\node at (26,4) {\lb};
\node at (28,0) {\lb};
\node at (29,1) {\lb};
\node at (30,2) {\lb};
\node at (30,1) {\lb};
\node at (30,0) {\lb};
\node at (34,0) {\lb};
\node at (34,1) {\lb};
\node at (34,2) {\lb};
\node at (34,3) {\lb};
\node at (36,4) {\lb};
\node at (37,5) {\lb};
\node at (38,5) {\lb};
\node at (38,4) {\lb};
\node at (38,3) {\lb};
\node at (38,2) {\lb};
\node at (38,1) {\lb};
\node at (39,2) {\lb};
\node at (40,3) {\lb};
\node at (42,4) {\lb};
\node at (42,5) {\lb};
\node at (46,5) {\lb};
\node at (3,-.6) {$67$};
\node at (7,-.6) {$71$};
\node at (11,-.6) {$75$};
\node at (17,-.6) {$35$};
\node at (18,0) {\lb};
\node at (19,1) {\lb};
\draw (18,0) -- (19,1);
\node at (21,-.6) {$39$};
\node at (25,-.6) {$43$};
\node at (30,-.6) {$67$};
\node at (34,-.6) {$71$};
\node at (38,-.6) {$75$};
\node at (42,-.6) {$79$};
\node at (46,-.6) {$83$};
\node at (7,-2.1) {$\ell=7,\,k=5$};
\node at (21,-2.1) {$\ell=7,\,k=4$};
\node at (36,-2.1) {$\ell=10,\,k=5$};

\end{\tz}
\end{center}
\end{fig}
\end{minipage}
\bigskip

In Figure \ref{B47}, we provide another example, $\Bt_{4,7}$. The portion from part (2) of Theorem \ref{ziBklthm} is in black, and parts (1), (3), and (4) are in red. The black part corresponds to the modified version of $\E_{4,7}$ and the edges under it. Compare with $\E'_{4,7}$ in Figure \ref{four}.

The low red part is $\Si^{32}C_{0,4}$ together with the edges under $\Si^{32}\E_{2,5}$, which are $\Si^{40}\E_{3,4}$, $\Si^{48}\E_{3,5}$, and $\Si^{56}\E_{4,5}$. Compare with Figure \ref{C04} and Figure \ref{four}. Note that elements in $\Si^{48}\E'_{3,5}$ support differentials killing elements in $\Si^{32}C_{0,4}$ in grading 59 and 60.

The top red part has part (1), which is the middle chart in Figure \ref{Sexpls}, and the three $\eta$ pairs, which are part (4).

\bigskip
\begin{minipage}{6in}
\begin{fig}\label{B47}

{\bf $\Bt_{4,7}$}

\begin{center}

\begin{\tz}[scale=.28]
\draw (9.5,0) -- (67,0);
\node at (12,-.7) {$12$};
\node at (16,-.7) {$16$};
\node at (20,-.7) {$20$};
\node at (24,-.7) {$24$};
\node at (28,-.7) {$28$};
\node at (32,-.7) {$32$};
\node at (36,-.7) {$36$};
\node  at (40,-.7) {$40$};
\node at (44,-.7) {$44$};
\node at (48,-.7) {$48$};
\node at (52,-.7) {$52$};
\node at (56,-.7) {$56$};
\node at (60,-.7) {$60$};
\node at (64,-.7) {$64$};
\draw (10,0) -- (11,1);
\draw (12,0) -- (12,1);
\draw (13,0) -- (14,1);
\draw (19,0) -- (20,1);
\draw (20,4) -- (20,3) -- (22,5);
\draw (27,0) -- (28,1) -- (28,0);
\draw (37,0) -- (38,1);
\draw (34,2) -- (36,4);
\draw (26,6) -- (27,7);
\draw (28,7) -- (30,9);
\node [red] at (46,0) {\llb};
\node [red] at (50,1) {\llb};
\node [red] at (56,0) {\llb};
\node [red] at (58,1) {\llb};
\node [red] at (59,1) {\llb};
\node [red] at (60,2) {\llb};
\node [red] at (62,3) {\llb};
\node [red] at (66,4) {\llb};
\node [red] at (66,0) {\llb};
\node [red] at (67,1) {\llb};
\draw [red] (59,1) -- (60,2);
\draw [red] (66,0) -- (67,1);
\draw [dashed] [red] (66,0) -- (66,4);

\draw [red] (36,11) -- (37,12);
\draw [red] (33,11) -- (35,13) -- (35,11);
\draw [red] (39,13) -- (39,14);
\draw [red] (41,15) -- (42,16);
\draw [red] (43,15) -- (43,14) -- (44,15);
\draw [red] (49,19) -- (50,20);
\draw [red] (57,23) -- (58,24);
\draw [red] (40,0) -- (40,1);
\draw [red] (42,1) -- (44,3) -- (44,2);
\draw [red] (48,3) -- (48,4);
\draw [red] (50,4) -- (52,6);
\node at (10,0) {\llb};
\node at (11,1) {\llb};
\node at (12,0) {\llb};
\node at (12,1) {\llb};
\node at (13,0) {\llb};
\node at (14,1) {\llb};
\node at (16,2) {\llb};
\node at (18,3) {\llb};
\node at (19,0) {\llb};
\node at (20,1) {\llb};
\node at (20,3) {\llb};
\node at (20,4) {\llb};
\node at (21,4) {\llb};
\node at (22,5) {\llb};
\node at (27,0) {\llb};
\node at (28,1) {\llb};
\node at (28,0) {\llb};
\node at (30,1) {\llb};
\node at (36,0) {\llb};
\node at (37,0) {\llb};
\node at (38,1) {\llb};
\node at (34,2) {\llb};
\node at (35,3) {\llb};
\node at (36,4) {\llb};
\node at (26,6) {\llb};
\node at (27,7) {\llb};
\node at (28,7) {\llb};
\node at (29,8) {\llb};
\node at (30,9) {\llb};
\node at (34,10) {\llb};
\node [red] at (36,11) {\llb};
\node [red] at (37,12) {\llb};
\node [red] at (36.2,0) {\llb};
\node [red] at (40,0) {\llb};
\node [red] at (40,1) {\llb};
\node [red] at (42,1) {\llb};
\node [red] at (43,2) {\llb};
\node [red] at (44,3) {\llb};
\node [red] at (44,2) {\llb};
\node [red] at (48,3) {\llb};
\node [red] at (48,4) {\llb};
\node [red] at (50,4) {\llb};
\node [red] at (51,5) {\llb};
\node [red] at (52,6) {\llb};
\node [red] at (56,7) {\llb};
\node [red] at (58,8) {\llb};
\node [red] at (64,11) {\llb};
\node [red] at (33,11) {\llb};
\node [red] at (34,12) {\llb};
\node [red] at (35,13) {\llb};
\node [red] at (35,12) {\llb};
\node [red] at (35,11) {\llb};
\node [red] at (39,13) {\llb};
\node [red] at (39,14) {\llb};
\node [red] at (41,15) {\llb};
\node [red] at (42,16) {\llb};
\node [red] at (43,15) {\llb};
\node [red] at (43,14) {\llb};
\node [red] at (44,15) {\llb};
\node [red] at (49,19) {\llb};
\node [red] at (50,20) {\llb};
\node [red] at (57,23) {\llb};
\node [red] at (58,24) {\llb};
\draw [dashed] (30,1) -- (30,9);
\draw [dashed] (34,2) -- (34,10);
\draw [dashed] [red] (50,4) -- (50,20);
\draw [dashed] [red] (58,8) -- (58,24);

\end{\tz}
\end{center}
\end{fig}
\end{minipage}
\bigskip

Although $z^i\Bt_{k,\ell}$ is built from various $z^iM^a_b$ ($=M^{a+i}_b$), and multiplying $M^a_b$ by $z^i$ just lowers all filtrations by $i$, the effect of $z^i$ on $\Bt_{k,\ell}$ is more complicated, due to differentials. In Figure \ref{ziB34}, we display $z^i\Bt_{k,\ell}$ for $1\le i\le 3$, for comparison with the case $i=0$ in Figure \ref{B34} (which is $\Si^{64}\Bt_{3,4}$).

\bigskip
\begin{minipage}{6in}
\begin{fig}\label{ziB34}

{\bf $z^i\Bt_{3,4}$}

\begin{center}

\begin{\tz}[scale=.45]
\draw (10,0) -- (38,0);
\draw (11,0) -- (12,1) -- (12,0);
\draw (13,0) -- (14,1);
\draw (18,2) -- (19,3);
\draw (20,3) -- (21,4);
\draw (25,5) -- (27,7) -- (27,6) -- (28,7);
\draw (27,0) -- (28,1) -- (28,0);
\draw (34,2) --  (35,3);
\draw (33,9) -- (34,10);
\draw [dashed] (34,2) -- (34,10);
\node at (11,0) {\lb};
\node at (12,1) {\lb};
\node at (12,0) {\lb};
\node at (13,0) {\lb};
\node at (14,1) {\lb};
\node at (18,2) {\lb};
\node at (19,3) {\lb};
\node at (20,3) {\lb};
\node at (21,4) {\lb};
\node at (20,0) {\lb};
\node at (22,1) {\lb};
\node at (27,0) {\lb};
\node at (28,1) {\lb};
\node at (28,0) {\lb};
\node at (32,1) {\lb};
\draw [dotted] (32,1) -- (33,9);
\node at (36,0) {\lb};
\node at (38,1) {\lb};
\node at (34,2) {\lb};
\node at (35,3) {\lb};
\node at (25,5) {\lb};
\node at (26,6) {\lb};
\node at (27,7) {\lb};
\node at (27,6) {\lb};
\node at (28,7) {\lb};
\node at (33,9) {\lb};
\node at (34,10) {\lb};
\node at (12,-.7) {$12$};
\node at (16,-.7) {$16$};
\node at (20,-.7) {$20$};
\node at (24,-.7) {$24$};
\node at (28,-.7) {$28$};
\node at (32,-.7) {$32$};
\node at (36,-.7) {$36$};
\node at (13,5) {$z^3\Bt_{3,4}$};

\draw (9,12) -- (36,12);
\draw (10,12) -- (11,13);
\draw (12,13) -- (12,12) -- (14,14);
\draw (19,12) -- (20,13);
\draw (20,16) -- (21,17);
\draw (26,12) -- (27,13);
\draw (25,18) -- (26,19);
\draw (35,12) -- (36,13);
\draw (33,22) -- (34,23);
\draw [dashed] (34,15) -- (34,23);
\node at (12,11.3) {$12$};
\node at (16,11.3) {$16$};
\node at (20,11.3) {$20$};
\node at (24,11.3) {$24$};
\node at (28,11.3) {$28$};
\node at (32,11.3) {$32$};
\node at (36,11.3) {$36$};
\node at (10,12) {\lb};
\node at (11,13) {\lb};
\node at (12,13) {\lb};
\node at (12,12) {\lb};
\node at (13,13) {\lb};
\node at (14,14) {\lb};
\node at (19,12) {\lb};
\node at (20,13) {\lb};
\node at (18,15) {\lb};
\node at (19,16) {\lb};
\node at (20,16) {\lb};
\node at (21,17) {\lb};
\node at (26,12) {\lb};
\node at (27,13) {\lb};
\node at (28,13) {\lb};
\node at (32,14) {\lb};
\node at (34,15) {\lb};
\node at (35,12) {\lb};
\node at (36,13) {\lb};
\node at (25,18) {\lb};
\node at (26,19) {\lb};
\node at (27,19) {\lb};
\node at (33,22) {\lb};
\node at (34,23) {\lb};
\node at (13,16) {$z^2\Bt_{3,4}$};
\draw [dotted] (32,14) -- (33,22);

\node at (16,24.3) {$16$};
\node at (20,24.3) {$20$};
\node at (24,24.3) {$24$};
\node at (28,24.3) {$28$};
\node at (32,24.3) {$32$};
\node at (8,25) {\lb};
\node at (10,26) {\lb};
\node at (11,26) {\lb};
\node at (12,27) {\lb};
\node at (12,26) {\lb};
\node at (13,27) {\lb};
\node at (14,28) {\lb};
\node at (17.8,25) {\lb};
\node at (19,26) {\lb};
\node at (18,29) {\lb};
\node at (17.8,29.2) {\lb};
\node at (19,30) {\lb};
\node at (19,29) {\lb};
\node at (20,30) {\lb};
\node at (21,31) {\lb};
\node at (25,32) {\lb};
\node at (26,33) {\lb};
\node at (24,25) {\lb};
\draw (7.5,25) -- (35,25);
\draw (11,26) -- (12,27) -- (12,26) -- (14,28);
\draw (18,29) -- (19,30) -- (19,29) -- (21,31);
\draw (17.8,25) -- (19,26);
\draw [dashed] (17.8,25) -- (17.8,29.2);
\draw (27,26) -- (28,27);
\draw (34,25) -- (35,26);
\draw (25,32) -- (26,33);
\draw [dotted] (32,28) -- (33,36);
\node at (8,24.3) {$8$};
\node at (12,24.3) {$12$};
\node at (26,26) {\lb};
\node at (27,26) {\lb};
\node at (28,27) {\lb};
\node at (34,25) {\lb};
\node at (35,26) {\lb};
\node at (32,28) {\lb};
\node at (34,29) {\lb};
\node at (33,36) {\lb};
\draw [dashed] (34,25) -- (34,29);
\draw (34,25) -- (35,26);
\node at (12,30) {$z\Bt_{3,4}$};
\end{\tz}
\end{center}
\end{fig}
\end{minipage}
\bigskip

We close this descriptive section by describing a way of visualizing how parts (1), (3), and (4) of Theorem \ref{ziBklthm} come about. For simplicity, we restrict to $i=0$. In Section \ref{Bklpfsec}, we will explain why this description is valid.

For this part of $\Bt_{k,\ell}$, we combine two charts. One chart has the resultant of
$$\Si^{2^{k+1}}(\Si M^2_{\ell+2}\lar\Si^8M_4^1\lar\cdots\Si^{2^i}M_4^1\lar\cdots\lar\Si^{2^{k+1}}M_4^1)$$
with the usual placement for $d_1$ differentials. The other chart has $\Si^{2^{k+1}}M^0_{\ell-k+3}$, with filtrations increased by $2^k-k-1$. See Figure \ref{M11} for $\Bt_{4,9}$.

Beginning in grading slightly greater than $2^{k+2}$, the two charts will have isomorphic forms, stably $M_{\ell-k+3}$, displaced by 1 horizontal unit. There must be a differential annihilating these; it will be $d_{2^k}$. We apply $(v_1^{-4})$ periodicity to these differentials whenever it applies. In Figure \ref{M11}, we indicate by small red dots the elements that are related by these differentials. For example, applying $v_1^{-4}$ to the $\Z/16$ in 72 and 71 shows that the bottom four elements in 64 support differentials, but the top one survives. The elements which survive are indicated by large black dots. Three $\eta$-pairs at the top are part (4) of Theorem \ref{ziBklthm}. The other black elements in the top half in filtration 0 to 4 (before the filtration shift) are type (1). The black elements in the bottom half are of type (3). They lie on or above the dashed line in Figure \ref{C04}. This description does not incorporate the type-(3) elements on subedges under $\Si^{2^{k+1}}\E_{2,k+1}$.

\bigskip
\begin{minipage}{6in}
\begin{fig}\label{M11}

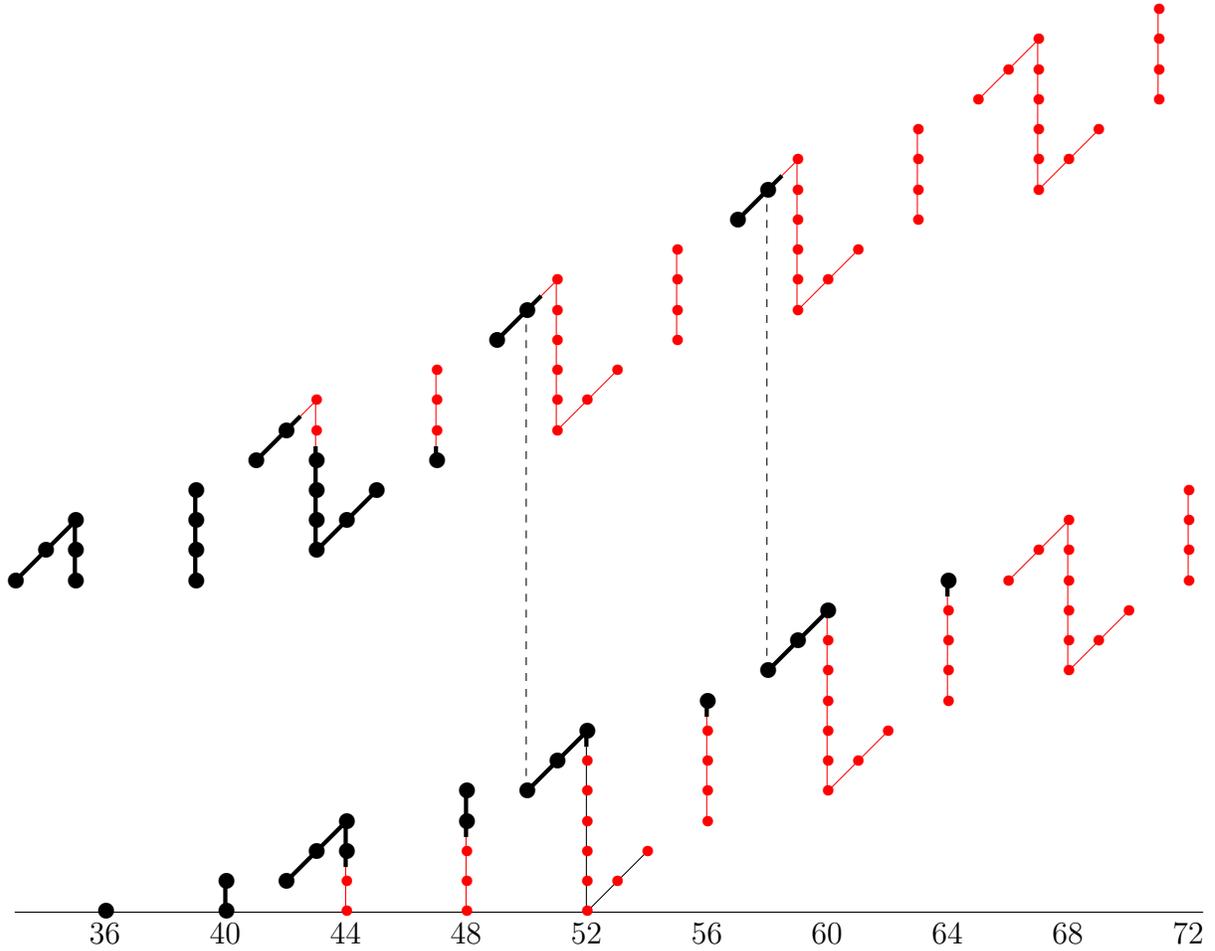
{\bf Some differentials (in red) in $\Bt_{4,9}$}

\begin{center}

\begin{\tz}[scale=.4]
\draw (33,0) -- (72.5,0);
\node at (72,-.7) {$72$};
\draw [ultra thick] (40,0) -- (40,1);
\draw [ultra thick] (42,1) -- (44,3) -- (44,1.5);
\draw [red] (44,1.5) -- (44,0);
\draw [red] (48,0) -- (48,2.5);
\draw [ultra thick] (48,2.5) -- (48,4);
\draw [ultra thick] (50,4) -- (52,6) -- (52,5.5);
\draw  (52,5.5) -- (52,0) -- (54,2);
\draw [red] (56,3) -- (56,6.5);
\draw [ultra thick] (56,6.5) -- (56,7);
\draw [red] (62,6) -- (60,4) -- (60,10);
\draw [ultra thick] (60,10) -- (58,8);
\node at (60,10) {\blb};
\node at (59,9) {\blb};
\draw [red] (64,7) -- (64,10.5);
\draw [ultra thick] (64,10.5) -- (64,11);
\draw [red] (66,11) -- (68,13) -- (68,8) -- (70,10);
\draw [red] (72,11) -- (72,14);
\draw [red] (71,27) -- (71,30);
\node [red] at (72,11) {\lb};
\node [red] at (72,12) {\lb};
\node [red] at (72,13) {\lb};
\node [red] at (72,14) {\lb};
\node [red] at (71,27) {\lb};
\node [red] at (71,28) {\lb};
\node [red] at (71,29) {\lb};
\node [red] at (71,30) {\lb};
\node at (36,0) {\blb};
\node at (40,0) {\blb};
\node at (40,1) {\blb};
\node at (42,1) {\blb};
\node at (43,2) {\blb};
\node at (44,3) {\blb};
\node at (44,2) {\blb};
\node at (48,3) {\blb};
\node at (48,4) {\blb};
\node at (50,4) {\blb};
\node at (51,5) {\blb};
\node at (52,6) {\blb};
\node at (56,7) {\blb};
\node at (58,8) {\blb};
\node at (64,11) {\blb};
\node [red] at (44,0) {\lb};
\node [red] at (44,1) {\lb};
\node [red] at (48,0) {\lb};
\node [red] at (48,1) {\lb};
\node [red] at (48,2) {\lb};
\node [red] at (52,0) {\lb};
\node [red] at (52,1) {\lb};
\node [red] at (52,2) {\lb};
\node [red] at (52,3) {\lb};
\node [red] at (52,4) {\lb};
\node [red] at (52,5) {\lb};
\node [red] at (53,1) {\lb};
\node [red] at (54,2) {\lb};
\node [red] at (56,3) {\lb};
\node [red] at (56,4) {\lb};
\node [red] at (56,5) {\lb};
\node [red] at (56,6) {\lb};
\node [red] at (60,4) {\lb};
\node [red] at (60,5) {\lb};
\node [red] at (60,6) {\lb};
\node [red] at (60,7) {\lb};
\node [red] at (60,8) {\lb};
\node [red] at (60,9) {\lb};
\node [red] at (61,5) {\lb};
\node [red] at (62,6) {\lb};
\node [red] at (64,7) {\lb};
\node [red] at (64,8) {\lb};
\node [red] at (64,9) {\lb};
\node [red] at (64,10) {\lb};
\node [red] at (66,11) {\lb};
\node [red] at (67,12) {\lb};
\node [red] at (68,13) {\lb};
\node [red] at (68,12) {\lb};
\node [red] at (68,11) {\lb};
\node [red] at (68,10) {\lb};
\node [red] at (68,9) {\lb};
\node [red] at (68,8) {\lb};
\node [red] at (69,9) {\lb};
\node [red] at (70,10) {\lb};
\draw [ultra thick] (33,11) -- (35,13) -- (35,11);
\draw [ultra thick]  (39,11) -- (39,14);
\draw [ultra thick] (41,15) -- (42.5,16.5);
\draw [ultra thick] (45,14) -- (43,12) -- (43,15.5);
\draw [ultra thick] (47,15) -- (47,15.5);
\draw [ultra thick] (49,19) -- (50.5,20.5);
\draw [ultra thick] (57,23) -- (58.5,24.5);
\draw [red] (43,15.5) -- (43,17) -- (42.5,16.5);
\draw [red] (47,15.5) -- (47,18);
\draw [red] (53,18) -- (51,16) -- (51,21) -- (50.5,20.5);
\draw [red] (55,19) -- (55,22);
\draw [red] (61,22) -- (59,20) -- (59,25) -- (58.5,24.5);
\draw [red] (63,23) -- (63,26);
\draw [red] (65,27) -- (67,29) -- (67,24) -- (69,26);
\node at (33,11) {\blb};
\node at (34,12) {\blb};
\node at (35,13) {\blb};
\node at (35,12) {\blb};
\node at (35,11) {\blb};
\node at (39,11) {\blb};
\node at (39,12) {\blb};
\node at (39,13) {\blb};
\node at (39,14) {\blb};
\node at (41,15) {\blb};
\node at (42,16) {\blb};
\node at (49,19) {\blb};
\node at (50,20) {\blb};
\node at (57,23) {\blb};
\node at (58,24) {\blb};
\node [red] at (43,16) {\lb};
\node [red] at (43,17) {\lb};
\node [red] at (47,16) {\lb};
\node [red] at (47,17) {\lb};
\node [red] at (47,18) {\lb};
\node [red] at (51,16) {\lb};
\node [red] at (51,17) {\lb};
\node [red] at (51,18) {\lb};
\node [red] at (51,19) {\lb};
\node [red] at (51,20) {\lb};
\node [red] at (51,21) {\lb};
\node [red] at (55,19) {\lb};
\node [red] at (55,20) {\lb};
\node [red] at (55,21) {\lb};
\node [red] at (55,22) {\lb};
\node [red] at (59,20) {\lb};
\node [red] at (59,21) {\lb};
\node [red] at (59,22) {\lb};
\node [red] at (59,23) {\lb};
\node [red] at (59,24) {\lb};
\node [red] at (59,25) {\lb};
\node [red] at (52,17) {\lb};
\node [red] at (53,18) {\lb};
\node [red] at (60,21) {\lb};
\node [red] at (61,22) {\lb};
\node [red] at (63,23) {\lb};
\node [red] at (63,24) {\lb};
\node [red] at (63,25) {\lb};
\node [red] at (63,26) {\lb};
\node [red] at (67,24) {\lb};
\node [red] at (67,25) {\lb};
\node [red] at (67,26) {\lb};
\node [red] at (67,27) {\lb};
\node [red] at (67,28) {\lb};
\node [red] at (67,29) {\lb};
\node [red] at (66,28) {\lb};
\node [red] at (65,27) {\lb};
\node [red] at (68,25) {\lb};
\node [red] at (69,26) {\lb};
\node at (36,-.7) {$36$};
\node at (40,-.7) {$40$};
\node at (44,-.7) {$44$};
\node at (48,-.7) {$48$};
\node at (52,-.7) {$52$};
\node at (56,-.7) {$56$};
\node at (60,-.7) {$60$};
\node at (64,-.7) {$64$};
\node at (68,-.7) {$68$};
\node at (43,12) {\blb};
\node at (43,13) {\blb};
\node at (43,14) {\blb};
\node at (43,15) {\blb};
\node at (44,13) {\blb};
\node at (45,14) {\blb};
\node at (47,15) {\blb};
\draw [dashed] (50,4) -- (50,20);
\draw [dashed] (58,8) -- (58,24);
\end{\tz}
\end{center}
\end{fig}
\end{minipage}
\bigskip

\section{Complete description through grading 42}\label{sec42}
We can easily depict $ko_*(K_2)$ for $*<42$, and much farther. In Figure \ref{42} we present it in three labeled rows, to avoid congestion. It is the sum of the three. We omit the trivial $ko_*$-submodule, which is enumerated through grading 24 at the end of Section \ref{kosummandsec}.

\bigskip
\begin{minipage}{6in}
\begin{fig}\label{42}

{\bf $ko_*(K_2)$, $*<42$}

\begin{center}

\begin{\tz}[scale=.39]
\draw (1.5,0) -- (42,0);
\node at (2,-.6) {$2$};
\node at (4,-.6) {$4$};
\node at (8,-.6) {$8$};
\node at (12,-.6) {$12$};
\node at (16,-.6) {$16$};
\node at (20,-.6) {$20$};
\node at (24,-.6) {$24$};
\node at (28,-.6) {$28$};
\node at (32,-.6) {$32$};
\node at (36,-.6) {$36$};
\node at (40,-.6) {$40$};
\node at (3,-1.6) {$\A_1$};
\node at (8.2,-1.6) {$\A_2$};
\node at (11.6,-1.6) {$\Si^8\A_1$};
\node at (19,-1.6) {$\Si^{16}\A_1$};
\node at (24.2,-1.6) {$\Si^{16}\A_2$};
\node at (27.6,-1.6) {$\Si^{24}\A_1$};
\node at (35,-1.6) {$\Si^{32}\A_1$};
\node at (41,-1.6) {$\Si^{32}\A_2$};
\draw (8,2) -- (8,0) -- (10,2);
\draw (11,2) -- (12,3);
\draw (24,2) -- (24,0) -- (26,2);
\draw (27,2) -- (28,3);
\draw (40,2) -- (40,0) -- (42,2);
\node at (2,0) {\lb};
\node at (4,1) {\lb};
\node at (8,0) {\lb};
\node at (8,1) {\lb};
\node at (8,2) {\lb};
\node at (9,1) {\lb};
\node at (10,2) {\lb};
\node at (10,0) {\lb};
\node at (12,1) {\lb};
\node at (11,2) {\lb};
\node at (12,3) {\lb};
\node at (18,0) {\lb};
\node at (20,1) {\lb};
\node at (24,2) {\lb};
\node at (24,1) {\lb};
\node at (24,0) {\lb};
\node at (25,1) {\lb};
\node at (26,2) {\lb};
\node at (26,0) {\lb};
\node at (28,1) {\lb};
\node at (27,2) {\lb};
\node at (28,3) {\lb};
\node at (34,0) {\lb};
\node at (36,1) {\lb};
\node at (40,0) {\lb};
\node at (40,1) {\lb};
\node at (40,2) {\lb};
\node at (41,1) {\lb};
\node at (42,2) {\lb};
\draw (2,7) -- (42,7);
\node at (16,6.4) {$16$};
\node at (20,6.4) {$20$};
\node at (24,6.4) {$24$};
\node at (28,6.4) {$28$};
\node at (32,6.4) {$32$};
\node at (36,6.4) {$36$};
\node at (40,6.4) {$40$};
\node at (21,5.4) {$\A_3$};
\node at (33,5.4) {$\A_4$};
\node at (37,5.4) {$\B_{1,3}$};
\draw (16,10) -- (16,7) -- (18,9);
\draw (19,10) -- (20,11) -- (20,10);
\draw (24,11) -- (26,13) -- (26,12) -- (27,13);
\draw (32,11) -- (32,7) -- (34,9);
\draw (36,10) -- (36,12) -- (35,11);
\draw (40,12) -- (40,11) -- (42,13);
\node at (16,7) {\lb};
\node at (16,8) {\lb};
\node at (16,9) {\lb};
\node at (16,10) {\lb};
\node at (17,8) {\lb};
\node at (18,9) {\lb};
\node at (19,10) {\lb};
\node at (20,11) {\lb};
\node at (20,10) {\lb};
\node at (24,11) {\lb};
\node at (25,12) {\lb};
\node at (26,13) {\lb};
\node at (26,12) {\lb};
\node at (27,13) {\lb};
\node at (28,7) {\lb};
\node at (30,8) {\lb};
\node at (32,7) {\lb};
\node at (32,8) {\lb};
\node at (32,9) {\lb};
\node at (32,10) {\lb};
\node at (32,11) {\lb};
\node at (33,8) {\lb};
\node at (34,9) {\lb};
\node at (35,7) {\lb};
\node at (36,10) {\lb};
\node at (36,11) {\lb};
\node at (36,12) {\lb};
\node at (35,11) {\lb};
\node at (37,7) {\lb};
\node at (39,8) {\lb};
\node at (40,11) {\lb};
\node at (40,12) {\lb};
\node at (41,12) {\lb};
\node at (42,13) {\lb};
\draw (2,17) -- (42,17);
\node at (20,16.4) {$20$};
\node at (24,16.4) {$24$};
\node at (28,16.4) {$28$};
\node at (32,16.4) {$32$};
\node at (36,16.4) {$36$};
\node at (40,16.4) {$40$};
\node at (20,15.4) {$\B_{1,2}$};
\node at (28,15.4) {$\Si^8\B_{1,2}$};
\node at (36,15.4) {$\Si^{16}\B_{1,2}$};
\node at (39,15.4) {$\B_{2,3}$};
\draw (20,17) -- (21,18);
\draw (28,17) -- (29,18);
\draw (36,17) -- (37,18);
\draw (37,17) -- (38,18);
\draw (41,18) -- (43,20) -- (43,19);
\node at (19,17) {\lb};
\node at (20,17) {\lb};
\node at (21,18) {\lb};
\node at (27,17) {\lb};
\node at (28,17) {\lb};
\node at (29,18) {\lb};
\node at (35,17) {\lb};
\node at (36,17) {\lb};
\node at (36.3,17) {\lb};
\node at (37,18) {\lb};
\node at (37,17) {\lb};
\node at (38,18) {\lb};
\node at (41,18) {\lb};
\node at (42,19) {\lb};
\end{\tz}
\end{center}
\end{fig}
\end{minipage}
\bigskip

\section{A formulaic description of edges}\label{expsec}
In this section, we give a more formulaic description of edges.

 Let $R=\Z_{(2)}[\eta]/(2\eta,\eta^3)$.

The upper edge $\E_{1,k}$ has a different form than the other edges.

\begin{defin}\label{Ldef} For $t\ge1$ and $0\le\eps\le2$, we define $R$-modules $L_{t,\eps}(x,y)$ to have generators $g$, $g'$, and $g''$ in positions $(x,y)$, $(x+2,y+t)$, and $(x+4,y+3)$, respectively, with relations $2^tg$, $\eta^2g'=2^{t-1}g''$,  $\eta^{3-\eps}g'$, $\eta g''$, and
$$2g'=\begin{cases} \eta^2g&t=1\\ 0&t>1.\end{cases}$$
Let $I_k$ denote the $R$-module with generators $G$, $G'$, and $G''$ in positions $(0,0)$, $(3,k)$, and $(4,3)$, respectively, with relations $2^{k+1}G$, $\eta G'=2^{k-2}G''$, $2G'$, and $\eta G''$.\end{defin}

For example, $I_4$ is the portion of Figure \ref{E14} in grading $\le7$.
The second subscript in $L_{t,\eps}$ indicates the number of elements killed in forming $L_{t,\eps}$ from $L_{t,0}$. In Figure \ref{Lchart} we show $L_{t,\eps}$ for $t\le 4$, omitting $(x,y)$, which is the position of the lower-left element.

\bigskip
\begin{minipage}{6in}
\begin{fig}\label{Lchart}

{\bf $L_{t,\eps}$}

\begin{center}

\begin{\tz}[scale=.5]
\draw (0,7) -- (2,9) -- (2,8) -- (4,10);
\draw (8,8) -- (8,7) -- (10,9);
\draw (10,9.3) -- (12,11) -- (12,10);
\draw (16,9) -- (16,7) -- (18,9);
\draw (18,10) -- (20,12) -- (20,10);
\draw (24,10) -- (24,7) -- (26,9);
\draw (26,11) -- (28,13) -- (28,10);
\draw (0,0) -- (2,2) -- (2,1) -- (3,2);
\draw (8,1) -- (8,0) -- (10,2);
\draw (16,2) -- (16,0) -- (18,2);
\draw (20,3) -- (20,4);
\draw (24,3) -- (24,0) -- (26,2);
\draw (28,3) -- (28,5);
\draw (0,-7) -- (2,-5) -- (2,-6);
\draw (8,-6) -- (8,-7) -- (10,-5);
\draw (16,-5) -- (16,-7) -- (18,-5);
\draw (20,-4) -- (20,-3);
\draw (24,-4) -- (24,-7) -- (26,-5);
\draw (28,-4) -- (28,-2);
\node at (3,7) {$L_{1,0}$};
\node at (3,0) {$L_{1,1}$};
\node at (3,-7) {$L_{1,2}$};
\node at (11,7) {$L_{2,0}$};
\node at (11,0) {$L_{2,1}$};
\node at (11,-7) {$L_{2,2}$};
\node at (19,7) {$L_{3,0}$};
\node at (19,0) {$L_{3,1}$};
\node at (19,-7) {$L_{3,2}$};
\node at (27,7) {$L_{4,0}$};
\node at (27,0) {$L_{4,1}$};
\node at (27,-7) {$L_{4,2}$};
\node at (0,0) {\lb};
\node at (1,1) {\lb};
\node at (2,2) {\lb};
\node at (2,1) {\lb};
\node at (3,2) {\lb};
\node at (8,0) {\lb};
\node at (8,1) {\lb};
\node at (9,1) {\lb};
\node at (10,2) {\lb};
\node at (10,2.3) {\lb};
\node at (11,3) {\lb};
\draw (10,2.3) -- (11,3);
\node at (12,3) {\lb};
\node at (16,0) {\lb};
\node at (16,1) {\lb};
\node at (16,2) {\lb};
\node at (17,1) {\lb};
\node at (18,2) {\lb};
\node at (18,3) {\lb};
\node at (19,4) {\lb};
\draw (18,3) -- (19,4);
\node at (20,3) {\lb};
\node at (20,4) {\lb};
\node at (24,0) {\lb};
\node at (24,1) {\lb};
\node at (24,2) {\lb};
\node at (24,3) {\lb};
\node at (25,1) {\lb};
\node at (26,2) {\lb};
\node at (26,4) {\lb};
\node at (27,5) {\lb};
\draw (26,4) -- (27,5);
\node at (28,3) {\lb};
\node at (28,4) {\lb};
\node at (28,5) {\lb};
\node at (0,-7) {\lb};
\node at (1,-6) {\lb};
\node at (2,-5) {\lb};
\node at (2,-6) {\lb};
\node at (8,-7) {\lb};
\node at (8,-6) {\lb};
\node at (9,-6) {\lb};
\node at (10,-5) {\lb};
\node at (10,-4.7) {\lb};
\node at (12,-4) {\lb};
\node at (16,-7) {\lb};
\node at (16,-6) {\lb};
\node at (16,-5) {\lb};
\node at (17,-6) {\lb};
\node at (18,-5) {\lb};
\node at (18,-4) {\lb};
\node at (20,-4) {\lb};
\node at (20,-3) {\lb};
\node at (24,-7) {\lb};
\node at (24,-6) {\lb};
\node at (24,-5) {\lb};
\node at (24,-4) {\lb};
\node at (25,-6) {\lb};
\node at (26,-5) {\lb};
\node at (26,-3) {\lb};
\node at (28,-4) {\lb};
\node at (28,-3) {\lb};
\node at (28,-2) {\lb};
\node at (0,7) {\lb};
\node at (1,8) {\lb};
\node at (2,9) {\lb};
\node at (2,8) {\lb};
\node at (3,9) {\lb};
\node at (4,10) {\lb};
\node at (8,7) {\lb};
\node at (8,8) {\lb};
\node at (9,8) {\lb};
\node at (10,9) {\lb};
\node at (10,9.3) {\lb};
\node at (11,10.15) {\lb};
\node at (12,11) {\lb};
\node at (12,10) {\lb};
\node at (16,7) {\lb};
\node at (16,8) {\lb};
\node at (16,9) {\lb};
\node at (17,8) {\lb};
\node at (18,9) {\lb};
\node at (18,10) {\lb};
\node at (19,11) {\lb};
\node at (20,12) {\lb};
\node at (20,11) {\lb};
\node at (20,10) {\lb};
\node at (24,7) {\lb};
\node at (24,8) {\lb};
\node at (24,9) {\lb};
\node at (24,10) {\lb};
\node at (25,8) {\lb};
\node at (26,9) {\lb};
\node at (26,11) {\lb};
\node at (27,12) {\lb};
\node at (28,13) {\lb};
\node at (28,12) {\lb};
\node at (28,11) {\lb};
\node at (28,10) {\lb};

\end{\tz}
\end{center}
\end{fig}
\end{minipage}
\bigskip

\begin{thm}\label{upperthm} Let $\lg(i)=[\log_2(i)]$. The upper edge, $\E_{1,k}$, of $\Sigma^{-2^{k+1}}\A_k$, as an $R$-module is
$$I_k\oplus\bigoplus_{i=1}^{2^{k-2}-1}L_{k-\lg(i)-2,f_1(i)}(8i,4i),$$
where $$f_1(i)=\begin{cases}0&2^j\le i\le2^j+[(j-1)/4]\\
1&i=2^j+\frac j4, j\equiv 0\ (4)\\
2&2^j+1+[j/4]\le i\le 2^{j+1}-1\end{cases}
$$
for some $j$.
\end{thm}

Many values of $f_1$ are presented in Table \ref{tbl}.
For example, the upper edge in Figure \ref{A5} satisfies
$$\E_{1,5}=I_5\oplus L_{3,1}\oplus L_{2,0}\oplus L_{2,2}\oplus L_{1,0}\oplus L_{1,2}\oplus L_{1,2}\oplus L_{1,2},$$
where the $i$th $L$ is in position $(8i,4i)$, $1\le i\le7$.

The other edges are formed from $R$-modules $K_{t,\eps}$, which we now define.
\begin{defin}\label{Kdef} For $t\ge1$ and $0\le\eps\le2$, the $R$-module $K_{t,\eps}(x,y)$ has generators
$g$, $g'$, $g''$, and $g'''$ in positions $(x,y)$, $(x-2,y+t-3)$, $(x+2,y+1)$, and $(x+4,y+2)$, respectively, with relations
$$2^{t-1}g=\eta^2g',\ 2g',\ \eta g,\ 2g'',\ \eta g'',\ 2^{t-1}g''',\ \eta g''',\ \eta^{3-\eps}g'.$$
\end{defin}

In Figure \ref{ch2}, we depict $K_{t,\eps}$ for $t\le4$. Position $(x,y)$ is circled. Note that $K_{1,2}$ has its only nonzero elements in $(x-2,y-2)$ and $(x+2,y+1)$. Again, the second subscript is the number of elements killed in forming $K_{t,\eps}$ from $K_{t,0}$.

\bigskip
\begin{minipage}{6in}
\begin{fig}\label{ch2}

{\bf $K_{t,\eps}$}

\begin{center}

\begin{\tz}[scale=.45]
\draw (2,0) -- (2,2);
\draw (6,2) -- (6,4);
\draw (12,1) -- (12,2);
\draw (2,7) -- (2,9);
\draw (1,9) -- (0,8);
\draw (6,9) -- (6,11);
\draw (10,8) -- (11,9); \draw (12,9) -- (12,8);
\draw (2,14) -- (2,17) -- (0,15);
\draw (6,16) -- (6,18);
\draw (10,15) -- (12,17) -- (12,15); 
\draw (16,3) -- (16,4);
\draw (16,10) -- (16,11);
\draw (20,8) -- (21,9);
\draw (30,8) -- (31,9);
\draw (16,17) -- (16,18);
\draw (20,15) -- (22,17) -- (22,16);
\draw (30,15) -- (32,17);
\node at (0,1) {\lb};
\node at (2,0) {\lb};
\node at (2,1) {\lb};
\node at (2,2) {\lb};
\node at (4,1) {\lb};
\node at (6,2) {\lb};
\node at (6,3) {\lb};
\node at (6,4) {\lb};
\node at (10,1) {\lb};
\node at (12,1) {\lb};
\node at (12,2) {\lb};
\node at (14,2) {\lb};
\node at (16,3) {\lb};
\node at (16,4) {\lb};
\node at (20,1) {\lb};
\node at (22,2) {\lb};
\node at (24,3) {\lb};
\node at (26,4) {\lb};
\node at (30,1) {\lb};
\node at (34,4) {\lb};
\node at (0,8) {\lb};
\node at (1,9) {\lb};
\node at (2,9) {\lb};
\node at (2,8) {\lb};
\node at (2,7) {\lb};
\node at (4,8) {\lb};
\node at (6,9) {\lb};
\node at (6,10) {\lb};
\node at (6,11) {\lb};
\node at (10,8) {\lb};
\node at (11,9) {\lb};
\node at (12,9) {\lb};
\node at (12,8) {\lb};
\node at (14,9) {\lb};
\node at (16,10) {\lb};
\node at (16,11) {\lb};
\node at (20,8) {\lb};
\node at (21,9) {\lb};
\node at (22,9) {\lb};
\node at (24,10) {\lb};
\node at (26,11) {\lb};
\node at (30,8) {\lb};
\node at (31,9) {\lb};
\node at (34,11) {\lb};
\node at (0,15) {\lb};
\node at (1,16) {\lb};
\node at (2,17) {\lb};
\node at (2,16) {\lb};
\node at (2,15) {\lb};
\node at (2,14) {\lb};
\node at (4,15) {\lb};
\node at (6,16) {\lb};
\node at (6,17) {\lb};
\node at (6,18) {\lb};
\node at (10,15) {\lb};
\node at (11,16) {\lb};
\node at (12,17) {\lb};
\node at (12,16) {\lb};
\node at (12,15) {\lb};
\node at (14,16) {\lb};
\node at (16,17) {\lb};
\node at (16,18) {\lb};
\node at (20,15) {\lb};
\node at (21,16) {\lb};
\node at (22,17) {\lb};
\node at (22,16) {\lb};
\node at (24,17) {\lb};
\node at (26,18) {\lb};
\node at (30,15) {\lb};
\node at (31,16) {\lb};
\node at (32,17) {\lb};
\node at (34,18) {\lb};
\draw (2,0) circle (.3);
\draw (12,1) circle (.3);
\draw (22,2) circle (.3);
\draw (2,7) circle (.3);
\draw (12,8) circle (.3);
\draw (22,9) circle (.3);
\draw (32,3) circle (.3);
\draw (32,10) circle (.3);
\draw (2,14) circle (.3);
\draw (12,15) circle (.3);
\draw (22,16) circle (.3);
\draw (32,17) circle (.3);
\node at (5,0) {$K_{4,2}$};
\node at (15,1) {$K_{3,2}$};
\node at (25,2) {$K_{2,2}$};
\node at (33.5,2) {$K_{1,2}$};
\node at (5,7) {$K_{4,1}$};
\node at (15,8) {$K_{3,1}$};
\node at (25,9) {$K_{2,1}$};
\node at (33.5,9) {$K_{1,1}$};
\node at (5,14) {$K_{4,0}$};
\node at (15,15) {$K_{3,0}$};
\node at (25,16) {$K_{2,0}$};
\node at (33,16) {$K_{1,0}$};

\end{\tz}
\end{center}
\end{fig}
\end{minipage}
\bigskip

Now we describe the remaining edges.

\begin{thm} \label{edges} For $0\le b\le3$ define $f_b$ by
$$f_b(0)=\begin{cases}1&b=0\\ 0&b=1,2,3,\end{cases}$$
and for $i\ge1$ and any $j$,
$$f_b(i)=\begin{cases}0&2^j\le i\le 2^j+[(b +j-2)/4]\\
1&i=2^j+(b+j-1)/4,\ b+j\equiv 1\ (4)\\
2&2^j+[(b +j-1)/4]+1\le i\le 2^{j+1}-1.\end{cases}$$
 For $2\le e\le k-1$, let $e=4a+b$ with $0\le b\le 3$. Then $$\E_{e,k}=\bigoplus_{i=0}^{2^{k-1-e}-1}K_{k-\lg(i)-1-e,f_b(i)}(8(i+a)+4,3-k+4(i+a)),$$
with $\lg(0)=-1$ and  the following modifications:
\begin{enumerate} 
\item If $b=3$ and $K_{t,\eps}(x,y)$ is the summand for $i$ a $2$-power, then the element $g'$ in   $(x-2,y+t-3)$ is replaced by an element in $(x-2,y+t-2)$;
\item If $b=0$ or $1$ and $K_{t,\eps}(x,y)$ is the summand for $i+1$ a $2$-power, and $t>1$, then the top element in grading $x+4$ is killed; i.e., $2^{t-2}g'''=0$;
\item If $b=0$ or $1$ and $i=2^{k-1-e}-1$, the element in grading $x+2$ in this $K_{1,\eps}(x,y)$ is killed.
\item If $b=3$, there is an additional element in $(2^{k+2-e}+8a+2,2^{k+1-e}-k+4a+1)$, which would be the first element if the string of $K_{1,2}$'s at the end was extended one farther.
\end{enumerate}
\end{thm}

In Table \ref{tbl}, we list a sample of values of $f_b(i)$. If $i<256$ is not included in the table, then $f_b(i)=2$ for all $b$.

\begin{table}[h]
\caption{Values of $f_b(i)$}
\label{tbl}

\begin{tabular}{c|ccccccccccccccccc}
$i$&0&1&2&4&5&8&9&16&17&32&33&64&65&66&128&129&130\\
\hline
$f_0(i)$&1&2&1&0&2&0&2&0&2&0&1&0&0&2&0&0&2\\
$f_1(i)$&0&1&0&0&2&0&2&0&1&0&0&0&0&2&0&0&2\\
$f_2(i)$&0&0&0&0&2&0&1&0&0&0&0&0&0&2&0&0&1\\
$f_3(i)$&0&0&0&0&1&0&0&0&0&0&0&0&0&1&0&0&0
\end{tabular}
\end{table}

For example, in Figure \ref{A6}, you can see $\Si^{32}\E_{2,5}$ with
\begin{equation}\label{E25}\E_{2,5}=K_{3,0}(4,-2)\oplus K_{2,0}(12,2)\oplus K_{1,0}(20,6)\oplus K_{1,2}(28,10),\end{equation}
and $\Si^{48}\E_{3,5}$ and $\Si^{80}\E_{3,5}$ with 
$$\E_{3,5}=K_{2,0}(4,-2)\oplus K'_{1,0}(12,2)\oplus\zt(18,3),$$ 
where $K'$ incorporates modification (1), and $\zt(18,3)$ is modification (3).
\section{Outline of proof}\label{outline}
In this section, we outline the proof, which occupies the seven sections which follow.

In \cite{DW}, we obtained the $E_2$ page of the ASS converging to $ku^*(K_2)$ and used a comparison with results about $k(1)^*(K_2)$ obtained in \cite{DRW} to obtain formulas for differentials in this ASS. In Section \ref{summandsec} here, we show that the result of \cite{DW}
{\em could have been grouped} into summands $A_k$ and $B_{k,\ell}$, multiplied by various coefficients, rather than the $A_k$, $B_k$, and $S_{k,\ell}$ used there. By this, we mean that the cohomology classes underlying these summands fill out $H^*(K_2)$ (Theorem \ref{everything}) and (Theorem \ref{closed}) they are closed under the differentials in the ASS obtained in \cite[Theorem 3.1]{DW}. 

These summands $A_k$ and $B_{k,\ell}$ have analogues for the ASS converging to $ku_*(K_2)$, and in Section \ref{kudifflsec}, we show how the duality between $ku_*(K_2)$ and $ku^*(K_2)$ obtained in \cite{DG} can be used to obtain differentials in the ASS converging to $ku_*(K_2)$. These will be used  in Section \ref{pfsec}.

In \cite{DW2}, we obtained the $ko$ analogues, $\M_k$, of the basic summands $M_k$ used in the $ku$ work in \cite{DW}. A complication for the $ko$ work is that whereas  coefficients $z_j$ for the $M_k$'s in $ku$ just resulted in suspensions, the $ko$ analogue causes a change in form, resulting in charts $M_k^i$ if there are $i$ $z_j$'s. In Section \ref{kosummandsec}, we describe the way in which the $ko_*$ analogues $\A_k$ of $A_k$, and $z^i\B_{k,\ell}$ of $B_{k,\ell}$  (if there are $i$ $z_j$'s) are built from a sequence of charts $M^a_b$ with differentials between them. The results of Section \ref{summandsec} imply that these are closed under differentials and fill out $ko_*(K_2)$.

In Section \ref{A4sec}, we use a small example, $\A_4$, to illustrate how the differentials combine to yield a nice picture for this summand of $ko_*(K_2)$. We also explain how the initial part of $\A_k$, which differs from everything else, is handled.

Our favored description of $\A_k$  in Section \ref{Akthmsec} is in terms of pre-edges and subedges. There is a nice pattern of differentials from the subedges to the pre-edges, turning pre-edges into edges. In Section \ref{derivationsec} we show how the $2^{k-1}$ charts $M_b^a$ which form $\A_k$ work together to form the pre-edges and subedges.  The upper edge, which involves the initial part considered in Section \ref{A4sec}, is slightly different than the others, and is discussed at the end of Section \ref{derivationsec}.

In Subsection \ref{131}, we describe how the summands that build $\At_k$ work together in the spectral sequence. The proof that the differentials are as claimed is given in Subsection \ref{132}, by comparison with the $ku_*$ differentials obtained in Sections \ref{kudifflsec} and \ref{summandsec}. In Subsection \ref{133}, we explain how the nice pattern of exotic extensions in the $ko_*$ ASS is established, using comparison with $ku_*$ extensions, Toda brackets, and Adams periodicity.

In Section \ref{Bklpfsec}, we discuss modifications required in the formation of $z^i\B_{k,\ell}$.

\section{Using the $ku^*$ differentials}\label{kudifflsec}
In this section, we make our first, very preliminary, step toward a proof. Our input is formulas
(\cite[Theorem 3.1]{DW}) for differentials in the ASS converging to $ku^*(K_2)$, which were derived in \cite{DW} by complicated ad hoc methods. In the final section of \cite{DW}, almost as an afterthought, we sketched an approach to $ku^*(K_2)$ and $ku_*(K_2)$ which involved summands in a splitting of $H^*(K_2)$\footnote{Coefficients of cohomology groups are $\zt$ unless noted to the contrary.} as a module over a subalgebra of the mod-2 Steenrod algebra. This summand approach will be the way that we compute $ko_*(K_2)$.

In this section, we focus on the summand $A_4$ of $ku^*(K_2)$ and show specifically how the formulas for the differentials in the ASS of $ku^*(K_2)$ can be applied to the summand approach to $ku^*(K_2)$ and then to $ku_*(K_2)$. There is an exterior subalgebra $E_1$ of the mod 2 Steenrod algebra with generators $Q_0$ and $Q_1$  of grading 1 and 3, respectively, such that the $E_2$ page of the ASS converging to $ku^*(K_2)$ is $\ext_{E_1}(\zt,H^*(K_2))$. We depict the spectral sequence using $ku^*$ gradings increasing from right to left, as was done in \cite{DW}. We draw $E_1$-modules using straight lines for the action of $Q_0$ and curved lines for the action of $Q_1$.

In \cite{DW}, we introduced classes $z_i\in H^{2^{i+2}+2}(K_2)$, $y_1\in H^4(K_2)$, and $q\in H^9(K_2)$, and $E_1$-submodules $M_k$ for $k\ge4$, which we picture in Figure \ref{modules} for $4\le k\le 6$, indicating the grading of the classes.

\bigskip
\begin{minipage}{6in}
\begin{fig}\label{modules}

{\bf  $E_1$-modules $M_k$.}

\begin{center}

\begin{\tz}[scale=.35]
\draw (0,0) -- (2,0);
\node at (0,0) {\lb};
\node at (2,0) {\lb};
\node at (0,-.7) {$17$};
\node at (2,-.7) {$18$};
\draw (10,0) -- (12,0);
\draw (14,0) -- (16,0);

\draw (10,0) to[out=45, in=135] (16,0);

\node at (10,-.7) {$33$};
\node at (14,-.7) {$35$};
\node at (12,-.7) {$34$};
\node at (16,-.7) {$36$};

\node at (10,0) {\lb};
\node at (12,0) {\lb};
\node at (14,0) {\lb};
\node at (16,0) {\lb};

\draw (20,0) -- (22,0);
\draw (24,0) -- (26,0);
\draw (28,0) -- (30,0);

\draw (20,0) to[out=45, in=135] (26,0);
\draw (24,0) to[out=315, in=225] (30,0);

\node at (20,-.7) {$65$};
\node at (24,-1) {$67$};
\node at (22,-.7) {$66$};
\node at (28,.8) {$69$};
\node at (26,.8) {$68$};
\node at (30,-1) {$70$};

\node at (20,0) {\lb};
\node at (22,0) {\lb};
\node at (24,0) {\lb};
\node at (26,0) {\lb};
\node at (28,0) {\lb};
\node at (30,0) {\lb};
\node at (1,-1.7) {$M_4$};
\node at (13,-1.7) {$M_5$};
\node at (25,-1.7) {$M_6$};
\end{\tz}
\end{center}
\end{fig}
\end{minipage}

\bigskip
The classes in Figure \ref{modules} with grading 18, 34, 36, 66, 68, and 70 are 
\begin{equation}\label{zs}z_2,\  z_3,\  z_2^2,\  z_4,\  z_3^2,\text{ and }z_2^2z_3,\end{equation}
respectively. Pictures of $\ext_{E_1}(\zt,M_k)$ are given in Figure \ref{Mchts}, where vertical lines indicate multiplication by $h_0$, and diagonal lines multiplication by $v\in ku^{-2}$.

\bigskip
\begin{minipage}{6in}
\begin{fig}\label{Mchts}

{\bf  $\ext_{E_1}(\zt,M_k)$.}

\begin{center}

\begin{\tz}[scale=.35]
\draw (-.5,0) -- (6,0);
\draw (9.5,0) -- (18,0);
\draw (23.5,0) -- (31,0);
\draw [->] (0,0) -- (6,3);
\draw [->] (10,0) -- (17,3.5);
\draw [->] (24,0) -- (32,4);
\draw [->] (12,0) -- (17,2.5);
\draw [->] (26,0) -- (32,3);
\draw [->] (28,0) -- (32,2);
\draw (12,0) -- (12,1);
\draw (14,1) -- (14,2);
\draw (16,2) -- (16,3);
\draw (26,0) -- (26,1);
\draw (28,0) -- (28,2);
\draw (30,1) -- (30,3);
\node at (0,0) {\lb};
\node at (2,1) {\lb};
\node at (4,2) {\lb};
\node at (10,0) {\lb};
\node at (12,1) {\lb};
\node at (14,2) {\lb};
\node at (16,3) {\lb};
\node at (12,0) {\lb};
\node at (14,1) {\lb};
\node at (16,2) {\lb};
\node at (24,0) {\lb};
\node at (26,1) {\lb};
\node at (28,2) {\lb};
\node at (30,3) {\lb};
\node at (26,0) {\lb};
\node at (28,1) {\lb};
\node at (30,2) {\lb};
\node at (28,0) {\lb};
\node at (30,1) {\lb};
\node at (0,-.7) {$18$};
\node at (2,-.7) {$16$};
\node at (10,-.7) {$36$};
\node at (12,-.7) {$34$};
\node at (14,-.7) {$32$};
\node at (24,-.7) {$70$};
\node at (26,-.7) {$68$};
\node at (28,-.7) {$66$};
\node at (2,-1.7) {$k=4$};
\node at (13,-1.7) {$k=5$};
\node at (27,-1.7) {$k=6$};
\end{\tz}
\end{center}
\end{fig}
\end{minipage}
\bigskip

\ni We call these $v$-towers. Their generators are the classes listed in (\ref{zs}).

Analogous to \cite[Figure 17]{DW} is the list in Figure \ref{71} of $E_1$-submodules which combine to form $A_4\subset ku^*(K_2)$. We list them in  order of decreasing grading to correspond to the order in the $ku^*(K_2)$ chart.

\bigskip
\begin{minipage}{6in}
\begin{fig}\label{71}

{\bf  $E_1$-modules building $A_4$.}

\begin{center}

\begin{\tz}[scale=.42]
\draw (0,0) -- (1,0);
\draw (2,0) -- (3,0);
\draw (4,0) -- (5,0);
\draw (5,-1) -- (6,-1);
\draw (10,0) -- (11,0);
\draw (13,0) -- (14,0);
\draw (15,0) -- (16,0);
\draw (18,0) -- (19,0);
\draw (20,0) -- (21,0);
\draw (23,0) -- (24,0);
\draw (28,0) -- (29,0);
\draw (32,0) -- (33,0);
\draw (34,0) -- (35,0);

\draw (0,0) to[out=25, in=155] (3,0);
\draw (2,0) to[out=335, in=205] (5,0);
\draw (32,0) to[out=25, in=155] (35,0);
\draw (13,0) to[out=335, in=205] (16,0);
\draw (18,0) to[out=335, in=205] (21,0);
\draw (34,0) to[out=335, in=205] (37,0);
\node at (2,1.75) {$M_6$};
\node at (5.5,1.75) {$y_1qz_3M_4$};
\node at (10.5,1.75) {$y_1^2z_3M_4$};
\node at (14.5,1.75) {$y_1^3qM_5$};
\node at (19.5,1.75) {$y_1^4M_5$};
\node at (23.5,1.75) {$y_1^5qM_4$};
\node at (28.5,1.75) {$y_1^6M_4$};
\node at (35.5,1.75) {$y_1^7N$};
\node at (38,1.75) {$y_1^8$};
\node at (0,.7) {$70$};
\node at (4,.7) {$66$};
\node at (6,-.2) {$64$};
\node at (10,.7) {$60$};
\node at (0,0) {\lb};
\node at (1,0) {\lb};
\node at (2,0) {\lb};
\node at (3,0) {\lb};
\node at (4,0) {\lb};
\node at (5,0) {\lb};
\node at (5,-1) {\lb};
\node at (6,-1) {\lb};
\node at (10,0) {\lb};
\node at (11,0) {\lb};
\node at (13,0) {\lb};
\node at (14,0) {\lb};
\node at (15,0) {\lb};
\node at (16,0) {\lb};
\node at (18,0) {\lb};
\node at (19,0) {\lb};
\node at (20,0) {\lb};
\node at (21,0) {\lb};
\node at (23,0) {\lb};
\node at (24,0) {\lb};
\node at (28,0) {\lb};
\node at (29,0) {\lb};
\node at (32,0) {\lb};
\node at (33,0) {\lb};
\node at (34,0) {\lb};
\node at (35,0) {\lb};
\node at (37,0) {\lb};
\node at (38,0) {\lb};
\node at (13,.7) {$57$};
\node at (18,.7) {$52$};
\node at (23,.7) {$47$};
\node at (28,.7) {$42$};
\node at (38,.7) {$32$};
\node at (35,.7) {$35$};
\node at (32,.7) {$38$};
\end{\tz}
\end{center}
\end{fig}
\end{minipage}
\bigskip

Here $N$ is an $E_1$-module pictured in Figure \ref{N}; $\ext_{E_1}(\zt,N)$ is also in Figure \ref{N}. The cohomology class in grading 9 is $q$, and, as shown in \cite[Figures 8 and 9]{DW}, the generator of the first infinite tower in Ext
corresponds to $v^2q$.

\bigskip
\begin{minipage}{6in}
\begin{fig}\label{N}

{\bf  $N$ and $\ext_{E_1}(\zt,N)$.}

\begin{center}

\begin{\tz}[scale=.5]
\node at (2,0) {\lb};
\node at (6,0) {\lb};
\node at (8,0) {\lb};
\node at (10,0) {\lb};
\node at (12,0) {\lb};
\draw (2,0) to[out=25, in=155] (8,0);
\draw (6,0) to[out=335, in=205] (12,0);
\node at (2,-.7) {$5$};
\node at (10,.7) {$9$};
\draw (19.5,-2) -- (28,-2);
\node at (20,-2) {\lb};
\node at (22,-2) {\lb};
\node at (22,-1) {\lb};
\node at (25,0) {\lb};
\node at (25,1) {\lb};
\node at (25,2) {\lb};
\node at (27,1) {\lb};
\node at (27,2) {\lb};
\draw [->] (25,0) -- (28,1.5);
\draw [->] (25,1) -- (28,2.5);
\draw (20,-2) -- (22,-1) -- (22,-2);
\draw [->] (25,0) -- (25,4);
\draw [->] (27,1) -- (27,4);
\node at (20,-2.7) {$10$};
\node at (22,-2.7) {$8$};
\node at (25,-2.7) {$5$};
\node at (29,4) {$\iddots$};
\draw (6,0) -- (8,0);
\draw (10,0) -- (12,0);
\end{\tz}
\end{center}
\end{fig}
\end{minipage}
\bigskip

We extract  from \cite[Theorem 3.1]{DW} the formulas for differentials in the ASS converging to $ku^*(K_2)$ relevant to this approach to $A_4$. Here $\L_{j+1}$ is an exterior algebra on classes $z_i$ with $i\ge j+1$, and $\ot\L_{j+1}$ means that the formula can be multiplied by any monomial in $\L_{j+1}$.
\begin{eqnarray*}d^5(y_1^8)&=&h_0^3v^2qy_1^7\\
d^2(qy_1^{2a+1}z_j)&=&v^2y_1^{2a}z_2z_j\qquad\ot\L_{j+1},\ j\ge2\\
d^5(qy_1^{4a+3}z_{j-1}^2)&=&v^5y_1^{4a}z_3z_j\qquad\ot\L_{j+1},\ j\ge3\\
d^{12}(qy_1^{8a+7}z_{j-2}^2z_{j-1})&=&v^{12}y_1^{8a}z_4z_j\qquad\ot\L_{j+1},\ j\ge4\\
d^2(v^2qy_1^{2a+1})&=&v^4y_1^{2a}z_2\\
d^5(h_0v^2qy_1^{4a+3})&=&v^8y_1^{4a}z_3\\
d^{12}(h_0^2v^2qy_1^{8a+7})&=&v^{16}y_1^{8a}z_4.\end{eqnarray*}

After the first formula is applied to $y_1^7N\oplus \langle y_1^8\rangle$, the surviving $v$-towers are $h_0^tv^2qy_1^7$ for $0\le t\le2$. The other formulas apply to  these $v$-towers and all those in the other summands in Figure \ref{71}. We show these in Figure \ref{alldiff}.

\bigskip
\begin{minipage}{6in}
\begin{fig}\label{alldiff}

{\bf  Differentials leading to $A_4$.}

\begin{center}

\begin{\tz}[scale=.36]
\node at (2,13) {$M_6$};
\node at (8,13) {$y_1qz_3M_4$};
\node at (14,13) {$y_1^2z_3M_4$};
\node at (20,13) {$y_1^3qM_5$};
\node at (26,13) {$y_1^4M_5$};
\node at (32,13) {$y_1^5qM_4$};
\node at (38,13) {$y_1^6M_4$};
\node at (44,13) {$y_1^7N$};
\node at (44,4) {$h_0^2y_1^7v^2q$};
\node at (44,2) {$h_0y_1^7v^2q$};
\node at (44,0) {$y_1^7v^2q$};
\node at (38,0) {$v^4y_1^6z_2$};
\draw [->] (42,0) -- (40,0);
\node at (41,.6) {$d^2$};
\node at (26,2) {$v^8y_1^4z_3$};
\draw [->] (42,2) -- (28,2);
\node at (35,2.6) {$d^5$};
\draw (0,12) -- (45,12);
\node at (2,10) {$v^2z_2^2z_3$};
\node at (8,10) {$y_1qz_3z_2$};
\draw [->] (6,10) -- (4,10);
\node at (5,10.6) {$d^2$};
\node at (26,10) {$v^2y_1^4z_2^2$};
\node at (32,10) {$y_1^5qz_2$};
\draw [->] (30,10) -- (28,10);
\node at (29,10.6) {$d^2$};
\node at (2,8) {$v^5z_3^2$};
\node at (20,8) {$y_1^3qz_2^2$};
\draw [->] (18,8) -- (4,8);
\node at (11,8.6) {$d^5$};
\node at (14,6) {$v^2y_1^2z_3z_2$};
\node at (20,6) {$y_1^3qz_3$};
\draw [->] (18,6) -- (16.4,6);
\node at (17,6.6) {$d^2$};
\node at (2,4) {$v^{16}z_4$};
\node at (44,4) {$h_0^2y_1^7v^2q$};
\draw [->] (42,4) -- (4,4);
\node at (26,4.6) {$d^{12}$};
\end{\tz}
\end{center}
\end{fig}
\end{minipage}
\bigskip

After these differentials, what remains in $A_4$  are truncated $v$-towers of the targets of the differentials in Figure \ref{alldiff}. We depict this in Figure \ref{kuA4}. The classes in the lower right corner come from Figure \ref{N} after multiplying by $y_1^7$. This chart agrees with the $A_4$ part of \cite[Figure 1]{DW} except that we do not address the exotic extensions here.

\bigskip
\begin{minipage}{6in}
\begin{fig}\label{kuA4}

{\bf  $v$-towers in $A_4$.}

\begin{center}

\begin{\tz}[scale=.45]
\draw (-.5,0) -- (34.5,0);
\draw (0,0) -- (2,1) -- (2,0) -- (10,4) -- (10,3);
\draw (4,1) -- (4,0) -- (34,15);
\draw (6,2) -- (6,1) -- (8,2) -- (8,3);
\draw (10,0) -- (12,1);
\draw (18,0) -- (20,1) -- (20,0) -- (34,7);
\draw (28,0) -- (34,3);
\draw (32,0) -- (34,1) -- (34,0);
\node at (0,-.7) {$70$};
\node at (4,-.7) {$66$};
\node at (10,-.7) {$60$};
\node at (18,-.7) {$52$};
\node at (28,-.7) {$42$};
\node at (34,-.7) {$36$};
\node at (0,0) {\lb};
\node at (2,0) {\lb};
\node at (2,1) {\lb};
\node at (4,0) {\lb};
\node at (4,1) {\lb};
\node at (6,1) {\lb};
\node at (6,2) {\lb};
\node at (8,2) {\lb};
\node at (8,3) {\lb};
\node at (10,0) {\lb};
\node at (10,3) {\lb};
\node at (10,4) {\lb};
\node at (12,4) {\lb};
\node at (14,5) {\lb};
\node at (16,6) {\lb};
\node at (18,7) {\lb};
\node at (20,8) {\lb};
\node at (22,9) {\lb};
\node at (24,10) {\lb};
\node at (26,11) {\lb};
\node at (28,12) {\lb};
\node at (30,13) {\lb};
\node at (32,14) {\lb};
\node at (34,15) {\lb};
\node at (12,1) {\lb};
\node at (18,0) {\lb};
\node at (20,1) {\lb};
\node at (20,0) {\lb};
\node at (22,1) {\lb};
\node at (24,2) {\lb};
\node at (26,3) {\lb};
\node at (28,4) {\lb};
\node at (30,5) {\lb};
\node at (32,6) {\lb};
\node at (34,7) {\lb};
\node at (28,0) {\lb};
\node at (30,1) {\lb};
\node at (32,2) {\lb};
\node at (34,3) {\lb};
\node at (32,0) {\lb};
\node at (34,0) {\lb};
\node at (34,1) {\lb};

\end{\tz}
\end{center}
\end{fig}
\end{minipage}
\bigskip

For $ku_*(K_2)$, the $E_2$-page is $\ext_{E_1}^{*,*}(H^*K_2,\zt)$. For its $A_4$, called $\as_4$, we use the summands of Figure \ref{71} but arrange them in the opposite order. The Ext groups corresponding to Figures \ref{Mchts} and \ref{N} are presented in Figure \ref{4chts}. See \cite[Figure 4.2]{DWSW}
for $N$.

\bigskip
\begin{minipage}{6in}
\begin{fig}\label{4chts}

{\bf  $\ext_{E_1}(M_k,\zt)$ and $\ext_{E_1}(N,\zt)$.}

\begin{center}

\begin{\tz}[scale=.32]
\draw (-.5,0) -- (6,0);
\draw (9.5,0) -- (18,0);
\draw (19.5,0) -- (27,0);
\draw [->] (0,0) -- (6,3);
\draw [->] (10,0) -- (17,3.5);
\draw [->] (20,0) -- (28,4);
\draw [->] (12,0) -- (17,2.5);
\draw [->] (22,0) -- (28,3);
\draw [->] (24,0) -- (28,2);
\draw (12,0) -- (12,1);
\draw (14,1) -- (14,2);
\draw (16,2) -- (16,3);
\draw (22,0) -- (22,1);
\draw (24,0) -- (24,2);
\draw (26,1) -- (26,3);
\node at (0,0) {\lb};
\node at (2,1) {\lb};
\node at (4,2) {\lb};
\node at (10,0) {\lb};
\node at (12,1) {\lb};
\node at (14,2) {\lb};
\node at (16,3) {\lb};
\node at (12,0) {\lb};
\node at (14,1) {\lb};
\node at (16,2) {\lb};
\node at (20,0) {\lb};
\node at (22,1) {\lb};
\node at (24,2) {\lb};
\node at (26,3) {\lb};
\node at (22,0) {\lb};
\node at (24,1) {\lb};
\node at (26,2) {\lb};
\node at (24,0) {\lb};
\node at (26,1) {\lb};
\node at (0,-.7) {$17$};
\node at (2,-.7) {$19$};
\node at (10,-.7) {$33$};
\node at (12,-.7) {$35$};
\node at (14,-.7) {$37$};
\node at (20,-.7) {$65$};
\node at (22,-.7) {$67$};
\node at (24,-.7) {$69$};
\node at (2,-1.7) {$k=4$};
\node at (13,-1.7) {$k=5$};
\node at (23,-1.7) {$k=6$};
\draw (31.5,0) -- (40.5,0);
\draw [->] (32,0) -- (32,4);
\draw [->] (34,0) -- (34,4);
\draw [->] (36,0) -- (36,4);
\draw [->] (38,1) -- (38,4);
\draw [->] (40,2) -- (40,4);
\draw [->] (36,0) -- (41,2.5);
\draw [->] (34,0) -- (41,3.5);
\draw [->] (32,0) -- (39,3.5);
\draw [->] (32,1) -- (37,3.5);
\draw [->] (32,2) -- (35,3.5);
\draw [->] (32,3) -- (33,3.5);

\node at (32,0) {\lb};
\node at (32,1) {\lb};
\node at (32,2) {\lb};
\node at (32,3) {\lb};
\node at (34,0) {\lb};
\node at (34,1) {\lb};
\node at (34,2) {\lb};
\node at (34,3) {\lb};
\node at (36,0) {\lb};
\node at (36,1) {\lb};
\node at (36,2) {\lb};
\node at (36,3) {\lb};
\node at (38,1) {\lb};
\node at (38,2) {\lb};
\node at (38,3) {\lb};
\node at (40,2) {\lb};
\node at (40,3) {\lb};
\node at (42,3) {$\iddots$};
\node at (32,-.7) {$5$};
\node at (36,-.7) {$9$};
\node at (40,-.7) {$13$};
\node at (38,-1.7) {$N$};

\end{\tz}
\end{center}
\end{fig}
\end{minipage}
\bigskip

We dualize the $ku^*$ differentials, similarly to \cite{DRW}. In Figure \ref{dual}, we show the dualization of a $d^2$ and $d^5$ differential from \ref{alldiff}.

\bigskip
\begin{minipage}{6in}
\begin{fig}\label{dual}

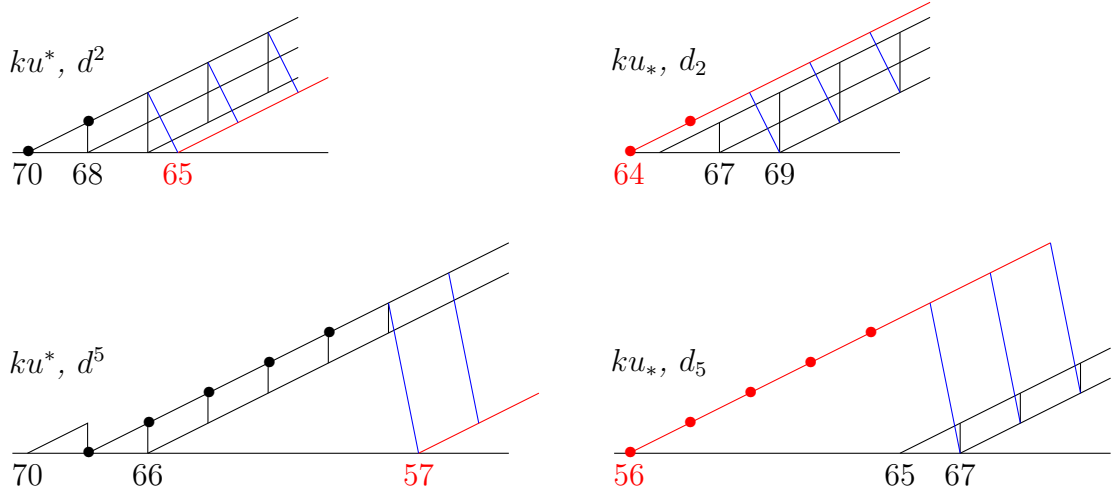
{\bf  Dualizing $d^2$ and $d^5$ differentials.}

\begin{center}

\begin{\tz}[scale=.4]
\draw (-.5,10) -- (10,10);
\draw (0,10) -- (9,14.5);
\draw (20,10) -- (29,10);
\draw (-.5,0) -- (16,0);
\draw (2,0) -- (16,7);
\draw (19.5,0) -- (36,0);
\draw [red] (20,0) -- (34,7);
\draw (2,10) -- (9,13.5);
\draw (4,10) -- (9,12.5);
\draw [red] (5,10) -- (10,12.5);
\draw (2,10) -- (2,11);
\draw (4,10) -- (4,12);
\draw (6,11) -- (6,13);
\draw (8,12) -- (8,14);
\draw [blue] (5,10) -- (4,12);
\draw [blue] (7,11) -- (6,13);
\draw [blue] (9,12) -- (8,14);
\draw [red] (20,10) -- (30,15);
\draw (21,10) -- (30,14.5);
\draw (23,10) -- (30,13.5);
\draw (25,10) -- (30,12.5);
\draw (23,10) -- (23,11);
\draw (25,10) -- (25,12);
\draw (27,11) -- (27,13);
\draw (29,12) -- (29,14);
\draw [blue] (25,10) -- (24,12);
\draw [blue] (27,11) -- (26,13);
\draw [blue] (29,12) -- (28,14);
\draw (0,0) -- (2,1) -- (2,0);
\draw (4,0) -- (16,6);
\draw (4,0) -- (4,1);
\draw (6,1) -- (6,2);
\draw (8,2) -- (8,3);
\draw (10,3) -- (10,4);
\draw (12,4) -- (12,5);
\draw [red] (13,0) -- (17,2);
\draw [blue] (13,0) -- (12,5);
\draw [blue] (15,1) -- (14,6);
\node at (0,9.3) {$70$};
\node [red] at (5,9.3) {$65$};
\node at (2,9.3) {$68$};
\node [red] at (20,9.3) {$64$};
\node at (23,9.3) {$67$};
\node at (25,9.3) {$69$};
\node at (1,13) {$ku^*$, $d^2$};
\node at (21,13) {$ku_*$, $d_2$};
\node at (0,10) {\lb};
\node at (2,11) {\lb};
\node [red] at (20,10) {\lb};
\node [red] at (22,11) {\lb};
\node at (2,0) {\lb};
\node at (4,1) {\lb};
\node at (6,2) {\lb};
\node at (8,3) {\lb};
\node at (10,4) {\lb};
\node at (0,-.7) {$70$};
\node at (4,-.7) {$66$};
\node [red] at (13,-.7) {$57$};
\node at (1,3) {$ku^*$, $d^5$};
\node at (21,3) {$ku_*$, $d_5$};
\node [red] at (20,0) {\lb};
\node [red] at (22,1) {\lb};
\node [red] at (24,2) {\lb};
\node [red] at (26,3) {\lb};
\node [red] at (28,4) {\lb};
\draw (29,0) -- (36,3.5);
\draw (31,0) -- (36,2.5);
\draw (31,0) -- (31,1);
\draw (33,1) -- (33,2);
\draw (35,2) -- (35,3);
\draw [blue] (31,0) -- (30,5);
\draw [blue] (33,1) -- (32,6);
\draw [blue] (35,2) -- (34,7);
\node [red] at (20,-.7) {$56$};
\node at (29,-.7) {$65$};
\node at (31,-.7) {$67$};
\end{\tz}
\end{center}
\end{fig}
\end{minipage}
\bigskip

The justification for this duality of differentials is the duality of groups,
\begin{equation}\label{kuduality}ku_*(K_2)\approx (ku^{*+4}(K_2))^\vee,\end{equation}
which was proved in \cite[Example 3.4]{DG} and restated in \cite[Theorem 1.20]{DW}. Here $M^\vee$ is the Pontryagin dual of $M$. The Pontryagin dual of the two $ku^*$ charts in Figure \ref{dual} have $v$-towers of height 2 and 5 on classes of grading 68 and 60, respectively, corresponding to $v$-towers on classes of grading 64 and 56, respectively, in the $ku_*$ charts. Note that in the Pontryagin dual, the action of $\cdot2$ and $v$ appears backwards from its usual appearance.

There is a $d_5$-differential from $\ext_{E_1}(y_1^7N,\zt)$ to $\ext_{E_1}(\langle y_1^8\rangle,\zt)$ pictured in Figure \ref{d5N}. This and the first differential in the above list of seven differential formulas in $ku^*(K_2)$ are derived from \cite{Bro} as  discussed in \cite[Section 4]{DWSW} and \cite[Section 3]{DW}. This truncates all the infinite vertical towers, leaving the classes indicated by dots.

\bigskip
\begin{minipage}{6in}
\begin{fig}\label{d5N}

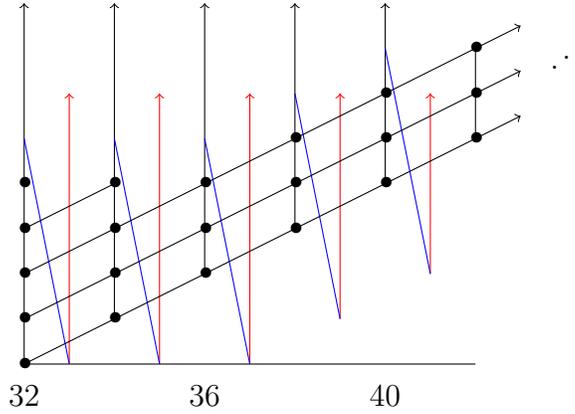
{\bf  Differential killing infinite towers in $\as_4$.}

\begin{center}

\begin{\tz}[scale=.6]
\draw (0,0) -- (10,0);
\draw [->] (0,0) -- (0,8);
\draw [->] (2,1) -- (2,8);
\draw [->] (4,2) -- (4,8);
\draw [->] (6,3) -- (6,8);
\draw [->] (8,4) -- (8,8);
\draw [->] [red] (1,0) -- (1,6);
\draw [->] [red] (3,0) -- (3,6);
\draw [->] [red] (5,0) -- (5,6);
\draw [->] [red] (7,1) -- (7,6);
\draw [->] [red] (9,2) -- (9,6);
\draw [blue] (1,0) -- (0,5);
\draw [blue] (3,0) -- (2,5);
\draw [blue] (5,0) -- (4,5);
\draw [blue] (7,1) -- (6,6);
\draw [blue] (9,2) -- (8,7);
\draw [->] (0,0) -- (11,5.5);
\draw [->] (0,1) -- (11,6.5);
\draw [->] (0,2) -- (11,7.5);
\draw (0,3) -- (2,4);
\node at (0,0) {\lb};
\node at (0,1) {\lb};
\node at (0,2) {\lb};
\node at (0,3) {\lb};
\node at (0,4) {\lb};
\node at (2,1) {\lb};
\node at (2,2) {\lb};
\node at (2,3) {\lb};
\node at (2,4) {\lb};
\node at (4,2) {\lb};
\node at (4,3) {\lb};
\node at (4,4) {\lb};
\node at (6,3) {\lb};
\node at (6,4) {\lb};
\node at (6,5) {\lb};
\node at (8,4) {\lb};
\node at (8,5) {\lb};
\node at (8,6) {\lb};
\node at (10,5) {\lb};
\node at (10,6) {\lb};
\node at (10,7) {\lb};
\node at (12,7) {$\iddots$};
\draw (10,5) -- (10,7);
\node at (0,-.7) {$32$};
\node at (4,-.7) {$36$};
\node at (8,-.7) {$40$};
\end{\tz}
\end{center}
\end{fig}
\end{minipage}
\bigskip

This leaves three infinite $v$-towers from the $\langle y_1^8\rangle\oplus y_1^7 N$ summand in the summand approach to determining  $\as_4$. In Figure \ref{homdiff} we dualize Figure \ref{alldiff}. We write the summands in the opposite order. The direction of the differentials is reversed. For example, the duals of the $d^2$ and $d^5$ differentials in the upper left part of Figure \ref{alldiff} are shown in Figure \ref{dual} and appear in the upper right corner of Figure \ref{homdiff}. The classes now have different names. What will be useful for us in passing to $ko_*$ is the type of differential between the summands. One difference is that the $d^2$, $d^5$, and $d^{12}$ differentials in the lower right part of Figure \ref{alldiff} become $d_6$, $d_9$, and $d_{16}$. We show this for the $d^2$ in Figure \ref{d2d6}. This is caused by the triangle of classes of height 4 in the $ku_*$ diagram.

\bigskip
\begin{minipage}{6in}
\begin{fig}\label{d2d6}

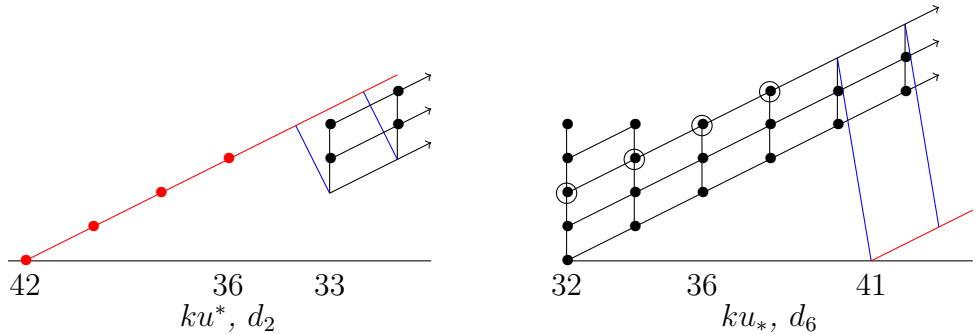
{\bf  Duality of a $d^2$ and $d_6$ differential.}

\begin{center}

\begin{\tz}[scale=.45]
\draw (-.5,0) -- (12,0);
\draw [red] (0,0) -- (11,5.5);
\draw (9,2) -- (9,4);
\draw (11,3) -- (11,5);
\draw [blue] (9,2) -- (8,4);
\draw [blue] (11,3) -- (10,5);
\node [red] at (0,0) {\lb};
\node [red] at (2,1) {\lb};
\node [red] at (4,2) {\lb};
\node [red] at (6,3) {\lb};
\node at (9,3) {\lb};
\node at (9,4) {\lb};
\node at (11,4) {\lb};
\node at (11,5) {\lb};
\draw (16,0) -- (28,0);
\draw (16,0) -- (16,4);
\draw (18,1) -- (18,4);
\draw (20,2) -- (20,4);
\draw (22,3) -- (22,5);
\draw (24,4) -- (24,6);
\draw (26,5) -- (26,7);
\draw [->] (9,2) -- (12,3.5);
\draw [->] (9,3) -- (12,4.5);
\draw [->] (9,4) -- (12,5.5);
\draw [->] (16,0) -- (27,5.5);
\draw [->] (16,1) -- (27,6.5);
\draw [->] (16,2) -- (27,7.5);
\draw (16,3) -- (18,4);
\draw [red] (25,0) -- (28,1.5);
\draw [blue] (25,0) -- (24,6);
\draw [blue] (27,1) -- (26,7);
\node at (0,-.7) {$42$};
\node at (6,-.7) {$36$};
\node at (9,-.7) {$33$};
\node at (6,-1.7) {$ku^*$, $d_2$};
\node at (22,-1.7) {$ku_*$, $d_6$};
\draw (22,5) circle (.3);
\draw (20,4) circle (.3);
\draw (18,3) circle (.3);
\draw (16,2) circle (.3);
\node at (16,-.7) {$32$};
\node at (20,-.7) {$36$};
\node at (25,-.7) {$41$};
\node at (16,0) {\lb};
\node at (16,1) {\lb};
\node at (16,2) {\lb};
\node at (16,3) {\lb};
\node at (16,4) {\lb};
\node at (18,1) {\lb};
\node at (18,2) {\lb};
\node at (18,3) {\lb};
\node at (18,4) {\lb};
\node at (20,2) {\lb};
\node at (20,3) {\lb};
\node at (20,4) {\lb};
\node at (22,3) {\lb};
\node at (22,4) {\lb};
\node at (22,5) {\lb};
\node at (24,4) {\lb};
\node at (24,5) {\lb};
\node at (26,5) {\lb};
\node at (26,6) {\lb};

\end{\tz}
\end{center}
\end{fig}
\end{minipage}
\bigskip

\ni The $v$-tower of height 4 in $ku^*(K_2)^\vee$ on a class of grading 36 corresponds to a $v$-tower of height 4 on a class of grading 32 in $ku_*(K_2)$.

In Figure \ref{homdiff}, we list the differentials for the summand approach to  $\as_4$. We write the gradings of the generators of the $v$ towers involved in the differentials. For example, the $ku_*$ differentials in Figure \ref{dual} are in the upper right corner of Figure \ref{homdiff}, while the $ku_*$ differential in Figure \ref{d2d6} is in the lower left.

\bigskip
\begin{minipage}{6in}
\begin{fig}\label{homdiff}

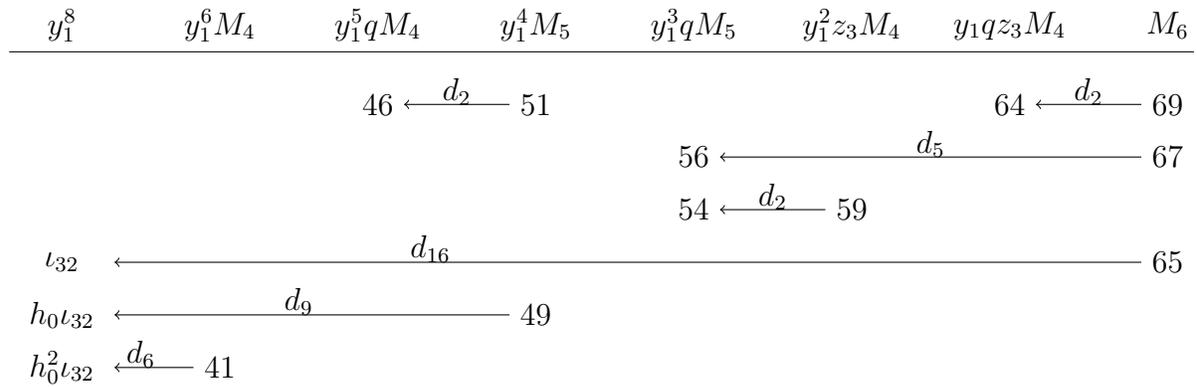
{\bf  Summand approach to  $\as_4$.}

\begin{center}

\begin{\tz}[scale=.35]
\node at (2,13) {$y_1^8$};
\node at (8,13) {$y_1^6M_4$};
\node at (14,13) {$y_1^5qM_4$};
\node at (20,13) {$y_1^4M_5$};
\node at (26,13) {$y_1^3qM_5$};
\node at (32,13) {$y_1^2z_3M_4$};
\node at (38,13) {$y_1qz_3M_4$};
\node at (44,13) {$M_6$};
\draw (0,12) -- (45,12);
\node at (14,10) {$46$};
\node at (20,10) {$51$};
\node at (38,10) {$64$};
\node at (44,10) {$69$};
\draw [->] (19,10) -- (15,10);
\draw [->] (43,10) -- (39,10);
\node at (17,10.5) {$d_2$};
\node at (41,10.5) {$d_2$};
\node at (26,8) {$56$};
\node at (44,8) {$67$};
\draw [->] (43,8) -- (27,8);
\node at (35,8.5) {$d_5$};
\node at (26,6) {$54$};
\node at (32,6) {$59$};
\draw [->] (31,6) -- (27,6);
\node at (29,6.5) {$d_2$};
\node at (44,4) {$65$};
\node at (2,4) {$\iota_{32}$};
\draw [->] (43,4) -- (4,4);
\node at (16,4.5) {$d_{16}$};
\node at (20,2) {$49$};
\node at (2,2) {$h_0\iota_{32}$};
\draw [->] (19,2) -- (4,2);
\node at (11,2.5) {$d_9$};
\node at (8,0) {$41$};
\node at (2,0) {$h_0^2\iota_{32}$};
\draw [->] (7,0) -- (4,0);
\node at (5,.5) {$d_6$};

\end{\tz}
\end{center}
\end{fig}
\end{minipage}
\bigskip

We close this section by displaying in Figure \ref{ku_*A4} the result of these differentials, the $A_4$ analogue of Figure \ref{kuA5} without the exotic extensions. 
For example, the $d_6$ differential to $h_0^2\io_{32}$ hits $v^4h_0^2\io_{32}$.

\bigskip
\begin{minipage}{6in}
\begin{fig}\label{ku_*A4}

{\bf  $\as_4$ without exotic extensions.}

\begin{center}

\begin{\tz}[scale=.4]
\draw (31.5,0) -- (66,0);
\draw (32,0) -- (32,4);
\draw (34,1) -- (34,4);
\draw (36,2) -- (36,4);
\draw (38,3) -- (38,5);
\draw (40,4) -- (40,5);
\draw (42,5) -- (42,6);
\draw (44,6) -- (44,7);
\draw (46,7) -- (46,8);
\draw (32,3) -- (34,4);
\draw (32,2) -- (38,5);
\draw (32,1) -- (46,8);
\draw (32,0) -- (62,15);
\draw (46,0) -- (48,1);
\draw (54,0) -- (56,1) -- (56,0) -- (64,4);
\draw (64,0) -- (66,1);
\node at (32,0) {\lb};
\node at (32,1) {\lb};
\node at (32,2) {\lb};
\node at (32,3) {\lb};
\node at (32,4) {\lb};
\node at (34,1) {\lb};
\node at (34,2) {\lb};
\node at (34,3) {\lb};
\node at (34,4) {\lb};
\node at (36,2) {\lb};
\node at (36,3) {\lb};
\node at (36,4) {\lb};
\node at (38,3) {\lb};
\node at (38,4) {\lb};
\node at (38,5) {\lb};
\node at (40,4) {\lb};
\node at (40,5) {\lb};
\node at (42,5) {\lb};
\node at (42,6) {\lb};
\node at (44,6) {\lb};
\node at (44,7) {\lb};
\node at (46,7) {\lb};
\node at (46,8) {\lb};
\node at (48,8) {\lb};
\node at (50,9) {\lb};
\node at (52,10) {\lb};
\node at (54,11) {\lb};
\node at (56,12) {\lb};
\node at (58,13) {\lb};
\node at (60,14) {\lb};
\node at (62,15) {\lb};
\node at (46,0) {\lb};
\node at (48,1) {\lb};
\node at (54,0) {\lb};
\node at (56,0) {\lb};
\node at (56,1) {\lb};
\node at (58,1) {\lb};
\node at (60,2) {\lb};
\node at (62,3) {\lb};
\node at (64,4) {\lb};
\node at (64,0) {\lb};
\node at (66,1) {\lb};
\node at (32,-.7) {$32$};
\node at (36,-.7) {$36$};
\node at (40,-.7) {$40$};
\node at (44,-.7) {$44$};
\node at (48,-.7) {$48$};
\node at (52,-.7) {$52$};
\node at (56,-.7) {$56$};
\node at (60,-.7) {$60$};
\node at (64,-.7) {$64$};
\end{\tz}
\end{center}
\end{fig}
\end{minipage}
\bigskip

\section{Summand approach to $ku^*(K_2)$}\label{summandsec}
In this section, we provide more details to the summand approach to $ku^*(K_2)$ initiated briefly in \cite[Section 7]{DW} and introduced here in Section \ref{kudifflsec}. In subsequent sections, we will use it to determine $ko_*(K_2)$. We use the notation of \cite{DW} as recalled early in Section \ref{kudifflsec}. 

First we define some $E_1$-submodules of $H^*(K_2)$ which generalize Figure \ref{71}. Here and elsewhere $\nu(m)$ is the 2-exponent of an integer $m$.

\begin{defin}\label{Ckdef} For $k\ge1$, let
$$C_k=(y_1q\oplus y_1^2)\cdot\bigoplus_{i=1}^{2^{k-2}-1}y_1^{2i-2}z(J(k,i))M_{\nu(i)+4},$$
where $J(k,i)= 2^k-4(2^{\nu(i)}+i)$ and $z(\sum 2^{j_t})=\prod z_{j_t}$ if $j_t$ are distinct integers.
\end{defin}

For example, $C_4$ is the sum of the modules in Figure \ref{71}, omitting $M_6$ and the last two.
Note that $C_1=C_2=0$. Also $M_3=0$. We recall notation from \cite{DW}, $Z_k^\ell=z_k\cdots z_{\ell-1}$.

\begin{defin}\label{ABdef} For $k\ge1$,
$$\Ah_k=M_{k+2}\oplus C_k\oplus y_1^{2^{k-1}-1}N\oplus \langle y_1^{2^{k-1}}\rangle.$$
For $1\le k<\ell$,
$$\Bh_{k,\ell}=z_\ell M_{k+2}\oplus z_\ell C_k\oplus y_1^{2^{k-1}-1}qM_{\ell+2}\oplus y_1^{2^{k-1}}M_{\ell+2}\oplus y_1^{2^{k-1}}Z_k^\ell C_k\oplus y_1^{2^k-1}qZ_{k+1}^\ell M_{k+2}.$$
\end{defin}

For example, $\Ah_4$ is in Figure \ref{71}. Since we are thinking here of the spectral sequence converging to $ku^*(K_2)$, and we depict that spectral sequence with gradings decreasing from left to right, we prefer to list our modules in that order, too.

In this section, we prove the following two key theorems. As in \cite{DW}, $\L_j$ denotes the exterior algebra generated by all $z_i$ with $i\ge j$, thought of as a set of monomials used as coefficients. We often use juxtaposition for tensor products.

\begin{thm}\label{everything} There is an isomorphism of $E_1$-modules,
$$H^*(K_2)=\bigoplus_{i\ge0,\,k\ge1}y_1^{2^k\,i}\Ah_k\oplus\bigoplus_{\substack{1\le k<\ell\\ i\ge0}}y_1^{2^k\,i}\L_{\ell+1}\Bh_{k,\ell}\oplus F,$$
where $F$ is a free $E_1$-module.
\end{thm}
\ni As in \cite{DW}, $F$ is computable but not of interest. It gives rise to a trivial $ku^*$-submodule of $ku^*(K_2)$.

\begin{thm}\label{closed} Each $y_1^{2^k\,i}\Ah_k$ and $y_1^{2^k\,i}\L_{\ell+1}\Bh_{k,\ell}$ is closed under the differentials in the ASS converging to $ku^*(K_2)$, as listed in \cite[Theorem 3.1]{DW}.
\end{thm}

The ideas in the following proofs will not be needed elsewhere,

\begin{proof}[Proof of Theorem \ref{everything}] By \cite[Proposition 2.11, (2.16)]{DW}, ignoring the free part, the $E_1$-summands of $H^*(K_2)$ are
\begin{equation}\label{E1spl}\{\langle y_1^i\rangle:i\ge1\}\cup\{y_1^iN:i\ge0\}\cup\{y_1^iM_j\L_{j-1},\,qy_1^iM_j\L_{j-1}:i\ge0,\,j\ge4\}.\end{equation}

From the last two summands of the $\Ah_k$'s, with their coefficients, we get
$$\{y_1^{2^{k-1}-1}N\cdot y_1^{2^k\,i}:k\ge1,i\ge0\}=\{y_1^iN:i\ge0\}$$ and $$\{\langle y_1^{2^{k-1}}\rangle y_1^{2^k\,i}:k\ge1,i\ge0\}=\{\langle y_1^i\rangle:i\ge1\}.$$
Since $C_1=0$ and $M_3=0$, all $\Bh_{1,\ell}$ with their coefficients give all $y_1^{2i+1}M_j\L_{j-1}$ and all $y_1^{2i}qM_j\L_{j-1}$.

For the rest, we consider the part without the $q$ factor. We will show that we obtain exactly all $y_1^{2i}M_j\L_{j-1}$ with $j\ge4$. The part with $q$ follows similarly.

As in \cite{DW}, let $y_k=y_1^{2^{k-1}}$. We let $\G_k=\{y_1^{2^{k-1}i}:i\ge0\}$, which we think of additively as an exterior algebra on $\{y_i:i\ge k\}$, analogous to $\L_k$ for the $z$'s.

We consider the coefficient of $M_k$. When we talk about $\Ah$'s and $\Bh$'s, we always include their coefficients as given in Theorem \ref{everything}. We consider first just the non-$C$ part in Definition \ref{ABdef}. Then $\Ah_{k-2}$ gives $\G_{k-1}$ as coefficient of $M_k$. Next, $\bigcup_{\ell\ge k-1}\Bh_{k-2,\ell}$ gives $\G_{k-1}\cdot\bigcup_{\ell\ge k-1}z_\ell\L_{\ell+1}=\G_{k-1}\overline \L_{k-1}$, where the bar denotes the augmentation ideal. Combining with the part obtained from $\Ah_{k-2}$ yields $\G_{k-1}\L_{k-1}$.

The $y_1^{2^{i-1}}M_k$ part of $\Bh_{i,k-2}$ for $2\le i\le k-3$ gives $\L_{k-1}\cdot\{y_1^t:1\le \nu(t)\le k-4\}$ as coefficient of $M_k$. The $\G_{k-1}\L_{k-1}$ obtained in the preceding paragraph gave the portion with $\nu(t)\ge k-2$. So the combination of the non-$C$ parts gave everything that we want except $y_{k-2}\G_{k-1}\L_{k-1}$.

One can show that, for $p\ge k-1$, the coefficient of $M_k$ in $C_p$ is 
\begin{equation}\label{Cp}y_{k-2}\prod_{i=k-1}^{p-1}\{z_i,y_i\},\end{equation}
where, as in \cite[(1.7)]{DW}, the product-of-sets notation means the set of all products containing one choice of the two elements in each factor. For example, if $p=k+1$, this product is $\{z_{k-1}y_{k-1},z_ky_{k-1},z_{k-1}y_k,z_ky_k\}$. An empty product equals 1. We verify the claim when $p=k+1$. The relevant values of $i$ in Definition \ref{Ckdef} are $2^{k-4}u$ for $u=1$, 3, 5, 7. For these, the $z(J(k+1,i))$ equal, respectively, $z_{k-1}z_k$, $z_k$, $z_{k-1}$, and 1, while $y_1^{2i}$ are $y_{k-2}$, $y_{k-2}y_{k-1}$, $y_{k-2}y_k$, and $y_{k-2}y_{k-1}y_k$. The case of arbitrary $p$ is similar.

Now we complete the proof by showing that the coefficient of $M_k$ coming from the occurrences of various $C_p$'s in $\Ah_p$'s and $\Bh_{p,\ell}$'s is $y_{k-2}\G_{k-1}\L_{k-1}$. Since $y_{k-2}$ is always a factor in (\ref{Cp}), we shall omit writing it in what follows. When we write ``From $\Ah_p$'' or ``from $\Bh_{p,\ell}$,'' we always mean to include their coefficients as stated in the theorem.

From $\Ah_{k-1}$, we get $\G_k$ (as coefficient of $M_k$), since the product in (\ref{Cp}) is empty. From $\bigcup_{\ell\ge k}\Bh_{k-1,\ell}$, we get
\begin{eqnarray*}&&\G_k\cdot\bigcup_{\ell\ge k}z_\ell\L_{\ell+1}\cup y_{k-1}\G_k\bigcup_{\ell\ge k}Z_{k-1}^\ell\L_{\ell+1}\\
&=&\G_k\overline{\L}_k\cup y_{k-1}\G_kz_{k-1}\L_k.\end{eqnarray*}
We combine with the part from $\Ah_{k-1}$ to remove the bar, yielding
\begin{equation}\label{one}\G_k\L_k\{1,y_{k-1}z_{k-1}\}.\end{equation}

From $\Ah_k$, using (\ref{Cp}), we get $\G_{k+1}\cdot\{z_{k-1},y_{k-1}\}$. From $\bigcup_{\ell\ge k+1}\Bh_{k,\ell}$, we get
\begin{eqnarray*} &&\{z_{k-1},y_{k-1}\}\G_{k+1}\bigcup_{\ell\ge k+1}\{z_\ell\L_{\ell+1},y_kZ_k^\ell\L_{\ell+1}\}\\
&=&\{z_{k-1},y_{k-1}\}\G_{k+1}\bigl(\overline{\L}_{k+1}\cup y_kz_k\L_{k+1}\bigr).\end{eqnarray*}
The part from $\Ah_k$ removes the bar, yielding
\begin{equation}\label{two}\G_{k+1}\L_{k+1}\{z_{k-1},y_{k-1}\}\cdot\{1,y_kz_k\}.\end{equation}

The part in (\ref{one}) is everything in $\G_{k-1}\L_{k-1}$ with an even number of factors with subscript $k-1$. The part in (\ref{two}) is everything with an odd (resp.~even) number of factors with subscript $k-1$ (resp.~$k$). Using $\Ah_{k+1}$ and $\Bh_{k+1,\ell}$, we will get those with odd (resp.~odd, even) number with subscript $k-1$ (resp.~$k$, $k+1$). This will continue, yielding all of $\G_{k-1}\L_{k-1}$, as desired.
\end{proof}

\begin{proof}[Proof of Theorem \ref{closed}]
We begin by considering the differentials within $C_k$. In Table \ref{C5tbl}, we list $C_5$ as a guide.

\renewcommand{\arraystretch}{1.5}
\begin{table}[h]
\caption{$C_5$}
\label{C5tbl}
\begin{tabular}{cccccccc}
&$y_1qz_3z_4M_4$&$y_1^2z_3z_4M_4$&$y_1^3qz_4M_5$&$y_1^4z_4M_5$&$y_1^5qz_4M_4$&$y_1^6z_4M_4$&$y_1^7qM_6$\\
$y_1^8M_6$&$y_1^9qz_3M_4$&$y_1^{10}z_3M_4$&$y_1^{11}qM_5$&$y_1^{12}M_5$&$y_1^{13}qM_4$&$y_1^{14}M_4$&

\end{tabular}
\end{table}

 Let $Y_{j,k}$ denote the term in $C_k$ containing $y_1^j$ as a factor. For example, $Y_{3,5}=y_1^3qz_4M_5$. We conflate $M_j$ with its chart, $\ext_{E_1}(\zt,M_j)$, as pictured in Figure \ref{Mchts}. Note that $M_j$ has $j-3$ $v$-towers. For example, the $v$-towers of $M_7$  are generated by $z_2^2z_3z_4$, $z_3^2z_4$, $z_4^2$, and $z_5$, in order of decreasing  grading.
  We let $b_i(M_r)$ denote the $i$th generator from the bottom of $M_r$, and $\tau_i(M_r)$ the $i$th generator from the top. For example, $b_1(M_7)=z_5$ and $\tau_2(M_7)=z_3^2z_4$. 

    We recall the notation of \cite[1.12)]{DW}:
 $$z_{i,j}=\begin{cases} z_i(z_i\cdots z_{j-1})&i<j\\ z_i&i=j\end{cases}$$
 Then $b_i(M_r)=z_{r-1-i,r-2}$ and $\tau_i(M_r)=z_{i+1,r-2}$.  
  If $Y_{r,k}=c\,M_j$, we let $b_i(Y_{r,k})=c\,b_i(M_j)$ and $\tau_i(Y_{r,k})=c\,\tau_i(M_j)$. Our claim is that the differentials of \cite[Theorem 3.1]{DW} are exactly
 \begin{equation}\label{bt}d^{2^{t+1}-(t+1)}(b_t(Y_{2^t(a+1)-1,k}))=v^{2^{t+1}-(t+1)}\tau_t(Y_{2^ta,k})\end{equation}
 for $1\le a\le2^{k-t-1}-2$ with $1\le t\le k-3$. This will remove or truncate all $v$-towers except $\tau_1(Y_{2^r-1,k})$ and $b_1(Y_{2^{k-1}-2^{r},k})$ for $1\le r\le k-2$. To justify this, we use $C_5$ as our example and consider $Y_{r,k}$ with $r$ even. The image of $d^2$ is all $\tau_1(Y_{2a,5})$ for $1\le a\le6$, while $\im(d^5)$ is $\{\tau_2(Y_{4,5}),\tau_2(Y_{8,5}\}$. Since $Y_{r,k}$ has $\nu(r)$ classes when $r$ is even, the elements not in $\im(d)$ are $\tau_3(Y_{8,5})$, $\tau_2(Y_{12,5})$, and $\tau_1(Y_{14,5})$, each of which is $b_1$.

Now we work toward proving (\ref{bt}).
 The formula in \cite[Theorem 3.1]{DW} is
 \begin{equation}\label{DWd}d^{2^{t+1}-(t+1)}(y_1^{2^t(a+1)-1}qz_{j-t+1,j})=v^{2^{t+1}-(t+1)}y_1^{2^t\,a}z_{t+1}z_j\end{equation} for $j\ge t+1$. This can be multiplied by $\L_{j+1}$.
 Our formula (\ref{bt}) says
 \begin{eqnarray}\nonumber&&d^{2^{t+1}-(t+1)}(y_1^{2^t(a+1)-1}qz(J(k,2^{t-1}(a+1)))b_t(M_{t+3+\nu(a+1)}))\\
 &=&v^{2^{t+1}-(t+1)}y_1^{2^ta}z(J(k,2^{t-1}a))\tau_t(M_{t+3+\nu(a)}).\label{our}\end{eqnarray}

 We have
 $$z(J(k,2^{t-1}(a+1)))=2^k-2^{t+1}(2^{\nu(a+1)}+a+1)$$
 and $$z(J(k,2^{t-1}a))=2^k-2^{t+1}(2^{\nu(a)}+a).$$
 If $a$ is odd, then
 $$z(J(k,2^{t-1}a))=z(J(k,2^{t-1}(a+1)))+2^{t+1+\nu(a+1)}$$
 with $\nu(z(J(k,2^{t-1}(a+1))))>t+1+\nu(a+1)$. If $a$ is even, then
 $$z(J(k,2^{t-1}(a+1)))=z(J(k,2^{t-1}a))+2^{t+1}(2^{\nu(a)}-2)$$
 with $\nu(z(J(k,2^{t-1}a)))>t+1+\nu(a)$.

 If $a$ is odd, then, with $\nu=\nu(a+1)$,
 $$b_t(M_{t+3+\nu})=z_{2+\nu,t+1+\nu}$$
 and
 $\tau_t(M_{t+3})=z_{t+1}$. In this case, our formula (\ref{our}) becomes
 $$d^{2^{t+1}-(t+1)}(y_1^{2^t(a+1)-1}qPz_{2+\nu,t+1+\nu})=v^{2^{t+1}-(t+1)}y_1^{2^ta}z_{t+1+\nu}Pz_{t+1}$$
 with $P\in \L_{t+2+\nu}$. This agrees with (\ref{DWd}) with $j=t+1+\nu$.

 If $a$ is even, then, with $\nu=\nu(j)$,
 $b_t(M_{t+3})=z_{2,t+1}$ and $\tau_t(M_{t+3+\nu})=z_{t+1,t+1+\nu}$. In this case, our formula (\ref{our}) becomes
 $$d^{2^{t+1}-(t+1)}(y_1^{2^t(a+1)-1}qz_{t+2}\cdots z_{t+\nu}Pz_{2,t+1})=y^{2^{t+1}-(t+1)}y_1^{2^ta}Pz_{t+1,t+1+\nu}$$
 with $P\in\L_{t+2+\nu}$. Since $z_{t+1,t+1+\nu}=z_{t+1}^2\cdot z_{t+2}\cdots z_{t+\nu}$, this agrees with (\ref{DWd}) with $j=t+1$. This completes the proof regarding differentials within $C_k$.

 Next we handle the differentials in $\Ah_k$ from $C_k$ to $M_{k+2}$. We claim that
 $$d^{2^{t+1}-(t+1)}(\tau_1(Y_{2^t-1,k}))=v^{2^{t+1}-(t+1)}\tau_t(M_{k+2})$$
 for $1\le t\le k-2$ follows from (\ref{DWd}). These can be seen in Figure \ref{alldiff}
 when $k=4$.
We have
\begin{eqnarray*}\tau_1(Y_{2^t-1,k})&=&y_1^{2^t-1}qz(J(k,2^{t-1}))\cdot\tau_1(M_{t+3})\\
&=&y_1^{2^t-1}qz_{t+2}\cdots z_{k-1}\cdot z_{2,t+1}.\end{eqnarray*}
By (\ref{DWd}), $d^{2^{t+1}-(t+1)}$ applied to this equals $$v^{2^{t+1}-(t+1)}z_{t+1}^2\cdot 
z_{t+2}\cdots z_{k-1}=v^{2^{t+1}-(t+1)}z_{t+1,k}=v^{2^{t+1}-(t+1)}\tau_t(M_{k+2}),$$
as claimed. The same argument applies to differentials in $\Bh_{k,\ell}$ from $z_\ell C_k$ to $z_\ell M_{k+2}$.

Next we explain the differentials in $\Ah_k$ from the $y_1^{2^{k-1}-1}N$ part. See Figure \ref{alldiff} for a depiction of the case $k=4$. Using the Ext calculation pictured in Figure \ref{N} and the differential
$$d^{\nu(i)+2}(y_1^i)=h_0^{\nu(i)}v^2qy_1^{i-1}$$
from \cite[Theorem 3.1]{DW}, the $v$-towers remaining in $y_1^{2^{k-1}-1}N$ are
$$h_0^sv^2qy_1^{2^{k-1}-1}\text{ for }0\le s\le k-2.$$
These satisfy
$$d^{2^{t+1}-(t+1)}(h_0^{t-1}v^2qy_1^{(a+1)2^t-1})=v^{2^{t+1}}y_1^{2^ta}z_{t+1}$$
by \cite[Theorem 3.1]{DW}. Here $1\le t\le k-1$. With $(a+1)2^t=2^{k-1}$, so $2^ta=2^{k-1}-2^t$, the target classes are the remaining classes $b_1(Y_{2^{k-1}-2^t,k})$ (see after (\ref{bt})), including also $b_1(M_{k+2})$.

Next we discuss the differentials in $\Bh_{k,\ell}$ from the $y_1^{2^k-1}qZ_{k+1}^\ell M_{k+2}$ part. We factor out the $Z_{k+1}^\ell$, so are hitting into $y_1^{2^{k-1}}z_kC_k$. We have, using \cite[Theorem 3.1]{DW}, if $1\le t\le k-2$,
\begin{eqnarray*}d^{2^{t+1}-(t+1)}(y_1^{2^k-1}qb_t(M_{k+2}))&=&d^{2^{t+1}-(t+1)}(y_1^{2^k-1}qz_{k+1-t,k})\\
&=&v^{2^{t+1}-(t+1)}y_1^{2^k-2^t}z_{t+1}z_k\\
&=&v^{2^{t+1}-(t+1)}y_1^{2^{k-1}}y_1^{2^{k-1}-2^t}z_{t+1}z_k\\
&=&v^{2^{t+1}-(t+1)}y_1^{2^{k-1}}b_1(Y_{2^{k-1}-2^t,k})z_k,\end{eqnarray*}
which were the remaining $v$-towers in $y_1^{2^{k-1}}z_kC_k$. If $t=k-1$, then
\begin{eqnarray*}&&d^{2^{t+1}-(t+1)}(y_1^{2^k-1}q\,Z_{k+1}^\ell b_t(M_{k+2}))\\&=&v^{2^{t+1}-(t+1)}y_1^{2^{k-1}}Z_{k+1}^\ell z_k^2\\&=&v^{2^{t+1}-(t+1)}y_1^{2^{k-1}}z_{k,\ell}\\
&=&v^{2^{t+1}-(t+1)}y_1^{2^{k-1}}\tau_{k-1}(M_{\ell+2})\end{eqnarray*}
in $y_1^{2^{k-1}}M_{\ell+2}$.

The argument for differentials in $\Bh_{k,\ell}$ from $y_1^{2^{k-1}-1}qM_{\ell+2}$ to $z_\ell C_k$ is almost identical.
\begin{eqnarray*}d^{2^{t+1}-(t+1)}(y_1^{2^{k-1}-1}qb_t(M_{\ell+2}))&=&d^{2^{t+1}-(t+1)}(y_1^{2^{k-1}-1}qz_{\ell+1-t,\ell})\\
&=&v^{2^{t+1}-(t+1)}y_1^{2^{k-1}-2^t}z_{t+1}z_\ell\\
&=&v^{2^{t+1}-(t+1)}b_1(Y_{2^{k-1}-2^t,k})z_\ell,\end{eqnarray*}
for $1\le t\le k-2$. If $t=k-1$, this is hitting $v^{2^{t+1}-(t+1)}z_\ell\tau_1(M_{k+2})$ in the $z_\ell M_{k+2}$ part of $\Bh_{k,\ell}$. Note that $v$-towers on $y_1^{2^{k-1}-1}qb_t(M_{\ell+2})$ for $k\le t\le \ell-1$ are not yet accounted for by these differentials.

Next we show that (\ref{DWd}) implies that the differentials from $y_1^{2^{k-1}}Z_k^\ell C_k$ into $y_1^{2^{k-1}}M_{\ell+2}$ hit $v^{2^{t+1}-(t+1)}y_1^{2^{k-1}}\tau_t(M_{\ell+2})$ for $1\le t\le k-2$. (The case $t=k-1$ was also hit, two paragraphs earlier.) We may ignore the $y_1^{2^{k-1}}$ factor. The $v$-towers in $C_k$ remaining to support differentials are
$$\tau_1(Y_{2^t-1,k})=y_1^{2^t-1}q\,z(J(k,2^{t-1}))\tau_1(M_{t+3})=y_1^{2^t-1}q\,Z_{t+2}^k\,z_{2,t+1}$$
for $1\le t\le k-2$. By (\ref{DWd})
$$d^{2^{t+1}-(t+1)}(Z_k^\ell\,\tau_1(Y_{2^t-1,k}))=Z_k^\ell\,Z_{t+2}^kv^{2^{t+1}-(t+1)}z_{t+1}^2=v^{2^{t+1}-(t+1)}z_{t+1,\ell}=v^{2^{t+1}-(t+1)}\tau_t(M_{\ell+2}).$$

What remains to consider in $\Bh_{k,\ell}$ are the $v$-towers on $y_1^{2^{k-1}}b_i(M_{\ell+2})$ and $y_1^{2^{k-1}-1}q\,b_{k-1+i}(M_{\ell+2})$ for $1\le i\le \ell-k$. (Here we have used $b_i(M_{\ell+2})=\tau_{\ell-i}(M_{\ell+2})$.) A different formula from \cite[Theorem 3.1]{DW} says
$$d^{k+1}(y_1^{2^{k-1}}z_\ell)=v^{k+1}y_1^{2^{k-1}-1}q\,z_{\ell-(k-1),\ell}$$
for $\ell\ge k+1$. This says
$$d^{k+1}(y_1^{2^{k-1}}b_1(M_{\ell+2}))=v^{k+1}y_1^{2^{k-1}-1}q\,b_k(M_{\ell+2}).$$
This differential is taking place in the sum of two charts of the form illustrated in Figure \ref{Mchts}. The charts are displaced by 1 unit. For $i\le\ell-k$, we have
\begin{eqnarray*}v^{i-1}d^{k+1}(y_1^{2^{k-1}}b_i(M_{\ell+2}))&=&d^{k+1}(y_1^{2^{k-1}}v^{i-1}b_i(M_{\ell+2}))\\
&=&d^{k+1}(y_1^{2^{k-1}}h_0^{i-1}b_1(M_{\ell+2}))\\
&=&v^{k+1}y_1^{2^{k-1}-1}q\,h_0^{i-1}b_k(M_{\ell+2})\\
&=&v^{k+1}y_1^{2^{k-1}-1}q\,v^{i-1}b_{k+i-1}(M_{\ell+2}).\end{eqnarray*}
Hence for $i\le\ell-k$, we have
$$d^{k+1}(y_1^{2^{k-1}}b_i(M_{\ell+2}))=v^{k+1}y_1^{2^{k-1}-1}q\,b_{k+i-1}(M_{\ell+2}).$$
\end{proof}

\section{The summands for $\A_k$ and $\B_{k,\ell}$}\label{kosummandsec}
As we illustrated for $A_4$ in Section \ref{kudifflsec}, if we apply $\ext_{E_1}(-,\zt)$ to the summands of $\Ah_k$ or $\Bh_{k,\ell}$, we obtain the $E_2$ page of the ASS converging to summands of $ku_*(K_2)$, and the differentials are dual to those in the ASS converging to $ku^*(K_2)$.

The ASS converging to $ko_*(K_2)$ has $E_2=\ext_{A(1)}(H^*K_2,\zt)$, where $A(1)$ is the subalgebra of the mod-2 Steenrod algebra generated by $\sq^1$ and $\sq^2$. In \cite{DW2}, we found $A(1)$-submodules of $H^*K_2$ corresponding to the $E_1$-submodules in the splitting (\ref{E1spl}). We also found $\ext_{A(1)}(-,\zt)$ for each of these summands. We review this now.

 Multiplying by $y_1^2$
just suspends by $\Si^8$. Corresponding to $y_1$ is an $A(1)$-module $U$. There is an $A(1)$-module $N$ which restricts to the $E_1$-module $N$. Corresponding to $y_1N$ is an $A(1)$-module $NU$, which satisfies $N\ot U=NU\oplus F$, with $F$ free. We defer discussion of the charts for these modules themselves until the next section.

We denote by $\M_k$ the $A(1)$-submodule of $H^*K_2$ found in \cite{DW2} whose $E_1$-module structure equals that of the $E_1$-module $M_k$ in (\ref{E1spl}) plus possibly a free $E_1$-module of rank 1, and let $\Mct_k=\Si^{-2^{k}}\M_k$. The chart for $\Mct_k$ is $M^0_k$, where $M^0_k$ is described in Figure \ref{Mcharts} and subsequent discussion. If $z_J=z_{j_1}\cdots z_{j_r}$ with $j_i$ distinct integers satisfying $j_i\ge k-1$, then $z_J\Mct_k$ is an $A(1)$-submodule of $H^*K_2$ with chart $\Si^{4D}M_k^r$, where $D=\sum 2^{j_i}$ and $M_k^r$ is as defined just before Figure \ref{M6}. It is formed from $M^0_k$ by decreasing filtrations by $r$. This follows from \cite[Proposition 3.5]{DW2}. 

Corresponding to multiplication by $q$ in (\ref{E1spl}) is just 9-fold suspension. Corresponding to $y_1z_JM_k$ in (\ref{E1spl}), with $J$ as above, is an $A(1)$-module $Uz_J\M_k$ described in \cite{DW2}, where its chart is shown to be $\Si^{2^k}\Si^{4D}M^{r+2}_k$, with $D$ as above. That is, multiplying $z_J\M_k$ by $U$ reduces filtrations by 2 without changing the number of suspensions.

We have now described the $A(1)$-submodules corresponding to all the $E_1$-submodules in (\ref{E1spl}), and \cite[Theorem 3.9]{DW2} says that a sum of them gives an $A(1)$-module splitting of $H^*K_2$. We define $\Ah_k^{\text{o}}$ and $z_J\Bh_{k,\ell}^{\text{o}}$ to be $A(1)$-modules formed from the $A(1)$-module analogues of the $E_1$-modules that formed  $\Ah_k$ and $z_J\Bh_{k,\ell}$ in Definition \ref{ABdef}. The proof of Theorem  \ref{everything} applies to $A(1)$-modules just as well as to $E_1$-modules, and so we have the $A(1)$ analogue of Theorem \ref{everything}.

\begin{thm}\label{koeverything} There is an isomorphism of $A(1)$-modules,
$$H^*(K_2)=\bigoplus_{i\ge0,\,k\ge1}y_1^{2^k\,i}\Ah_k^o\oplus\bigoplus_{\substack{1\le k<\ell\\ i\ge0}}y_1^{2^k\,i}\L_{\ell+1}\Bh_{k,\ell}^o\oplus F,$$
where $F$ is a free $A(1)$-module.
\end{thm}

We will see in subsequent sections how the differentials among the summands of $\Ah_k$ or $z_J\Bh_{k,\ell}$ imply differentials among the corresponding summands of $\Ah_k^{\text{o}}$ or $\Bh_{k,\ell}^{\text{o}}$, and these are closed under differentials by comparison with the $ku$ result,  Theorem \ref{closed}.

Incorporating this observation, we can prove Theorem \ref{main}.

\begin{proof}[Proof of Theorem \ref{main}] Since $y_1^2=\Si^8$ and $\A_k$ is the resultant of differentials and extensions in the ASS of $\Ah^o_k$, the first half follows immediately from the first half of Theorem \ref{koeverything}. There is a 1-1 correspondence between integers $2^{\ell+1}i$ and monomials in $\L_{\ell+1}$, given by $2^{\ell+1}i=\sum 2^{j_i}\leftrightarrow \prod z_{j_i}=z_J$ with $j_i\ge\ell+1$ distinct. Since $\a(i)$ equals the number of $2^j$'s, we have $z_J\B_{k,\ell}=\Si^{4\cdot2^{\ell+1}i}z^{\a(i)}\B_{k,\ell}$, and the second half of the theorem follows from this and the second half of Theorem \ref{koeverything}.\end{proof}

We will be working with charts rather than modules, and so we wish to take advantage of the $M_k^r$ notation. Recall from Sections 
\ref{Akexpls} and \ref{Bklsec} that $\At_k=\Si^{-2^{k+1}}\A_k$ and $\Bt_{k,\ell}=\Si^{-2^{\ell+2}}\B_{k,\ell}$. We define a tableau for $\At_k$ or $\Bt_{k,\ell}$ to be a list or direct sum of the charts which combine to give their $E_2$ page, but also usually include arrows for the $d_2$-differentials. These are listed in order of increasing grading. The first two summands of $\At_k$ are different and will be discussed in the next section; they involve higher differentials. The first summand, called $V_k$, was defined already in Definition \ref{E1'def}.

In either the $ku^*$ or $ku_*$ or $ko_*$ context, we associate to $M_k$ or $\M_k$ the grading $2^k$ (even though it starts in grading $2^k+1$.) Then $M_k^r$ will have associated grading 0, even though it starts in some positive grading. A summand $y_1^{2i}z(J(k,i))\M_{\nu(i)+4}$ appearing in $C_k$ in Definition \ref{Ckdef} has associated grading $$8i+4(2^k-4(2^{\nu(i)}+i))+2^{\nu(i)+4}=2^{k+2}-8i.$$ For $\At_k$, we subtract $2^{k+1}$ from the grading, and so, listing them in order of increasing grading, opposite to the order in Section \ref{kudifflsec}, associated to grading $8t$ in $\At_k$ is $y_1^{2^{k-1}-2t}z(J(k,2^{k-2}-t))\M_{\nu(t)+4}$.
For $1\le t<2^{k-2}$, $J(k,2^{k-2}-t)=4(t-\nu(t))$, so there are $\a(t)-1$ $z$'s in $z(J(k,2^{k-2}-t))$. Hence the chart will be $\Si^{8t}M_{\nu(t)+4}^{\a(t)-1}$. The next term in the list will be the one with the same $i$ value in Definition \ref{Ckdef} but with $y_1^2$ replaced by $y_1q$, so $\Si^8$ replaced by $\Si^9U$, and $U$ adds 2 to the superscript of $M$. Thus the tableau for $\At_k$ is 
\begin{equation}\label{Aktableau}(V_k\lar\Si^8M_4^0)\oplus\bigoplus_{i=1}^{2^{k-2}-1}(\Si^{8i+1}M_{4+\nu(i)}^{\a(i)+1}\lar \Si^{8i+8}M_{4+\nu(i+1)}^{\a(i+1)-1}).\end{equation}
For example, the tableaux for $\At_4$ and $\At_5$ are in (\ref{A4tableau}) and (\ref{A5tableau}).
\begin{equation}\label{A4tableau}V_4\lar  \Si^8M_4^0\quad \Si^9M_4^2\lar \Si^{16}M_5^0\quad \Si^{17}M_5^2\lar \Si^{24}M_4^1\quad \Si^{25}M_4^3\lar \Si^{32}M_6^0\end{equation}
\begin{align}
\nonumber V_5\lar  \Si^8M_4^0\quad \Si^9M_4^2\lar \Si^{16}M_5^0&\quad \Si^{17}M_5^2\lar \Si^{24}M_4^1\quad \Si^{25}M_4^3\lar \Si^{32}M_6^0\\
\label{A5tableau}\Si^{33}M_6^2\lar \Si^{40}M_4^1\quad \Si^{41}M_4^3\lar \Si^{48}M_5^1&\quad \Si^{49}M_5^3\lar\Si^{56}M_4^2\quad\Si^{57}M_4^4\lar \Si^{64}M_7^0.\end{align}
If the reader wishes to compare with Table \ref{C5tbl}, note that the order of terms is reversed. For example the entry $\Si^{33}M_6^2$ at the beginning of the second row of (\ref{A5tableau}) corresponds to the $y_1^7qM_6$ at the end of the first row of Table \ref{C5tbl}. The grading associated to the latter term in $\At_5$ is $24+9+64-64$. Recall that multiplying by a single $y_1$ doesn't change the number of suspensions; it lowers filtrations by 2 since its $A(1)$ analogue is $U$.

The tableau for $\Bt_{k,\ell}$ when $k\ge3$ is in (\ref{Bkltableau}). For $k\le2$, they are in (\ref{B2}). These are seen by comparison with Definition \ref{ABdef}.
\begin{gather}\nonumber (\Si M_{k+2}^{\ell-k+1}\lar\Si^8M_4^{\ell-k})\oplus\bigoplus_{i=1}^{2^{k-2}-2}(\Si^{8i+1}M_{4+\nu(i)}^{\ell-k+\a(i)+1}\lar\Si^{8i+8}M_{4+\nu(i+1)}^{\ell-k+\a(i+1)-1})\\
\nonumber \oplus(\Si^{2^{k+1}-7}M_4^{\ell-1}\lar\Si^{2^{k+1}}M_{\ell+2}^0)\oplus(\Si^{2^{k+1}+1}M_{\ell+2}^2\lar\Si^{2^{k+1}+8}M_4^1)\\
\label{Bkltableau}\oplus\bigoplus_{i=2^{k-2}+1}^{2^{k-1}-2}(\Si^{8i+1}M_{4+\nu(i)}^{\a(i)+1}\lar\Si^{8i+8}M_{4+\nu(i+1)}^{\a(i+1)-1})\oplus(\Si^{2^{k+2}-7}M_4^k\lar\Si^{2^{k+2}}M_{k+2}^1).\end{gather}

\begin{align}\Bt_{1,\ell}:&\ M^2_{\ell+2}\lar\Si^9M_{\ell+2}^0\label{B2}\\
\nonumber \Bt_{2,\ell}:&\ \Si M_4^{\ell-1}\lar\Si^8M_{\ell+2}^0\quad \Si^9M_{\ell+2}^2\lar\Si^{16}M_4^1.
\end{align}
\noindent The surprising change in position of superscripts when $k=1$ is due to the $2^{k-1}$ in Definition \ref{ABdef}.

We will compare these tableaux with their $ku$-homology analogues. In the $ku$ context, the role of the superscript is $\Si^2$. This is because the role of $z_j$ in $ku$-homology or $ku$-cohomology charts is $\Si^{2^{j+2}+2}$. Recalling that reducing filtration in $M$ charts by 4 is equivalent to an 8-fold suspension lends credence to the correspondence between reducing filtration by 1 (for $ko$) and suspending twice (for $ku$).

The tableau for $z^i\Bt_{k,\ell}$ is obtained from (\ref{Bkltableau}) by increasing all superscripts by $i$. As mentioned in the introduction, $z^i$ refers to any product of $i$ distinct $\Si^{-2^{j+2}}z_j$'s.

We close this section by discussing briefly the trivial $ko_*$-submodule of $ko_*(K_2)$.
The $A(1)$-module $\tilde H^*(K_2)$ can be decomposed as $I\oplus F$, where $I$ has no free summands, and $F$ is free. The trivial $ko_*$-submodule corresponds to a basis of $F$. By \cite[Theorem 3.9]{DW2}, $I$ can be written as
$$P[y_1^2]\ot(U\oplus\langle y_1^2\rangle\oplus N\oplus NU\oplus(\zt\oplus\Si^9)\bigoplus_{k\ge4}(M_k\L_{k-1}\oplus M_k\L_{k-1}U)).$$
Using results of \cite{DW2}, one can write a Poincar\'e series for $I$. For summands such as $z_JM_k$, we know the Ext chart, and  can determine from that the $A(1)$-module structure. Since things become quite complicated, here we just present through grading 24. In this range, $I$ is just $P[y_1^2]\ot(U\oplus\langle y_1^2\rangle\oplus N\oplus NU\oplus M_4\oplus M_4U)$, and the Poincar\'e series $P_I$ is \begin{gather}x^2+x^3+x^4+2x^5+x^6+x^7+2x^8+2x^9+2x^{10}+3x^{11}+3x^{12}+3x^{13}+3x^{14}\\+2x^{15}+2x^{16}+3x^{17}+3x^{18}+4x^{19}+5x^{20}+5x^{21}+5x^{22}+4x^{23}+3x^{24}.\end{gather}
Unreduced $H^*(K_2)$ is a polynomial algebra on classes of grading $2^j+1$, $j\ge0$, from which its Poincar\'e series $P_H$ is easily written. The Poincar\'e series for $A(1)$ is $P_{A(1)}=1+x+x^2+2x^3+x^4+x^5+x^6$. The Poincar\'e series $P_F$ for $F$ must satisfy
$$P_H-1=P_I+ (P_{A(1)}\cdot P_F).$$
We compute $P_F$ to begin  $$P_F=x^6+x^{10}+x^{12}+3x^{14}+x^{15}+x^{16}+3x^{18}+x^{19}+3x^{20}+x^{21}+5x^{22}+3x^{23}+5x^{24},$$
 the start of the Poincar\'e series for the trivial $ko_*$-submodule of $ko_*(K_2)$ through grading 24.

\section{Computation of $\A_4$ and $\A_1$}\label{A4sec}
In this section, we present a detailed argument for the $ko$-homology chart $\A_4$ given in Figure \ref{A4}, setting the stage for more general arguments to follow. At the end of the section, we also do $\A_1$, which does not follow some formulas for larger $k$.

The first part of $\Ah_4^o$ is the $A(1)$-module $\langle y_1^8\rangle\oplus y_1^6NU$, where $NU$ is an $A(1)$-module discussed in \cite{DW2}. Its $A(1)$-module structure and the chart $\ext_{A(1)}(NU,\zt)$ appear in Figure \ref{NU}.

\bigskip
\begin{minipage}{6in}
\begin{fig}\label{NU}

{\bf  $NU$ and its chart.}

\begin{center}

\begin{\tz}[scale=.4]
\draw (18,-12) -- (20,-10.5);
\node at (10,-10.5) {\lb};
\node at (14,-10.5) {\lb};
\draw (14,-10.5) -- (16,-10.5);
\draw (10,-10.5) to[out=30, in=150] (14,-10.5);

\node at (16,-10.5) {\lb};
\node at (20,-10.5) {\lb};
\node at (14,-12) {\lb};
\node at (16,-12) {\lb};
\node at (18,-12) {\lb};
\node at (20,-12) {\lb};
\node at (22,-12) {\lb};
\draw (14,-12) -- (16,-12);
\draw (20,-12) -- (22,-12);

\draw (27,-12) -- (39,-12);
\draw (28,-12) -- (28,-6);
\draw (32,-11) -- (32,-6);
\draw (36,-6) -- (36,-10) -- (38,-8);
\node at (10,-11.2) {$9$};
\node at (14,-12.7) {$11$};
\node at (30,-12) {\lb};
\node at (28,-12.7) {$9$};
\node at (32,-12.7) {$13$};
\node at (36,-12.7) {$17$};
\node at (40,-6) {$\iddots$};
\draw (16,-12) to[out=330, in=210] (20,-12);
\draw (16,-10.5) to[out=30, in=150] (20,-10.5);
\draw (14,-12) to[out=30, in=150] (18,-12);
\draw (18,-12) to[out=30, in=150] (22,-12);
\end{\tz}
\end{center}
\end{fig}
\end{minipage}
\bigskip

The chart for $\langle y_1^8\rangle\oplus y_1^6NU$ appears in Figure \ref{UNpic}. There is a $d_5$-differential from grading 33 to 32 implied by comparison with Browder's result (\cite{Bro}) that $H_{32}(K_2)\approx\Z/32$. See \cite[p.124]{W}, \cite[near Figure 4.2]{DWSW}, or \cite[Figure 11]{DW} for similar discussions. The $d_5$-differential is promulgated by the action of $ko_*$ on the spectral sequence.

There is an exotic $\eta$ extension from $(35,0)$ to $(36,5)$, and this is promulgated by $v_1^4$-periodicity. This follows from a result in \cite[p.228]{DM}, which says that
\begin{equation}\label{bracket}\text{if }d_r(\langle \a,\eta,2\rangle)=2\b,\text{ then }\a\eta=\b.\end{equation}
Here $\langle \a,\eta,2\rangle$ is a Toda bracket which gives $v_1\a$. The bracket relation can be deduced from Figure \ref{NU}, where the class in $(13,1)$ comes from the relation $\sq^3g_{11}=\sq^2\sq^3g_9$. That $v_1$ times the class in $(43,4)$ in Figure \ref{UNpic} equals the class in $(45,5)$ can be deduced from the morphism $ku_*\mapright{r}ko_*$. We prefer to elevate filtrations of classes to make $\eta$ extensions easier to see. In effect, whenever we have a $d_r$-differential, we elevate the classes supporting the differential, and most of those already related to them via $h_0$ or $h_1$ or $v_1^4$, so that the differential looks like a $d_1$. When this is done to the chart on the left side of Figure \ref{UNpic}, the chart on the right hand side is obtained.

\bigskip
\begin{minipage}{6in}
\begin{fig}\label{UNpic}

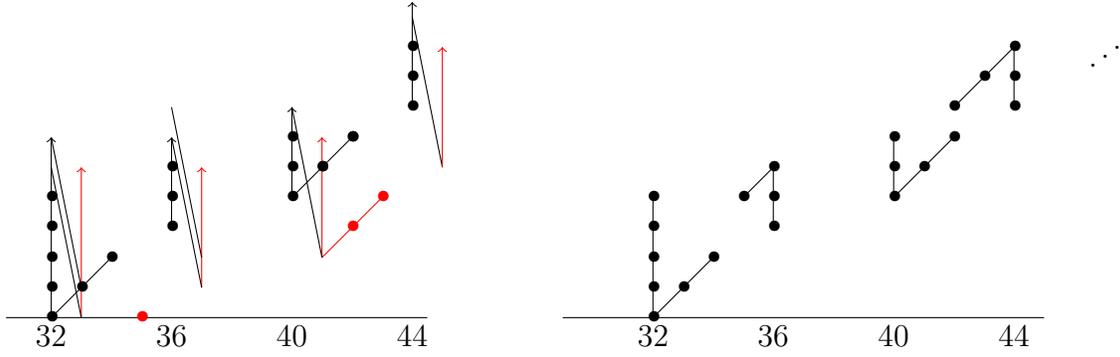
{\bf $\langle y_1^8\rangle\oplus y_1^6N U$ with differentials, and redrawing.}

\begin{center}

\begin{\tz}[scale=.4]

\draw (26.5,0) -- (40.5,0);

\draw [->] (28,0) -- (28,6);
\draw (28,0) -- (30,2);
\draw [->] (32,3) -- (32,6);
\draw [->] (36,4) -- (36,7);
\draw (36,4) -- (38,6);
\draw [->] (40,7) -- (40,10.5);

\draw [red] [->] (29,0) -- (29,5);
\draw [red] [->] (33,1) -- (33,5);
\draw [red] [->] (37,2) -- (37,6);
\draw [red] (37,2) -- (39,4);

\draw (29,0) -- (28,5);
\draw (29,1) -- (28,6);
\draw (33,1) -- (32,6);
\draw (33,2) -- (32,7);
\draw (37,2) -- (36,7);
\draw (41,5) -- (40,10);
\draw [red] [->] (41,5) -- (41,9);

\node at (28,-.6) {$32$};
\node at (32,-.6) {$36$};
\node at (36,-.6) {$40$};
\node at (40,-.6) {$44$};

\node at (28,0) {\lb};
\node at (28,1) {\lb};
\node at (28,2) {\lb};
\node at (28,3) {\lb};
\node at (28,4) {\lb};
\node at (32,3) {\lb};
\node at (32,4) {\lb};
\node at (32,5) {\lb};
\node at (36,4) {\lb};
\node at (36,5) {\lb};
\node at (36,6) {\lb};
\node at (37,5) {\lb};
\node at (38,6) {\lb};
\node at (40,8) {\lb};
\node at (40,7) {\lb};
\node at (40,9) {\lb};

\node [red] at (31,0) {\lb};

\node at (29,1) {\lb};
\node at (30,2) {\lb};
\node at (37,5) {\lb};
\node at (38,6) {\lb};
\node [red] at (38,3) {\lb};
\node [red] at (39,4) {\lb};
\draw (45,0) -- (61,0);
\draw (48,4) -- (48,0) -- (50,2);
\draw (52,3) -- (52,5) -- (51,4);
\draw (56,6) -- (56,4) -- (58,6);
\draw (58,7) -- (60,9) -- (60,7);
\node at (48,0) {\lb};
\node at (48,1) {\lb};
\node at (48,2) {\lb};
\node at (48,3) {\lb};
\node at (48,4) {\lb};
\node at (49,1) {\lb};
\node at (50,2) {\lb};
\node at (52,3) {\lb};
\node at (52,4) {\lb};
\node at (52,5) {\lb};
\node at (51,4) {\lb};
\node at (56,4) {\lb};
\node at (56,5) {\lb};
\node at (56,6) {\lb};
\node at (57,5) {\lb};
\node at (58,6) {\lb};
\node at (58,7) {\lb};
\node at (59,8) {\lb};
\node at (60,9) {\lb};
\node at (60,8) {\lb};
\node at (60,7) {\lb};
\node at (48,-.6) {$32$};
\node at (52,-.6) {$36$};
\node at (56,-.6) {$40$};
\node at (60,-.6) {$44$};
\node at (63,9) {$\iddots$};
\end{\tz}
\end{center}
\end{fig}
\end{minipage}
\bigskip

The right side of Figure \ref{UNpic} is what we called $\Si^{32}V_4$ in (\ref{Aktableau}).
It was useful to incorporate the $\Si^{32}$ at the outset because it affected the differential in $V_4$. But now we wish to switch to $\At_4$, by applying $\Si^{-32}$.
Now we consider
$(V_4\lar \Si^8M_4^0)$ in (\ref{Aktableau}). The chart is on the left hand side of Figure \ref{V4d}.
The differential follows by comparison with the $ku$ case, which is pictured in the right side of Figure \ref{d2d6}. The differential agrees with the first morphism in the 
 exact sequence
$$ko_*(M(2))\to ko_*(M(8))\to ko_*(M(4))\to ko_*(\Si M(2)),$$
where $M(n)$ is the mod-$n$ Moore spectrum. So the resultant must be $ko_*(M(4))$. In the middle part of Figure \ref{V4d}, we show the differential interpreted as a $d_1$, and in the right hand side we show the result, with the $\eta$ extension incorporated. The middle and right parts omit the portion in grading less than 8. It should be included, and several of its classes are promulgated by $v_1^4$-periodicity. The $\eta$ extension from 11 to 12 in the right hand chart can also be deduced from (\ref{bracket}).

\bigskip
\begin{minipage}{6in}
\begin{fig}\label{V4d}

{\bf $V_4\lar\Si^8M_4^0$}

\begin{center}

\begin{\tz}[scale=.38]
\node at (62,8) {$\iddots$};
\draw (45,0) -- (61,0);
\draw (48,4) -- (48,0) -- (50,2);
\draw (52,3) -- (52,5) -- (51,4);
\draw (56,6) -- (56,4) -- (58,6);
\draw (58,7) -- (60,9) -- (60,7);
\node at (48,0) {\lb};
\node at (48,1) {\lb};
\node at (48,2) {\lb};
\node at (48,3) {\lb};
\node at (48,4) {\lb};
\node at (49,1) {\lb};
\node at (50,2) {\lb};
\node at (52,3) {\lb};
\node at (52,4) {\lb};
\node at (52,5) {\lb};
\node at (51,4) {\lb};
\node at (56,4) {\lb};
\node at (56,5) {\lb};
\node at (56,6) {\lb};
\node at (57,5) {\lb};
\node at (58,6) {\lb};
\node at (58,7) {\lb};
\node at (59,8) {\lb};
\node at (60,9) {\lb};
\node at (60,8) {\lb};
\node at (60,7) {\lb};
\node at (48,-.6) {$0$};
\node at (52,-.6) {$4$};
\node at (56,-.6) {$8$};
\node at (60,-.6) {$12$};
\draw [red] (57,0) -- (59,2) -- (59,1) -- (61,3);
\node [red] at (57,0) {\lb};
\node [red] at (58,1) {\lb};
\node [red] at (59,2) {\lb};
\node [red] at (59,1) {\lb};
\node [red] at (60,2) {\lb};
\node [red] at (61,3) {\lb};
\draw [blue] (57,0) -- (56,6);
\draw [blue] (59,1) -- (58,7);
\draw [blue] (60,2) -- (59,8);
\draw [blue] (61,3) -- (60,9);
\draw (68,6) -- (68,4) -- (70,6);
\draw (70,7) -- (72,9) -- (72,7);
\draw [red] (69,5.3) -- (71,7) -- (71,6) -- (73,7.6);
\draw [blue] (69,5.3) -- (68,6);
\draw [blue] (71,6) -- (70,7);
\draw [blue] (71.9,6.92) -- (71,8);
\draw [blue] (73,7.6) -- (72,9);
\draw (67,0) -- (73,0);
\node at (74,9) {$\iddots$};
\node at (68,6) {\lb};
\node at (68,5) {\lb};
\node at (68,4) {\lb};
\node at (69,5) {\lb};
\node at (70,6) {\lb};
\node at (70,7) {\lb};
\node at (71,8) {\lb};
\node at (72,9) {\lb};
\node at (72,8) {\lb};
\node at (72,7) {\lb};
\node [red] at (69,5.3) {\lb};
\node [red] at (70,6.15) {\lb};
\node [red] at (71,7) {\lb};
\node [red] at (71,6) {\lb};
\node [red] at (73,7.6) {\lb};
\node [red] at (71.9,6.92) {\lb};
\node at (68,-.6) {$8$};
\node at (72,-.6) {$12$};
\draw (80,5) -- (80,4) -- (82,6);
\draw (82,6.3) -- (84,8) -- (84,7);
\draw (79,0) -- (85,0);
\node at (80,-.6) {$8$};
\node at (84,-.6) {$12$};
\node at (80,5) {\lb};
\node at (80,4) {\lb};
\node at (81,5) {\lb};
\node at (82,6) {\lb};
\node at (84,7) {\lb};
\node at (84,8) {\lb};
\node at (82,6.3) {\lb};
\node at (83.1,7.25) {\lb};
\node at (86,9) {$\iddots$};
\end{\tz}
\end{center}
\end{fig}
\end{minipage}
\bigskip

In Figure \ref{kokud2}, we show how the next $d_2$ differential in the tableau (\ref{A4tableau}),
$\Si^9M_4^2\lar \Si^{16}M_5^0$, is deduced by naturality from the $d_2$ in the $ku$ spectral sequence. The $ku$ differential can be seen as the $46\lar51$ arrow in Figure \ref{homdiff}. In Figure \ref{kokud2}, we are using gradings for $\At_4$, which are 32 greater than the $A_4$ gradings in Figure \ref{homdiff}. The first $ko$ differential is implied by naturality, and the next two are implied by the action of $\eta$. The fourth one can be thought of as being implied by the action of $v_1$, or can be deduced from the map $ku_{*+2}(-)\to ko_*(-)$. In the figure, we use big dots for the ones that map across under $c$.

\bigskip
\begin{minipage}{6in}
\begin{fig}\label{kokud2}

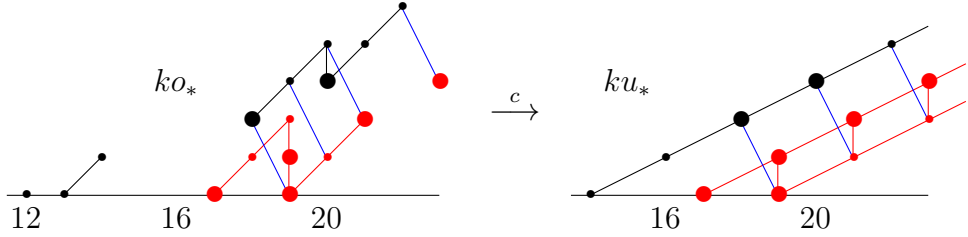
{\bf Deducing a $ko$ differential from $ku$}

\begin{center}

\begin{\tz}[scale=.5]
\draw (11.5,0) -- (23,0);
\draw (26.5,0) -- (36,0);
\draw (13,0) -- (14,1); 
\draw (18,2) -- (20,4) -- (20,3) -- (22,5);
\draw (27,0) -- (36,4.5);
\draw [red] (17,0) -- (19,2) -- (19,0) -- (21,2);
\draw [blue] (19,0) -- (18,2);
\draw [blue] (20,1) -- (19,3);
\draw [blue] (21,2) -- (20,4);
\draw [blue] (23,3) -- (22,5);
\draw [red] (30,0) -- (37,3.5);
\draw [red] (32,0) -- (37,2.5);
\draw [blue] (32,0) -- (31,2);
\draw [blue] (34,1) -- (33,3);
\draw [blue] (36,2) -- (35,4);
\node at (12,-.6) {$12$};
\node at (16,-.6) {$16$};
\node at (20,-.6) {$20$};
\node at (29,-.6) {$16$};
\node at (33,-.6) {$20$};
\node at (12,0) {\llb};
\node at (13,0) {\llb};
\node at (14,1) {\llb};
\node at (19,3) {\llb};
\node at (20,4) {\llb};
\node at (21,4) {\llb};
\node at (22,5) {\llb};
\node at (27,0) {\llb};
\node at (29,1) {\llb};
\node at (35,4) {\llb};
\node at (18,2) {\blb};
\node at (20,3) {\blb};
\node at (31,2) {\blb};
\node at (33,3) {\blb};
\node [red] at (17,0) {\blb};
\node [red] at (19,0) {\blb};
\node [red] at (19,1) {\blb};
\node [red] at (21,2) {\blb};
\node [red] at (23,3) {\blb};
\node [red] at (30,0) {\blb};
\node [red] at (32,0) {\blb};
\node [red] at (32,1) {\blb};
\node [red] at (34,2) {\blb};
\node [red] at (36,3) {\blb};
\node [red] at (18,1) {\llb};
\node [red] at (19,2) {\llb};
\node [red] at (20,1) {\llb};
\node [red] at (34,1) {\llb};
\node [red] at (36,2) {\llb};
\node at (25,2) {$\mapright{c}$};
\node at (16,3) {$ko_*$};
\node at (28,3) {$ku_*$};
\draw [red] (32,0) -- (32,1);
\draw [red] (34,1) -- (34,2);
\draw [red] (36,2) -- (36,3);

\end{\tz}
\end{center}
\end{fig}
\end{minipage}
\bigskip

In the next three figures, we show the three $d_2$-differentials involved in forming $\At_4$. 
(See (\ref{A4tableau}).) 
In each case, the domain chart has its filtrations increased by 1, so that the $d_2$ looks like a $d_1$. In the right side of each chart, we show the resultant. Of the early classes in the target chart, some are $v_1^4$-periodic in the resultant charts, while others are not. The ones that are not are indicated with open circles in the resultant chart. They will ``stay behind'' when subsequent differentials are considered, and will form the subedges.

\bigskip
\begin{minipage}{6in}
\begin{fig}\label{d2a}

{\bf $\Si^9M_4^2\lar \Si^{16}M_5^0$}

\begin{center}

\begin{\tz}[scale=.5]
\node at (12,0) {\lb};
\node at (13,0) {\lb};
\node at (14,1) {\lb};
\node at (18,2) {\lb};
\node at (19,3) {\lb};
\node at (20,4) {\lb};
\node at (20,3) {\lb};
\node at (21,4) {\lb};
\node at (22,5) {\lb};
\node at (31,0) {\lb};
\node at (35,1) {\lb};
\node at (36,2) {\lb};
\node at (37,3) {\lb};
\node at (37,2) {\lb};
\node at (38,3) {\lb};
\node at (39,4) {\lb};
\draw (30,0) circle (.15);
\draw (32,1) circle (.15);
\node [red] at (17.1,.8) {\lb};
\node [red] at (18.1,1.8) {\lb};
\node [red] at (19.1,2.8) {\lb};
\node [red] at (19.1,2) {\lb};
\node [red] at (19.1,1) {\lb};
\node [red] at (20.1,2) {\lb};
\node [red] at (21.1,3) {\lb};
\node [red] at (23,4) {\lb};
\draw (13,0) -- (14,1);
\draw (18,2) -- (20,4) -- (20,3) -- (22,5);
\draw [red] (17.1,.8) -- (19.1,2.8) -- (19.1,1) -- (21.1,3);
\draw [blue] (19.1,1) -- (18,2);
\draw [blue] (20.1,2) -- (19,3);
\draw [blue] (21.1,3) -- (20,4);
\draw [blue] (23,4) -- (22,5);
\draw (11,0) -- (23,0);
\node at (12,-.6) {$12$};
\node at (16,-.6) {$16$};
\node at (20,-.6) {$20$};
\node at (23,-.6) {$23$};
\draw (29,0) -- (39,0);
\node at (30,-.6) {$12$};
\node at (34,-.6) {$16$};
\node at (38,-.6) {$20$};
\draw (31.1,.1) -- (32,1);
\draw (35,1) -- (37,3) -- (37,2) -- (39,4);
\node at (27,3) {$=$};
\node at (41,5) {$\iddots$};
\end{\tz}
\end{center}
\end{fig}
\end{minipage}
\bigskip

\bigskip
\begin{minipage}{6in}
\begin{fig}\label{d2b}

{\bf $\Si^{17}M_5^2\lar \Si^{24}M_4^1$}

\begin{center}

\begin{\tz}[scale=.55]
\node at (20,0) {\lb};
\node at (22,0) {\lb};
\node at (24,1) {\lb};
\node at (26,2) {\lb};
\node at (27,3) {\lb};
\node at (28,4) {\lb};
\node at (28,3) {\lb};
\node at (28,2) {\lb};
\node at (29,3) {\lb};
\node at (30,4) {\lb};
\node at (38,0) {\lb};
\node at (42,1) {\lb};
\node at (43,2) {\lb};
\node at (44,3) {\lb};
\node at (44,2) {\lb};
\node at (45,3) {\lb};
\node at (46,4) {\lb};
\node at (48,5) {$\iddots$};
\draw (36,0) circle (.16);
\draw (40,1) circle (.16);
\draw (26,2) -- (28,4) -- (28,2) -- (30,4);
\draw (42,1) -- (44,3) -- (44,2) -- (46,4);
\draw [red] (26,1) -- (27,2) -- (27,1) -- (29,3.3);
\draw [blue] (27,1) -- (26,2);
\draw [blue] (29,3.3) -- (28,4);
\draw [blue] (27.8,1.92) -- (27,3);
\node [red] at (26,1) {\lb};
\node [red] at (27,2) {\lb};
\node [red] at (27,1) {\lb};
\node [red] at (27.8,1.92) {\lb};
\node [red] at (29,3.3) {\lb};
\draw (19,0) -- (30,0);
\draw (35,0) -- (46,0);
\node at (32.5,3) {$=$};
\node at (20,-.6) {$20$};
\node at (24,-.6) {$24$};
\node at (28,-.6) {$28$};
\node at (36,-.6) {$20$};
\node at (40,-.6) {$24$};
\node at (44,-.6) {$28$};

\end{\tz}
\end{center}
\end{fig}
\end{minipage}
\bigskip

\bigskip
\begin{minipage}{6in}
\begin{fig}\label{d2c}

{\bf $\Si^{25}M_4^3\lar \Si^{32}M_6^0$}

\begin{center}

\begin{\tz}[scale=.55]
\draw [red] (33,1) -- (35,3) -- (35,1);
\draw [red] (36,1) -- (37,2);
\draw [red] (39,4) -- (39,3);
\node at (30,0) {\lb};
\node at (34,1) {\lb};
\node at (35,2) {\lb};
\node at (36,3) {\lb};
\node at (36,2) {\lb};
\node at (37,3) {\lb};
\node at (38,4) {\lb};
\node at (48,1) {\lb};
\node at (49,2) {\lb};
\node at (50,3) {\lb};
\node at (50,2) {\lb};
\node at (50,1) {\lb};
\node at (51,2) {\lb};
\node at (52,3) {\lb};
\node at (54,4) {\lb};
\node at (56,5) {$\iddots$};
\draw (45,0) circle (.16);
\draw (49,1) circle (.16);
\draw (34,1) -- (36,3) -- (36,2) -- (38,4);
\draw (48,1) -- (50,3) -- (50,1) -- (52,3);
\draw (29,0) -- (39,0);
\draw (45,0) -- (54,0);
\node at (30,-.6) {$30$};
\node at (35,-.6) {$35$};
\node at (39,-.6) {$39$};
\node at (45,-.6) {$30$};
\node at (50,-.6) {$35$};
\node at (54,-.6) {$39$};
\node at (42,3) {$=$};
\node [red] at (33,1) {\lb};
\node [red] at (34,2) {\lb};
\node [red] at (35,3) {\lb};
\node [red] at (35,2.3) {\lb};
\node [red] at (35,1) {\lb};
\node [red] at (36,1) {\lb};
\node [red] at (37,2) {\lb};
\node [red] at (39,3) {\lb};
\node [red] at (39,4) {\lb};
\draw [blue] (36,1) -- (35,2);
\draw [blue] (37,2) -- (36,3);
\draw [blue] (39,3) -- (38,4);

\end{\tz}
\end{center}
\end{fig}
\end{minipage}
\bigskip

In Figure \ref{new4}, we combine Figures \ref{V4d}, \ref{d2a}, \ref{d2b}, and 
\ref{d2c}. There is a $d_8$ from the periodic part of Figure \ref{d2a} to \ref{V4d}, and a $d_4$ from the periodic part of Figure \ref{d2c} to \ref{d2b}. These are deduced by naturality from  the $d_5$ and $d_9$ in Figure \ref{homdiff}. The indices of the differentials are decreased by 1 because of the filtration shift incorporated in the $d_2$'s above. In Figure \ref{new4}, we increase filtrations of periodic classes so that the differentials look like $d_1$'s. This allows for nicer pictures showing $\cdot2$'s and $\eta$'s more clearly, which is useful in consideration of subsequent differentials.

We use double circles for the previously-circled classes, and use single circles for classes in Figure \ref{d2b} which survive but do not lead to periodic classes, since $v_1^{4i}$ times them are hit by differentials. These circled classes form another edge. Note the significance of the periodic class in grading 13 in Figure \ref{d2a}; it supports a differential in Figure \ref{new4} because $v_1^4$ times it does.

\bigskip
\begin{minipage}{6in}
\begin{fig}\label{new4}

{\bf Forming $\At_4$}

\begin{center}

\begin{\tz}[scale=.38]
\draw (0,4) -- (0,0) -- (2,2);
\draw (4,3) -- (4,5) -- (3,4);
\draw (8,5) -- (8,4) -- (10,5.9);
\draw (10,6.1) -- (12,8) -- (12,7);
\draw (16,9) -- (16,8) -- (18,9.9);
\draw (18,10.1) -- (20,12) -- (20,11);
\draw (24,13) -- (24,12) -- (26,13.9);
\draw (26,14.1) -- (28,16) -- (28,15);
\draw (32,17) -- (32,16) -- (34,17.9);
\draw (34,18.1) -- (36,20) -- (36,19);
\draw [red] (17,8) -- (19,10) -- (19,9) -- (21,11);
\draw [red] (25,12) -- (27,14) -- (27,13) -- (29,15);
\draw [red] (33,16) -- (35,18) -- (35,17) -- (37,19);
\draw (26,1) -- (28,3) -- (28,2) -- (30,4);
\draw (34,5) -- (36,7) -- (36,6) -- (38,8);
\draw [red] (33,4) -- (35,5.7) -- (35,4) -- (37,6);
\draw (-.5,0) -- (38,0);
\node at (0,-.6) {$0$};
\node at (4,-.6) {$4$};
\node at (8,-.6) {$8$};
\node at (12,-.6) {$12$};
\node at (16,-.6) {$16$};
\node at (20,-.6) {$20$};
\node at (24,-.6) {$24$};
\node at (28,-.6) {$28$};
\node at (32,-.6) {$32$};
\node at (36,-.6) {$36$};
\node at (0,0) {\lb};
\node at (0,1) {\lb};
\node at (0,2) {\lb};
\node at (0,3) {\lb};
\node at (0,4) {\lb};
\node at (1,1) {\lb};
\node at (2,2) {\lb};
\node at (4,3) {\lb};
\node at (4,4) {\lb};
\node at (4,5) {\lb};
\node at (3,4) {\lb};
\node at (8,4) {\lb};
\node at (8,5) {\lb};
\node at (9,5) {\lb};
\node at (10,5.9) {\lb};
\node at (10,6.1) {\lb};
\node at (11,7) {\lb};
\node at (12,8) {\lb};
\node at (12,7) {\lb};
\node at (24,12) {\lb};
\node at (24,13) {\lb};
\node at (25,13) {\lb};
\node at (26,13.9) {\lb};
\node at (26,14.1) {\lb};
\node at (27,15) {\lb};
\node at (28,16) {\lb};
\node at (28,15) {\lb};
\node at (16,8) {\lb};
\node at (16,9) {\lb};
\node at (17,9) {\lb};
\node at (18,9.9) {\lb};
\node at (18,10.1) {\lb};
\node at (19,11) {\lb};
\node at (20,12) {\lb};
\node at (20,11) {\lb};
\node at (32,16) {\lb};
\node at (32,17) {\lb};
\node at (33,17) {\lb};
\node at (34,17.9) {\lb};
\node at (34,18.1) {\lb};
\node at (35,19) {\lb};
\node at (36,20) {\lb};
\node at (36,19) {\lb};
\node at (28,2) {\lb};
\node at (29,3) {\lb};
\node at (34,5) {\lb};
\node at (35,6) {\lb};
\node at (36,7) {\lb};
\node at (36,6) {\lb};
\node at (37,7) {\lb};
\node at (38,8) {\lb};
\draw (22,0) circle (.16);
\draw (26,1) circle (.16);
\draw (27,2) circle (.16);
\draw (28,3) circle (.16);
\draw (30,4) circle (.16);
\draw (12,0) circle (.12);
\draw (12,0) circle (.22);
\draw (14,1) circle (.12);
\draw (14,1) circle (.22);
\draw (20,0) circle (.12);
\draw (20,0) circle (.22);
\draw (24,1) circle (.12);
\draw (24,1) circle (.22);
\draw (30,0) circle (.12);
\draw (30,0) circle (.22);
\draw (34,1) circle (.12);
\draw (34,1) circle (.22);
\node [red] at (13,7) {\lb};
\node [red] at (17,8) {\lb};
\node [red] at (18,9) {\lb};
\node [red] at (19,10) {\lb};
\node [red] at (19,9) {\lb};
\node [red] at (20,10) {\lb};
\node [red] at (21,11) {\lb};
\node [red] at (25,12) {\lb};
\node [red] at (26,13) {\lb};
\node [red] at (27,14) {\lb};
\node [red] at (27,13) {\lb};
\node [red] at (28,14) {\lb};
\node [red] at (29,15) {\lb};
\node [red] at (33,16) {\lb};
\node [red] at (34,17) {\lb};
\node [red] at (35,18) {\lb};
\node [red] at (35,17) {\lb};
\node [red] at (36,18) {\lb};
\node [red] at (37,19) {\lb};
\node [red] at (33,4) {\lb};
\node [red] at (34,4.85) {\lb};
\node [red] at (35,5.7) {\lb};
\node [red] at (35,5) {\lb};
\node [red] at (35,4) {\lb};
\node [red] at (36,5) {\lb};
\node [red] at (37,6) {\lb};
\node [red] at (39,7) {\lb};
\draw [blue] (13,7) -- (12,8);
\draw [blue] (17,8) -- (16,9);
\draw [blue] (19,9) -- (18,10.1);
\draw [blue] (20,10) -- (19,11.05);
\draw [blue] (21,11) -- (20,12);
\draw [blue] (25,12) -- (24,13);
\draw [blue] (27,13) -- (26,14.1);
\draw [blue] (28,14) -- (27,15.05);
\draw [blue] (29,15) -- (28,16);
\draw [blue] (33,16) -- (32,17);
\draw [blue] (35,17) -- (34,18.1);
\draw [blue] (36,18) -- (35,19.05);
\draw [blue] (37,19) -- (36,20);
\draw [blue] (35,4) -- (34,5);
\draw [blue] (36,5) -- (35,6);
\draw [blue] (37,6) -- (36,7);
\draw [blue] (39,7) -- (38,8);

\end{\tz}
\end{center}
\end{fig}
\end{minipage}
\bigskip

The result of each of the three sets of differentials pictured along the upper edge of Figure \ref{new4} is a lightning flash, starting in grading 16, 24, and 32, continuing indefinitely. The result of the group of differentials in the lower edge is a lightning flash, beginning in 33, promulgated by $v_1^4$-periodicity. There is a $d_{12}$-differential from this lightning flash in 33 to the one in 32, deduced from the morphism of ASS's for $ko_*(K_2)\to ku_*(K_2)$. The differential in the $ku_*$ spectral sequence is the $d_{16}$ in Figure \ref{homdiff}. (This would have been a $d_{12}$ if we had employed similar filtration shifts in the $ku_*$ spectral sequence. Note also that gradings in Figure \ref{new4} differ by 32 from those in Figures \ref{homdiff}, \ref{A4}, and \ref{ku_*A4} because those are for $A_4$, while here we have been considering $\At_4$.) The $d_{12}$ totally annihilates both lightning flashes and all subsequent ones, turning what had been an infinite picture into a finite one.

Note that $v_1^4$ times the classes in $(28,2)$ and $(29,3)$ in Figure \ref{new4} are in the lightning flash supporting differentials, and so they support $d_{12}$ differentials into the end of the lightning flash that started in 24. The result of this is Figure \ref{A4}.

The relationship of the above analysis with the approach to edges given in Section \ref{Akthmsec} will be explained in the next section. As a preview, the circled classes in Figure \ref{new4} are $\Si^8\E_{2,3}$ (grading 12 and 14), $\Si^{16}\E_{2,4}$ (gradings 20 to 30(top class)), and $\Si^{24}\E_{3,4}$ (30 and 34). (See Figure \ref{four}.)

The exact sequence
\begin{equation}\label{exact}ko_{*-1}(K_2)\mapright{\eta}ko_{*}(K_2)\mapright{c} ku_*(K_2)\to ko_{*-2}(K_2)\mapright{\eta}ko_{*-1}(K_2)\end{equation}
can be used to deduce extensions. In \cite{DW}, we established that there are nontrivial extensions in Figure \ref{ku_*A4} in grading 58, 60, 62, and 64.
The classes in $ko_{58}(K_2)$ which are not in $\im(\eta)$ must map to $ku_{58}(K_2)$; the extension in the latter implies one in the former.
The two classes in $ko_{62}(K_2)$
in Figure \ref{A4} must be in the image from $ku_{64}(K_2)$, where the extension is present, implying the extension in $ko_{62}(K_2)$. 

We conclude this section by deriving the chart for $\A_1$, since it is slightly special. Its chart is formed from the $A(1)$-module $U\oplus N$. These modules are defined in \cite[Section 3]{DW2} and their charts are in \cite[Figure 3.10]{DW2}. We reproduce them here in Figure \ref{U+N} with $N$ in red.

\bigskip
\begin{minipage}{6in}
\begin{fig}\label{U+N}

{\bf Forming $\A_1$}
\begin{center}

\begin{\tz}[scale=.6]
\draw (1,0) -- (14,0);
\draw (4,1) -- (4,6);
\draw [red] (5,0) -- (5,6);
\draw [blue] (5,0) -- (4,2);
\draw [blue] (5,1) -- (4,3);
\draw [blue] (5,2) -- (4,4);
\draw (10,4) -- (8,2) -- (8,6);
\draw [red] (11,2) -- (9,0) -- (9,6);
\draw [blue] (11,2) -- (10,4);
\draw [blue] (10,1) -- (9.2,3.2);
\draw [blue] (9,0) -- (8,2);
\draw [blue] (9,1) -- (8,3);
\draw (12,5) -- (12,8);
\draw [red] (13,3) -- (13,7);
\draw [blue] (13,3) -- (12,5);
\draw [blue] (13,4) -- (12,6);
\draw [blue] (13,5) -- (12,7);
\node at (14,7) {$\iddots$};
\node at (2,-.5) {$2$};
\node at (4,-.5) {$4$};
\node at (8,-.5) {$8$};
\node at (12,-.5) {$12$};
\node at (2,0) {\lb};
\node at (4,1) {\lb};
\end{\tz}
\end{center}
\end{fig}
\end{minipage}
\bigskip

The $d_2$-differential was first noted in \cite[p.124]{W} and was shown in \cite[Figure 3.1]{DWSW}.
It follows from \cite{Bro}. The only classes in $\A_1$ are the $\zt$'s in grading 2 and 4.

\section{Derivation of edge description}\label{derivationsec}

The derivation of Definition \ref{E'def} requires a careful study of the tableau (\ref{Aktableau}) for $\At_k$  and the ways that edges (or pre-edges) are associated to arrows in the tableau. For every arrow with target $\Si^{8i+1}M_\ell^s$, there is a corresponding pre-edge  $\Si^{8i}\E'_{s,s+\ell-3}$, which involves $2^{\ell-4}$ arrows beginning with the specified one. We have found it very useful to list, as an example, the tableau for $\At_7$ with arrows labeled by their corresponding pre-edge. This is done in (\ref{A7tableau}). The labeled tableau for $\At_k$ for $k=4$, 5, and 6 can be read off from this by taking the first 1, 2, or 4 rows, respectively, with $V_7$ replaced by $V_k$.

\begin{align}\nonumber V_7\ml{1,7}\Si^8 M_4^0\quad \Si^9M_4^2\ml{2,3}\Si^{16}M_5^0&\quad \Si^{17}M_5^2\ml{2,4}\Si^{24}M_4^1\quad \Si^{25}M_4^3\ml{3,4}\Si^{32}M_6^0\\
\Si^{33}M_6^2\ml{2,5}\Si^{40}M_4^1\quad \Si^{41}M_4^3\ml{3,4}\Si^{48}M_5^1&\quad \Si^{49}M_5^3\ml{3,5}\Si^{56}M_4^2\quad\Si^{57}M_4^4\ml{4,5}\Si^{64}M_7^0\nonumber\\
\Si^{65}M_7^2\ml{2,6}\Si^{72}M_4^1\quad\Si^{73}M_4^3\ml{3,4}\Si^{80}M_5^1&\quad\Si^{81}M_5^3\ml{3,5}\Si^{88}M_4^2\quad \Si^{89}M_4^4\ml{4,5}\Si^{96}M_6^1\nonumber\\
\Si^{97}M_6^3\ml{3,6}\Si^{104}M_4^2\quad\Si^{105}M_4^4\ml{4,5}\Si^{112}M_5^2&\quad\Si^{113}M_5^4\ml{4,6}\Si^{120}M_4^3\quad\Si^{121}M_4^5\ml{5,6}\Si^{128}M_8^0\nonumber\\
\Si^{129}M_8^2\ml{2,7}\Si^{136}M_4^1\quad\Si^{137}M_4^3\ml{3,4}\Si^{144}M_5^1&\quad\Si^{145}M_5^3\ml{3,5}\Si^{152}M_4^2\quad\Si^{153}M_4^4\ml{4,5}\Si^{160}M_6^1\nonumber\\
\Si^{161}M_6^3\ml{3,6}\Si^{168}M_4^2\quad \Si^{169}M_4^4\ml{4,5}\Si^{176}M_5^2&\quad\Si^{177}M_5^4\ml{4,6}\Si^{184}M_4^3\quad\Si^{185}M_4^5\ml{5,6}\Si^{192}M_7^1\nonumber\\
\Si^{193}M_7^3\ml{3,7}\Si^{200}M_4^2\quad \Si^{201}M_4^4\ml{4,5}\Si^{208}M_5^2&\quad \Si^{209}M_5^4\ml{4,6}\Si^{216}M_4^3\quad \Si^{217}M_4^5\ml{5,6}\Si^{224}M_6^2\nonumber\\
\label{A7tableau}\Si^{225}M_6^4\ml{4,7}\Si^{232}M_4^3\quad \Si^{233}M_4^5\ml{5,6}\Si^{240}M_5^3&\quad \Si^{241}M_5^5\ml{5,7}\Si^{248}M_4^4\quad \Si^{249}M_4^6\ml{6,7}\Si^{256} M_9^0\end{align}

We begin our derivation of Definition \ref{E'def} by considering $\E'_{5,6}$. As can be seen in (\ref{A7tableau}), it can come from any of $\Si M_4^5\lar \Si^8M_t^{8-t}$ for $5\le t\le8$. In Figure \ref{E56} we consider the cases $t=6$ and $t=8$. The other two are similar. The left side depicts the arrow as a $d_1$ differential, with $\Si M_4^5$ in red, while the right side shows the result, with $\E'_{5,6}$ in red, and $\Si^8 M_t^{7-t}$ in black.

\bigskip
\begin{minipage}{6in}
\begin{fig}\label{E56}

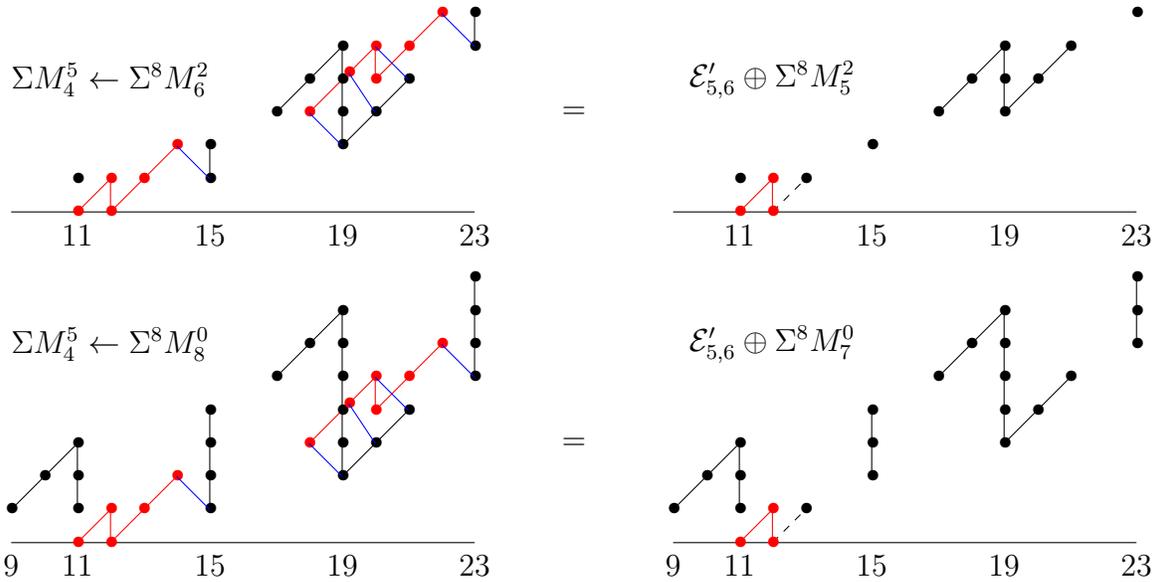
{\bf Two ways of forming $\E'_{5,6}$}

\begin{center}

\begin{\tz}[scale=.44]
\node at (26,13) {$=$};
\node at (26,3) {$=$};
\draw (9,1) -- (11,3) -- (11,1);
\draw (15,1) -- (15,4);
\draw (17,5) -- (19,7) -- (19,2) -- (21,4);
\draw (23,5) -- (23,8);
\draw (29,1) -- (31,3) -- (31,1);
\draw (35,2) -- (35,4);
\draw (37,5) -- (39,7) -- (39,3) -- (41,5);
\draw (43,6) -- (43,8);
\draw [red] (11,0) -- (12,1) -- (12,0) -- (14,2);
\draw [red] (18,3) -- (20,5) -- (20,4) -- (22,6);
\draw [red] (31,0) -- (32,1) -- (32,0);
\draw [dashed] (32,0) -- (33,1);
\draw (9,0) -- (23,0);
\draw (29,0) -- (43,0);
\draw (9,10) -- (23,10);
\draw (29,10) -- (43,10);
\draw (15,11) -- (15,12);
\draw (17,13) -- (19,15) -- (19,12) -- (21,14);
\draw (23,15) -- (23,16);
\draw (37,13) -- (39,15) -- (39,13) -- (41,15);
\draw [red] (11,10) -- (12,11) -- (12,10) -- (14,12);
\draw [red] (18,13) -- (20,15) -- (20,14) -- (22,16);
\draw [red] (31,10) -- (32,11) -- (32,10);
\draw [dashed] (32,10) -- (33,11);
\node at (9,-.7) {$9$};
\node at (11,-.7) {$11$};
\node at (15,-.7) {$15$};
\node at (19,-.7) {$19$};
\node at (23,-.7) {$23$};
\node at (29,-.7) {$9$};
\node at (31,-.7) {$11$};
\node at (35,-.7) {$15$};
\node at (39,-.7) {$19$};
\node at (43,-.7) {$23$};
\node at (11,9.3) {$11$};
\node at (15,9.3) {$15$};
\node at (19,9.3) {$19$};
\node at (23,9.3) {$23$};
\node at (31,9.3) {$11$};
\node at (35,9.3) {$15$};
\node at (39,9.3) {$19$};
\node at (43,9.3) {$23$};
\node at (12,14) {$\Si M_4^5\lar\Si^8M_6^2$};
\node at (32,14) {$\E'_{5,6}\oplus \Si^8M_5^2$};
\node at (12,6) {$\Si M_4^5\lar \Si^8M_8^0$};
\node at (32,6) {$\E'_{5,6}\oplus \Si^8M_7^0$};
\node at (9,1) {\lb};
\node at (10,2) {\lb};
\node at (11,3) {\lb};
\node at (11,2) {\lb};
\node at (11,1) {\lb};
\node at (15,1) {\lb};
\node at (15,2) {\lb};
\node at (15,3) {\lb};
\node at (15,4) {\lb};
\node at (17,5) {\lb};
\node at (18,6) {\lb};
\node at (19,7) {\lb};
\node at (19,6) {\lb};
\node at (19,5) {\lb};
\node at (19,4) {\lb};
\node at (19,3) {\lb};
\node at (19,2) {\lb};
\node at (20,3) {\lb};
\node at (21,4) {\lb};
\node at (23,5) {\lb};
\node at (23,6) {\lb};
\node at (23,7) {\lb};
\node at (23,8) {\lb};
\node at (29,1) {\lb};
\node at (30,2) {\lb};
\node at (31,3) {\lb};
\node at (31,2) {\lb};
\node at (31,1) {\lb};
\node at (35,2) {\lb};
\node at (35,3) {\lb};
\node at (35,4) {\lb};
\node at (37,5) {\lb};
\node at (38,6) {\lb};
\node at (39,7) {\lb};
\node at (39,6) {\lb};
\node at (39,5) {\lb};
\node at (39,4) {\lb};
\node at (39,3) {\lb};
\node at (40,4) {\lb};
\node at (41,5) {\lb};
\node at (43,6) {\lb};
\node at (43,7) {\lb};
\node at (43,8) {\lb};
\node at (11,11) {\lb};
\node at (15,11) {\lb};
\node at (15,12) {\lb};
\node at (17,13) {\lb};
\node at (18,14) {\lb};
\node at (19,15) {\lb};
\node at (19,14) {\lb};
\node at (19,13) {\lb};
\node at (19,12) {\lb};
\node at (20,13) {\lb};
\node at (21,14) {\lb};
\node at (23,15) {\lb};
\node at (23,16) {\lb};
\node at (31,11) {\lb};
\node at (33,11) {\lb};
\node at (35,12) {\lb};
\node at (37,13) {\lb};
\node at (38,14) {\lb};
\node at (39,15) {\lb};
\node at (39,14) {\lb};
\node at (39,13) {\lb};
\node at (40,14) {\lb};
\node at (41,15) {\lb};
\node at (43,16) {\lb};
\node [red] at (11,0) {\lb};
\node [red] at (12,1) {\lb};
\node [red] at (12,0) {\lb};
\node [red] at (13,1) {\lb};
\node [red] at (14,2) {\lb};
\node [red] at (18,3) {\lb};
\node [red] at (19.2,4.2) {\lb};
\node [red] at (20,5) {\lb};
\node [red] at (20,4) {\lb};
\node [red] at (21,5) {\lb};
\node [red] at (22,6) {\lb};
\node [red] at (31,0) {\lb};
\node [red] at (32,1) {\lb};
\node [red] at (32,0) {\lb};
\node [red] at (11,10) {\lb};
\node [red] at (12,11) {\lb};
\node [red] at (12,10) {\lb};
\node [red] at (13,11) {\lb};
\node [red] at (14,12) {\lb};
\node [red] at (18,13) {\lb};
\node [red] at (19.2,14.2) {\lb};
\node [red] at (20,15) {\lb};
\node [red] at (20,14) {\lb};
\node [red] at (21,15) {\lb};
\node [red] at (22,16) {\lb};
\node [red] at (31,10) {\lb};
\node [red] at (32,11) {\lb};
\node [red] at (32,10) {\lb};
\draw [blue] (15,1) -- (14,2);
\draw [blue] (19,2) -- (18,3);
\draw [blue] (20,3) -- (19.2,4.2);
\draw [blue] (21,4) -- (20,5);
\draw [blue] (23,5) -- (22,6);
\draw [blue] (15,11) -- (14,12);
\draw [blue] (19,12) -- (18,13);
\draw [blue] (20,13) -- (19.2,14.2);
\draw [blue] (21,14) -- (20,15);
\draw [blue] (23,15) -- (22,16);
\node at (33,1) {\lb};

\end{\tz}
\end{center}
\end{fig}
\end{minipage}
\bigskip

The job of this $d_2$ pair is to change $\Si^8M_t^{8-t}$ to $\Si^8M_{t-1}^{8-t}$. The $\E'_{5,6}$ is the extra part, and is not involved in subsequent arrows until we consider the differentials which change $\E'_{e,\ell}$ to $\E_{e,\ell}$. Its class in $(12,0)$ will support such a differential. In Figure \ref{E'56}, we show that the result of $\Si M_4^5\lar \Si^{10}\Mh_4^3$ is exactly the $\E'_{5,6}$ just obtained, as claimed in Definition \ref{E'def}. Here $\Si^{10}\Mh_4^3$ is in red.

\bigskip
\begin{minipage}{6in}
\begin{fig}\label{E'56}

{\bf $\Si M_4^5\lar \Si^{10}\Mh_4^3$}

\begin{center}

\begin{\tz}[scale=.6]
\draw (11,0) -- (12,1) -- (12,0) -- (14,2);
\draw (18,3) -- (20,5) -- (20,4) -- (22,6);
\draw [red] (14,0) -- (15,1);
\draw [red] (19,2) -- (21,4) -- (21,3) -- (23,5);
\draw [blue] (14,0) -- (13,1);
\draw [blue] (15,1) -- (14,2);
\draw [blue] (19,2) -- (18,3);
\draw [blue] (20,3) -- (19,4);
\draw [blue] (21,4) -- (20,5);
\draw [blue] (21,3) -- (20,4);
\draw [blue] (22,4) -- (21,5);
\draw [blue] (23,5) -- (22,6);
\draw (10,0) -- (23,0);
\node at (11,-.7) {$11$};
\node at (14,-.7) {$14$};
\node at (19,-.7) {$19$};
\node at (23,-.7) {$23$};
\node at (11,0) {\blb};
\node at (12,1) {\blb};
\node at (12,0) {\blb};
\node at (13,1) {\lb};
\node at (14,2) {\lb};
\node at (18,3) {\lb};
\node at (19,4) {\lb};
\node at (20,5) {\lb};
\node at (20,4) {\lb};
\node at (21,5) {\lb};
\node at (22,6) {\lb};
\node [red] at (14,0) {\lb};
\node [red] at (15,1) {\lb};
\node [red] at (19,2) {\lb};
\node [red] at (20,3) {\lb};
\node [red] at (21,4) {\lb};
\node [red] at (21,3) {\lb};
\node [red] at (22,4) {\lb};
\node [red] at (23,5) {\lb};

\end{\tz}
\end{center}
\end{fig}
\end{minipage}
\bigskip

We now give a chart-theoretic explanation of why this occurs, which is useful in seeing how it generalizes. We illustrate with $M_8^0$. The cofiber of $M_7^0\to M_8^0$ is $\Si^2 \Mh_4^3$. In Figure \ref{M78} we show that $M_8^0$ can be formed from $M_7^0+\Si^2 \Mh_4^3$. 

\bigskip
\begin{minipage}{6in}
\begin{fig}\label{M78}

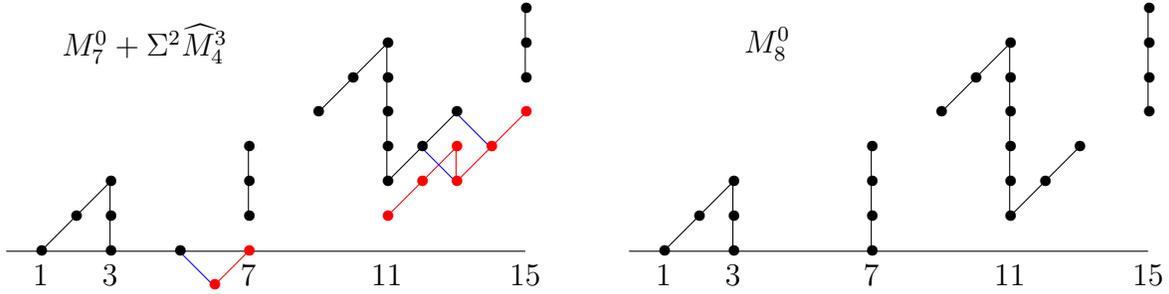
{\bf $M_7^0+\Si^2\Mh_4^3=M_8^0$}

\begin{center}

\begin{\tz}[scale=.46]
\draw (1,0) -- (3,2) -- (3,0);
\draw (7,1) -- (7,3);
\draw (9,4) -- (11,6) -- (11,2) -- (13,4);
\draw (15,5) -- (15,7);
\draw (19,0) -- (21,2) -- (21,0);
\draw (25,0) -- (25,3);
\draw (27,4) -- (29,6) -- (29,1) -- (31,3);
\draw (33,4) -- (33,7);
\draw [red] (6,-1) -- (7,0);
\draw [red] (11,1) -- (13,3) -- (13,2) -- (15,4);
\draw [blue] (6,-1) -- (5,0);
\draw [blue] (13,2) -- (12,3);
\draw [blue] (14,3) -- (13,4);
\draw (0,0) -- (15,0);
\draw (18,0) -- (33,0);
\node at (4,6) {$M_7^0+ \Si^2\Mh_4^3$};
\node at (22,6) {$M_8^0$};
\node at (1,-.7) {$1$};
\node at (3,-.7) {$3$};
\node at (7,-.7) {$7$};
\node at (11,-.7) {$11$};
\node at (15,-.7) {$15$};
\node at (19,-.7) {$1$};
\node at (21,-.7) {$3$};
\node at (25,-.7) {$7$};
\node at (29,-.7) {$11$};
\node at (33,-.7) {$15$};
\node at (1,0) {\lb};
\node at (2,1) {\lb};
\node at (3,2) {\lb};
\node at (3,1) {\lb};
\node at (3,0) {\lb};
\node at (5,0) {\lb};
\node at (7,1) {\lb};
\node at (7,2) {\lb};
\node at (7,3) {\lb};
\node at (9,4) {\lb};
\node at (10,5) {\lb};
\node at (11,6) {\lb};
\node at (11,5) {\lb};
\node at (11,4) {\lb};
\node at (11,3) {\lb};
\node at (11,2) {\lb};
\node at (12,3) {\lb};
\node at (13,4) {\lb};
\node at (15,5) {\lb};
\node at (15,6) {\lb};
\node at (15,7) {\lb};
\node at (19,0) {\lb};
\node at (20,1) {\lb};
\node at (21,2) {\lb};
\node at (21,1) {\lb};
\node at (21,0) {\lb};
\node at (25,0) {\lb};
\node at (25,1) {\lb};
\node at (25,2) {\lb};
\node at (25,3) {\lb};
\node at (27,4) {\lb};
\node at (28,5) {\lb};
\node at (29,6) {\lb};
\node at (29,5) {\lb};
\node at (29,4) {\lb};
\node at (29,3) {\lb};
\node at (29,2) {\lb};
\node at (29,1) {\lb};
\node at (30,2) {\lb};
\node at (31,3) {\lb};
\node at (33,4) {\lb};
\node at (33,5) {\lb};
\node at (33,6) {\lb};
\node at (33,7) {\lb};
\node [red] at (6,-1) {\lb};
\node [red] at (7,0) {\lb};
\node [red] at (11,1) {\lb};
\node [red] at (12,2) {\lb};
\node [red] at (13,3) {\lb};
\node [red] at (13,2) {\lb};
\node [red] at (14,3) {\lb};
\node [red] at (15,4) {\lb};

\end{\tz}
\end{center}
\end{fig}
\end{minipage}
\bigskip

In Figure \ref{samp} we modify
 the lower left part of Figure \ref{E56} by replacing $\Si^8M_8^0$ by $\Si^8M_7^0+ \Si^{10}\Mh_4^3$. The result is the desired $\E'_{5,6}\oplus\Si^8M_7^0$. Indeed, the $\Si M_4^5$ and $\Si^{10}\Mh_4^3$ combine to yield $\E'_{5,6}$ as in Figure \ref{E'56}, and the $\Si^8M_7^0$ remains unchanged. In chart arithmetic, we have
 $$(\Si M_4^5\lar M_7^0\oplus\Si^2\Mh_4^3)=(\E'_{5,6}\oplus M_7^0),$$
 then cancel the $M_7^0$ to obtain the $\E'_{5,6}$ case of Definition \ref{E'def}.
 This may seem roundabout, but will be useful.

 \bigskip
\begin{minipage}{6in}
\begin{fig}\label{samp}

{\bf Figure \ref{E56} after replacement}

\begin{center}

\begin{\tz}[scale=.66]
\draw (9,0) -- (23,0);
\draw (9,1) -- (11,3) -- (11,1);
\draw (14,0) -- (15,1);
\draw (15,2) -- (15,4);
\draw (17,5) -- (19,7) -- (19,3) -- (21,5);
\draw (19,2) -- (21,4) -- (21,3) -- (23,5);
\draw (23,6) -- (23,8);
\draw [red] (11,0) -- (12,1) -- (12,0) -- (14,2);
\draw [red] (18,3) -- (20,5.1) -- (20,4.3) -- (22,6.3);
\draw [blue] (14,0) -- (13,1);
\draw [blue] (15,1) -- (14,2);
\draw [blue] (19,2) -- (18,3);
\draw [blue] (20,3) -- (19.15,4.2);
\draw [blue] (21,4) -- (20,5.1);
\draw [blue] (21,3) -- (20,4.3);
\draw [blue] (22,4) -- (21.2,5.5);
\draw [blue] (23,5) -- (22,6.3);
\node at (9,1) {\lb};
\node at (10,2) {\lb};
\node at (11,3) {\lb};
\node at (11,2) {\lb};
\node at (11,1) {\lb};
\node at (12.8,1.2) {\lb};
\node at (14,0) {\lb};
\node at (15,1) {\lb};
\node at (15,2) {\lb};
\node at (15,3) {\lb};
\node at (15,4) {\lb};
\node at (17,5) {\lb};
\node at (18,6) {\lb};
\node at (19,7) {\lb};
\node at (19,6) {\lb};
\node at (19,5) {\lb};
\node at (19,3.8) {\lb};
\node at (19,3) {\lb};
\node at (19.8,3.8) {\lb};
\node at (21,5) {\lb};
\node at (19,2) {\lb};
\node at (20,3) {\lb};
\node at (21,4) {\lb};
\node at (21,3) {\lb};
\node at (22,4) {\lb};
\node at (23,5) {\lb};
\node at (23,6) {\lb};
\node at (23,7) {\lb};
\node at (23,8) {\lb};
\node [red] at (11,0) {\lb};
\node [red] at (12,1) {\lb};
\node [red] at (12,0) {\lb};
\node [red] at (13,1) {\lb};
\node [red] at (14,2) {\lb};
\node [red] at (18,3) {\lb};
\node [red] at (19.15,4.2) {\lb};
\node [red] at (20,5.1) {\lb};
\node [red] at (20,4.3) {\lb};
\node [red] at (21.2,5.5) {\lb};
\node [red] at (22,6.3) {\lb};
\node at (9,-.7) {$9$};
\node at (11,-.7) {$11$};
\node at (15,-.7) {$15$};
\node at (19,-.7) {$19$};
\node at (23,-.7) {$23$};

\end{\tz}
\end{center}
\end{fig}
\end{minipage}
\bigskip

The following proposition is the generalization of the phenomenon observed in Figure \ref{M78}.
\begin{prop} \label{gen} For $k\ge4$ and $i\ge0$, we have
$$M_k^i+\Si^2\Mh^{k+i-4}_4=M_{k+1}^i$$
in the sense that if $\Si^2\Mh_4^{k+i-4}$ is placed beneath $M_k^i$ so that there are $d_1$ differentials as in Figure \ref{d1}, then the result is $M_{k+1}^i$.
\end{prop}
\begin{proof} First note that all $M_k^i$ with the same mod 4 value of $k+i$ have the same general form as indicated in Figure \ref{same}. This refers to the way in which they leave filtration 0.

\bigskip
\begin{minipage}{6in}
\begin{fig}\label{same}

{\bf $M_k^i$}

\begin{center}

\begin{\tz}[scale=.27]
\draw (0,0) -- (13,0);
\draw (16,0) -- (28,0);
\draw (31,0) -- (43,0);
\draw (46,0) -- (58,0);
\draw (3,0) -- (3,1);
\draw [dashed] (3,1) -- (3,2.5);
\draw (3,2.5) -- (3,4) -- (1,2);
\draw (7,0) -- (7,2);
\draw [dashed] (7,2) -- (7,3.5);
\draw (7,3.5) -- (7,5);
\draw (13,2) -- (11,0) -- (11,2.5);
\draw [dashed] (11,2.5) -- (11,4.5);
\draw (11,4.5) -- (11,8) -- (9,6);
\draw [red] (14,-1) -- (14,9);
\node at (1,2) {\llb};
\node at (2,3) {\llb};
\node at (3,4) {\llb};
\node at (3,3) {\llb};
\node at (3,0) {\llb};
\node at (7,0) {\llb};
\node at (7,1) {\llb};
\node at (7,5) {\llb};
\node at (7,4) {\llb};
\node at (9,6) {\llb};
\node at (11,5) {\llb};
\node at (11,6) {\llb};
\node at (11,7) {\llb};
\node at (11,8) {\llb};
\node at (10,7) {\llb};
\node at (11,0) {\llb};
\node at (11,1) {\llb};
\node at (11,2) {\llb};
\node at (12,1) {\llb};
\node at (13,2) {\llb};
\node at (0,1) {$\iddots$};
\node at (13,5) {$\iddots$};
\node at (3,6) {$k+i\equiv1$};
\draw (16,2) -- (18,4) -- (18,2.5);
\draw (18,0) -- (18,1);
\draw [dashed] (18,1) -- (18,2.5);
\node at (15,1) {$\iddots$};
\node at (28,6) {$\iddots$};
\draw (22,0) -- (22,2);
\draw [dashed] (22,2) -- (22,3.5);
\draw (22,3.5) -- (22,5);
\draw (24,6) -- (26,8) -- (26,6);
\draw (28,3) -- (26,1) -- (26,3);
\draw [dashed] (26,3) -- (26,6);
\node at (18,6) {$k+i\equiv 0$};
\node at (16,2) {\llb};
\node at (17,3) {\llb};
\node at (18,4) {\llb};
\node at (18,3) {\llb};
\node at (18,1) {\llb};
\node at (18,0) {\llb};
\node at (22,0) {\llb};
\node at (22,1) {\llb};
\node at (22,2) {\llb};
\node at (22,4) {\llb};
\node at (22,5) {\llb};
\node at (28,3) {\llb};
\node at (27,2) {\llb};
\node at (26,1) {\llb};
\node at (26,2) {\llb};
\node at (26,3) {\llb};
\node at (24,6) {\llb};
\node at (25,7) {\llb};
\node at (26,8) {\llb};
\node at (26,7) {\llb};
\node at (26,6) {\llb};
\draw [red] (29,-1) -- (29,9);
\draw (31,2) -- (33,4) -- (33,2.5);
\draw [dashed] (33,1) -- (33,2.5);
\draw (33,0) -- (33,1);
\draw (37,1) -- (37,2.5);
\draw [dashed] (37,2.5) -- (37,4);
\draw (37,4) -- (37,5);
\draw (39,6) -- (41,8) -- (41,6);
\draw [dashed] (41,6) -- (41,4);
\draw (41,4) -- (41,2) -- (43,4);
\node at (31,2) {\llb};
\node at (32,3) {\llb};
\node at (33,4) {\llb};
\node at (33,3) {\llb};
\node at (33,0) {\llb};
\node at (33,1) {\llb};
\node at (35,0) {\llb};
\node at (37,1) {\llb};
\node at (37,2) {\llb};
\node at (37,4) {\llb};
\node at (37,5) {\llb};
\node at (39,6) {\llb};
\node at (40,7) {\llb};
\node at (41,8) {\llb};
\node at (41,7) {\llb};
\node at (41,6) {\llb};
\node at (41,2) {\llb};
\node at (41,3) {\llb};
\node at (42,3) {\llb};
\node at (43,4) {\llb};
\node at (30,1) {$\iddots$};
\node at (43,6) {$\iddots$};
\draw [red] (44,-1) -- (44,9);
\draw (46,2) -- (48,4) -- (48,2.5);
\draw (48,0) -- (48,1);
\draw [dashed] (48,1) -- (48,2.5);
\node at (33,6) {$k+i\equiv3$};
\node at (48,6) {$k+i\equiv2$};
\node at (41,4) {\llb};
\node at (52,3) {\llb};
\draw (49,0) -- (50,1);
\draw (52,2) -- (52,2.5);
\draw (52,5) -- (52,4);
\draw [dashed] (52,2.5) -- (52,4);
\draw (54,6) -- (56,8) -- (56,6);
\draw [dashed] (56,6) -- (56,4);
\draw (56,4) -- (56,3) -- (58,5);
\node at (45,1) {$\iddots$};
\node at (58,7) {$\iddots$};
\node at (46,2) {\llb};
\node at (47,3) {\llb};
\node at (48,4) {\llb};
\node at (48,3) {\llb};
\node at (48,0) {\llb};
\node at (48,1) {\llb};
\node at (49,0) {\llb};
\node at (50,1) {\llb};
\node at (52,2) {\llb};
\node at (52,5) {\llb};
\node at (52,4) {\llb};
\node at (58,5) {\llb};
\node at (57,4) {\llb};
\node at (56,3) {\llb};
\node at (56,4) {\llb};
\node at (56,6) {\llb};
\node at (56,7) {\llb};
\node at (56,8) {\llb};
\node at (55,7) {\llb};
\node at (54,6) {\llb};
\end{\tz}
\end{center}
\end{fig}
\end{minipage}
\bigskip

Then note that $k$ and $k+i$ can each be increased by 1 by the additions in Figure \ref{d1}, each of which illustrates $M_k^i+\Si^2\Mh_4^{k+i-4}$. Here we use that $\Mh_4^{4t+i}=\Si^{8t}\Mh_4^i$.
\end{proof}

\bigskip
\begin{minipage}{6in}
\begin{fig}\label{d1}

{\bf $M_k^i+\Si^2\Mh_4^{k+i-4}$}

\begin{center}

\begin{\tz}[scale=.33]
\draw (-.5,0) -- (8,0);
\draw (11.5,0) -- (20,0);
\draw (23.5,0) -- (34,0);
\draw (37.5,0) -- (49,0);
\draw (0,0) -- (0,1);
\draw [dashed] (0,1) -- (0,3);
\draw [dashed] (0,3) -- (0,4);
\draw (6,2) -- (4,0) -- (4,2);
\draw (2,5) -- (4,7) -- (4,5);
\draw [dashed] (4,2) -- (4,5);
\draw (8,3) -- (8,4);
\draw [dashed] (8,4) -- (8,7);
\draw (8,7) -- (8,8);
\draw [red] (5,0) -- (6,1) -- (6,0) -- (8,2);
\draw [blue] (6,0) -- (5,1);
\draw [blue] (7,1) -- (6,2);
\node at (0,0) {\llb};
\node at (0,1) {\llb};
\node at (0,3) {\llb};
\node at (0,4) {\llb};
\node at (2,5) {\llb};
\node at (3,6) {\llb};
\node at (4,7) {\llb};
\node at (4,6) {\llb};
\node at (4,5) {\llb};
\node at (4,2) {\llb};
\node at (4,1) {\llb};
\node at (4,0) {\llb};
\node at (5,1) {\llb};
\node at (6,2) {\llb};
\node at (8,3) {\llb};
\node at (8,4) {\llb};
\node at (8,8) {\llb};
\node at (8,7) {\llb};
\node at (4,-.8) {$8p+3$};
\node at (4,9) {$k+i=4p+5$};
\node [red] at (5,0) {\llb};
\node [red] at (6,1) {\llb};
\node [red] at (6,0) {\llb};
\node [red] at (7,1) {\llb};
\node [red] at (8,2) {\llb};
\draw [red] (10,-1) -- (10,10);
\draw (12,0) -- (12,1);
\draw [dashed] (12,1) -- (12,3);
\draw (12,3) -- (12,4);
\draw (14,5) -- (16,7) -- (16,5);
\draw (16,3) -- (16,1) -- (18,3);
\draw (20,4) -- (20,5);
\draw [dashed] (20,5) -- (20,7);
\draw (20,7) -- (20,8);
\draw [red] (16,0) -- (18,2) -- (18,1) -- (20,3);
\draw [blue] (18,1) -- (17,2);
\draw [blue] (19,2) -- (18,3);
\node at (16,-.8) {$8p+3$};
\node at (16,9) {$k+i=4p+4$};
\node at (12,0) {\llb};
\node at (12,1) {\llb};
\node at (12,3) {\llb};
\node at (12,4) {\llb};
\node at (14,5) {\llb};
\node at (15,6) {\llb};
\node at (16,7) {\llb};
\node at (16,6) {\llb};
\node at (16,5) {\llb};
\node at (18,3) {\llb};
\node at (17,2) {\llb};
\node at (16,1) {\llb};
\node at (16,2) {\llb};
\node at (16,3) {\llb};
\node at (20,4) {\llb};
\node at (20,5) {\llb};
\node at (20,8) {\llb};
\node at (20,7) {\llb};
\draw [dashed] (16,3) -- (16,5);
\node [red] at (16,0) {\llb};
\node [red] at (17,1) {\llb};
\node [red] at (18,2) {\llb};
\node [red] at (18,1) {\llb};
\node [red] at (19,2) {\llb};
\node [red] at (20,3) {\llb};
\draw (26,1) -- (26,2);
\draw [dashed] (26,2) -- (26,3);
\draw (26,3) -- (26,4);
\node at (26,2) {\llb};
\draw (28,5) -- (30,7) -- (30,5);
\draw (32,4) -- (30,2) -- (30,3.5);
\draw [dashed] (30,3.5) -- (30,5);
\draw [red] (25,-1) -- (26,0);
\draw [red] (30,1) -- (32,3) -- (32,2) -- (34,4);
\draw [blue] (25,-1) -- (24,0);
\draw [blue] (32,2) -- (31,3);
\draw [blue] (33,3) -- (32,4);
\node at (24,0) {\llb};
\node at (26,1) {\llb};
\node at (26,4) {\llb};
\node at (26,3) {\llb};
\node at (28,5) {\llb};
\node at (29,6) {\llb};
\node at (30,7) {\llb};
\node at (30,6) {\llb};
\node at (30,5) {\llb};
\node at (32,4) {\llb};
\node at (31,3) {\llb};
\node at (30,2) {\llb};
\node at (30,3) {\llb};
\node [red] at (25,-1) {\llb};
\node [red] at (26,0) {\llb};
\node [red] at (30,1) {\llb};
\node [red] at (31,2) {\llb};
\node [red] at (32,3) {\llb};
\node [red] at (32,2) {\llb};
\node [red] at (33,3) {\llb};
\node [red] at (34,4) {\llb};
\node at (30,-.8) {$8p+3$};
\node at (30,9) {$k+i=4p+3$};
\draw (34,5) -- (34,6);
\draw [dashed] (34,6) -- (34,7);
\draw (34,7) -- (34,8);
\node at (34,7) {\llb};
\node at (34,5) {\llb};
\node at (34,6) {\llb};
\node at (34,8) {\llb};
\draw (38,0) -- (39,1);
\draw [red] (22,-1) -- (22,10);
\draw [red] (36,-1) -- (36,10);
\draw (41,2) -- (41,3);
\draw [dashed] (41,3) -- (41,4);
\node at (41,3) {\llb};
\draw (43,5) -- (45,7) -- (45,5.5);
\draw [dashed] (45,5.5) -- (45,4);
\draw (45,4) -- (45,3) -- (47,5);
\draw (49,6) -- (49,7);
\draw [dashed] (49,7) -- (49,8);
\node at (49,7) {\llb};
\draw [red] (39,0) -- (39,-1) -- (41,1);
\draw [red] (45,2) -- (47,4) -- (47,3) -- (49,5);
\draw [blue] (39,-1) -- (38,0);
\draw [blue] (40,0) -- (39,1);
\draw [blue] (47,3) -- (46,4);
\draw [blue] (48,4) -- (47,5);
\node at (45,-.8) {$8p+3$};
\node at (44,9) {$k+i=4p+2$};
\node at (38,0) {\llb};
\node at (39,1) {\llb};
\node at (41,2) {\llb};
\node at (41,4) {\llb};
\node at (43,5) {\llb};
\node at (44,6) {\llb};
\node at (45,7) {\llb};
\node at (45,6) {\llb};
\node at (45,3) {\llb};
\node at (45,4) {\llb};
\node at (46,4) {\llb};
\node at (47,5) {\llb};
\node at (49,6) {\llb};
\node at (49,8) {\llb};
\node [red] at (39,0) {\llb};
\node [red] at (39,-1) {\llb};
\node [red] at (40,0) {\llb};
\node [red] at (41,1) {\llb};
\node [red] at (45,2) {\llb};
\node [red] at (46,3) {\llb};
\node [red] at (47,4) {\llb};
\node [red] at (47,3) {\llb};
\node [red] at (48,4) {\llb};
\node [red] at (49,5) {\llb};

\end{\tz}
\end{center}
\end{fig}
\end{minipage}
\bigskip

We saw in Figure \ref{E56} that $\E'_{5,6}$ is what was left over after $\Si M_4^5\lar\Si^8M_t^{8-t}$ is used to change $\Si^8M_t^{8-t}$ to $\Si^8M_{t-1}^{8-t}$, and then since $\Si^8M_t^{8-t}$ can be obtained from $\Si^8M_{t-1}^{8-t}+\Si^{10}\Mh_4^3$, we can obtain $\E'_{5,6}$ itself from $\Si M_4^5\lar\Si^{10}\Mh_4^3$. The same argument, together with Proposition \ref{gen}, yields that, for any $e>1$, $\E'_{e,e+1}$, the part left over when $\Si M_4^e\lar\Si^8M_t^{e+3-t}$ is used to change $\Si^8M_t^{e+3-t}$ to $\Si^8M_{t-1}^{e+3-t}$, can be obtained from $\Si M_4^e\lar\Si^2\Mh_4^{e-2}$, as claimed in Definition \ref{E'def}.

Next we consider $\E'_{4,6}$, which (\ref{A7tableau}) shows to be related to
$$\Si M_5^4\ml{4,6} \Si^8M_4^3\quad \Si^9M_4^5\ml{5,6}\Si^{16}M_t^{8-t} $$
for $t\in\{6,7,8\}$. We have already seen that the resultant of the second arrow is $\Si^8\E'_{5,6}+\Si^{16}M_{t-1}^{8-t}$. The $\E'_{5,6}$ is not involved in finding $\E'_{4,6}$. (It may be involved in differentials which change $\E'$ to $\E$; this will be considered in Section \ref{pfsec}.) We will show in the next paragraph that the resultant of $\Si M_5^4\lar\Si^8M_4^3\lar\Si^{16}M_{t-1}^{8-t}$ is $\Si^{16}M_{t-2}^{8-t}$ plus a remainder term, which is defined to be $\E'_{4,6}$. Then we can use Proposition \ref{gen} to replace $\Si^{16}M_{t-1}^{8-t}$ by $\Si^{16}M_{t-2}^{8-t}+\Si^{18}\Mh_4^2$, and deduce, similarly to what we did before, that $\E'_{4,6}$ itself can be obtained from
$$\Si M_5^4\lar \Si^8M_4^3\lar\Si^{18}\Mh_4^2.$$

In this paragraph, we use $t=8$ in the above discussion, but $t=6$ or 7 work similarly. In the left side of Figure \ref{E46}, we show $\Si M_5^4\lar\Si^8M_4^3$, with $\Si^8M_4^3$ in red. On the right side, we place the resultant of the left side in black and $\Si^{16}M_7^0$ in red, placed so that the differential appears as a $d_1$. (The justification for the differentials will be given in Section \ref{pfsec}.) The differential reduces the $\Si^{16}M_7^0$ to $\Si^{16}M_6^0$, as claimed in the preceding paragraph. The classes in 20 and 21 are part of $\Si^{16}M_6^0$. The black classes in grading $\le19$ form $\E'_{4,6}$. A similar thing happens in consideration of any $\E'_{e,e+2}$. It will come from
$$\Si M_5^e\lar \Si^8M_4^{e-1}\quad \Si^9M_4^{e+1}\lar\Si^{16}M_t^{e+4-t}.$$
The second arrow produces $\Si^8\E'_{e+1,e+2}\oplus \Si^{16}M_{t-1}^{e+4-t}$.
After some initial deviation, the first arrow will produce pure lightning flashes in time to work with $\Si^{16}M_{t-1}^{e+4-t}$ to form $\Si^{16}M_{t-2}^{e+4-t}$.

\bigskip
\begin{minipage}{6in}
\begin{fig}\label{E46}

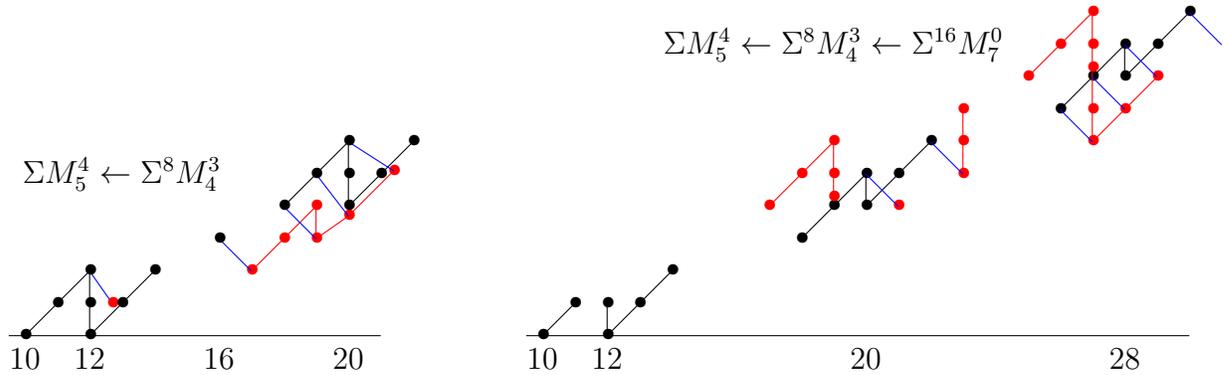
{\bf Forming $\E'_{4,6}+\Si^{16}M_6^0$}

\begin{center}

\begin{\tz}[scale=.43]
\node at (13,5) {$\Si M_5^4\lar \Si^8M_4^3$};
\node at (35,9) {$\Si M_5^4\lar \Si^8M_4^3\lar \Si^{16}M_7^0$};
\draw (10,0) -- (12,2) -- (12,0) -- (14,2);
\draw (18,4) -- (20,6) -- (20,4) -- (22,6);
\draw (9.5,0) -- (21,0);
\draw (25.5,0) -- (46,0);
\draw (26,0) -- (27,1);
\draw (28,1) -- (28,0) -- (30,2);
\draw (34,3) -- (36,5) -- (36,4) -- (38,6);
\draw (42,7) -- (44,9) -- (44,8) -- (46,10);
\draw [red] (17,2) -- (19,4) -- (19,3) -- (20,3.7) -- (21.4,5.1);
\draw [blue] (12.7,1) -- (12,2);
\node [red] at (12.7,1) {\lb};
\node [red] at (21.4,5.1) {\lb};
\node [red] at (20,3.7) {\lb};
\node [red] at (17,2) {\lb};
\node [red] at (18,3) {\lb};
\node [red] at (19,4) {\lb};
\node [red] at (19,3) {\lb};
\draw [blue] (17,2) -- (16,3);
\draw [blue] (19,3) -- (18,4);
\draw [blue] (20,3.7) -- (19,5);
\draw [blue] (21.4,5.1) -- (20,6);
\node at (10,0) {\lb};
\node at (11,1) {\lb};
\node at (12,2) {\lb};
\node at (12,1) {\lb};
\node at (12,0) {\lb};
\node at (13,1) {\lb};
\node at (14,2) {\lb};
\node at (16,3) {\lb};
\node at (18,4) {\lb};
\node at (19,5) {\lb};
\node at (20,6) {\lb};
\node at (20,5) {\lb};
\node at (20,4) {\lb};
\node at (21,5) {\lb};
\node at (22,6) {\lb};
\node at (10,-.7) {$10$};
\node at (12,-.7) {$12$};
\node at (16,-.7) {$16$};
\node at (20,-.7) {$20$};
\node at (26,-.7) {$10$};
\node at (28,-.7) {$12$};
\node at (36,-.7) {$20$};
\node at (44,-.7) {$28$};
\node at (26,0) {\lb};
\node at (27,1) {\lb};
\node at (28,1) {\lb};
\node at (28,0) {\lb};
\node at (29,1) {\lb};
\node at (30,2) {\lb};
\node at (34,3) {\lb};
\node at (35,4) {\lb};
\node at (36,5) {\lb};
\node at (36,4) {\lb};
\node at (37,5) {\lb};
\node at (38,6) {\lb};
\node at (42,7) {\lb};
\node at (43,8) {\lb};
\node at (44,9) {\lb};
\node at (44,8) {\lb};
\node at (45,9) {\lb};
\node at (46,10) {\lb};
\draw [red] (33,4) -- (35,6) -- (35,4.3);
\draw [red] (39,5) -- (39,7);
\draw [red] (41,8) -- (43,10) -- (43,6) -- (45,8);
\node [red] at (33,4) {\lb};
\node [red] at (34,5) {\lb};
\node [red] at (35,6) {\lb};
\node [red] at (35,5) {\lb};
\node [red] at (35,4.3) {\lb};
\node [red] at (37,4) {\lb};
\node [red] at (39,5) {\lb};
\node [red] at (39,6) {\lb};
\node [red] at (39,7) {\lb};
\node [red] at (41,8) {\lb};
\node [red] at (42,9) {\lb};
\node [red] at (43,10) {\lb};
\node [red] at (43,9) {\lb};
\node [red] at (43,8.3) {\lb};
\node [red] at (43,7) {\lb};
\node [red] at (43,6) {\lb};
\node [red] at (44,7) {\lb};
\node [red] at (45,8) {\lb};
\draw [blue] (39,5) -- (38,6);
\draw [blue] (37,4) -- (36,5);
\draw [blue] (43,6) -- (42,7);
\draw [blue] (44,7) -- (43,8);
\draw [blue] (45,8) -- (44,9);
\draw [blue] (46,10) -- (47,9);

\end{\tz}

\end{center}
\end{fig}
\end{minipage}
\bigskip

A similar argument works for any $\E'_{e,e+d}$. We are deriving a definition, Definition \ref{E'def}. The tableau (\ref{A7tableau}) is very helpful in getting started. For $e\ge2$ and $d\ge1$, there are $\binom{k-2-d}{e-2}$ arrows labeled $(e,e+d)$ in the tableau for $\At_k$. Each gives rise to an occurrence of $\E'_{e,e+d}$. The structure of $\E'_{e,e+d}$ is determined by $2^{d-1}$ arrows beginning with  one labeled $(e,e+d)$. The sequence of these arrow labels is the same for every occurrence of $(e,e+d)$, and the sequence of the $2^d$ $M$ charts is also the same except for the very last one, which can be $\Si^{2^{d+2}}M_t^{e+d+2-t}$ for $t=e+d+2\le k+2$ or $d+4\le t\le e+d+1\le k$.

We take a major step toward justifying Definition \ref{E'def} by considering $\E'_{e,e+3}$. This is determined by the four arrows in (\ref{e+3}),
\begin{equation}\label{e+3}\overset{e,e+3}{\Si M^e_6\lar \Si^8M_4^{e-1}}\quad\overset{e+1,e+2}{\Si^9M_4^{e+1}\lar\Si^{16}M_5^{e-1}}\quad\overset{e+1,e+3}{\Si^{17}M_5^{e+1}\lar\Si^{24}M_4^e}\quad\overset{e+2,e+3}{\Si^{25}M_4^{e+2}\lar\Si^{32}M_t^{e+5-t}}\end{equation}
where either $t=e+5\le k+2$ or $7\le t\le e+4\le k$. This can be seen in (\ref{Aktableau}), where it will be suspended by $\Si^{64b+32}$
for some $b$. Then (\ref{e+3}) corresponds to $i=8b+(4,5,6,7)$ in (\ref{Aktableau}). Also $e=\a(b)+2$ and $M_t^{e+5-t}=M_{7+\nu(b+1)}^{\a(b+1)-1}$.

Similarly to our analysis of $\E'_{5,6}$ and $\E'_{4,6}$, the resultant of the second arrow in (\ref{e+3}) is $\Si^8\E'_{e+1,e+2}\oplus \Si^{16}M_4^{e-1}$, and the resultant of the last two arrows is $\Si^{16}\E'_{e+1,e+3}+\Si^{24}\E'_{e+2,e+3}\oplus \Si^{32}M_{t-2}^{e+5-t}$. We leave the $\Si^8\E'_{e+1,e+2}$, $\Si^{16}\E'_{e+1,e+3}$, and 
$\Si^{24}\E'_{e+2,e+3}$ behind, and are left with
\begin{equation}\label{for}\Si M_6^e\lar \Si^8M_4^{e-1}\lar \Si^{16}M_4^{e-1}\lar\Si^{32}M_{t-2}^{e+5-t}.\end{equation} In Figure \ref{thr}, which is the case $e=4$, we show the result of the first two arrows.

\bigskip
\begin{minipage}{6in}
\begin{fig}\label{thr}

{\bf First two arrows of (\ref{for})}

\begin{center}

\begin{\tz}[scale=.4]
\draw (9.5,0) -- (22,0);
\draw (27.5,0) -- (48,0);
\draw (10,0) -- (12,2) -- (12,0);
\draw (13,0) -- (14,1);
\draw (16,2) -- (16,3);
\draw (18,4) -- (20,6) -- (20,3) -- (22,5);
\draw [red] (17,2) -- (19,4) -- (19,3) -- (21,5);
\draw [blue] (13,1) -- (12,2);
\draw [blue] (17,2) -- (16,3);
\draw [blue] (19,3) -- (18,4);
\draw [blue] (19.8,3.8) -- (19,5);
\draw [blue] (21,5) -- (20,6);
\node at (10,-.7) {$10$};
\node at (12,-.7) {$12$};
\node at (16,-.7) {$16$};
\node at (20,-.7) {$20$};
\node at (10,0) {\lb};
\node at (11,1) {\lb};
\node at (12,2) {\lb};
\node at (12,1) {\lb};
\node at (12,0) {\lb};
\node at (13,0) {\lb};
\node at (14,1) {\lb};
\node at (16,2) {\lb};
\node at (16,3) {\lb};
\node at (18,4) {\lb};
\node at (19,5) {\lb};
\node at (20,6) {\lb};
\node at (20,5) {\lb};
\node at (20,4.1) {\lb};
\node at (20,3) {\lb};
\node at (21,4) {\lb};
\node at (22,5) {\lb};
\node [red] at (13,1) {\lb};
\node [red] at (17,2) {\lb};
\node [red] at (18,3) {\lb};
\node [red] at (19,4) {\lb};
\node [red] at (19,3) {\lb};
\node [red] at (19.8,3.8) {\lb};
\node [red] at (21,5) {\lb};
\draw (28,0) -- (29,1);
\draw (30,0) -- (30,1);
\draw (31,0) -- (32,1);
\draw (36,3) -- (38,5) -- (38,3) -- (40,5);
\draw (44,7) -- (46,9) -- (46,7) -- (48,9);
\draw [red] (43,5) -- (45,7) -- (45,6) -- (46,6.7) -- (47.4,8.1);
\draw [blue] (38.8,4.2) -- (38,5);
\draw [blue] (43,5) -- (42,6);
\draw [blue] (45,6) -- (44,7);
\draw [blue] (46,6.7) -- (45,8);
\draw [blue] (47.4,8.1) -- (46,9);
\node at (28,0) {\lb};
\node at (29,1) {\lb};
\node at (30,0) {\lb};
\node at (30,1) {\lb};
\node at (31,0) {\lb};
\node at (32,1) {\lb};
\node at (34,2) {\lb};
\node at (36,3) {\lb};
\node at (37,4) {\lb};
\node at (38,5) {\lb};
\node at (38,4) {\lb};
\node at (38,3) {\lb};
\node at (39,4) {\lb};
\node at (40,5) {\lb};
\node at (42,6) {\lb};
\node at (44,7) {\lb};
\node at (45,8) {\lb};
\node at (46,9) {\lb};
\node at (46,8) {\lb};
\node at (46,7) {\lb};
\node at (47,8) {\lb};
\node at (48,9) {\lb};
\node [red] at (38.8,4.2) {\lb};
\node [red] at (43,5) {\lb};
\node [red] at (44,6) {\lb};
\node [red] at (45,7) {\lb};
\node [red] at (45,6) {\lb};
\node [red] at (46,6.7) {\lb};
\node [red] at (47.4,8.1) {\lb};
\node at (28,-.7) {$10$};
\node at (30,-.7) {$12$};
\node at (38,-.7) {$20$};
\node at (46,-.7) {$28$};
\node at (14,6) {$\Si M_6^4\lar \Si^8M_4^3$};
\node at (37,7) {$\Si M_6^4\lar \Si^8M_4^3\lar \Si^{16}M_4^3$};
\end{\tz}
\end{center}
\end{fig}
\end{minipage}
\bigskip

Then in Figure \ref{E3}, which is the case $t=9$, we place $\Si^{32}M_{t-2}^{e+5-t}$ so that it has a $d_1$ differential into the stable lightning flash beyond the end of Figure \ref{thr}. This reduces it to $\Si^{32}M_{t-3}^{e+5-t}$, and the leftover part is defined to be $\E'_{e,e+3}$. Note that the classes in the lightning flash in grading 36 and 37 become part of $\Si^{32}M_6^0$, while the classes in 34 and 35 are the last two elements in $\E'_{4,7}$. So the result of (\ref{for}) is $\Si^{32}M_{t-3}^{e+5-t}\oplus \E'_{e,e+3}$.

\bigskip
\begin{minipage}{6in}
\begin{fig}\label{E3}

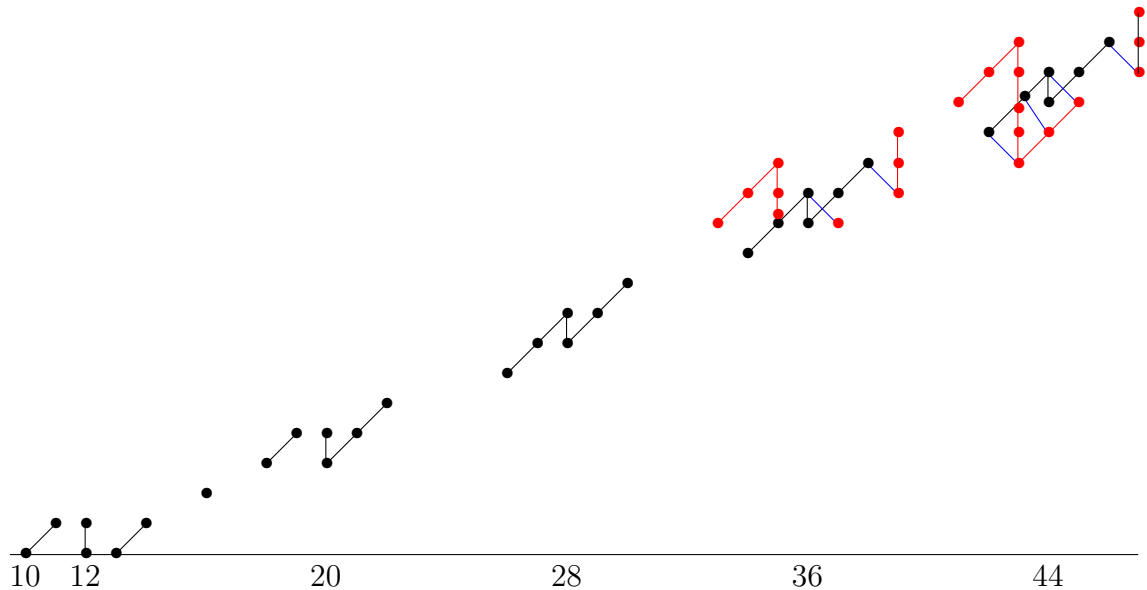
{\bf Forming $\Si^{32}M_{t-3}^{e+5-t}\oplus \E'_{e,e+3}$}

\begin{center}

\begin{\tz}[scale=.4]
\draw (9.5,0) -- (47,0);
\draw (10,0) -- (11,1);
\draw (12,0) -- (12,1);
\draw (13,0) -- (14,1);
\draw (18,3) -- (19,4);
\draw (20,4) -- (20,3) -- (22,5);
\draw (26,6) -- (28,8) -- (28,7) -- (30,9);
\draw (34,10) -- (36,12) -- (36,11) -- (38,13);
\draw (42,14) -- (44,16) -- (44,15) -- (46,17);
\draw [red] (33,11) -- (35,13) -- (35,11.3);
\draw [red] (39,12) -- (39,14);
\draw [red] (41,15) -- (43,17) -- (43,13) -- (45,15);
\draw [blue] (37,11) -- (36,12);
\draw [blue] (39,12) -- (38,13);
\draw [blue] (43,13) -- (42,14);
\draw [blue] (45,15) -- (44,16);
\draw [blue] (44,14) -- (43.2,15.2);
\draw [blue] (47,16) -- (46,17);
\node at (10,0) {\lb};
\node at (11,1) {\lb};
\node at (12,0) {\lb};
\node at (12,1) {\lb};
\node at (13,0) {\lb};
\node at (14,1) {\lb};
\node at (16,2) {\lb};
\node at (18,3) {\lb};
\node at (19,4) {\lb};
\node at (20,4) {\lb};
\node at (20,3) {\lb};
\node at (21,4) {\lb};
\node at (22,5) {\lb};
\node at (26,6) {\lb};
\node at (27,7) {\lb};
\node at (28,8) {\lb};
\node at (28,7) {\lb};
\node at (29,8) {\lb};
\node at (30,9) {\lb};
\node at (34,10) {\lb};
\node at (35,11) {\lb};
\node at (36,12) {\lb};
\node at (36,11) {\lb};
\node at (37,12) {\lb};
\node at (38,13) {\lb};
\node at (42,14) {\lb};
\node at (43.2,15.2) {\lb};
\node at (44,16) {\lb};
\node at (44,15) {\lb};
\node at (45,16) {\lb};
\node at (46,17) {\lb};
\node [red] at (33,11) {\lb};
\node [red] at (34,12) {\lb};
\node [red] at (35,13) {\lb};
\node [red] at (35,12) {\lb};
\node [red] at (35,11.3) {\lb};
\node [red] at (37,11) {\lb};
\node [red] at (39,12) {\lb};
\node [red] at (39,13) {\lb};
\node [red] at (39,14) {\lb};
\node [red] at (47,16) {\lb};
\node [red] at (47,17) {\lb};
\node [red] at (47,18) {\lb};
\node [red] at (41,15) {\lb};
\node [red] at (42,16) {\lb};
\node [red] at (43,17) {\lb};
\node [red] at (43,16) {\lb};
\node [red] at (43,14.8) {\lb};
\node [red] at (43,14) {\lb};
\node [red] at (43,13) {\lb};
\node [red] at (44,14) {\lb};
\node [red] at (45,15) {\lb};
\node at (10,-.7) {$10$};
\node at (12,-.7) {$12$};
\node at (20,-.7) {$20$};
\node at (28,-.7) {$28$};
\node at (36,-.7) {$36$};
\node at (44,-.7) {$44$};
\draw (47,16) -- (47,18);
\end{\tz}
\end{center}
\end{fig}
\end{minipage}
\bigskip

We use Proposition \ref{gen} to replace $\Si^{32}M_{t-2}^{e+5-t}$ in (\ref{for}) by $\Si^{32}M_{t-3}^{e+5-t}+\Si^{34}\Mh_4^{e-2}$, similarly to Figure \ref{samp}. Since the result of
$$\Si M_6^e\lar\Si^8M_4^{e-1}\lar\Si^{16}M_4^{e-1}\lar(\Si^{32}M_{t-3}^{e+5-t}+\Si^{34}\Mh_4^{e-2})$$
equals $\Si^{32}M_{t-3}^{e+5-t}\oplus\E'_{e,e+3}$, and, as in Figure \ref{samp}, the two copies of $\Si^{32}M_{t-3}^{e+5-t}$ split off, we derive the $\ell=e+3$ case of Definition \ref{E'def}.

The general case $\E'_{e,e+d}$ does not differ significantly from the case $d=3$ just considered. The arrows involve $i=2^db+(2^{d-1},\ldots,2^d-1)$ in (\ref{Aktableau}) with $e=\a(b)+2$. Rather than constructing a full-blown induction argument, we will explain how the result for $d=4$ follows from the argument for $d<4$. The general case $\E'_{e,e+4}$ is embodied in the last two rows of (\ref{A7tableau}), which are $\Si^{192}$ times the case with $e=3$. For arbitrary $e$, increase all superscripts of $M$ by $e-3$. We subtract 192 from all the suspension parameters and discuss the eight arrows. By the work we have already performed, the result of the second arrow is $\E'_{4,5}\oplus\Si^{16}M_4^2$, the result of the third and fourth arrows is $\E'_{4,6}+\E'_{5,6}\oplus\Si^{32}M_4^2$, and that of the final four arrows is $\E'_{4,7}+\E'_{5,6}+\E'_{5,7}+\E'_{6,7}\oplus\Si^{32}M_6^0$. (We distinguish between $+$ and $\oplus$ because the $\E'$'s can have differentials into one another, but not into the $M$.) We split off all the $\E'$'s and are left with
\begin{equation}\label{37}\Si M_7^3\lar\Si^8M_4^2\lar \Si^{16}M_4^2\lar \Si^{32}M_4^2\lar\Si^{64}M_6^0.\end{equation}
The resultant of the first three arrows will, after initial deviations, stabilize into a sequence of lightning flashes well before the $\Si^{64}M_6^0$ has started. We claim that the lightning flashes will change it to $\Si^{64}M_5^0$, and $\E'_{3,7}$ is defined to be everything left over.

Considering the slightly subtle phenomenon observed in Figure \ref{E3} regarding the classes in 34, 35, 36, and 37, scrutiny is warranted here. We claim that if $M_k^i$ is placed into a sequence of lightning flashes so that its generators of (full) order $2^{k-2}$ hit generators of lightning flashes with a $d_1$, then the result is $M_{k-1}^i$ plus the first $(k+i-1)$ mod 4 gradings of the last lightning flash before total incorporation. We illustrate this in Figure \ref{k+i}, in which bigger red dots are incorporated into the $M$.

Recall from Figure \ref{same} that the way $M_k^i$ leaves filtration 0 depends only on $(k+i)$ mod 4. In Figure \ref{k+i}, we see that these differentials decrease the $(k+i)$ mod 4 type by 1 and decrease $k$ by 1, so they leave $i$ unchanged.

Now we explain why the ``$(k+i-1)$ mod 4 gradings'' noted above is consistent with everything else. As discussed just before (\ref{e+3}), the sequence of $2^d$ charts which define $\E'_{e,e+d}$ ends with $\Si^{2^{d+2}}M_t^{e+d+2-t}$.
When this sequence is reduced to a sequence of the type illustrated by (\ref{for}), the subscript of $M$ will be reduced $d-1$ times, so that the sequence of type (\ref{for}) will end with $M_k^i$ satisfying $k+i=e+3$. Thus the first $(e+2)$ mod 4 gradings of the last lightning flash before total incorporation will be left behind for the end of $\E'_{e,e+d}$. This is consistent with the endings of Figures \ref{four}, \ref{2,6}, \ref{3,7}, \ref{4,8}, and \ref{5,9}.

\bigskip
\begin{minipage}{6in}
\begin{fig}\label{k+i}

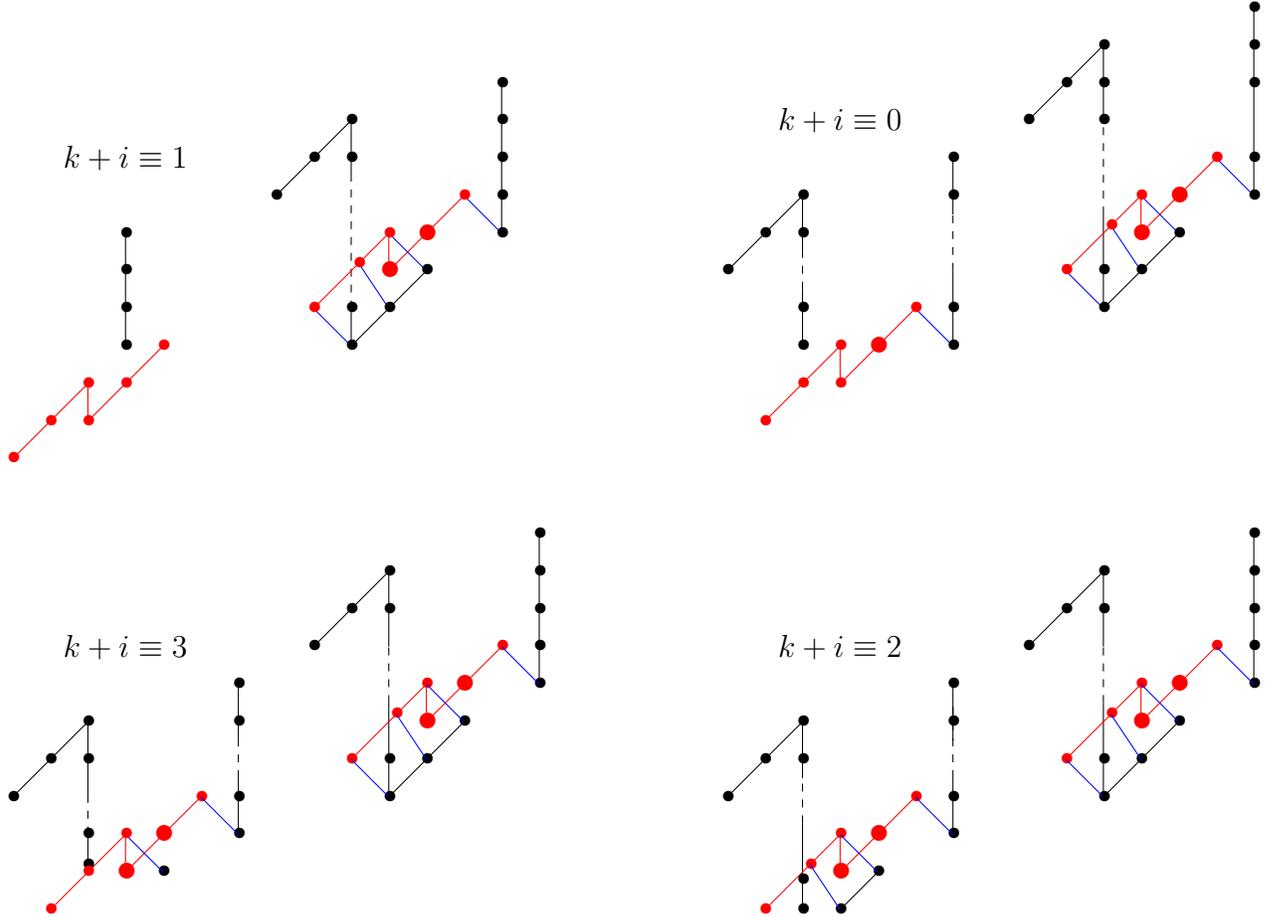
{\bf Cutting off $M_k^i$ with lightning flashes}

\begin{center}

\begin{\tz}[scale=.5]
\draw (0,3) -- (2,5) -- (2,3);
\draw (2,1.2) -- (2,2);
\draw [dashed] (2,2) -- (2,3);
\draw (6,2) -- (6,3.5);
\draw (6,4.5) -- (6,6);
\draw [dashed] (6,3.5) -- (6,4.5);
\draw (8,7) -- (10,9) -- (10,7);
\draw (12,5) -- (10,3) -- (10,5.5);
\draw [dashed] (10,5.5) -- (10,7);
\draw [red] (1,0) -- (3,2) -- (3,1) -- (5,3);
\draw [red] (9,4) -- (11,6) -- (11,5) -- (13,7);
\draw (14,6) -- (14,10);
\node at (0,3) {\lb};
\node at (1,4) {\lb};
\node at (2,5) {\lb};
\node at (2,4) {\lb};
\node at (2,1.2) {\lb};
\node at (2,2) {\lb};
\node at (4,1) {\lb};
\node at (6,2) {\lb};
\node at (6,3) {\lb};
\node at (6,6) {\lb};
\node at (6,5) {\lb};
\node at (10,3) {\lb};
\node at (10,4) {\lb};
\node at (11,4) {\lb};
\node at (12,5) {\lb};
\node at (8,7) {\lb};
\node at (9,8) {\lb};
\node at (10,9) {\lb};
\node at (10,8) {\lb};
\node at (14,6) {\lb};
\node at (14,7) {\lb};
\node at (14,8) {\lb};
\node at (14,9) {\lb};
\node at (14,10) {\lb};
\draw [blue] (4,1) -- (3,2);
\node [red] at (1,0) {\lb};
\node [red] at (2,1) {\lb};
\node [red] at (3,2) {\lb};
\node [red] at (5,3) {\lb};
\node [red] at (9,4) {\lb};
\node [red] at (10.2,5.2) {\lb};
\node [red] at (11,6) {\lb};
\node [red] at (13,7) {\lb};
\draw [blue] (6,2) -- (5,3);
\draw [blue] (10,3) -- (9,4);
\draw [blue] (11,4) -- (10.2,5.2);
\draw [blue] (12,5) -- (11,6);
\draw [blue] (14,6) -- (13,7);
\node [red] at (3,1) {\blb};
\node [red] at (4,2) {\blb};
\node [red] at (11,5) {\blb};
\node [red] at (12,6) {\blb};
\node at (3,7) {$k+i\equiv3$};
\draw (19,3) -- (21,5) -- (21,3.5);
\draw (21,0) -- (21,2.3);
\draw [dashed] (21,2.3) -- (21,3.5);
\draw (22,0) -- (23,1);
\draw (25,2) -- (25,3.5);
\draw (25,6) -- (25,4.5);
\draw [dashed] (25,3.5) -- (25,5.5);
\draw (27,7) -- (29,9) -- (29,7);
\draw (31,5) -- (29,3) -- (29,5.5);
\draw [dashed] (29,5.5) -- (29,7);
\draw (33,6) -- (33,10);
\draw [red] (20,0) -- (22,2) -- (22,1) -- (24,3);
\draw [red] (28,4) -- (30,6) -- (30,5) -- (32,7);
\node at (19,3) {\lb};
\node at (20,4) {\lb};
\node at (21,5) {\lb};
\node at (21,4) {\lb};
\node at (21,0) {\lb};
\node at (21,.8) {\lb};
\node at (22,0) {\lb};
\node at (23,1) {\lb};
\node at (25,2) {\lb};
\node at (25,3) {\lb};
\node at (25,6) {\lb};
\node at (25,5) {\lb};
\node at (27,7) {\lb};
\node at (28,8) {\lb};
\node at (29,9) {\lb};
\node at (29,8) {\lb};
\node at (31,5) {\lb};
\node at (30,4) {\lb};
\node at (29,3) {\lb};
\node at (29,4) {\lb};
\node at (33,6) {\lb};
\node at (33,7) {\lb};
\node at (33,8) {\lb};
\node at (33,9) {\lb};
\node at (33,10) {\lb};
\node [red] at (20,0) {\lb};
\node [red] at (21.2,1.2) {\lb};
\node [red] at (22,2) {\lb};
\node [red] at (24,3) {\lb};
\node [red] at (28,4) {\lb};
\node [red] at (29.2,5.2) {\lb};
\node [red] at (30,6) {\lb};
\node [red] at (32,7) {\lb};
\node [red] at (22,1) {\blb};
\node [red] at (23,2) {\blb};
\node [red] at (30,5) {\blb};
\node [red] at (31,6) {\blb};
\draw [blue] (22,0) -- (21.2,1.2);
\draw [blue] (23,1) -- (22,2);
\draw [blue] (25,2) -- (24,3);
\draw [blue] (29,3) -- (28,4);
\draw [blue] (30,4) -- (29.2,5.2);
\draw [blue] (31,5) -- (30,6);
\draw [blue] (33,6) -- (32,7);
\node at (22,7) {$k+i\equiv2$};
\draw (3,15) -- (3,18);
\draw (7,19) -- (9,21) -- (9,19.5);
\draw (11,17) -- (9,15) -- (9,16);
\draw (13,18) -- (13,22);
\draw [dashed] (9,16) -- (9,19.5);
\draw [red] (0,12) -- (2,14) -- (2,13) -- (4,15);
\draw [red] (8,16) -- (10,18) -- (10,17) -- (12,19);
\draw [blue] (9,15) -- (8,16);
\draw [blue] (10,16) -- (9.2,17.2);
\draw [blue] (11,17) -- (10,18);
\draw [blue] (13,18) -- (12,19);
\node at (3,15) {\lb};
\node at (3,16) {\lb};
\node at (3,17) {\lb};
\node at (3,18) {\lb};
\node at (7,19) {\lb};
\node at (8,20) {\lb};
\node at (9,21) {\lb};
\node at (9,20) {\lb};
\node at (9,15) {\lb};
\node at (9,16) {\lb};
\node at (10,16) {\lb};
\node at (11,17) {\lb};
\node at (13,18) {\lb};
\node at (13,19) {\lb};
\node at (13,20) {\lb};
\node at (13,21) {\lb};
\node at (13,22) {\lb};
\node [red] at (0,12) {\lb};
\node [red] at (1,13) {\lb};
\node [red] at (2,14) {\lb};
\node [red] at (2,13) {\lb};
\node [red] at (3,14) {\lb};
\node [red] at (4,15) {\lb};
\node [red] at (8,16) {\lb};
\node [red] at (9.2,17.2) {\lb};
\node [red] at (10,18) {\lb};
\node [red] at (12,19) {\lb};
\node [red] at (10,17) {\blb};
\node [red] at (11,18) {\blb};
\node at (3,20) {$k+i\equiv1$};
\draw (19,17) -- (21,19) -- (21,17.5);
\draw (21,15) -- (21,16.5);
\draw (25,15) -- (25,17);
\draw (25,20) -- (25,18.5);
\draw [dashed] (25,17) -- (25,18.5);
\draw [dashed] (21,16.5) -- (21,17.5);
\draw (27,21) -- (29,23) -- (29,21);
\draw (31,18) -- (29,16) -- (29,18.5);
\draw [dashed] (29,18.5) -- (29,21);
\draw (33,19) -- (33,24);
\draw [red] (20,13) -- (22,15) -- (22,14) -- (24,16);
\draw [red] (28,17) -- (30,19) -- (30,18) -- (32,20);
\draw [blue] (25,15) -- (24,16);
\draw [blue] (29,16) -- (28,17);
\draw [blue] (30,17) -- (29.2,18.2);
\draw [blue] (31,18) -- (30,19);
\draw [blue] (33,19) -- (32,20);
\node at (19,17) {\lb};
\node at (20,18) {\lb};
\node at (21,19) {\lb};
\node at (21,18) {\lb};
\node at (21,15) {\lb};
\node at (21,16) {\lb};
\node at (25,15) {\lb};
\node at (25,16) {\lb};
\node at (25,20) {\lb};
\node at (25,19) {\lb};
\node at (27,21) {\lb};
\node at (28,22) {\lb};
\node at (29,23) {\lb};
\node at (29,22) {\lb};
\node at (29,21) {\lb};
\node at (31,18) {\lb};
\node at (30,17) {\lb};
\node at (29,16) {\lb};
\node at (29,17) {\lb};
\node at (33,19) {\lb};
\node at (33,20) {\lb};
\node at (33,24) {\lb};
\node at (33,22) {\lb};
\node at (33,23) {\lb};
\node [red] at (20,13) {\lb};
\node [red] at (21,14) {\lb};
\node [red] at (22,15) {\lb};
\node [red] at (22,14) {\lb};
\node [red] at (24,16) {\lb};
\node [red] at (28,17) {\lb};
\node [red] at (29.2,18.2) {\lb};
\node [red] at (30,19) {\lb};
\node [red] at (32,20) {\lb};
\node [red] at (23,15) {\blb};
\node [red] at (30,18) {\blb};
\node [red] at (31,19) {\blb};
\node at (22,21) {$k+i\equiv0$};

\end{\tz}
\end{center}
\end{fig}
\end{minipage}
\bigskip

Returning to $\E'_{3,7}$, (\ref{37}) results in $\Si^{64}M_5^0\oplus\E'_{3,7}$. Replace $\Si^{64}M_6^0$ by $\Si^{64}M_5^0\oplus\Si^{66}\Mh_4^2$ in (\ref{37}), using Proposition \ref{gen}. The $\Si^{64}M_5^0$ splits off, similarly to Figure \ref{samp}, and we deduce that the result of
$$\Si M_7^3\lar\Si^8M_4^2\lar \Si^{16}M_4^2\lar \Si^{32}M_4^2\lar\Si^{66}\Mh_4^2$$
is $\E'_{3,7}$. The same argument applies if the $M_9^0$ at the end of (\ref{A7tableau}) is replaced by $M_8^1$, the other way that the sequence of arrows defining $\E'_{3,7}$ can end. Then the $M_6^0$ in (\ref{37}) becomes $M_5^1$, with the same value of $k+i$. If $e$ is changed, then all superscripts in (\ref{37}) are changed by the same amount, and the replacement argument using Proposition \ref{gen} still applies.

This completes the induction step for the case $d=4$. The same argument works for arbitrary $d$.
The conclusion of this section is the following theorem.
\begin{thm}\label{edgthm} A sequence
$$\Si M^e_{\ell-e+3}\lar \Si^8M_4^{e-1}\lar\Si^{16}M_4^{e-1}\lar\cdots\lar\Si^{2^{\ell-e+1}}M_4^{e-1}\lar\Si^{2^{\ell-e+2}}M^i_{k+1}$$
with all arrows being $d_1$ on generators in grading $1$, $3$, $4$, and $5$ mod $8$, and $\eta$ inserted, and satisfying $k+i=e+2$ equals $\E'_{e,\ell}\oplus \Si^{2^{\ell-e+2}}M_k^i$.
\end{thm}
\ni The theorem also holds if the sequence and the conclusion are suspended by any amount.

Now we justify Definition \ref{E1'def}. The description of $V_k$ in Definition \ref{E1'def} is a straightforward generalization of Figure \ref{UNpic}. For the sequence  (\ref{E1seq}) in Definition \ref{E1'def}, we refer to (\ref{A7tableau}) for the case $k=7$. We have seen in the analysis above that the $\ml{2,3}$ yields $\E'_{2,3}\oplus \Si^{16}M_4^0$, and that the last two arrows in the first row of (\ref{A7tableau}) yield $\E'_{2,4}+\E'_{3,4}\oplus\Si^{32}M_4^0$. The second row of (\ref{A7tableau}) yields $\E'_{2,5}+\E'_{3,4}+\E'_{3,5}+\E'_{4,5}\oplus \Si^{64}M_4^0$, as the $M_7^0$ at the end will have its subscript decreased three times. A similar analysis applies to the third and fourth rows of (\ref{A7tableau}) together, yielding eight $\E'$'s $\oplus\Si^{128}M_4^0$, and to the last four rows of (\ref{A7tableau}) combined, yielding sixteen $\E'$'s $\oplus\Si^{256}M_4^0$. We split off all the $\E'$'s, and the various $\Si^{2^i}M_4^0$'s combine to give the sequence in Definition \ref{E1'def}. The placement of the $\Si^{2^i}M_4^0$'s in the sequence will be justified in the paragraph preceding Figure \ref{koku4}.

We can use (\ref{A7tableau})  similarly to justify the subedge structure of $\At_k$ described in Theorem \ref{Akthm}. Under $\E'_{1,7}$, we have $\E'_{2,3}$ and $\E'_{2,4}$ in the first row of (\ref{A7tableau}), $\E'_{2,5}$ in the second, $\E'_{2,6}$ in the third, and $\E'_{2,7}$ in the fifth. Then, under the $\E'_{2,4}$ is the $\E'_{3,4}$ just after it, and under the $\E'_{2,5}$ are the $\E'_{3,4}$ and $\E'_{3,5}$ in the second row. Also under that $\E'_{3,5}$ is the $\E'_{4,5}$ at the end of the second row. That exhausts the second row of (\ref{A7tableau}).

Next we consider the third and fourth rows together. The $\E'_{3,4}$, $\E'_{3,5}$, and $\E'_{3,6}$ under $\E'_{2,6}$ are apparent. Then under the $\E'_{3,5}$ is the $\E'_{4,5}$ which follows it, and under the $\E'_{3,6}$ are the $\E'_{4,5}$ and $\E'_{4,6}$ in the fourth row. We also have $\E'_{5,6}$ under the $\E'_{4,6}$.

A similar analysis applies to the last four rows of (\ref{A7tableau}) to give the subedge structure under $\E'_{2,7}$. This procedure generalizes to any $\At_k$. The discussion following (\ref{e+3}) suggests how the formulas in (\ref{Aktableau}) enable one to see how the pattern of edges in (\ref{A7tableau}) generalizes.

\section{The spectral sequence for $\A_k$, with proofs}\label{pfsec}
In this section, we show how the spectral sequence converging to $\At_k$ works, and then justify that the differentials and extensions work as claimed, by comparing with the $ku$ analysis.
Referring to the tableau for $\At_7$ in (\ref{A7tableau}) is strongly recommended. We will eventually specialize to $\At_7$.

\subsection{Description of spectral sequence}\label{131}

The following definition contains several useful notations.
\begin{defin}\label{Wdef} We will use the following notations.
\begin{itemize}
\item For $n>0$ and $n\equiv 0$ or $7$ mod $8$, let $W^n$ denote the $\Si^n$ term in the tableau for $\At_k$, which is independent of $k$ except that $n$ must be $\le 2^{k+1}$. If $W^n=\Si^nM^s_t$, let $W^n_i=\Si^nM^s_{t-i}$.
\item If $C$ is a chart and $M$ is an $M^s_t$ chart, positioned so that there are $d_1$ differentials from the stable lower edge of $M^s_t$ into $C$, we denote by $C\cup_{d_1}M$ the resulting chart, including $\eta$ extensions from $2\a$ to $v_1\b$ inserted whenever $d_1(\a)=\b$, following (\ref{bracket}).
\item Let $L(8a)$ denote the label on the arrow $W^{8a-7}\lar W^{8a}$.
\end{itemize}
\end{defin}

Now we describe the way the ASS for $\At_k$ works.
We make frequent use of Theorem \ref{edgthm} to see the increasing subscripts of $W$ (corresponding to decreased subscripts of the associated $M$).
After the first two steps, we specialize to $\At_7$. All claims will be justified later in this section.

{\bf Step 1} says that for $a\ge1$,
\begin{equation}\label{step1}W^{16a-7}\cup_{d_1}\Phi^1W^{16a}=\Phi^1W^{16a}_1\oplus \Si^{16a-8}\E'_{L(16a)}.\end{equation}
Then {\bf Step 2} says that for $a\ge1$,
\begin{equation}\label{step2}(W^{32a-15}\cup_{d_1}\Phi^1W^{32a-8})\cup_{d_1}\Phi^4W_1^{32a}=\Phi^4W_2^{32a}\oplus\Si^{32a-16}\E'_{L(32a-8)}.\end{equation}
At each step, we leave all $\E'$'s behind, until consideration of differentials among them at the very end.

The second row of (\ref{A7tableau}) now, which applies to any $\At_k$ with $k\ge5$, is 
\begin{equation}\label{70}((\Si^{33}M_6^2\cup_{d_1}\Phi^1\Si^{40}M_4^1)\lar\Phi^1\Si^{48}M_4^1)\lar\Phi^4\Si^{64}M_5^0,\end{equation}
where filtrations of $\Phi^1\Si^{48}M_4^1$ and then $\Phi^4\Si^{64}M_5^0$ will have to be increased (to $\Phi^4$ and $\Phi^{11}$, respectively) in order to get the $d_1$ differentials. We prefer to write it as 
$$\Si^{33}M_6^2\ml{2,5}\Si^{40}M_4^1\lar \Si^{48}M_4^1\lar \Si^{64}M_5^0,$$
ignoring the filtration increases, which we think of as being implicit. In fact, the corresponding $\Phi$'s are $\Phi^{2^{t+1}-t-2}$, $t\ge1$ starting with the domain of the first arrow. By Theorem \ref{edgthm}, this sequence equals $\Si^{64}M_4^0\oplus\Si^{32}\E'_{2,5}$.

We now restrict to $\At_7$, as in (\ref{A7tableau}). 
Rows 4, 6, and 8 in (\ref{A7tableau}) simplify similarly to row 2 just described. With previous $\E'$'s omitted, they become $\Si^{128}M_5^0\oplus\Si^{96}\E'_{3,6}$, $\Si^{192}M_4^1\oplus \Si^{160}\E'_{3,6}$, and $\Si^{256}M_6^0\oplus \Si^{224}\E'_{4,7}$. In each case, the $M$ term is $W^{64a}_3$.
This doesn't work quite so well in rows 3, 5, and 7. The requirement in Theorem \ref{edgthm} that $k+1+i=e+3$ is crucial.

Now (\ref{A7tableau}) has simplified to the following, with all $\E'$'s omitted.
\begin{gather*}V_7\ml{1,7}\Si^8M_4^0\lar\Si^{16}M_4^0\lar\Si^{32}M_4^0\\
\lar\Si^{64}M_4^0\\
\Si^{65}M^2_7\ml{2,6}\Si^{72}M_4^1\lar\Si^{80}M_4^1\lar\Si^{96}M_4^1\\
\lar\Si^{128}M_5^0\\
\Si^{129}M_5^2\ml{2,7}\Si^{136}M_4^1\lar \Si^{144}M_4^1\lar\Si^{160}M_4^1\\
\lar\Si^{192}M_4^1\\
\Si^{193}M_7^3\ml{3,7}\Si^{200}M_4^2\lar \Si^{208}M_4^2\lar\Si^{224}M_4^2\\
\lar\Si^{256}M_6^0\end{gather*}
Theorem \ref{edgthm} applies to the combination of rows 3 and 4, yielding $\Si^{128}M_4^0\oplus \Si^{64}\E'_{2,6}$, and to the combination of rows 7 and 8, yielding $\Si^{256}M_5^0\oplus \Si^{192}\E'_{3,7}$. Now, ignoring $\E'$'s and combining some rows, we obtain
\begin{gather}V_7\ml{1,7}\Si^8M_4^0\lar\Si^{16}M_4^0\lar\Si^{32}M_4^0\lar\Si^{64}M_4^0\lar\Si^{128}M_4^0\\
\Si^{129}M^2_8\ml{2,7}\Si^{136}M_4^1\lar\Si^{144}M_4^1\lar\Si^{160}M_4^1\lar\Si^{192}M_4^1\lar\Si^{256}M_5^0.\end{gather}
Theorem \ref{edgthm} applies to the second row, yielding $\Si^{256}M_4^0\oplus\Si^{128}\E'_{2,7}$.
Appending $\Si^{256}M_4^0$ to the first row, we obtain $\E'_{1,7}$ of Definition \ref{E1'def}. As noted there,  we obtain that $\E'_{1,7}$ equals
$$I_7\oplus\Phi^4\Si^8M(32)\bigr|_{16}\oplus\Phi^8\Si^{16}M(16)\bigr|_{32}\oplus\Phi^{16}\Si^{32}M(8)\bigr|_{64}\oplus\Phi^{32}\Si^{64}M(4)\bigr|_{128}\oplus\Phi^{64}\Si^{128}M(2)\bigr|_{256},$$
where $M(2^i)=ko_*(M(2^i))$ and $M\bigr|_x$ is the portion of $M$ in grading less than $x$. Here $I_k$ is as in Definition \ref{Ldef}. Figure \ref{E14} illustrates the case $k=4$.

The spectral sequence behaves for any $k$ in a manner which generalizes what we displayed here for $k=7$. The general tableau is in (\ref{Aktableau}), and an arrow with target $\Si^pM_\ell^i$ is labeled with $(i,i+\ell-3)$. The sequences which form all $\E'_{e,\ell}$ with $e\ge2$ in $\At_k$ are described in the following result.
\begin{thm}\label{genthm} For every odd positive integer $p$ and $t\ge0$ satisfying $2^{t+3}(p+1)\le 2^{k+1}$, there is a sequence satisfying Theorem \ref{edgthm}
\begin{equation}\label{Wtp} W^{2^{t+3}p+1}\lar W_0^{2^{t+3}p+8}\lar W_1^{2^{t+3}p+16}\lar\cdots\lar W_{t-1}^{2^{t+3}p+2^{t+2}}\lar W_t^{2^{t+3}(p+1)},\end{equation}
which forms $W^{2^{t+3}(p+1)}_{t+1}\oplus\Si^{2^{t+3}p}\E'_{\a(p)+1,\a(p)+t+2}$. Every $W^{8n}_i$ with $$0\le i\le\begin{cases}\nu(n)&\a(n)>1\\ \nu(n)-1&\a(n)=1\end{cases}$$
fits into exactly one of these sequences. If $i=\nu(n)$, then $t=\nu(n-2^i)$ and $p=(n-2^i)/2^t$, while if $i<\nu(n)$, then $t=i$ and $p=\frac n{2^i}-1$.\end{thm}

\begin{expl} If $n=20$, then
\begin{itemize}
    \item $i=0$ has $t=0$, $p=19$, and the sequence is $W^{153}\lar W_0^{160}$, forming $W_1^{160}\oplus \Si^{152}\E'_{4,5}$,
    \item $i=1$ has $t=1$, $p=9$, and the sequence is
    $W^{145}\lar W_0^{152}\lar W_1^{160}$, forming $W_2^{160}\oplus \Si^{144}\E'_{3,5}$, and
    \item $i=2$ has $t=4$, $p=1$, and the sequence is
    $$W^{129}\lar W_0^{136}\lar W_1^{144} \lar W_2^{160}\lar W_3^{192}\lar W_4^{256},$$
    forming $W^{256}_5\oplus \Si^{128}\E'_{2,7}$.
    \end{itemize}
(\ref{Aktableau}) can be used to see that these sequences are as follows. They are obtained from (\ref{A7tableau}) using Theorem \ref{edgthm}.
\begin{align*}
    &\Si^{153}M_4^4\ml{4,5}\Si^{160}M_6^1\\
&\Si^{145}M_5^3 \ml{3,5}\Si^{152}M_4^2\lar \Si^{160}M_5^1\\
&\Si^{129}M_7^2\ml{2,7}\Si^{136}M_4^1\lar \Si^{144}M_4^1\lar \Si^{160}M_4^1\lar\Si^{192}M_4^1\lar\Si^{256}M_5^0.
\end{align*}
    \end{expl}
    
\ni This example illustrates how the subscripts of $W^{8n}$    are successively increased until it becomes an $M_4^i$, at which point it is melded into another sequence, unless $n$ is a 2-power, in which case it melds into the sequence (\ref{E1seq}) which forms $\E'_{1,k}$. We are left with just $\E'_{1,k}$ and all the $\E'$'s which were left behind. Finally, there are the differentials among the $\E'$'s described in Theorem \ref{third}, which will be justified in the second paragraph after Figure \ref{last}.

\begin{proof}[Proof of Theorem \ref{genthm}] Using (\ref{Aktableau}), the sequence becomes, omitting $\Si^{2^{t+3}p}$,
$$\Si M^{\a(p)+1}_{4+t}\lar \Si^8M_4^{\a(p)}\lar\cdots\lar\Si^{2^{t+2}}M_4^{\a(p)}\lar\Si^{2^{t+3}}M_{\nu(p+1)+4}^{\a(p+1)-1}.$$
This satisfies the requirement in Theorem \ref{edgthm} because $\a(p+1)+\nu(p+1)=\a(p)+1$.\end{proof}

\subsection{Justification of differentials}\label{132}

The main ingredient in justifying the differentials is the generalization of what was done in Section \ref{kudifflsec} and Theorem \ref{closed} to pass from differentials in $ku^*(K_2)$ obtained in \cite[Theorem 3.1]{DW} to differentials among $ku_*$ summands as illustrated in Figure \ref{homdiff}. We will postpone temporarily consideration of differentials involving the summand $V_k$ of lowest grading.
In (\ref{bt}), we showed how the theorem of \cite{DW} could be interpreted in $ku^*$ summands as
$$d^{2^{t+1}-(t+1)}(b_t(Y_{2^t(a+1)-1,k}))=v^{2^{t+1}-(t+1)}\tau_t(Y_{2^ta,k}).$$

If $\sq^1(\a)=\b$ in $H^*(K_2)$, a $v$-tower in the formation of $ku^*(K_2)$ arises from $\b$, while a $v$-tower in the formation of $ku_*(K_2)$ arises from $\a$. The $ku^*$ differentials in Figure \ref{dual} are $d^2(b_1(Y_{1,4}))=v^2\tau_1(Y_{0,4})$ and $d^5(b_2(Y_{3,4}))=v^5\tau_2(Y_{0,4})$. 
The dual differentials in $ku_*$ go in the opposite direction on the $\a$ classes associated to the corresponding  $\b$ classes.
For example, dual to the $d^5$ from 57 to 68 is a $d_5$ from 67 to 56. See Figure \ref{71}.
This is consistent with (\ref{kuduality}) since it gives a $v$-tower of height 5 in $(ku^*(K_2))^\vee$ on a class of grading 60, and on a class in $ku_{56}(K_2)$.

This always works, leading to a generalization of Figure \ref{homdiff}. Suppose $\sq^1(\a_2)=\b_2$ and $\sq^1(\a_1)=\b_1$ with $|\a_2|=2x_2-1$ and $|\a_1|=2x_1$, and there is a $ku^*$ differential from the $\b_1$ $v$-tower to the $\b_2$ $v$-tower. The differential hits in grading $2x_1+2$, giving a $v$-tower of height $x_2-x_1-1$ on a class in $(ku^*(K_2))^\vee$ of grading $2x_1+4$. The differential in $ku_*(K_2)$ gives a $v$-tower of height $x_2-x_1-1$ on a class of grading $2x_1$, consistent with (\ref{kuduality})..

We can translate from the $ku$ notation used to label the columns of Figure \ref{homdiff} to the $ko$ notation explained prior to (\ref{Aktableau}). The contributions to $\Si$ of $y_1^2$, $y_1q$, $z_j$, and $M_j$ are 8, 9, $2^{j+2}$, and $2^j$, respectively. Also, $v_1q$ and $z_j$ add 2 and 1, respectively, to the superscript of $M$. In the notation of (\ref{bt}) and Definition \ref{Wdef}, $Y_{2a,k}\leftrightarrow W^{2^{k+1}-8a}$ and $Y_{2a+1,k}\leftrightarrow W^{2^{k+1}-8a+1}$. Then (\ref{bt}) implies the following proposition, adapting to $ku$ the notation in (\ref{A7tableau}) and Definition \ref{Wdef}. For example, the columns of Figure \ref{homdiff} are now labeled with the $W$ version of (\ref{A4tableau}).
\begin{prop}\label{ku*diff} A $d_{2^{t+1}-(t+1)}$ $ku_*$-differential will go from a $W^{2^{t+2}a}$ summand to the $W^{2^{t+2}(a-1)+1}$ summand.\end{prop}

We begin illustrating how the description of the spectral sequence, especially the differentials, in Subsection \ref{131} follows from Proposition \ref{ku*diff} using $\E'_{3,6}$ in  the sixth row of (\ref{A7tableau}). After desuspending, this is
\begin{equation}\label{6}\Si M_6^3\ml{3,6}\Si^8M_4^2\quad \Si^9M_4^4\ml{4,5}\Si^{16}M_5^2\quad \Si^{17}M_5^4\ml{4,6}\Si^{24}M_4^3\quad \Si^{25}M_4^5\ml{5,6}\Si^{32}M_7^1.\end{equation}

In Figure \ref{d2s}, we show, for each of the four arrows in (\ref{6}), the two charts combined, without increasing filtration, in both the $ko$ and $ku$ context. For the $ku$ context, $M_k^i$ refers to $\ext_{E_1}(\Si^{2i-2^k}M_k,\zt)$, analogous to those in Figure \ref{4chts}. We use black for the $W^{16a-7}$ chart and red for $W^{16a}$. We use big dots to indicate the classes that map across under $ko_*(K_2)\mapright{c} ku_*(K_2)$. We insert differentials in the $ku$ spectral sequence implied by Proposition \ref{ku*diff}, and then the differentials in the $ko$ spectral sequence implied by those. After presenting the figure, we will comment on some of the implications.

\bigskip
\begin{minipage}{6in}
\begin{fig}\label{d2s}

{\bf $d_2$ differentials ($ko_*\mapright{c}ku_*$)}

\begin{center}

\begin{\tz}[scale=.49]
\draw (1,0) -- (2,1) -- (2,0) -- (4,2);
\draw (8,3) -- (10,5) -- (10,4) -- (12,6);
\draw [red] (0,0) -- (1,1) -- (1,.3);
\draw [red] (5,0) -- (5,2);
\draw [red] (7,3) -- (9,5) -- (9,1) -- (11,3);
\draw [->] (17,0) -- (28,5.5);
\draw [red] [->] (16,0) -- (28,6);
\draw [red] [->] (18,0) -- (28,5);
\draw [red] [->] (20,0) -- (28,4);
\draw [red] [->] (22,0) -- (28,3);
\draw [red] (18,0) -- (18,1);
\draw [red] (20,0) -- (20,2);
\draw [red] (22,0) -- (22,3);
\draw [red] (24,1) -- (24,4);
\draw [red] (26,2) -- (26,5);
\draw (0,0) -- (11,0);
\draw (15,0) -- (27,0);
\node at (1,0) {\llb};
\node at (2,1) {\llb};
\node at (3,1) {\llb};
\node at (4,2) {\llb};
\node at (9.2,4.2) {\llb};
\node at (10,5) {\llb};
\node at (11,5) {\llb};
\node at (12,6) {\llb};
\node at (19,1) {\llb};
\node at (21,2) {\llb};
\node at (27,5) {\llb};
\node [red] at (0,0) {\llb};
\node [red] at  (1,1) {\llb};
\node [red] at  (5,0) {\blb};
\node [red] at  (5,1) {\blb};
\node [red] at  (5,2) {\blb};
\node [red] at  (8,4) {\llb};
\node [red] at  (9,5) {\llb};
\node [red] at  (10,2) {\llb};
\node [red] at  (11,3) {\llb};
\node [red] at  (18,0) {\llb};
\node [red] at  (18,1) {\llb};
\node [red] at  (22,0) {\llb};
\node [red] at  (22,1) {\llb};
\node [red] at  (22,2) {\llb};
\node [red] at  (26,2) {\llb};
\node [red] at  (26,3) {\llb};
\node [red] at  (26,4) {\llb};
\node [red] at  (26,5) {\llb};
\node at (0,-.7) {$34$};
\node at (2,-.7) {$36$};
\node at (6,-.7) {$40$};
\node at (10,-.7) {$44$};
\node at (17,-.7) {$36$};
\node at (21,-.7) {$40$};
\node at (25,-.7) {$44$};
\node at (2,0) {\blb};
\node at (8,3) {\blb};
\node at (10,4) {\blb};
\node at (17,0) {\blb};
\node at (23,3) {\blb};
\node at (25,4) {\blb};
\node [red] at (16,0) {\blb};
\node [red] at (20,0) {\blb};
\node [red] at (20,1) {\blb};
\node [red] at (20,2) {\blb};
\node [red] at (22,3) {\blb};
\node [red] at (24,1) {\blb};
\node [red] at (24,2) {\blb};
\node [red] at (24,3) {\blb};
\node [red] at (24,4) {\blb};
\node [red] at (1,.3) {\blb};
\node [red] at (7,3) {\blb};
\node [red] at (9,1) {\blb};
\node [red] at (9,2) {\blb};
\node [red] at (9,3) {\blb};
\node [red] at (9,3.8) {\blb};
\draw [blue] (5,0) -- (4,2);
\draw [blue] [ultra thick] (9,1) -- (8,3);
\draw [blue] (10,2) -- (9.2,4.2);
\draw [blue] (11,3) -- (10,5);
\draw [blue] (22,0) -- (21,2);
\draw [blue] [ultra thick] (24,1) -- (23,3);
\draw [blue] (26,2) -- (25,4);
\node at (4,5) {$\Si^{25}M_4^5\lar \Si^{32}M_7^1$};
\node at (14,2) {$\mapright{c}$};

\draw (0,8) -- (2,10) -- (2,8) -- (4,10);
\draw (8,12) -- (10,14) -- (10,12) -- (12,14);
\draw [->] (16,8) -- (27,13.5);
\draw [red] (7,9) -- (9,11) -- (9,10) -- (11,12);
\draw  [->] (18,8) -- (27,12.5);
\draw [red] [->] (21,8) -- (27,11);
\draw [blue] (12,6) -- (13,4);
\draw (-.5,8) -- (12,8);
\draw (15.5,8) -- (27,8);
\node at (2,7.3) {$28$};
\node at (6,7.3) {$32$};
\node at (10,7.3) {$36$};
\node at (18,7.3) {$28$};
\node at (22,7.3) {$32$};
\node at (26,7.3) {$36$};
\node at (1,9) {\llb};
\node at (2,10) {\llb};
\node at (3,9) {\llb};
\node at (4,10) {\llb};
\node at (9,13) {\llb};
\node at (10,14) {\llb};
\node at (11,13) {\llb};
\node at (12,14) {\llb};
\node at (20,10) {\llb};
\node at (20,9) {\llb};
\node at (22,10) {\llb};
\node at (24,11) {\llb};
\node at (0,8) {\blb};
\node at (2,8) {\blb};
\node at (2,9) {\blb};
\node at (6,11) {\blb};
\node at (8,12) {\blb};
\node at (10,12) {\blb};
\node at (10,13) {\blb};
\node at (16,8) {\blb};
\node at (18,8) {\blb};
\node at (18,9) {\blb};
\node at (22,11) {\blb};
\node at (24,12) {\blb};
\node at (26,13) {\blb};
\node at (26,12) {\blb};
\node [red] at (3,8) {\llb};
\node [red] at (8,10) {\llb};
\node [red] at (9,11) {\llb};
\node [red] at (10,11) {\llb};
\node [red] at (11,12) {\llb};
\node [red] at (21,8) {\llb};
\node [red] at (7,9) {\blb};
\node [red] at (9,10) {\blb};
\node [red] at (23,9) {\blb};
\node [red] at (25,10) {\blb};
\draw [blue] (3,8) -- (2,10);
\draw [blue] [ultra thick] (7,9) -- (6,11);
\draw [blue] [ultra thick] (9,10) -- (8,12);
\draw [blue] (10,11) -- (9,13);
\draw [blue] (11,12) -- (10,14);
\draw [blue] (21,8) -- (20,10);
\draw [blue] [ultra thick] (23,9) -- (22,11);
\draw [blue] [ultra thick] (25,10) -- (24,12);
\node at (14,11) {$\mapright{c}$};
\node at (3,13) {$\Si^{17}M_5^4\lar  \Si^{24}M_4^3$};
\draw (18,8) -- (18,9);
\draw (20,9) -- (20,10);
\draw (22,10) -- (22,11);
\draw (24,11) -- (24,12);
\draw (26,12) -- (26,13);

\draw (0,16) -- (2,18) -- (2,17) -- (4,19);
\draw (8,20) -- (10,22) -- (10,21) -- (12,23);
\draw [->] (16,16) -- (28,22);
\draw [red] [->] (19,16) -- (28,20.5);
\draw [red] [->] (21,16) -- (28,19.5);
\draw [red] (21,16) -- (21,17);
\draw [red] (23,17) -- (23,18);
\draw [red] (25,18) -- (25,19);
\draw [red] (27,19) -- (27,20);
\draw (-.5,16) -- (12,16);
\draw (17.5,16) -- (27,16);
\node at (1,17) {\llb};
\node at (2,18) {\llb};
\node at (3,18) {\llb};
\node at (4,19) {\llb};
\node at (9,21) {\llb};
\node at (10,22) {\llb};
\node at (11,22) {\llb};
\node at (12,23) {\llb};
\node at (20,18) {\llb};
\node at (22,19) {\llb};
\node at (0,16) {\blb};
\node at (2,17) {\blb};
\node at (8,20) {\blb};
\node at (10,21) {\blb};
\node at (16,16) {\blb};
\node at (18,17) {\blb};
\node at (24,20) {\blb};
\node at (26,21) {\blb};
\node [red] at (1,16) {\llb};
\node [red] at (8,19) {\llb};
\node [red] at (9,20) {\llb};
\node [red] at (10,19) {\llb};
\node [red] at (21,16) {\llb};
\node [red] at (23,17) {\llb};
\node [red] at (27,19) {\llb};
\node [red] at (3,16) {\blb};
\node [red] at (5,17) {\blb};
\node [red] at (7,18) {\blb};
\node [red] at (9,19) {\blb};
\node [red] at (11,20) {\blb};
\node [red] at (13,21) {\blb};
\node [red] at (9,18) {\blb};
\node [red] at (19,16) {\blb};
\node [red] at (21,17) {\blb};
\node [red] at (23,18) {\blb};
\node [red] at (25,19) {\blb};
\node [red] at (27,20) {\blb};
\node [red] at (25,18) {\blb};
\draw [blue] (3,16) -- (2,18);
\draw [blue] (5,17) -- (4,19);
\draw [blue] [ultra thick] (9,18) -- (8,20);
\draw [blue] (10,19) -- (9,21);
\draw [blue] (11,20) -- (10,22);
\draw [blue] (21,16) -- (20,18);
\draw [blue] (23,17) -- (22,19);
\draw [blue] [ultra thick] (25,18) -- (24,20);
\draw [blue] (27,19) -- (26,21);
\node at (2,15.3) {$20$};
\node at (6,15.3) {$24$};
\node at (10,15.3) {$28$};
\node at (18,15.3) {$20$};
\node at (22,15.3) {$24$};
\node at (26,15.3) {$28$};
\node at (14,19) {$\mapright{c}$};
\draw [red] (7,18) -- (9,20) -- (9,18) -- (11,20);
\draw [blue] (13,21) -- (12,23);
\draw (15.5,16) -- (28,16);
\node at (3,21) {$\Si^9M_4^4\lar \Si^{16}M_5^2$};

\draw (2,26) -- (4,28) -- (4,25) -- (6,27);
\node [red] at (3,25) {\llb};
\draw (8,28) -- (8,29);
\draw (10,30) -- (12,32) -- (12,29) -- (14,31);
\draw [red] (3.8,25.2) -- (4.8,26.2);
\draw [red] (9,27) -- (11,29) -- (11,28) -- (11.8,29.2) -- (12.8,30.2);
\draw [blue] (3.8,25.2) -- (3,27);
\draw [blue] (4.8,26.2) -- (4,28);
\draw [blue] [ultra thick] (9,27) -- (8,29);
\draw [blue] [ultra thick] (11,28) -- (10,30);
\draw [blue] (11.8,29.2) -- (11,31);
\draw [blue] (12.8,30.2) -- (12,32);
\node at (3,27) {\llb};
\node at (4,28) {\llb};
\node [red] at  (4.8,26.2) {\llb};
\node [red] at  (3.8,25.2) {\llb};
\node at (5,26) {\llb};
\node at (6,27) {\llb};
\node at (8,28) {\llb};
\node  [red] at (10,28) {\llb};
\node  [red] at (11,29) {\llb};
\node at (11,31) {\llb};
\node at (12,32) {\llb};
\node at (13,30) {\llb};
\node at (14,31) {\llb};
\node [red] at (11.8,29.2) {\llb};
\node [red] at (12.8,30.2) {\llb};
\node at (0,25) {\blb};
\node at (2,26) {\blb};
\node at (4,25) {\blb};
\node at (4,26) {\blb};
\node at (4,27) {\blb};
\node at (8,29) {\blb};
\node  [red] at (9,27) {\blb};
\node at (10,30) {\blb};
\node at (12,29) {\blb};
\node at (12,30) {\blb};
\node at (12,31) {\blb};
\node [red] at (11,28) {\blb};
\draw [->] (18,25) -- (30,31);
\draw [->] (20,25) -- (30,30);
\draw [->] (22,25) -- (30,29);
\draw [->] [red] (23,25) -- (30,28.5);
\draw (22,25) -- (22,27);
\draw (24,26) -- (24,28);
\draw (26,27) -- (26,29);
\draw (28,28) -- (28,30);
\draw [blue] (23,25) -- (22,27);
\draw [blue] (25,26) -- (24,28);
\draw [blue] [ultra thick] (27,27) -- (26,29);
\draw [blue] [ultra thick] (29,28) -- (28,30);
\node at (20,25) {\llb};
\node at (24,26) {\llb};
\node at (24,27) {\llb};
\node at (24,28) {\llb};
\node at (26,27) {\llb};
\node at (26,28) {\llb};
\node at (28,28) {\llb};
\node at (28,29) {\llb};
\node at (18,25) {\blb};
\node at (20,26) {\blb};
\node at (22,25) {\blb};
\node at (22,26) {\blb};
\node at (22,27) {\blb};
\node at (26,29) {\blb};
\node at (28,30) {\blb};
\node [red] at (23,25) {\llb};
\node [red] at (25,26) {\llb};
\node [red] at (27,27) {\blb};
\node [red] at (29,28) {\blb};
\draw (-.5,25) -- (13,25);
\draw (17.5,25) -- (29,25);
\node at (0,24.3) {$8$};
\node at (4,24.3) {$12$};
\node at (8,24.3) {$16$};
\node at (12,24.3) {$20$};
\node at (18,24.3) {$8$};
\node at (22,24.3) {$12$};
\node at (26,24.3) {$16$};
\node at (14,28) {$\mapright{c}$};
\node at (3,30) {$\Si M_6^3\lar \Si^8M_4^2$};

\end{\tz}
\end{center}
\end{fig}
\end{minipage}
\bigskip

First we note that the classes which map across are exactly those not in $\im(\eta)$. This follows from (\ref{exact}). In each of the four situations, there is at least one differential in the $ku$ spectral sequence in which both classes involved in the differential are in $\im(c)$. That implies the corresponding differential in the $ko$ spectral sequence. We indicate these differentials with thicker lines. Other differentials in the $ko$ spectral sequence are implied by naturality with respect to the action of $\eta$, $v_1^4$, or the generator of $ko_4$.

Now, in each of the four cases, we elevate by 1 the filtration of all classes in the $W^{8a}$ chart (the red classes) so that the $d_2$ differentials look like $d_1$'s. We did not do this in our $ku$ work in \cite{DW}, but it will be useful for subsequent comparisons with $ko$ charts. Insert $\eta$ extensions from $\a$ to $\b$ if $d_1(v_1\a)=2\b$. This follows from (\ref{bracket}).
This comparison with $ku$ is the justification for the filtration shift in (\ref{step1}), which applies to the second and fourth situations in Figure \ref{d2s}.

In Figure \ref{d5}, the top left chart contains the result of the first $ko$ chart of Figure \ref{d2s} in black and the result of the second $ko$ chart of Figure \ref{d2s} in red.
The top right chart in Figure \ref{d5} is the analogue of this for the first two $ku$ charts in Figure \ref{d2s}. We use big dots for classes that map across under $\mapright{c}$. We insert differentials in the $ku$ chart implied by Proposition \ref{ku*diff}. They will appear as a $d_4$ rather than a $d_5$ because we increased the filtration of the source classes. Then we insert differentials in the $ko$ chart forced by the morphism $c$ and the $ku$ differentials, again indicating with darker lines the $ku$ differentials that pull back to $ko$ differentials. Then we insert differentials implied by the inserted differentials and naturality with respect to $\eta$, $v_1^4$, and the generator of $ko_4$. The bottom half of Figure \ref{d5} does the analogous thing for the third and fourth charts of Figure \ref{d2s}.

\bigskip
\begin{minipage}{6in}
\begin{fig}\label{d5}

{\bf $d_5$ differentials (look like $d_4$)}

\begin{center}

\begin{\tz}[scale=.33]
\draw (8,0) -- (27,0);
\draw (30,0) -- (49,0);
\draw (8,0) -- (9,1);
\draw (10,1) -- (10,0) -- (12,2);
\draw (16,3) -- (18,5) -- (18,4) -- (20,6);
\draw (24,7) -- (26,9) -- (26,8) -- (28,10);
\draw (30,0) -- (32,1) -- (32,0);
\draw [->] (32,0) -- (50,9);
\node at (8,-.7) {$26$};
\node at (10,-.7) {$28$};
\node at (14,-.7) {$32$};
\node at (18,-.7) {$36$};
\node at (22,-.7) {$40$};
\node at (26,-.7) {$44$};
\node at (30,-.7) {$26$};
\node at (32,-.7) {$28$};
\node at (36,-.7) {$32$};
\node at (40,-.7) {$36$};
\node at (44,-.7) {$40$};
\node at (48,-.7) {$44$};
\node at (9,1) {\llb};
\node at (11,1) {\llb};
\node at (12,2) {\llb};
\node at (17,4) {\llb};
\node at (18,5) {\llb};
\node at (19,5) {\llb};
\node at (20,6) {\llb};
\node at (25,8) {\llb};
\node at (26,9) {\llb};
\node at (27,9) {\llb};
\node at (28,10) {\llb};
\node at (34,1) {\llb};
\node at (36,2) {\llb};
\node at (42,5) {\llb};
\node at (44,6) {\llb};
\node at (8,0) {\blb};
\node at (10,0) {\blb};
\node at (10,1) {\blb};
\node at (16,3) {\blb};
\node at (18,4) {\blb};
\node at (24,7) {\blb};
\node at (26,8) {\blb};
\node at (30,0) {\blb};
\node at (32,0) {\blb};
\node at (32,1) {\blb};
\node at (38,3) {\blb};
\node at (40,4) {\blb};
\node at (46,7) {\blb};
\node at (48,8) {\blb};
\draw [red] (16,1) -- (17,2) -- (17,1);
\draw [red] (17,0) -- (18,1) -- (18,0);
\draw [red] (21,2) -- (21,3);
\draw [red] (23,4) -- (25,6) -- (25,3) -- (27,5);
\draw [red] (29,6) -- (29,7);
\draw [red] [->] (39,1) -- (50,6.5);
\draw [red] [->] (41,1) -- (50,5.5);
\draw [red] [->] (43,1) -- (50,4.5);
\draw [red] (40,0) -- (42,1);
\draw [red] (41,1) -- (41,2);
\draw [red] (43,1) -- (43,3);
\draw [red] (45,2) -- (45,4);
\draw [red] (47,3) -- (47,5);
\draw [red] (49,4) -- (49,6);
\node [red] at (16,1) {\llb};
\node [red] at (17,2) {\llb};
\node [red] at (17,0) {\llb};
\node [red] at (18,1) {\llb};
\node [red] at (19,1) {\llb};
\node [red] at (24,5) {\llb};
\node [red] at (25,6) {\llb};
\node [red] at (26,4) {\llb};
\node [red] at (27,5) {\llb};
\node [red] at (41,1) {\llb};
\node [red] at (41,2) {\llb};
\node [red] at (43,1) {\llb};
\node [red] at (45,2) {\llb};
\node [red] at (45,3) {\llb};
\node [red] at (49,4) {\llb};
\node [red] at (49,5) {\llb};
\node [red] at (49,6) {\llb};
\node [red] at (17,1) {\blb};
\node [red] at (18,0) {\blb};
\node [red] at (21,2) {\blb};
\node [red] at (21,3) {\blb};
\node [red] at (23,4) {\blb};
\node [red] at (25,3) {\blb};
\node [red] at (25,4) {\blb};
\node [red] at (25,5) {\blb};
\node [red] at (29,6) {\blb};
\node [red] at (29,7) {\blb};
\node [red] at (39,1) {\blb};
\node [red] at (40,0) {\blb};
\node [red] at (43,2) {\blb};
\node [red] at (43,3) {\blb};
\node [red] at (45,4) {\blb};
\node [red] at (47,3) {\blb};
\node [red] at (47,4) {\blb};
\node [red] at (47,5) {\blb};
\draw [blue] (21,2) -- (20,6);
\draw [blue] [ultra thick] (25,3) -- (24,7);
\draw [blue] (26,4) -- (25,8);
\draw [blue] (27,5) -- (26,9);
\draw [blue] (29,6) -- (28,10);
\draw [blue] (43,1) -- (42,5);
\draw [blue] (45,2) -- (44,6);
\draw [blue] [ultra thick] (47,3) -- (46,7);
\draw [blue] (49,4) -- (48,8);
\draw [red] (18,0) -- (19,1);
\node [red] at (42,1) {\llb};
\node at (29,2) {$\mapright{c}$};
\node at (14,8.7) {$\Si^{17}M_5^4\ml{4,6}\Si^{24}M_4^3\quad \Phi^1\Si^{32}M_6^1\oplus \Si^{24}\E'_{5,6}$};

\draw (8,12) -- (27,12);
\draw (30,12) -- (49,12);
\draw (11,13) -- (12,14) -- (12,12) -- (14,14);
\draw (18,16) -- (20,18) -- (20,16) -- (22,18);
\draw (26,20) -- (28,22) -- (28,20) -- (30,22);
\draw [red] (18,12) -- (19,13);
\draw [red] (20,13) -- (21,14);
\draw [red] (25,15) -- (27,17) -- (27,16) -- (29,18);
\draw (30,12) -- (32,13) -- (32,12);
\draw [->] (32,12) -- (50,21);
\draw [->] (34,12) -- (50,20);
\draw (34,12) -- (34,13);
\draw (36,13) -- (36,14);
\draw (38,14) -- (38,15);
\draw (40,15) -- (40,16);
\draw (42,16) -- (42,17);
\draw (44,17) -- (44,18);
\draw (46,18) -- (46,19);
\draw (48,19) -- (48,20);
\draw [red] (40,12) -- (42,13);
\draw [red] [->] (43,13) -- (50,16.5);
\node at (8,11.2) {$8$};
\node at (12,11.2) {$12$};
\node at (16,11.2) {$16$};
\node at (20,11.2) {$20$};
\node at (24,11.2) {$24$};
\node at (30,11.2) {$8$};
\node at (34,11.2) {$12$};
\node at (38,11.2) {$16$};
\node at (42,11.2) {$20$};
\node at (46,11.2) {$24$};
\node at (11,13) {\llb};
\node at (12,14) {\llb};
\node at (13,13) {\llb};
\node at (14,14) {\llb};
\node at (19,17) {\llb};
\node at (20,18) {\llb};
\node at (21,17) {\llb};
\node at (22,18) {\llb};
\node at (27,21) {\llb};
\node at (28,22) {\llb};
\node at (29,21) {\llb};
\node at (30,22) {\llb};
\node at (32,12) {\llb};
\node at (34,13) {\llb};
\node at (36,13) {\llb};
\node at (38,14) {\llb};
\node at (40,15) {\llb};
\node at (44,17) {\llb};
\node at (46,18) {\llb};
\node at (48,19) {\llb};
\node at (8,12) {\blb};
\node at (10,13) {\blb};
\node at (12,12) {\blb};
\node at (12,13) {\blb};
\node at (16,15) {\blb};
\node at (18,16) {\blb};
\node at (20,16) {\blb};
\node at (20,17) {\blb};
\node at (24,19) {\blb};
\node at (26,20) {\blb};
\node at (28,20) {\blb};
\node at (28,21) {\blb};
\node at (30,12) {\blb};
\node at (32,13) {\blb};
\node at (34,12) {\blb};
\node at (36,14) {\blb};
\node at (38,15) {\blb};
\node at (40,16) {\blb};
\node at (42,16) {\blb};
\node at (42,17) {\blb};
\node at (46,19) {\blb};
\node at (48,20) {\blb};
\node [red] at (19,13) {\llb};
\node [red] at (18.7,13.3) {\llb};
\node [red] at (21,14) {\llb};
\node [red] at (26,16) {\llb};
\node [red] at (27,17) {\llb};
\node [red] at (28,17) {\llb};
\node [red] at (29,18) {\llb};
\node [red] at (43,13) {\llb};
\node [red] at (45,14) {\llb};
\node [red] at (18,12) {\blb};
\node [red] at (20,13) {\blb};
\node [red] at (25,15) {\blb};
\node [red] at (27,16) {\blb};
\node [red] at (40,12) {\blb};
\node [red] at (42,13) {\blb};
\node [red] at (47,15) {\blb};
\node [red] at (49,16) {\blb};
\draw [blue] (20,13) -- (19,17);
\draw [blue] (21,14) -- (20,18);
\draw [blue] [ultra thick] (25,15) -- (24,19);
\draw [blue] [ultra thick] (27,16) -- (26,20);
\draw [blue] (28,17) -- (27,21);
\draw [blue] (29,18) -- (28,22);
\draw [blue] (43,13) -- (42,17);
\draw [blue] (45,14) -- (44,18);
\draw [blue] [ultra thick] (47,15) -- (46,19);
\draw [blue] [ultra thick] (49,16) -- (48,20);
\draw [blue] (18,0) -- (17,4);
\draw [blue] (19,1) -- (18,5);
\node at (14,20) {$\Si M_6^3\ml{3,6}\Si^8M_4^2\quad \Phi^1\Si^{16}M_4^2\oplus \Si^8\E'_{4,5}$};
\node at (44,18) {\llb};
\node at (30,15) {$\mapright{c}$};
\end{\tz}
\end{center}
\end{fig}
\end{minipage}
\bigskip

We increase filtrations of $v_1^4$-periodic red classes in the $ko$ charts of Figure \ref{d5}, remove classes involved in differentials, and insert $\eta$ extensions implied by (\ref{bracket}), obtaining Figure \ref{d12}, in which the top (resp.~bottom) half of Figure \ref{d5} is in black (resp.~red). Step 2 (\ref{step2}) says that the lower chart in Figure \ref{d5}, after simplification, is $\Phi^4\Si^{32}M_5^1\oplus\Si^{16}\E'_{4,6}$. The $\Phi^4\Si^{32}M_5^1$ appears in Figure \ref{d12} as all red elements in filtration $\ge4$. Our charts also include $\Si^{24}\E'_{5,6}$, which was dropped in the transition from (\ref{step1}) to (\ref{step2}).

If we did the same for the $ku$ charts, we would have, among other things, a $v$-tower arising from $(12,0)$ and $v$-towers connected by $h_0$ ($\cdot2$) arising from $(35,1)$ and $(37,1)$ (or $(35,4)$ and $(37,4)$, if we increased filtrations analogously to $ko$). See Figure \ref{d5}. The classes in the $ko$ chart which map across to these $v$-towers are indicated by big dots in Figure \ref{d12}. There are $d_{11}$ differentials\footnote{They are $d_{12}$ in Proposition \ref{ku*diff}, but $d_{11}$ here due to the filtration shift.} in the $ku$ spectral sequence from the $(37,1)$ $v$-tower to the $(12,0)$ $v$-tower implied by Proposition \ref{ku*diff}. The classes in the $ko$ chart in $(43,7)$ and $(42,15)$, being not in $\im(\eta)$, map across to classes involved in a differential in the $ku$ spectral sequence. Thus there is a differential in the $ko$ spectral sequence from $(43,7)$ to $(42,15)$, and the other differentials in Figure \ref{d12} are implied by naturality, as before.

\bigskip
\begin{minipage}{6in}
\begin{fig}\label{d12}

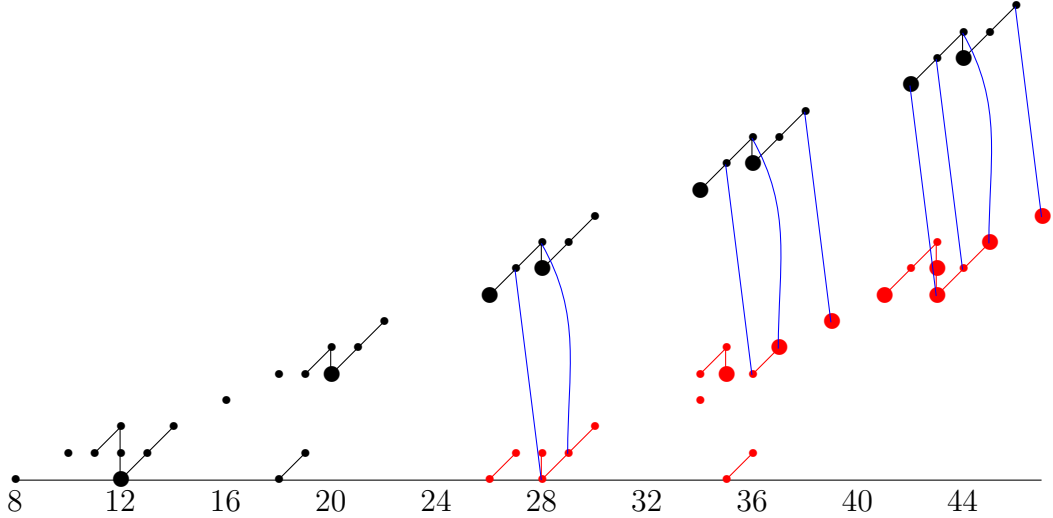
{\bf Forming the result of (\ref{6})}

\begin{center}

\begin{\tz}[scale=.35]
\draw (8,0) -- (47,0);
\draw (11,1) -- (12,2) -- (12,0) -- (14,2);
\draw (18,0) -- (19,1);
\draw (19,4) -- (20,5) -- (20,4) -- (22,6);
\draw (26,7) -- (28,9) -- (28,8) -- (30,10);
\draw (34,11) -- (36,13) -- (36,12) -- (38,14);
\draw (42,15) -- (44,17) -- (44,16) -- (46,18);
\draw [red] (26,0) -- (27,1);
\draw [red] (28,1) -- (28,0) -- (30,2);
\draw [red] (35,0) -- (36,1);
\draw [red] (34,4) -- (35,5) -- (35,4);
\draw [red] (36,4) -- (37,5);
\draw [red] (41,7) -- (43,9) -- (43,7) -- (45,9);
\node at (8,-.8) {$8$};
\node at (12,-.8) {$12$};
\node at (16,-.8) {$16$};
\node at (20,-.8) {$20$};
\node at (24,-.8) {$24$};
\node at (28,-.8) {$28$};
\node at (32,-.8) {$32$};
\node at (36,-.8) {$36$};
\node at (40,-.8) {$40$};
\node at (44,-.8) {$44$};
\node at (8,0) {\llb};
\node at (10,1) {\llb};
\node at (11,1) {\llb};
\node at (12,2) {\llb};
\node at (12,1) {\llb};
\node at (13,1) {\llb};
\node at (14,2) {\llb};
\node at (16,3) {\llb};
\node at (18,4) {\llb};
\node at (18,0) {\llb};
\node at (19,1) {\llb};
\node at (19,4) {\llb};
\node at (20,5) {\llb};
\node at (21,5) {\llb};
\node at (22,6) {\llb};
\node at (27,8) {\llb};
\node at (28,9) {\llb};
\node at (29,9) {\llb};
\node at (30,10) {\llb};
\node at (35,12) {\llb};
\node at (36,13) {\llb};
\node at (37,13) {\llb};
\node at (38,14) {\llb};
\node at (43,16) {\llb};
\node at (44,17) {\llb};
\node at (45,17) {\llb};
\node at (46,18) {\llb};
\node at (12,0) {\blb};
\node at (20,4) {\blb};
\node at (26,7) {\blb};
\node at (28,8) {\blb};
\node at (34,11) {\blb};
\node at (36,12) {\blb};
\node at (42,15) {\blb};
\node at (44,16) {\blb};
\node [red] at (26,0) {\llb};
\node [red] at (27,1) {\llb};
\node [red] at (28,0) {\llb};
\node [red] at (28,1) {\llb};
\node [red] at (29,1) {\llb};
\node [red] at (30,2) {\llb};
\node [red] at (35,0) {\llb};
\node [red] at (36,1) {\llb};
\node [red] at (34,3) {\llb};
\node [red] at (34,4) {\llb};
\node [red] at (35,5) {\llb};
\node [red] at (36,4) {\llb};
\node [red] at (42,8) {\llb};
\node [red] at (43,9) {\llb};
\node [red] at (44,8) {\llb};
\node [red] at (35,4) {\blb};
\node [red] at (37,5) {\blb};
\node [red] at (39,6) {\blb};
\node [red] at (41,7) {\blb};
\node [red] at (43,8) {\blb};
\node [red] at (43,7) {\blb};
\node [red] at (45,9) {\blb};
\node [red] at (47,10) {\blb};
\draw [blue] (28,0) -- (27,8);
\draw [blue] (29,1) to[out=90, in=300] (28,9);
\draw [blue] (36,4) -- (35,12);
\draw [blue] (37,5) to[out=90, in=300] (36,13);
\draw [blue] (39,6) -- (38,14);
\draw [blue] (43,7) -- (42,15);
\draw [blue] (44,8) -- (43,16);
\draw [blue] (45,9) to[out=90, in=300] (44,17);
\draw [blue] (47,10) -- (46,18);

\end{\tz}
\end{center}
\end{fig}
\end{minipage}
\bigskip

After increasing filtrations of all remaining red classes of filtration $\ge4$ by 7, and removing classes involved in differentials, we obtain Figure \ref{last}. This is just $\E_{3,6}$, $\Si^8\E_{4,5}$, $\Si^{16}\E_{4,6}$, $\Si^{24}\E_{5,6}$, and $\Phi^{11}\Si^{32}M_4^1$, with $\E$'s as in Figure \ref{four}, except that the classes in $\E'_{3,6}$ which support differentials (into $\E'_{2,6}$) have not been removed. Figure \ref{d12} is the following analogue of (\ref{70}) except that we have included all the $\eps'$'s. 
$$((\Si M_6^3\cup_{d_1}\Phi^1\Si^8M_4^2)\lar\Phi^4\Si^{16}M_4^2)\lar\Phi^4\Si^{32}M_5^1$$
It corresponds to the case $e=3$, $\ell=6$, $k=4$, $i=1$ of Theorem \ref{edgthm}, and, if the $\Si^{160}$ is included, the case $t=2$, $p=5$ case of (\ref{Wtp}).
 The main point of the work we have been doing here is how the differentials are implied by the $ku$ differentials in Proposition \ref{ku*diff}.

\bigskip
\begin{minipage}{6in}
\begin{fig}\label{last}

{\bf Final result of (\ref{6})}

\begin{center}

\begin{\tz}[scale=.35]
\draw (8,0) -- (44,0);
\draw (11,1) -- (12,2) -- (12,0) -- (14,2);
\draw (18,0) -- (19,1);
\draw (19,4) -- (20,5) -- (20,4) -- (22,6);
\draw (26,0) -- (27,1);
\draw (35,0) -- (36,1);
\draw (28,8) -- (30,10);
\draw (41,14) -- (43,16) -- (43,15) -- (45,17);
\node at (46,18) {$\iddots$};
\draw (34,11) -- (35,12) -- (35,11) -- (37,13);
\node at (8,-.7) {$8$};
\node at (12,-.7) {$12$};
\node at (16,-.7) {$16$};
\node at (20,-.7) {$20$};
\node at (24,-.7) {$24$};
\node at (28,-.7) {$28$};
\node at (32,-.7) {$32$};
\node at (36,-.7) {$36$};
\node at (40,-.7) {$40$};
\node at (44,-.7) {$44$};
\node at (8,0) {\lb};
\node at (10,1) {\lb};
\node at (11,1) {\lb};
\node at (12,2) {\lb};
\node at (12,1) {\lb};
\node at (12,0) {\lb};
\node at (13,1) {\lb};
\node at (14,2) {\lb};
\node at (16,3) {\lb};
\node at (18,4) {\lb};
\node at (18,0) {\lb};
\node at (19,1) {\lb};
\node at (19,4) {\lb};
\node at (20,5) {\lb};
\node at (20,4) {\lb};
\node at (21,5) {\lb};
\node at (22,6) {\lb};
\node at (26,7) {\lb};
\node at (26,0) {\lb};
\node at (27,1) {\lb};
\node at (28,1) {\lb};
\node at (30,2) {\lb};
\node at (34,3) {\lb};
\node at (35,0) {\lb};
\node at (36,1) {\lb};
\node at (28,8) {\lb};
\node at (29,9) {\lb};
\node at (30,10) {\lb};
\node at (34,11) {\lb};
\node at (35,12) {\lb};
\node at (35,11) {\lb};
\node at (36,12) {\lb};
\node at (37,13) {\lb};
\node at (33.7,11.1) {\lb};
\node at (41,14) {\lb};
\node at (42,15) {\lb};
\node at (43,16) {\lb};
\node at (43,15) {\lb};
\node at (44,16) {\lb};
\node at (45,17) {\lb};

\end{\tz}
\end{center}
\end{fig}
\end{minipage}
\bigskip

The behavior described in the detailed example just completed will continue regarding how differentials in the $ko_*$ spectral sequence are implied by the $d_{2^{t+1}-(t+1)}$ $ku_*$ differentials of Theorem \ref{ku*diff}, as $t$ increases. It will always be the case, as in nicely illustrated in Figure \ref{d5}, that the potential differentials go stably from $M_k$ into $\Si M_{k'}$ with one of $k$ and $k'$ equal to 4 and the other greater than 4. The generator of $ko_3(M_k)$ is not in $\im(\eta)$ and hence maps across to $ku_3(M_k)$, and similarly for the generator of $ko_2(\Si M_{k'})$. Proposition \ref{ku*diff} says that there will be a differential between the $ku$ classes, and hence there will be a differential between the $ko$ classes. Other differentials in the $ko$ spectral sequence  are implied by the action of $\eta$, $v_1^4$, and the generator of $ko_4$.

This completes the proof that the description of the spectral sequence given in Subsection \ref{131}, which implies Definition \ref{E'def}, is implied by Proposition \ref{ku*diff}. Also, the differentials in Theorem \ref{third} follow since, as is nicely illustrated in Figure \ref{d12}, the classes in grading 5 mod 8, after applying $v_1^{4i}$,  yield stable elements in $ko_5(M_k)$, in the notation of the preceding paragraph, which supported differentials thanks to the action of $\eta^2$. 

Additional differentials can be ruled out using the action of $\eta$ and $v_1^4$. Edges $\E_{e,\ell}$ tend to have less $v_1^4$ periodicity than edges above them. For example, in $\At_6$ as shown in Figure \ref{A6}, a conceivable differential on the class in $(74,1)$ is ruled out by the action of $v_1^4$. The only possible short edges are shown in Figure \ref{four}, and the structure of longer edges is already suggested in Figures \ref{2,6}, \ref{3,7}, \ref{4,8}, and \ref{5,9}.
Note that the lower edge of the chart of an edge, such as the circled elements in Figure \ref{2,6},  is $v_1^4$ periodic throughout the range of the chart, but the upper edge of the chart, such as the big dots in Figure \ref{2,6}, ceases to be $v_1^4$ periodic at each 2-power.

The upper edge $\E_{1,k}$ was discussed near the end of Section \ref{derivationsec}. It was stated there that we would ``justify the placement of the $\Si^{2^i}M_4^0$'s in the sequence of Definition \ref{E1'def}.'' In Sections \ref{kudifflsec} and \ref{A4sec}, we gave a fairly complete treatment of the obtaining of the $\A_4$ chart, including the upper edge. Everything generalizes, but two things lacking there were the implication from $ku_*$ differentials to $ko_*$ differentials, as we have been doing for the lower edges in this section, and a thorough approach to the arrangement of summands, as in (\ref{A7tableau}). The summands for $\At_4$ are the first row of (\ref{A7tableau}) with $V_7$ changed to $V_4$. The $ku$ analogue of $V_4$ was derived in Figure \ref{d5N} and appears in the right half of Figure \ref{koku4}, which shows the morphism $c$ from $V_4\lar\Si^8M_4^0$ to its $ku$ analogue. We use big dots for classes that map across under $c$, the ones not in $\im(\eta)$. The differential in the $ku$ chart was derived in Figure \ref{d2d6}, and is reproduced in Figure \ref{koku4}, which incorporates the implied differentials in the $ko$ chart.

\bigskip
\begin{minipage}{6in}
\begin{fig}\label{koku4}

{\bf $c$ for early part of $\At_4$}

\begin{center}

\begin{\tz}[scale=.45]
\draw (-.5,0) -- (13,0);
\draw (17.5,0) -- (30,0);
\draw (0,4) -- (0,0) -- (2,2);
\draw (3,4) -- (4,5) -- (4,3);
\draw (8,6) -- (8,4) -- (10,6);
\draw (10,7) -- (12,9) -- (12,7);
\draw (9,0) --(11,2) --(11,1) -- (13,3);
\draw [blue] (9,0) -- (8,6);
\draw [blue] (11,1) -- (10,7);
\draw [->] (18,0) -- (30,6);
\draw [->] (18,1) -- (30,7);
\draw [->] (18,2) -- (30,8);
\draw [->] (27,0) -- (30,1.5);
\draw (18,3) -- (20,4);
\draw (18,0) -- (18,4);
\draw (20,1) -- (20,4);
\draw (22,2) -- (22,4);
\draw (24,3) -- (24,5);
\draw (26,4) -- (26,6);
\draw (28,5) -- (28,7);
\node at (0,-.7) {$0$};
\node at (4,-.7) {$4$};
\node at (8,-.7) {$8$};
\node at (12,-.7) {$12$};
\node at (18,-.7) {$0$};
\node at (22,-.7) {$4$};
\node at (26,-.7) {$8$};
\node at (1,1) {\llb};
\node at (2,2) {\llb};
\node at (3,4) {\llb};
\node at (4,5) {\llb};
\node at (9,5) {\llb};
\node at (10,6) {\llb};
\node at (10,1) {\llb};
\node at (11,2) {\llb};
\node at (12,2) {\llb};
\node at (13,3) {\llb};
\node at (11,8) {\llb};
\node at (12,9) {\llb};
\node at (20,1) {\llb};
\node at (20,2) {\llb};
\node at (20,3) {\llb};
\node at (20,4) {\llb};
\node at (24,3) {\llb};
\node at (24,4) {\llb};
\node at (24,5) {\llb};
\node at (28,5) {\llb};
\node at (28,6) {\llb};
\node at (0,1) {\blb};
\node at (0,2) {\blb};
\node at (0,3) {\blb};
\node at (0,4) {\blb};
\node at (18,4) {\blb};
\node at (4,3) {\blb};
\node at (4,4) {\blb};
\node at (8,4) {\blb};
\node at (8,5) {\blb};
\node at (8,6) {\blb};
\node at (9,0) {\blb};
\node at (11,1) {\blb};
\node at (10,7) {\blb};
\node at (12,7) {\blb};
\node at (12,8) {\blb};
\node at (0,0) {\blb};
\node at (18,3) {\blb};
\node at (18,2) {\blb};
\node at (18,1) {\blb};
\node at (18,0) {\blb};
\node at (22,3) {\blb};
\node at (22,4) {\blb};
\node at (27,0) {\blb};
\node at (29,1) {\blb};
\node at (26,4) {\blb};
\node at (26,5) {\blb};
\node at (26,6) {\blb};
\node at (28,7) {\blb};
\node at (15,2) {$\mapright{c}$};
\draw [blue] (12,2) -- (11,8);
\draw [blue] (13,3) -- (12,9);
\draw [blue] (27,0) -- (26,6);
\draw [blue] (29,1) -- (28,7);
\end{\tz}
\end{center}
\end{fig}
\end{minipage}
\bigskip

There is an $\eta$ extension from grading 11 in the resulting chart, which is shown in the right side of Figure \ref{V4d}, in which filtrations have been increased to improve the picture. The second arrow of (\ref{A7tableau}) results in $\Si^8\E'_{2,3}\oplus \Si^{16}M_4^0$ pictured in Figure \ref{d2a} and in Figure \ref{new4} with a filtration shift. The differential from 17 to 16 in forming $\E_{1,4}$ in the $ko$ spectral sequence, pictured in Figure \ref{new4}, is implied by the $d_9$ in Figure \ref{homdiff}. This all generalizes, showing how the upper edge is formed to yield Definition \ref{E1'def}.

\subsection{Justifying extensions}\label{133}

Finally, we prove Theorem \ref{extnthm}, the exotic extensions. Let $\as_k$ denote the $ku_*$ version of $\A_k$, and $\Ast_k=\Si^{-2^{k+1}}\as_k$.
The extensions from $\Sigma^{2^\ell}\E_{2,\ell}$ into $\E_{1,k}$ stated in Theorem \ref{extnthm}(1) follow from the exact sequence (\ref{exact}) and the extensions in $\Ast_k$. The first extension for a fixed $\ell$ occurs in grading $3\cdot2^{\ell-1}+2$, so grading $2^{\ell-1}+2$ for $\E_{2.\ell}$. This is where the first lightning flash occurs for $\E'_{2,\ell}$. This can be seen in grading 26 and 50 in Figure \ref{A5} and in Figures \ref{four} and \ref{2,6}. One can compare with the $ku$ chart in Figure \ref{kuA5}. We illustrate in Figure \ref{5026}, in which the groups along the top are in $\E_{1,5}$, and those along the bottom are in $\Si^{2^\ell}\E_{2,\ell}$. The extension in $ko$ is implied by the exact sequence.

\bigskip
\begin{minipage}{6in}
\begin{fig}\label{5026}

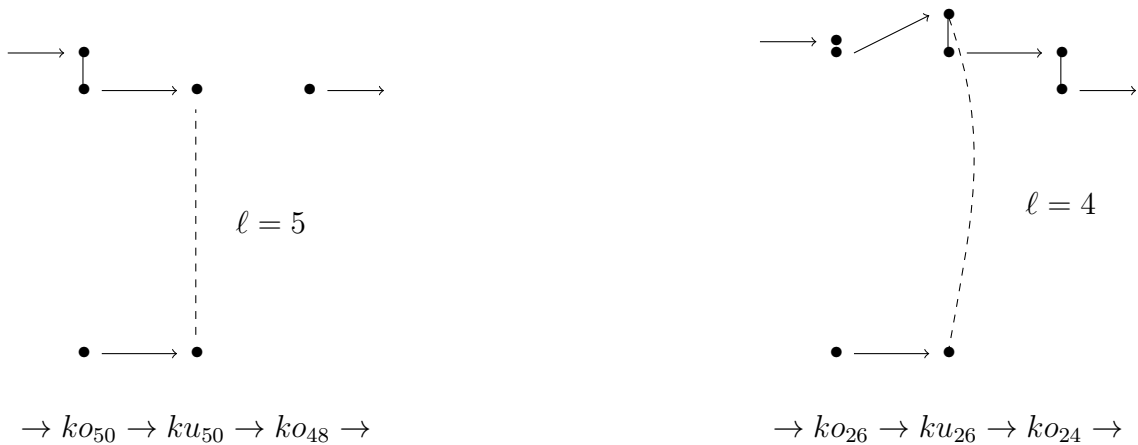
{\bf Extensions from $\Si^{2^\ell}\E_{2,\ell}$ into $\E_{1,5}$}

\begin{center}

\begin{\tz}[scale=.5]
\node at (0,0) {\lb};
\node at (3,0) {\lb};
\node at (0,7) {\lb};
\node at (0,8) {\lb};
\node at (3,7) {\lb};
\node at (6,7) {\lb};
\node at (20,0) {\lb};
\node at (23,0) {\lb};
\node at (20,8) {\lb};
\node at (20,8.3) {\lb};
\node at (23,8) {\lb};
\node at (23,9) {\lb};
\node at (26,8) {\lb};
\node at (26,7) {\lb};
\draw [->] (-2,8) -- (-.5,8);
\draw [->] (.5,7) -- (2.5,7); 
\draw [->] (6.5,7) -- (8,7);
\draw [->] (18,8.3) -- (19.5,8.3);
\draw [->] (20.5,8) -- (22.5,9);
\draw [->] (23.5,8) -- (25.5,8);
\draw [->] (26.5,7) -- (28,7);
\draw [->] (0.5,0) -- (2.5,0);
\node at (5,3.5) {$\ell=5$};
\node at (26,4) {$\ell=4$};
\draw [dashed] (3,.5) -- (3,6.5);
\draw [dashed] (23,0) to[out=80, in=290] (23,9);
\node at (3,-2) {$\to ko_{50}\to ku_{50} \to ko_{48}\to$};
\node at (23,-2) {$\to ko_{26}\to ku_{26} \to ko_{24}\to$};
\draw (0,7) -- (0,8);
\draw (23,8) -- (23,9);
\draw (26,7) -- (26,8);
\draw [->] (20.5,0) -- (22.5,0);
\end{\tz}
\end{center}
\end{fig}
\end{minipage}
\bigskip

It always happens this way: the only two possibilities for the groups in the exact sequence are illustrated in Figure \ref{5026}, and the $ku$ extension is implied by \cite[Definition 1.5 and Theorem 1.23]{DW}. On the other hand, classes $x$ in $\E_{2,\ell}$ in grading $2^{\ell-1}-8\Delta+2$ with $\Delta>0$ have $v_1^{4\Delta}x$ of filtration greater than that of the class in $2^{\ell-1}+2$ which supported the extension, so they cannot support an extension.
That there are no other extensions into the upper edge is true because 
extensions cannot hit into elements $x$ satisfying $\eta x\ne0$, nor into elements already divisible by 2 (since this could be accounted for by renaming), nor into periodic elements along the lower portion of the upper edge. Because of the very nice form of the upper edge, one can verify that this rules out everything except for a possible extension from $\Si^{2^k}\io_4$, where $\io_4$ is the bottom class of $\E_{2,k}$, into the last element of the only complete lightning flash in $\At_k$, an element $x$ such that $v_1^4x$ is hit by a differential. But both of these elements in $ko_{2^k+4}(K_2)$ are in the image from $ku_{2^k+6}(K_2)$, where there is no extension.

To prove the extensions in part (2) of Theorem \ref{extnthm}, we first note that for $e\ge2$ all nonzero groups in $\E_{e,\ell}$ in grading $j\equiv6$ mod 8 are $\zt$ groups, annihilated by $\eta$, and related to one another under $v_1^4$ periodicity. They pull back to $ku_{j+2}(K_2)$, where there are extensions between the pull-back elements, as can be seen from the inductive structure of the $\Ast$'s, implying the extension from $\Si^{2^{\ell+1-e}}\E_{e+1,\ell}$ to $\E_{e,\ell}$. 

There are two differences between the $\zt$'s in 2 mod 8 and those in 6 mod 8. In 6 mod 8, all classes are related by $v_1^4$ and are in $\ker(\eta)$, so they pull back to $ku_{*+2}(K_2)$ in (\ref{exact}). In grading 2 mod 8, they are related by $v_1^4$ until the next 2-power in Definition \ref{E'def} or Theorem \ref{edges}. The break points can also be thought of as the position a fraction $1/2^t$ along the edge for some $t\ge1$. Also, in 2 mod 8 the classes are not in $\im(\eta)$, so they map across to $ku_*(K_2)$.

In grading 2 mod 8, the only things that can possibly extend into the last half of $\E_{e,k}$ are (the last half of) $\E_{e+d,k}$ for some $d>0$. The only things that can extend into the second quarter of it are (the last half of) $\E_{e+d,k-1}$ for some $d>0$. Similarly, for second eighth, $\E_{e+d,k-2}$, etc. However, if $d>1$, this portion of $\E_{e+d,\ell}$ will have already extended into $\E_{e+d-1,\ell}$. We illustrate this in a schematic Figure \ref{lines}. It shows the edges beneath $\E_{2,7}$ and the key places where they drop their filtration in grading 2 mod 8 (halfway, one quarter, etc.). Extensions in grading 2 mod 8 occur from the second half of edges (where the lightning flashes start in $\E'$) into the edge immediately above them.  Other portions of edges have their periodicity in 2 mod 8 end prior to that of the portion of the edge above them. The extensions are implied by the extensions in $ku_*(K_2)$ between the corresponding classes.

\bigskip
\begin{minipage}{6in}
\begin{fig}\label{lines}

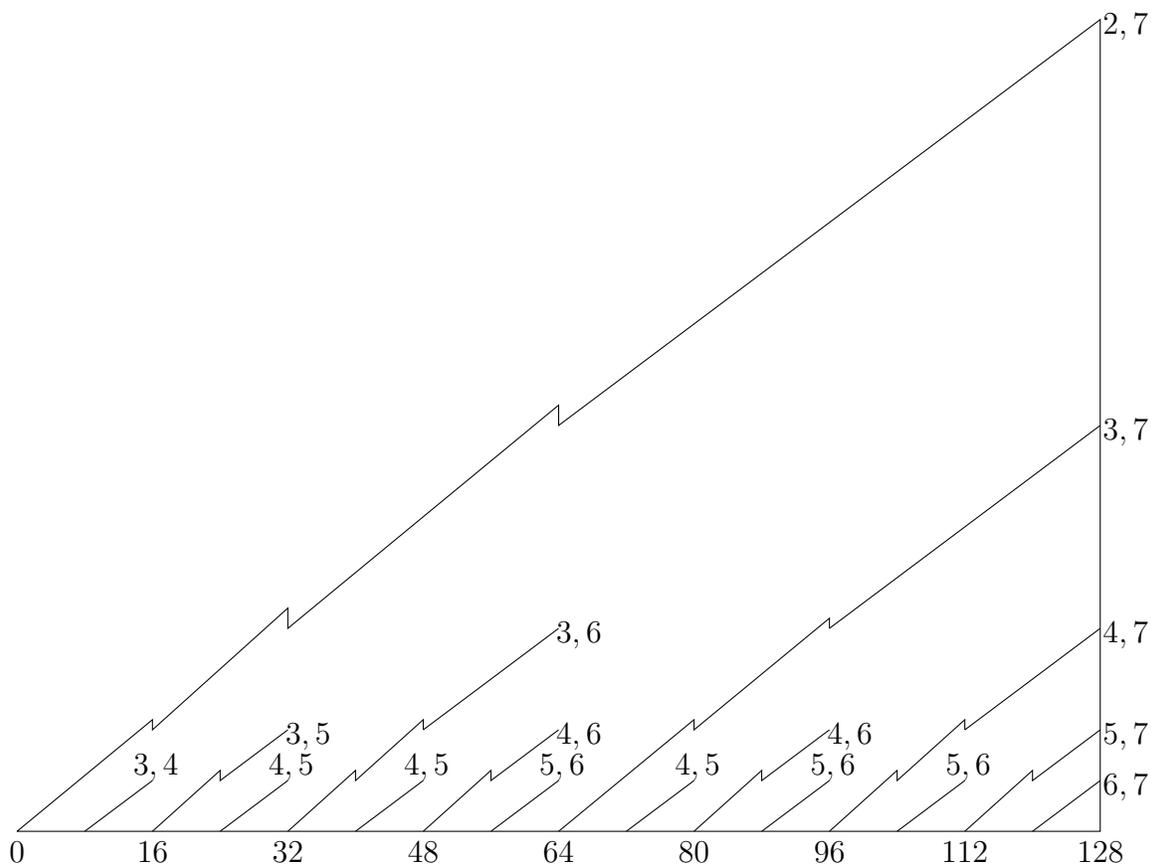
{\bf Schematic of edges under $\E_{2,7}$ in grading 2 mod 8}

\begin{center}

\begin{\tz}[scale=.45]
\draw (0,0) -- (4,3.3) -- (4,3) -- (8,6.6) -- (8,6) -- (16,12.6) -- (16,12) -- (32,24) -- (32,0) -- (0,0);
\draw  (4,0) -- (6,1.8) -- (6,1.5) -- (8,3);
\node at (32.75,23.8) {$2,7$};
\node at (8.6,2.8) {$3,5$};
\node at (0,-.6) {$0$};
\node at (4,-.6) {$16$};
\node at (8,-.6) {$32$};
\node at (12,-.6) {$48$};
\node at (16,-.6) {$64$};
\node at (20,-.6) {$80$};
\node at (24,-.6) {$96$};
\node at (28,-.6) {$112$};
\node at (32,-.6) {$128$};
\draw (8,0) --(10,1.8) -- (10,1.5) -- (12,3.3) -- (12,3);
\draw  (12,3) -- (16,6);
\node at (16.6,5.8) {$3,6$};
\draw  (12,0) -- (14,1.8) -- (14,1.5) -- (16,3);
\node at (16.6,2.8) {$4,6$};
\draw (16,0) -- (20,3.3) -- (20,3) -- (24,6.3) -- (24,6);
\draw (24,6) -- (32,12);
\node at (32.75,11.8) {$3,7$};
\draw  (20,0) -- (22,1.8) -- (22,1.5) -- (24,3);
\node at (24.6,2.8) {$4,6$};
\draw (24,0) --(26,1.8) -- (26,1.5) -- (28,3.3) -- (28,3);
\draw  (28,3) -- (32,6);
\node at (32.75,5.8) {$4,7$};
\draw (28,0) -- (30,1.8) -- (30,1.5) -- (32,3);
\node at (32.75,2.8) {$5,7$};
\draw (2,0) -- (4,1.5);
\draw (6,0) -- (8,1.5);
\draw (10,0) -- (12,1.5);
\draw (14,0) -- (16,1.5);
\draw (18,0) -- (20,1.5);
\draw (22,0) -- (24,1.5);
\draw (26,0) -- (28,1.5);
\draw (30,0) -- (32,1.5);
\node at (4.1,1.9) {$3,4$};
\node at (8.1,1.9) {$4,5$};
\node at (12.1,1.9) {$4,5$};
\node at (16.1,1.9) {$5,6$};
\node at (20.1,1.9) {$4,5$};
\node at (24.1,1.9) {$5,6$};
\node at (28.1,1.9) {$5,6$};
\node at (32.75,1.3) {$6,7$};
\end{\tz}
\end{center}
\end{fig}
\end{minipage}
\bigskip

There is no extension from the last class in $\E_{e,e+d}$ when $e\equiv3$ mod 4 because either there is nothing into which to extend, or else it would extending into an element $x$ satisfying $\eta x\ne0$.

Many extensions also follow from the following fact (\cite[near (2.2)]{DM}) about Toda brackets  and differentials in the Adams spectral sequence:
\begin{equation}\label{Toda} \text{if }d_r(\langle \a,2,\eta\rangle)=\b\eta,\text{ then }2\a=\b.\end{equation}
This bracket represents $v_1\a$; the generator of the $\Z/4$ in a lightning flash is obtained from the bottom class via this bracket. The situation in Figure \ref{bkt} appears frequently, implying $2\a=\b$. 
 Often the classes $\a\eta$ and $\a\eta^2$ will be hit by differentials, but that doesn't matter.

\bigskip
\begin{minipage}{6in}
\begin{fig}\label{bkt}

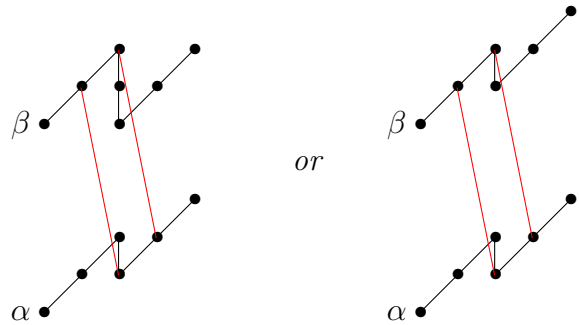
{\bf Charts implying $2\a=\b$}

\begin{center}

\begin{\tz}[scale=.5]
\draw (0,0) -- (2,2) -- (2,1) -- (4,3);
\draw (0,5) -- (2,7) -- (2,5) -- (4,7);
\node at (0,0) {\lb};
\node at (1,1) {\lb};
\node at (2,2) {\lb};
\node at (2,1) {\lb};
\node at (3,2) {\lb};
\node at (4,3) {\lb};
\node at (0,5) {\lb};
\node at (1,6) {\lb};
\node at (2,7) {\lb};
\node at (2,6) {\lb};
\node at (2,5) {\lb};
\node at (3,6) {\lb};
\node at (4,7) {\lb};
\draw [red] (2,1) -- (1,6);
\draw [red] (3,2) -- (2,7);
\node at (-.6,0) {$\a$};
\node at (-.6,5) {$\b$};
\draw (10,0) -- (12,2) -- (12,1) -- (14,3);
\draw (10,5) -- (12,7) -- (12,6) -- (14,8);
\node at (10,0) {\lb};
\node at (11,1) {\lb};
\node at (12,2) {\lb};
\node at (12,1) {\lb};
\node at (13,2) {\lb};
\node at (14,3) {\lb};
\node at (10,5) {\lb};
\node at (11,6) {\lb};
\node at (12,7) {\lb};
\node at (12,6) {\lb};
\node at (7,4) {or};
\node at (13,7) {\lb};
\node at (14,8) {\lb};
\draw [red] (12,1) -- (11,6);
\draw [red] (13,2) -- (12,7);
\node at (9.4,0) {$\a$};
\node at (9.4,5) {$\b$};

\end{\tz}
\end{center}
\end{fig}
\end{minipage}
\bigskip

All other possible extensions can be  ruled out by $v_1^4$ periodicity, as we have been using, the fact that $2\eta=0$, and comparing with $ku$. Probably the most subtle of these is illustrated by the situation involving $\E_{3,8}$ and $\Si^{64}\E_{4,8}$ in grading $64+12$. 
Using Figure \ref{4,8}, we illustrate the problem in Figure \ref{problem}, in which $\E'_{3,8}$ has the lightning flashes, which should have much higher filtration, and the worrisome extension question involves classes $A$ and $B$, indicated by big dots. The problem is that $v_1^4B=0$ (due to the differential) and $A$ is not in $\im(\eta)$, so we cannot use $v_1^4$ or $\eta$ to detect the non-extension..

\bigskip
\begin{minipage}{6in}
\begin{fig}\label{problem}

{\bf $\E'_{3,8}$ and $\Si^{64}\E'_{4,8}$}

\begin{center}

\begin{\tz}[scale=.5]
\draw (9.5,0) -- (22,0);
\draw (10,0) -- (11,1);
\draw (12,0) -- (12,1);
\draw (16,1) -- (16,2);
\draw (18,3) -- (19,4);
\draw (20,4) -- (20,2) -- (22,4);
\draw (10,7) -- (12,9) -- (12,8) -- (14,10);
\draw (18,11) -- (20,13) -- (20,12) -- (22,14);
\draw [blue] [->] (12,8) -- (11,11);
\draw [blue] [->] (13,9) -- (12,12);
\draw [blue] [->] (20,12) -- (19,15);
\draw [blue] [->] (21,13) -- (20,16);
\draw [blue] (20,2) -- (19,12);
\draw [blue] (21,3) to[out=95, in=295] (20,13);
\node at (10,0) {\llb};
\node at (11,1) {\llb};
\node at (12,0) {\llb};
\node at (14,0) {\llb};
\node at (16,1) {\llb};
\node at (16,2) {\llb};
\node at (18,3) {\llb};
\node at (19,4) {\llb};
\node at (20,2) {\llb};
\node at (20,3) {\llb};
\node at (20,4) {\llb};
\node at (21,3) {\llb};
\node at (22,4) {\llb};
\node at (10,7) {\llb};
\node at (11,8) {\llb};
\node at (22,14) {\llb};
\node at (12,8) {\llb};
\node at (13,9) {\llb};
\node at (14,10) {\llb};
\node at (18,11) {\llb};
\node at (19,12) {\llb};
\node at (20,13) {\llb};
\node at (20,12) {\llb};
\node at (21,13) {\llb};
\node at (12,1) {\blb};
\node at (12,9) {\blb};
\node at (12.5,1) {$A$};
\node at (12,9.6) {$B$};
\node at (8,-.6) {$64+$};
\node at (10,-.6) {$10$};
\node at (12,-.6) {$12$};
\node at (20,-.6) {$20$};
\end{\tz}
\end{center}
\end{fig}
\end{minipage}
\bigskip

However, $A$ pulls back to $ku_{*+2}(K_2)$, since $\eta A=0$, and there is no extension on the pullback element, as the $ku$ analogue of $\E_{4,8}$ begins as in Figure \ref{kuE}, with the first extension occurring on the element in grading 16. So $2A=0$.

\bigskip
\begin{minipage}{6in}
\begin{fig}\label{kuE}

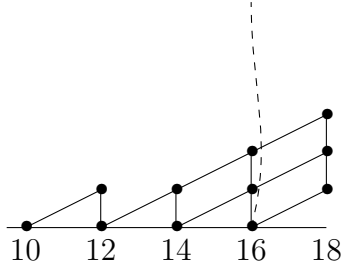
{\bf Beginning of $ku$ analogue of $\E_{4,k}$}

\begin{center}

\begin{\tz}[scale=.5]
\draw (9.5,0) -- (18,0);
\draw (10,0) -- (12,1) -- (12,0) -- (18,3) -- (18,2) -- (14,0) -- (14,1);
\draw (18,1) -- (16,0) -- (16,2);
\draw (18,1) -- (18,2);
\node at (10,0) {\lb};
\node at (12,0) {\lb};
\node at (12,1) {\lb};
\node at (14,0) {\lb};
\node at (14,1) {\lb};
\node at (16,0) {\lb};
\node at (16,1) {\lb};
\node at (16,2) {\lb};
\node at (18,1) {\lb};
\node at (18,2) {\lb};
\node at (18,3) {\lb};
\draw [dashed] (16,0) to[out=75, in=270] (16,6);
\node at (10,-.6) {$10$};
\node at (12,-.6) {$12$};
\node at (14,-.6) {$14$};
\node at (16,-.6) {$16$};
\node at (18,-.6) {$18$};

\end{\tz}
\end{center}
\end{fig}
\end{minipage}
\bigskip

\section{Proof for $z^i\B_{k,\ell}$}\label{Bklpfsec}
In this section, we prove Theorem \ref{ziBklthm}. We have defined a $ku$ version, $B_{k,\ell}$, of $\B_{k,\ell}$ via summands in Definition \ref{ABdef} and shown in Theorem \ref{everything} that these summands are what are required along with the $A_k$'s to fill out the $E_1$-module $H^*(K_2)$, and we showed in Theorem \ref{closed} that each $B_{k,\ell}$ is closed under the differentials of \cite[Theorem 3.1]{DW}. It is worth noting, but perhaps not necessary, to point out that $B_{k,\ell}$ corresponds to $B_kz_\ell\oplus y_1^{2^{k-1}-1}qS_{k,\ell}\oplus y_kB_kZ_k^\ell$ in \cite{DW}.

Most of our discussion will deal with the case $i=0$ in $z^i\B_{k,\ell}$. We will comment briefly on the effect of $i$ at the end of this section.
In (\ref{Bkltableau}), we gave a general description of the corresponding $ko$ summands, but we feel that the following explicit example, the summands for $\Bt_{5,9}$ is useful in understanding the proof.
\begin{align}\Si M_7^5\lar \Si^8M_4^4\quad \Si^9M_4^6\lar \Si^{16}M_5^4&\quad \Si^{17}M_5^6\lar\Si^{24}M_4^5\quad\Si^{25}M_4^7\lar\Si^{32}M_6^4\nonumber\\
\Si^{33}M_6^6\lar\Si^{40}M_4^5\quad\Si^{41}M_4^7\lar\Si^{48}M_5^5&\quad\Si^{49}M_5^7\lar\Si^{56}M_4^6\quad \Si^{57}M_4^8\lar\Si^{64}M_{11}^0\nonumber\\
\Si^{65}M_{11}^2\lar \Si^{72}M_4^1\quad\Si^{73}M_4^3\lar\Si^{80}M_5^1&\quad\Si^{81}M_5^3\lar \Si^{88}M_4^2\quad\Si^{89}M_4^4\lar\Si^{96}M_6^1\nonumber\\
\Si^{97}M_6^3\lar\Si^{104}M_4^2\quad\Si^{105}M_4^4\lar\Si^{112}M_5^2&\quad\Si^{113}M_5^4\lar\Si^{120}M_4^3\quad\Si^{121}M_4^5\lar\Si^{128}M_7^1\label{B59}\end{align}

\noindent If $\ell$ is increased by 1, the superscripts in the first two rows are increased by 1, except for the last term, whose subscript is increased by 1, while the last two rows remain unchanged except that the subscript of the first term (of the second half) is increased by 1.

What can be observed here, and is true in general using (\ref{Bkltableau}), is that $\Bt_{k,\ell}$ has $2^k$ summands of which (a) the first half agree with those that form $\E_{\ell-k+1,\ell}$ and everything under it in $\At_\ell$, and (b) the second half agree with the  second half of $\At_{k+1}$ except that the first summand is $\Si^{2^{k+1}+1}M^2_{\ell+2}$ instead of $\Si^{2^{k+1}+1}M^2_{k+2}$.\footnote{There is a difference in the last summand, too, $\Si^{2^{k+2}}M^1_{k+2}$ versus $\Si^{2^{k+2}}M^0_{k+3}$, which does not affect this part of the analysis.} (Comparison of (\ref{B59}) and (\ref{A7tableau}): first half of (\ref{B59}) with superscripts of all but last summand decreased by 2 equals $\Si^{-192}$ of last two rows of (\ref{A7tableau}), while (b)  last two rows of (\ref{B59})  equals rows 3 and 4 of (\ref{A7tableau}) except for the first (and last) summand.)

Part (2) of Theorem \ref{ziBklthm} follows readily from observation (a) above. The differentials among the summands which form $\E_{\ell-k+1,\ell}$ are the same regardless of whether the summands are appearing in $\At_\ell$ or in $\Bt_{k,\ell}$. There are two differences. One is that in $\Bt_{k,\ell}$ there is not an edge $\E_{\ell-k,\ell}$ to be hit by the differentials in grading 4 and 5 mod 8 that occurred in $\At_\ell$. The other has to do with what happens to the $\Si^{2^{k+1}}M_{\ell+2}^0$ at the end of the sequence of summands being considered. In neither case does it contribute to $\E_{\ell-k+1,\ell}$. In both cases, it has its subscript decreased $k-1$ times by differentials into lightning flashes in this batch of summands. In  $\At_\ell$, it continues to have its subscript decreased by differentials into lightning flashes formed by earlier edges. In $\Bt_{k,\ell}$, it combines with part of the $\Si^{2^{k+1}+1}M_{\ell+2}^2$ which occurs at the beginning of the second half of the list of summands forming $\Bt_{k,\ell}$ to form part (3) of Theorem \ref{ziBklthm} in a way that we shall discuss shortly.

We explain now why the filtration of  $\Si^{2^{k+1}}M_{\ell-k+3}$, which is formed after the last summand of the first half of $\Bt_{k,\ell}$ has its subscript reduced $k-1$ times, is increased by $2^k-k-1$, as claimed in part (1) of the theorem. The first quarter of the summands of $\Bt_{k,\ell}$ will have truncated the top part of the initial $\Si M_{k+2}^{\ell-k+1}$ chart $k-2$ times resulting stably in lightning flashes with initial classes in position $(8i+2,4i-\ell+1)$. [\![The initial class of $\Si M_{k+2}^0$ is in $(2,0)$, so stably the initial classes of $\Si M_{k+2}^{\ell-k+1}$ are in $(8i+2,4i-(\ell-k+1))$, and the $k-2$ truncations of the top reduce the second component by $k-2$.]\!] Stably the last step in forming  $\Si^{2^{k+1}}M_{\ell-k+3}$ will be as depicted in Figure \ref{laststep}, where the lightning flash is the one just described, and the generator of the $\Z/2^{\ell-k+1}$'s in $\Si^{2^{k+1}}M_{\ell-k+3}$ will be the big dot. Since the generator of the $\Z/2^{\ell-k+1}$'s in $\Si^{2^{k+1}}M_{\ell-k+3}$ are in position $(2^{k+1}+3+8j,2-\ell+k+4j)$, which must equal $(8i+3,4i-\ell+1)$, we obtain that the difference in filtrations is $4i-\ell+1-(2-\ell+k+4j)=2^k-k-1$, since $2^{k+1}+8j=8i$.

\bigskip
\begin{minipage}{6in}
\begin{fig}\label{laststep}

{\bf Last step in forming $\Si^{2^{k+1}}M_{\ell-k+3}$}

\begin{center}

\begin{\tz}[scale=.45]
\draw (5,2) -- (3,0) -- (3,7) -- (1,5);
\draw [red] (2,1) -- (4,3) -- (4,2) -- (6,4);
\draw [blue] (3,0) -- (2,1);
\draw [blue] (4,1) -- (3.2,2.2);
\draw [blue] (5,2) -- (4,3);
\node at (3,1) {\blb};
\end{\tz}
\end{center}
\end{fig}
\end{minipage}
\bigskip

As noted earlier, the second half of the summands of $\Bt_{k,\ell}$ agree with those of the second half of $\At_{k+1}$ except that the first summand is $\Si^{2^{k+1}+1}M_{\ell+2}^2$ instead of $\Si^{2^{k+1}+1}M_{k+2}^2$.\footnote{The difference in the last summand noted in the previous footnote affects $\E_{1,k+1}$
but not $\E_{2,k+1}$.}
The edge $\Si^{2^{k+1}}\E_{2,k+1}$ is formed from $\Si^{2^{k+1}+1}M_{k+2}^2$ by interactions with $\Si^{2^i}M_4^1$'s formed from the summands which follow it in $\At_{k+1}$. The leftover parts of those summands form the edges under $\Si^{2^{k+1}}\E_{2,k+1}$. The same interactions will occur with $\Si^{2^{k+1}+1}M^2_{\ell+2}$, forming something that will contribute to $\Si^{2^{k+1}}C_{0,k}$ in part (3) of Theorem \ref{ziBklthm}. The leftover parts from the subsequent summands will be the same here as they were in $\Si^{2^{k+1}}\E_{2,k+1}$, as stated in part (3) of the theorem.

The $\Si^{2^{k+1}+1}M_{\ell+2}^2$ is interacted with by $k-1$ $\Si^{2^i}M_4^1$'s, so that stably it is $\Si^{2^{k+1}+1}M_{\ell-k+3}$. Stably there must be differentials from this isomorphically to the $\Si^{2^{k+1}}M_{\ell-k+3}$ formed in the paragraph preceding Figure \ref{laststep}. As described in the paragraph preceding Figure \ref{M11} and illustrated there, we apply $(v_1^4)^{-1}$ to obtain additional differentials. To see that these differentials are $d_{2^k}$, we observe that the bottoms of the towers in the lower half of Figure \ref{M11} are like those of $M_{\ell+2}^2$, while those of the top half are those of $M^0_{\ell-k+3}$ with filtrations increased by $2^k-k-1$. Noting that larger subscript of $M$ causes smaller filtration of generators, we obtain that the difference in filtrations of generators is
$$2^k-k-1-(\ell-k+3)-(-(\ell+4))=2^k.$$

Since the differentials are $d_{2^k}$ and $M_{\ell-k+3}$ has its filtrations increased by $2^k-k-1$, all classes in $M_{\ell-k+3}$ of filtration $\le k$ will not be hit. Nor will the $\eta$-pairs along the upper edge until grading $2^{k+2}+2$, when the final truncation in the second half of the summands of $\Bt_{k,\ell}$ will have occurred.

This covers parts (1), (3), and (4) of Theorem \ref{ziBklthm} when $i=0$, completing the proof in this case, except for some fine tuning.
\begin{itemize}
\item We explain the last sentence of part (3) by means of an example. In $\Bt_{4,9}$, we have  $\Si^{32}C_{0,4}$, which is the big dots in the bottom half of Figure \ref{M11}, and the edges under $\Si^{32}\E_{2,5}$ together with differentials involving those subedges. You should compare the chart for $\E_{2,5}$ in Figure \ref{four} with $C_{0,4}$. Upper edge elements of $\E_{2,5}$ inject. However, the classes in grading 27 and 28 indicated by big dots in Figure \ref{four} are hit by differentials from $\Si^{16}\E'_{3,5}$. The meaning of the last sentence of part (3) of Theorem \ref{ziBklthm} is that the corresponding elements of $\Si^{32}C_{0,4}$ are also hit by differentials from $\Si^{48}\E'_{3,5}$.
\item We explain Remark \ref{rem} via an example, $\Bt_{3,7}$. The first half of the summands are
$$\Si M_5^5\lar\Si^8M_4^4\quad \Si^9M_4^6\lar\Si^{16}M_9^0.$$
After the indicated $d_2$ differentials and a filtration shift of the second half, the picture will be as in Figure \ref{B37}, where the portion from the first $d_2$ is in black and the second in red. The class in $(21,4)$ is in the $\Si^{16}M_6^0$ part of $\Si^{16}C_{0,3}$, while the class in $(20,3)$ is in the $\E'_{5,7}$ part of $\Bt_{3,7}$.
\item The exotic extensions in part (4) of the theorem follow from (\ref{Toda}) and the differential from the next-to-top element of one tower to the top of another, like the ones from grading 52 and 60 in Figure \ref{M11}. The  extension and the  exotic $\eta$ extension of Remark \ref{rem2} are nicely visualized in Figure \ref{P3}. With $P_3$ the usual stunted real projective space, Figure \ref{P3} illustrates the easily-proved result that the cofiber of $P_3\w bo\mapright{2}P_3\w bo$ equals $\Si^4 M\w bo\vee H$, where $M$ is the mod-2 Moore spectrum and $H$ a mod-2 Eilenberg MacLane spectrum. The diagram shows just the first few gradings of the cofiber sequence, and indicates the gradings both in the cofiber sequence and the corresponding elements in Figure \ref{M11} if filtrations of the bottom part were increased to make the $d_{16}$ differential look like a $d_1$.
\end{itemize}

\bigskip
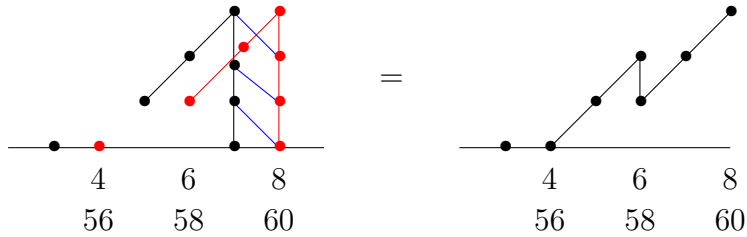
\begin{minipage}{6in}
\begin{fig}\label{B37}

{\bf Result from first half of summands of $\Bt_{3,7}$}

\begin{center}

\begin{\tz}[scale=.5]
\draw (10.5,0) -- (30,0);
\draw (11,0) -- (12,1) -- (12,0);
\draw (13,0) -- (14,1);
\draw (18,2) -- (20,4) -- (20,3) -- (22,5);
\draw (26,6) -- (28,8) -- (28,7) -- (30,9);
\draw [red] (17,4) -- (19,6) -- (19,4);
\draw [red] (23,4) -- (23,7);
\draw [red] (25,8) -- (27,10) -- (27,5) -- (29,7);
\draw [blue] (21,3) -- (20,4);
\draw [blue] (23,4) -- (22,5);
\draw [blue] (27,5) -- (26,6);
\draw [blue] (28,6) -- (27.2,7.2);
\draw [blue] (29,7) -- (28,8);
\draw [blue] (30,9) -- (31,8);
\node at (12,-.7) {$12$};
\node at (16,-.7) {$16$};
\node at (20,-.7) {$20$};
\node at (24,-.7) {$24$};
\node at (28,-.7) {$28$};
\node at (11,0) {\lb};
\node at (12,1) {\lb};
\node at (12,0) {\lb};
\node at (13,0) {\lb};
\node at (14,1) {\lb};
\node at (18,2) {\lb};
\node at (19,3) {\lb};
\node at (20,4) {\lb};
\node at (20,3) {\lb};
\node at (21,4) {\lb};
\node at (22,5) {\lb};
\node at (26,6) {\lb};
\node at (27.2,7.2) {\lb};
\node at (28,8) {\lb};
\node at (28,7) {\lb};
\node at (29,8) {\lb};
\node at (30,9) {\lb};
\node [red] at (20,0) {\lb};
\node [red] at (22,1) {\lb};
\node [red] at (17,4) {\lb};
\node [red] at (18,5) {\lb};
\node [red] at (19,6) {\lb};
\node [red] at (19,5) {\lb};
\node [red] at (19,4) {\lb};
\node [red] at (23,4) {\lb};
\node [red] at (23,5) {\lb};
\node [red] at (23,6) {\lb};
\node [red] at (23,7) {\lb};
\node [red] at (25,8) {\lb};
\node [red] at (26,9) {\lb};
\node [red] at (27,10) {\lb};
\node [red] at (27,9) {\lb};
\node [red] at (27,8) {\lb};
\node [red] at (27,7) {\lb};
\node [red] at (27,6) {\lb};
\node [red] at (27,5) {\lb};
\node [red] at (28,6) {\lb};
\node [red] at (29,7) {\lb};
\node [red] at (21,3) {\lb};

\end{\tz}
\end{center}
\end{fig}
\end{minipage}
\bigskip

\bigskip
\begin{minipage}{6in}
\begin{fig}\label{P3}

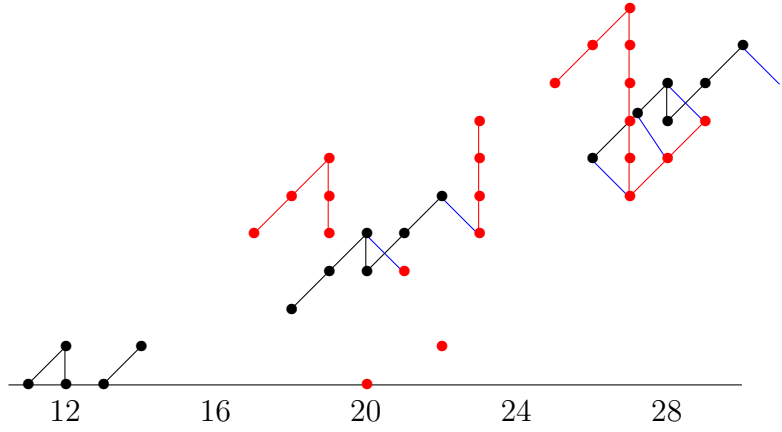
{\bf Cofiber of 2 on $P_3\w bo$ equals $\Si^4M\w bo\vee H$}

\begin{center}

\begin{\tz}[scale=.6]
\draw (2,0) -- (9,0);
\draw (5,1) -- (7,3) -- (7,0);
\draw [red] (6,1) -- (8,3) -- (8,0);
\draw [blue] (8,0) -- (7,1);
\draw [blue] (8,1) -- (7,1.8);
\draw [blue] (8,2) -- (7,3);
\draw (12,0) -- (18,0);
\draw (14,0) -- (16,2) -- (16,1) -- (18,3);
\node at (4,-.7) {$4$};
\node at (6,-.7) {$6$};
\node at (8,-.7) {$8$};
\node at (4,-1.6) {$56$};
\node at (6,-1.6) {$58$};
\node at (8,-1.6) {$60$};
\node at (14,-.7) {$4$};
\node at (16,-.7) {$6$};
\node at (18,-.7) {$8$};
\node at (14,-1.6) {$56$};
\node at (16,-1.6) {$58$};
\node at (18,-1.6) {$60$};
\node at (10.5,1.5) {$=$};
\node at (3,0) {\lb};
\node at (5,1) {\lb};
\node at (6,2) {\lb};
\node at (7,3) {\lb};
\node at (7,0) {\lb};
\node at (7,1) {\lb};
\node at (7,1.8) {\lb};
\node at (13,0) {\lb};
\node at (14,0) {\lb};
\node at (15,1) {\lb};
\node at (16,2) {\lb};
\node at (16,1) {\lb};
\node at (17,2) {\lb};
\node at (18,3) {\lb};
\node [red] at (6,1) {\lb};
\node [red] at (7.2,2.2) {\lb};
\node [red] at (8,3) {\lb};
\node [red] at (8,2) {\lb};
\node [red] at (8,1) {\lb};
\node [red] at (8,0) {\lb};
\node [red] at (4,0) {\lb};

\end{\tz}
\end{center}
\end{fig}
\end{minipage}
\bigskip

So far in this section, we have been dealing with $z^0\Bt_{k,\ell}$. For $z^i\Bt_{k,\ell}$, all summands are obtained from the corresponding summands of $\Bt_{k,\ell}$ by increasing superscripts by $i$. The effect on edges $\E_{e,\ell}$ with $e>1$ (which is all that is relevant for $\Bt_{k,\ell}$) is to increase both subscripts by $i$. To see this, we first observe that all summands, except the last, in the tableau for $\E_{e+1,\ell+1}$ are obtained from those of $\E_{e,\ell}$ by increasing all superscripts by 1. This is nicely illustrated in (\ref{A7tableau}), in which rows 3 and 4 are the summands for $\E_{2,6}$, while rows 7 and 8 are the summands for $\E_{3,7}$. That the different possibilities for the last summand do not affect the edge is discussed in and around Figure \ref{E56} and (\ref{e+3}).

The $\Si^{2^{k+1}}M^0_{\ell+2}$ and $\Si^{2^{k+1}+1}M^2_{\ell+2}$ in the middle of the tableau for $\Bt_{k,\ell}$ are responsible for parts (1), (3), and (4) of Theorem \ref{ziBklthm}. If both superscripts are increased by $i$, the filtration shift ($2^k-k-1$) required in forming Figure \ref{laststep} is unchanged since both charts have their filtrations changed by the same amount.

The paragraph following Figure \ref{laststep} explains that $C_{0,k}$ is formed from $\E'_{2,\infty}$ because the second half of the summands of $\Bt_{k,\ell}$
(except the first) are the same as the first $2^{k-1}$ summands of $\Si^{2^{k+1}}\E_{2,\infty}$. Increasing the superscripts to form $z^i\Bt_{k,\ell}$ causes the summands to be the same as those of $\Si^{2^{k+1}}\E_{2+i,\infty}$, which is relevant in Definition \ref{hdef}. The other thing relevant there is the function $h_k$, which is determined by the upper edge of the stable chart, which will be $i$ units lower.

\section{The exact $ko_*$-$ku_*$ sequence} \label{exactsec}
 In this section, we discuss how the exact sequence (\ref{exact}) relating $ko_*(K_2)$ and $ku_*(K_2)$ can be used to justify the consistency of our results, and also to obtain additional information about $ko_*(K_2)$.  There are three things that make these comparisons somewhat cumbersome.

\begin{enumerate}
\item $ku_*(K_2)$ contains filtration-0 elements due to free $E_1$-module summands that are not part of free $A(1)$-module summands. These played a significant role in \cite{DW}, and also play a significant role in analyzing the exact sequence (\ref{exact}). In Proposition \ref{E1prop}, we explain exactly where these classes occur in terms of the $A(1)$-module summands with chart $M_k^i$, which we also denote here as $M_k^i$.
\item We perform much shifting of filtrations in our determination of $ko_*(K_2)$. This is done to improve the appearance of the charts and to make $\eta$ extensions easier to see. This causes the $ko$- and $ku$-charts to not match up as nicely as one might like.
\item When elements of $ko_*(K_2)$ are involved in the homomorphism $\eta$ and also in a differential in the ASS of $ko_*(K_2)$, that can cause ambiguity regarding filtration of classes.
\end{enumerate}

\begin{prop}\label{E1prop} The free $E_1$-module summands of the $A(1)$-module $M_k^i$ have generators in grading
\begin{eqnarray*}8j+2&\text{if}&i=4j+1\\
8j+3&\text{if}&i=4j+2\\
8j+4&\text{if}&k+i=4j+6\\
8j+5&\text{if}&k+i=4j+7.\end{eqnarray*}
The $A(1)$-module $NU$ in Figure \ref{NU}, which is used in forming the charts $V_k$, has a single free $E_1$-module summand, with generator in grading 11.\end{prop}
\begin{proof}
The result for $M_k^0$ is immediate from \cite[Figure 2.3]{DW2}. For each $M_k^i$, it is just a matter of drawing the $A(1)$-module that matches the Ext chart. We illustrate with $M_9^1$. The Ext charts for $M^0_9$ and $M_9^1$ are shown in Figure \ref{Ext19}, and then the $A(1)$-module yielding the Ext chart for $M_9^1$ is in Figure \ref{M19}.

\bigskip
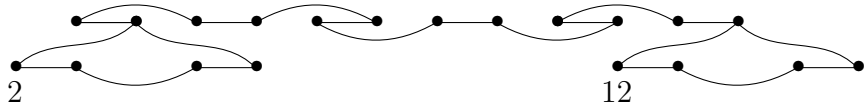
\begin{minipage}{6in}
\begin{fig}\label{Ext19}

{\bf Ext charts for $M^0_9$ and $M_9^1$.}

\begin{center}

\begin{\tz}[scale=.45]
\draw (0,0) -- (16,0);
\draw (1,0) -- (3,2) -- (3,0);
\draw (7,0) -- (7,3);
\draw (9,4) -- (11,6) -- (11,0) -- (13,2);
\draw (15,3) -- (15,7);
\draw (21,0) -- (36,0);
\draw (22,0) -- (23,1) -- (23,0);
\draw (27,0) -- (27,2);
\draw (29,3) -- (31,5) -- (31,0);
\draw (32,0) -- (33,1); 
\draw (35,2) -- (35,6);
\node at (1,-.7) {$1$};
\node at (3,-.7) {$3$};
\node at (7,-.7) {$7$};
\node at (11,-.7) {$11$};
\node at (22,-.7) {$2$};
\node at (27,-.7) {$7$};
\node at (32,-.7) {$12$};
\node at (1,0) {\lb};
\node at (2,1) {\lb};
\node at (3,2) {\lb};
\node at (3,1) {\lb};
\node at (3,0) {\lb};
\node at (7,0) {\lb};
\node at (7,1) {\lb};
\node at (7,2) {\lb};
\node at (7,3) {\lb};
\node at (9,4) {\lb};
\node at (10,5) {\lb};
\node at (11,6) {\lb};
\node at (11,5) {\lb};
\node at (11,4) {\lb};
\node at (11,3) {\lb};
\node at (11,2) {\lb};
\node at (11,1) {\lb};
\node at (11,0) {\lb};
\node at (12,1) {\lb};
\node at (13,2) {\lb};
\node at (15,3) {\lb};
\node at (15,4) {\lb};
\node at (15,5) {\lb};
\node at (15,6) {\lb};
\node at (15,7) {\lb};
\node at (22,0) {\lb};
\node at (23,1) {\lb};
\node at (23,0) {\lb};
\node at (27,0) {\lb};
\node at (27,1) {\lb};
\node at (27,2) {\lb};
\node at (29,3) {\lb};
\node at (30,4) {\lb};
\node at (31,5) {\lb};
\node at (31,4) {\lb};
\node at (31,3) {\lb};
\node at (31,2) {\lb};
\node at (31,1) {\lb};
\node at (31,0) {\lb};
\node at (32,0) {\lb};
\node at (33,1) {\lb};
\node at (35,2) {\lb};
\node at (35,3) {\lb};
\node at (35,4) {\lb};
\node at (35,5) {\lb};
\node at (35,6) {\lb};
\node at (5,5) {$\M^0_9$};
\node at (25,4) {$M_9^1$};

\end{\tz}
\end{center}
\end{fig}
\end{minipage}

\bigskip
\begin{minipage}{6in}
\begin{fig}\label{M19}

{\bf $M_9^1$.}

\begin{center}

\begin{\tz}[scale=.4]

\node at (-4,-6.5) {\lb};
\node at (-2,-6.5) {\lb};

\draw (-4,-6.5) -- (-2,-6.5);

\draw (-4,-6.5) to[out=30, in=150] (0,-6.5);


\node at (0,-6.5) {\lb};

\node at (4,-6.5) {\lb};
\node at (6,-6.5) {\lb};
\node at (8,-6.5) {\lb};
\node at (2,-6.5) {\lb};
\node at (12,-6.5) {\lb};
\node at (14,-8 ) {\lb};
\node at ( 22,-8) {\lb};
\node at ( 20,-8) {\lb};
\node at ( 16,-8) {\lb};
\node at (12,-6.5) {\lb};
\node at (14,-6.5) {\lb};
\node at (16,-6.5) {\lb};

\draw (0,-6.5) -- (2,-6.5);
\draw (4,-6.5) -- (6,-6.5);
\draw (8,-6.5) -- (10,-6.5);
\draw (14,-8) -- (16,-8);
\draw (20,-8) -- (22,-8);

\draw (4,-6.5) to[out=330, in=210] (8,-6.5);
\draw (2,-6.5) to[out=30, in=150] (6,-6.5);
\draw (10,-6.5) to[out=330, in=210] (14,-6.5);
\draw (12,-6.5) to[out=30, in=150] (16,-6.5);

\draw (16,-8) to[out=330, in=210] (20,-8);
\draw (14,-8) to[out=45, in=225] (18,-6.5);
\draw (18,-6.5) to[out=315, in=135] (22,-8);
\node at (18,-6.5) {\lb};
\draw (16,-6.5) -- (18,-6.5);
\node at (10,-6.5) {\lb};
\draw (12,-6.5) -- (14,-6.5);
\node at (-6,-8) {\lb};
\node at (-4,-8) {\lb};
\node at (2,-8) {\lb};
\node at (0,-8) {\lb};
\draw (-6,-8) -- (-4,-8);
\draw (2,-8) -- (0,-8);
\node at (-6,-8.8) {$2$};
\node at (14,-8.8) {$12$};
\draw (-4,-8) to[out=330, in=210] (0,-8);
\draw (-6,-8) to[out=45, in=225] (-2,-6.5);
\draw (-2,-6.5) to[out=315, in=135] (2,-8);
\node at (-2,-6.5) {\lb};
\draw (-4,-6.5) -- (-2,-6.5);
\end{\tz}
\end{center}
\end{fig}
\end{minipage}

\bigskip

\end{proof}
\bigskip

Our main example of the exact sequence (\ref{exact}) will be for $\At_3$. In Figure \ref{A3}, we present the $ko$ and $ku$ versions from Figure \ref{A3chart} and \cite[Figure 2]{DG}. 
The tableau for $\At_3$, from (\ref{Aktableau}), is
$$V_3\lar \Si^{8}M_4^0\quad \Si^{9}M_4^2\lar \Si^{16}M_5^0.$$
By Proposition \ref{E1prop}, the $ku$ version  $\Ast_3$ will have additional filtration-0 classes in gradings 3, 12, and 13, which have been incorporated into Figure \ref{A3}.

\bigskip
\begin{minipage}{6in}
\begin{fig}\label{A3}

{\bf $\At_3$ and $\Ast_3$.}

\begin{center}

\begin{\tz}[scale=.5]

\draw (33.5,0) -- (49,0);
\draw (34,3) -- (34,0) -- (36,2);
\node at (37,3) {\lb};
\draw (38,3) -- (38,4) -- (37,3);
\draw (42,4) -- (44,6) -- (44,5) -- (45,6);
\node at (46,0) {\lb};
\node at (48,1) {\lb};
\node at (48,-.6) {$14$};
\node at (34,-.6) {$0$};
\node at (42,-.6) {$8$};
\node at (34,0) {\lb};
\node at (34,1) {\lb};
\node at (34,2) {\lb};
\node at (34,3) {\lb};
\node at (35,1) {\lb};
\node at (36,2) {\lb};
\node at (38,3) {\lb};
\node at (38,4) {\lb};
\node at (42,4) {\lb};
\node at (43,5) {\lb};
\node at (44,6) {\lb};
\node at (44,5) {\lb};
\node at (45,6) {\lb};
\draw (52,0) -- (52,3);
\draw (52,1) -- (58,4);
\draw (52,2) -- (54,3);
\draw (54,1) -- (54,3);
\draw (56,2) -- (56,3);
\draw (58,3) -- (58,4);
\draw (52,0) -- (66,7);
\draw (66,0) -- (66,7);
\draw (66,0) -- (68,1);
\node at (55,0) {\lb};
\node at (64,0) {\lb};
\node at (65,0) {\lb};
\node at (52,0) {\lb};
\node at (52,1) {\lb};
\node at (52,2) {\lb};
\node at (52,3) {\lb};
\node at (54,1) {\lb};
\node at (54,2) {\lb};
\node at (54,3) {\lb};
\node at (56,2) {\lb};
\node at (56,3) {\lb};
\node at (58,3) {\lb};
\node at (58,4) {\lb};
\node at (60,4) {\lb};
\node at (62,5) {\lb};
\node at (64,6) {\lb};
\node at (66,7) {\lb};
\node at (66,0) {\lb};
\node at (68,1) {\lb};
\draw (51.5,0) -- (68,0);
\node at (52,-.6) {$0$};
\node at (60,-.6) {$8$};
\node at (66,-.6) {$14$};
\node at (42,-1.6) {$ko_*$-version};
\node at (58,-1.6) {$ku_*$-version};
\end{\tz}
\end{center}
\end{fig}
\end{minipage}
\bigskip

One convenient way of seeing how the exact sequence (\ref{exact}) can be used to relate $ko_*(X)$ and $ku_*(X)$ is to draw  $ko_*(X)\oplus ko_*(\Si^2 X)$, with differentials for the action of $\eta$. What is left after the differentials should be $ku_*(X)$. We have done this for $\At_3$ in Figure \ref{koku}. In order to illustrate difficulty (3) listed at the beginning of this section, we have also included the classes involved in the $d_7$-differential from $\Si^8\E'_{2,3}$ to $\E'_{1,3}$. 

\bigskip
\begin{minipage}{6in}
\begin{fig}\label{koku}

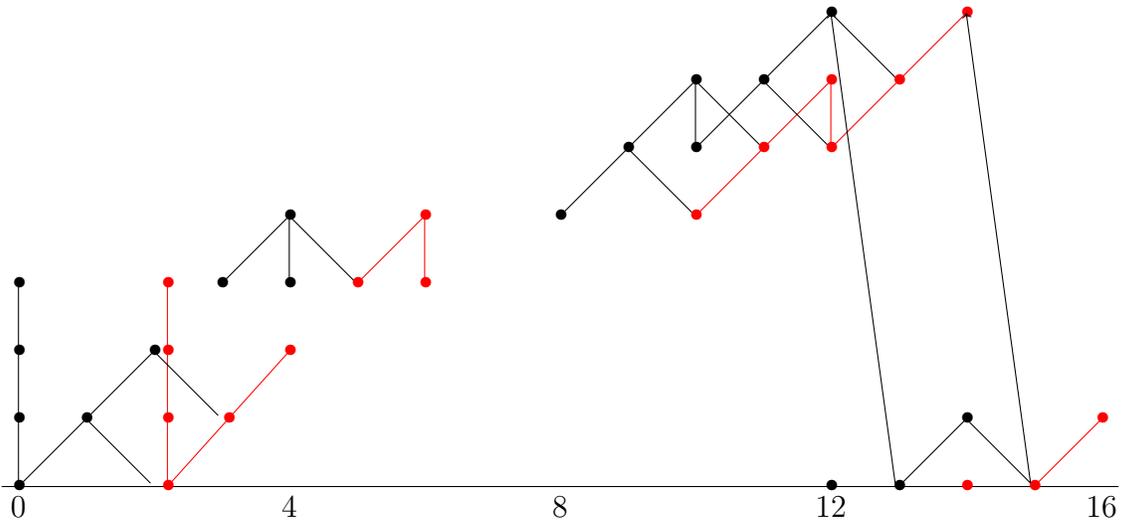
{\bf Forming $\Ast_3$ from $\At_3\oplus\Sigma^2\At_3$.}

\begin{center}

\begin{\tz}[scale=.45]
\draw (31.5,0) -- (64.5,0);
\draw (36,4) -- (32,0) -- (32,6);
\draw [red] (36.4,6) -- (36.4,0) -- (40,4);
\node at (32,0) {\lb};
\node at (32,2) {\lb};
\node at (32,4) {\lb};
\node at (32,6) {\lb};
\node at (34,2) {\lb};
\node at (36,4) {\lb};
\node[red] at (36.4,0) {\lb};
\node[red] at (36.4,2) {\lb};
\node [red] at (36.4,4) {\lb};
\node [red] at (36.4,6) {\lb};
\node [red] at (38.2,2) {\lb};
\node [red] at (40,4) {\lb};
\node [red] at (42,6) {\lb};
\node [red] at (44,6) {\lb};
\node [red] at (44,8) {\lb};
\node at (38,6) {\lb};
\draw (40,6) -- (40,8) -- (38,6);
\draw [red] (44,6) -- (44,8) -- (42,6);
\draw (35.9,.1) -- (34,2);
\draw (37.9,2.1) -- (36,4);
\draw (41.9,6.1) -- (40,8);
\node at (40,6) {\lb};
\node at (40,8) {\lb};
\draw (48,8) -- (52,12) -- (52,10) -- (56,14);
\draw [red] (52,8) -- (56,12) -- (56,10) -- (60,14);
\draw (51.9,8.1) -- (50,10);
\draw (53.9,10.1) -- (52,12);
\draw (55.9,10.1) -- (54,12);
\draw (57.9,12.1) -- (56,14);
\draw (58,0) -- (60,2);
\draw [red] (62,0) -- (64,2);
\draw (61.9,.1) -- (60,2);
\node [red] at (52,8) {\lb};
\node [red] at (54,10) {\lb};
\node [red] at (56,12) {\lb};
\node [red] at (56,10) {\lb};
\node [red] at (58,12) {\lb};
\node [red] at (60,14) {\lb};
\node [red] at (60,0) {\lb};
\node [red] at (62,0) {\lb};
\node [red] at (64,2) {\lb};
\node at (48,8) {\lb};
\node at (50,10) {\lb};
\node at (52,12) {\lb};
\node at (52,10) {\lb};
\node at (54,12) {\lb};
\node at (56,14) {\lb};
\node at (56,0) {\lb};
\node at (58,0) {\lb};
\node at (60,2) {\lb};
\draw [->] (57.9,.1) -- (56,14);
\draw [->] (61.9,.1) -- (60,14);
\node at (32,-.6) {$0$};
\node at (40,-.6) {$4$};
\node at (48,-.6) {$8$};
\node at (56,-.6) {$12$};
\node at (64,-.6) {$16$};

\end{\tz}
\end{center}
\end{fig}
\end{minipage}
\bigskip

The classes which remain in the left half of Figure \ref{koku} agree with the left half of the $ku$ chart in Figure \ref{A3} except for the filtration of the class in grading 3. The class in $(3,3)$ in $\At_3$ was initially in filtration 0, but had its filtration increased to 3 for reasons discussed prior to Figure \ref{UNpic}, so this is an example of difficulty (2) in the $ko$-$ku$ comparison. 

In the right side of Figure \ref{koku}, note that two classes in grading 13 are hitting the top class in 12, one due to an Adams differential and the other due to the $\eta$ homomorphism. If we had removed the classes involved in the Adams differential, it would have looked like the $ku$ class in grading 13 should have filtration 6. But if we had removed the Ext classes involved in the $\eta$ homomorphism first, then a filtration-0 class in 13 would be left, which is where it should be in the $ku$ chart. Similarly, one class in grading 14 can be considered to have either filtration 1 or 7, depending on whether the $\eta$ is considered before or after the Adams differential.

So you can see that the exact sequence works, but care is required.
For our final example, we consider $\B_{3,4}$, but in much less detail. We compare Figures \ref{kuB34} and \ref{B34}. Proposition \ref{E1prop} will add filtration-0 $\zt$'s to Figure \ref{kuB34} in grading 68, 70, 74, 78, 84, 90, 94, 98, and 100.

We make the following observations:
\begin{itemize}
    \item The lightning flash in the middle of Figure \ref{B34} plus its double suspension with $\eta$ differential exactly gives the $v$-tower from $(81,0)$ in Figure \ref{kuB34}, but with filtration 4 larger. 
    \item One justification for the filtration increases of the middle lightning flash and the high $\eta$ pair in Figure \ref{B34} is that the last two classes in the middle lightning flash and the high $\eta$ pair are both in $\im(v_1^4)$. This is not the case in the $ku$ situation for parity reasons.
    \item One way to see that $\eta$ is nonzero on the class in $(88,0)$ in Figure \ref{B34} is that the element in grading 89 must be in $\im(\eta)$ because it must map trivially in (\ref{exact}) since the $ku$ group is 0 in grading 89.
    \item If the high $\eta$ pair in 89 and 90 are lowered to filtration 1 and 2 and the $\eta$ inserted on the class in $(88,0)$, the analogue of the analysis like that in Figure \ref{koku} is straightforward.
    
\end{itemize}

\section{Consequences for Spin manifolds}\label{SWsec}
In this section, we review the relationship of $ko_*(K_2)$ to Stiefel-Whitney classes and immersions of Spin manifolds, and prove Theorem \ref{SWthm}, building on work done by
 the author and W.S.Wilson in \cite{DWSW}.

The generalized homology theory associated to the Thom spectrum MSpin has the property that ${\rm MSpin}_n(K(\zt,k))$ is the set of cobordism classes of pairs $(M,x)$, where $M$ is an $n$-dimensional Spin manifold and $x\in H^k(M,\zt)$. (\cite{At}) By \cite{ABP}, localized at 2, $bo$ is a split summand of MSpin, and so $ko_n(K(\zt,k))$ is a direct summand of ${\rm MSpin}_n(K(\zt,k))$. In particular, all of the elements in our calculation of $ko_n(K_2)$ give cobordism classes of pairs, an $n$-dimensional Spin manifold $M$ together with an element of $H^2(M;\zt)$.

The following  result was proved in \cite[Theorem 1.2]{DW}.
\begin{thm} \label{1.2} Let $h:ko_*(X)\to H_*(X;\zt)$ denote the Hurewicz homomorphism. There exists an $n$-dimensional {\rm Spin}-manifold $M$ with nonzero dual Stiefel-Whitney class $\overline{w}_{n-k}(M)$ if and only if there exists an element $\a\in ko_n(K(\zt,k))$ such that $\langle\chi\sq^{n-k}\io_k,h_*(\a)\rangle\ne0$.\end{thm}
\ni Here $\chi$ is the antiautomorphism of the Steenrod algebra, and $\io_k$ is the fundamental class in $H^k(K(\zt,k);\zt)$.

Dual Stiefel-Whitney classes are important because if $\overline{w}_c(M)\ne0$ for an $n$-manifold $M$, then $M$ cannot be immersed in $\R^{n+c-1}$. Thus we have the following corollary of our new result, Theorem \ref{SWthm}.
\begin{cor} The values of $n$ for which
 there exists an $n$-dimensional {\rm Spin} manifold that does not immerse in $\R^{2n-3}$, detected by Stiefel-Whitney classes, are exactly all $2$-powers $\ge8$.\end{cor}

\begin{proof} [Proof of Theorem \ref{SWthm}]
It was shown in \cite[Theorem 1.2]{DWSW} that a necessary condition for existence of an $n$-dimensional {\rm Spin}-manifold $M$ with $\overline{w}_{n-k}(M)\ne0$ is $\chi\sq^{n-k}\io_k\not\in\im(\sq^1,\sq^2)$, and in \cite[Theorem 1.3]{DWSW} it was shown that $\chi\sq^{n-2}\io_2\not\in\im(\sq^1,\sq^2)$ if and only if $n$ or $n-1$ is a 2-power $\ge8$.
It was shown in \cite{anti} that $\chi\sq^{2^e-1}=\sq^{2^{e-1}}\sq^{2^{e-2}}\cdots\sq^1$ 
and $\chi\sq^{2^e-2}=\sq^{2^{e-1}}\sq^{2^{e-2}}\cdots\sq^2$ .

Thus $\chi\sq^{2^e-2}\io_2=\io_2^{2^{e-1}}$. The bottom class of $\A_k$ is dual to $\io^{2^k}$. This follows from Definition \ref{ABdef} and the definition of $\Ah^o_k$ in Section \ref{kosummandsec}. Since the bottom class of $\A_k$ supports a nonzero subgroup of $ko_{2^{k+1}}(K_2)$ for $k\ge2$, its generator gives the desired element $\a$  in Theorem \ref{1.2}.

On the other hand, the  class of lowest grading in our $A(1)$-module $\M_{k+2}$ equals $\chi\sq^{2^{k+2}-1}\io_2$ mod decomposables. This follows from \cite[Theorem 2.2]{DW2} since $$\sq^{2^{k+1}}\cdots\sq^1\io_2=u_{k+2},$$ a multiplicative generator of $H^*(K_2)$. By Definition \ref{ABdef} and its adaptation to $\Ah^o_k$, $\M_{k+2}$ occurs at the end of  $\Ah^o_k$. Its  class of lowest grading supports a differential. To see this, note that in Definition \ref{E1'def}, $\M_{k+2}$ becomes the $\Si^{2^{k+1}}M_4^0$ at the end of the sequence defining $\E_{1,k}'$. ($\M_{k+2}$ has an additional $\Si^{2^{k+1}}$ since $\E_{1,k}'$ is for $\At_k$, while $\Ah^o_k$ is for $\A_k$.)
In forming (\ref{E1seq}), $\Si^{-2^{k+1}}\M_{k+2}$ started as $\Si^{2^{k+1}}M_{k+2}^0$ but had its bottom cut down by differentials $k-2$ times, turning it into  $\Si^{2^{k+1}}M_4^0$.
Because of the differential on the class arising from $\chi\sq^{2^{k+2}-1}\io_2$, we deduce that when $n$ is of the form $2^e+1$ and $k=2$, there is no class that works as $\a$ in Theorem \ref{1.2}.

 Using $\At_4$ as an example, in Figure \ref{E14} the class in grading 33 is the class supporting the differential discussed here. It started as the $\Si^{32}M_6^0$ at the end of the first row of (\ref{A7tableau}) but became $M_4^0$ after  differentials into the $\Si^{25}M_4^3$  preceding it and into the resultant of the arrow preceding that. Its filtration was increased several times.
\end{proof}

\section{$ko^*(K_2)$}\label{cohsec}
In this section, we state results regarding $ko^*(K_2)$ and several duality relationships. As we did for $ku^*(K_2)$ in \cite{DW}, we depict charts with $ko$-cohomology grading increasing from right to left.

\begin{defin} \label{Ak*def}We define $\At^*_k$ for $k\ge2$ and $(z^i\Bt_{k,\ell})^*$ for $1\le k<\ell$ to be the $ko$-cohomology charts obtained by applying
$\ext_{A(1)}(\zt,-)$ to the $A(1)$-modules
 $\Ah^o_k$ and $z^i\Bh^o_{k,\ell}$ defined in Section \ref{kosummandsec}, and incorporating differentials and extensions in the ASS, and filtration increases of the sort used  above. Also $\A_k^*=\Si^{2^{k+1}}\At_k^*$ and $(z^i\B_{k,\ell})^*=\Si^{2^{\ell+2}}(z^i\Bt_{k,\ell})^*$. We define $\A_1^*$ by incorporating $d_2$ differentials into $\ext_{A(1)}(\zt,U\oplus N)$, similarly to Figure \ref{U+N}.\end{defin}

 The proof of the following analogue of Theorem \ref{main} is identical to that of Theorem \ref{main}, which appeared in Section \ref{kosummandsec}.
 \begin{thm} There is an isomorphism of $ko^*$-modules
$$ko^*(K_2)\approx\bigoplus_{k\ge1} \bigoplus_{i\ge0}\Si^{2^{k+2}i}\A_k^*\oplus\ \bigoplus_{1\le k<\ell}\bigoplus_{i,j\ge0}\Si^{2^{k+2}i+2^{\ell+3}j}(z^{\a(j)}\B_{k,\ell})^*$$
plus a trivial $ko^*$-module.\label{ko^*thm}
\end{thm}

The next two results, which, with Theorem \ref{ko^*thm}, determine $ko^*(K_2)$, are proved in Section \ref{cohpfsec}.

\begin{thm}\label{B^*thm} For any $i$, there is an isomorphism of $ko^*$-modules
\begin{equation}\label{Biso}(z^t\Bt_{k,\ell})^*\approx(z^{4i-\ell-t}\Bt_{k,\ell})_{2^{k+2}+8i+4-*},\end{equation}
where $z^i\Bt_{k,\ell}$ is as determined in Theorem \ref{ziBklthm}.\end{thm}

The value of $i$ is irrelevant since $(z^{r+4}\Bt_{k,\ell})_{*+8}\approx(z^r\Bt)_*$. One would usually choose $i$ so that $0\le 4i-\ell-t\le3$. Since $ko^*=ko_{-*}$, the right hand side of (\ref{Biso}) is a $ko^*$-module. Since we picture $ko^*(-)$ charts with gradings increasing from right to left, the charts on both sides of (\ref{Biso}) will look the same. For example, the charts in Figure \ref{ziB34} for $i=1$, 2, 3  can be interpreted as charts of $(z^t\B_{3,4})^*$ for $t=3$, 2, 1 if the indicated gradings $x$ are replaced by $52-x$.

\begin{thm}\label{Ak^*thm} For $t\ge1$ and $\eps\in\{1,2\}$, $\At^*_{2t+\eps}$ 
is $(2^{2t+\eps+1}+8t+4)$-dual to the chart described as follows. 
Let $e=2t+2-\eps$ and $\ell=4t+1$, and $T=2^{2t+\eps+1}+8t-2$. For $e\le e'\le4t$, let $D=2^{\ell-e+2}-2^{\ell-e'+2}$ and let $\Si^D\E^s_{e',\ell}$ be a modified version of $\Si^D\E'_{e',\ell}$ which extends the sequence of lightning flashes or portions thereof through grading $T$. The chart is formed from these together with all edges under them with second subscript $<\ell$, and including differentials and extensions among these edges and subedges as described in Theorems \ref{third} and \ref{extnthm}. There are additional differentials on the added lightning flashes with $e'<e$ on classes in grading 4 or 5 mod 8. There are additional classes $x\in(T-1,0)$, $y\in(T,0)$, and $\eta x=2y$, with $d_2(x)\ne0$. Finally, all of the new classes (except $y$ and those in $\E^s_{e,\ell}$) have exotic extensions into the classes above them.

\end{thm}

For $k\le2$, $\At_k^*$ does not quite fit the theorem. In Figures \ref{ko^*A3}, \ref{ko^*A4}, and \ref{ko^*A5}, we display $\A_1^*$ and $\At^*_k$ for $2\le k\le5$. Keep in mind that $\At_k^*$ is always suspended by odd multiples of $2^{k+1}$ in $ko^*(K_2)$ if $k\ge2$.

We illustrate the theorem with $\At_5^*$. It is 84-dual to a chart which we now describe. For comparison, use Figure \ref{ko^*A5} with grading $x$ replaced by $84-x$. Start with $\E'_{5,9}$ in Figure \ref{5,9} with its last lightning flash completed and followed by one more (in grading 74 to 78). The classes in Figure \ref{5,9} with big dots are hit by differentials, but the ones with circles do not support differentials. Next we add $\Si^{32}\E'_{6,9}$, $\Si^{48}\E'_{7,9}$, and $\Si^{56}\E'_{8,9}$ (use Figure \ref{four} with these replaced by $\Si^{40}\E'_{2,5}$, $\Si^{56}\E'_{3,5}$, and $\Si^{64}\E'_{4,5}$)\footnote{The extra $\Si^8$ is due to increasing subscripts by 4.} extended through grading 78, which means that the latter two must have lightning flashes completed, and one additional lightning flash added to $\Si^{40}\E'_{2,5}$ in grading 74 to 78. Big dots and circles in these will all take effect. Also add subedges of all of these with second subscript $<9$, but do not extend them. Add classes $x$ and $y$ as in Theorem \ref{Ak^*thm}, and insert differentials and extensions as in Figure \ref{78}, in which the $\E_{5,9}$ part should be 11 higher.

\bigskip
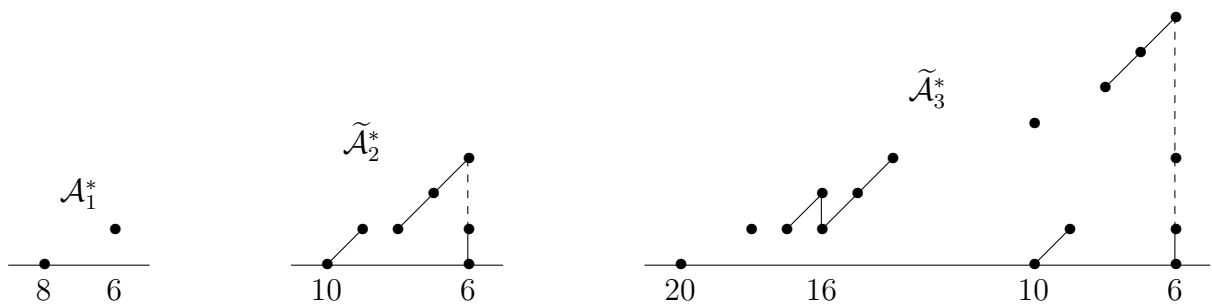
\begin{minipage}{6in}
\begin{fig}\label{78}

{\bf Forming the end of (the 84-dual of) $\At_5^*$}

\begin{center}

\begin{\tz}[scale=.57]
\draw (66,0) -- (79,0);
\draw (67,1) -- (68,2) -- (68,1) -- (70,3);
\draw (66,8) -- (68,10) -- (68,9) -- (70,11);
\draw (66,13) -- (68,15) -- (68,14) -- (70,16);
\draw (74,0) -- (76,2) -- (76,1) -- (78,3);
\draw (74,4) -- (76,6) -- (76,5) -- (78,7);
\draw (74,12) -- (76,14) -- (76,13) -- (78,15);
\draw (74,17) -- (76,19) -- (76,18) -- (78,20);
\draw [dashed] (66,8) -- (66,13);
\draw [dashed] (70,3) -- (70,16);
\draw [dashed] (74,0) -- (74,17);
\draw [dashed] (78,1) -- (78,20);
\node at (66,-.5) {$66$};
\node at (70,-.5) {$70$};
\node at (74,-.5) {$74$};
\node at (78,-.5) {$78$};
\node at (77,-.4) {$x$};
\node at (78.4,.3) {$y$};
\node at (78.7,3) {$8,9$};
\node at (78.7,7) {$7,9$};
\node at (78.7,15) {$6,9$};
\node at (78.7,20) {$5,9$};
\draw (77,0) -- (78,1) -- (78,0);
\draw [red] (68,1) to[out=110, in=270] (67,9);
\node at (67,1) {\lb};
\node at (68,2) {\lb};
\node at (70,3) {\lb};
\node at (70,11) {\lb};
\node at (70,16) {\lb};
\node at (69,15) {\lb};
\node at (68,14) {\lb};
\node at (66,13) {\lb};
\node at (66,8) {\lb};
\node at (74,0) {\lb};
\node at (75,1) {\lb};
\node at (78,0) {\lb};
\node at (78,1) {\lb};
\node at (78,3) {\lb};
\node at (78,7) {\lb};
\node at (78,15) {\lb};
\node at (78,20) {\lb};
\node at (77,19) {\lb};
\node at (76,18) {\lb};
\node at (74,17) {\lb};
\node at (74,4) {\lb};
\node at (74,12) {\lb};
\draw [red] (69,2) to[out=90, in=300] (68,10);
\draw [red] (68,9) to[out=110, in=270] (67,14);
\draw [red] (69,10) to[out=90, in=300] (68,15);
\draw [red] (77,0) -- (76,2);
\draw [red] (76,1) to[out=110, in=270] (75,5);
\draw [red] (77,2) to[out=90, in=300] (76,6);
\draw [red] (76,5) to[out=110, in=270] (75,13);3
\draw [red] (77,6) to[out=90, in=300] (76,14);
\draw [red] (76,13) to[out=110, in=270] (75,18);
\draw [red] (77,14) to[out=90, in=300] (76,19);

\end{\tz}
\end{center}
\end{fig}
\end{minipage}

\bigskip
\begin{minipage}{6in}
\begin{fig}\label{ko^*A3}

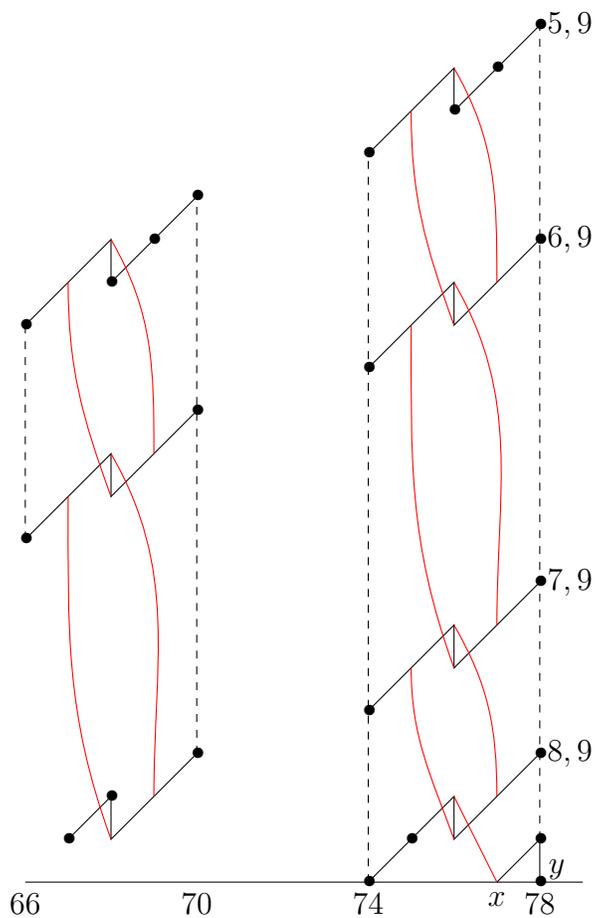
{\bf Charts for $k\le3$}

\begin{center}

\begin{\tz}[scale=.47]
\draw (-1,0) -- (3,0); 
\draw (7,0) -- (13,0);
\draw (17,0) -- (33,0);
\node at (0,-.7) {$8$};
\node at (2,-.7) {$6$};
\node at (8,-.7) {$10$};
\node at (12,-.7) {$6$};
\node at (18,-.7) {$20$};
\node at (22,-.7) {$16$};
\node at (28,-.7) {$10$};
\node at (32,-.7) {$6$};
\node at (0,0) {\lb};
\node at (2,1) {\lb};
\node at (8,0) {\lb};
\node at (9,1) {\lb};
\node at (10,1) {\lb};
\node at (11,2) {\lb};
\node at (12,3) {\lb};
\node at (12,1) {\lb};
\node at (12,0) {\lb};
\node at (18,0) {\lb};
\node at (20,1) {\lb};
\node at (21,1) {\lb};
\node at (22,2) {\lb};
\node at (22,1) {\lb};
\node at (23,2) {\lb};
\node at (24,3) {\lb};
\node at (28,4) {\lb};
\node at (30,5) {\lb};
\node at (31,6) {\lb};
\node at (32,7) {\lb};
\node at (32,3) {\lb};
\node at (32,1) {\lb};
\node at (32,0) {\lb};
\node at (28,0) {\lb};
\node at (29,1) {\lb};
\draw (8,0) -- (9,1);
\draw (10,1) -- (12,3);
\draw (12,0) -- (12,1);
\draw (21,1) -- (22,2) -- (22,1) -- (24,3);
\draw (28,0) -- (29,1);
\draw (32,0) -- (32,1);
\draw [dashed] (32,1) -- (32,7);
\draw (32,7) -- (30,5);
\draw [dashed] (12,1) -- (12,3);
\node at (1,2) {$\A_1^*$};
\node at (9,3.5) {$\At_2^*$};
\node at (25,5) {$\At_3^*$};
\end{\tz}
\end{center}
\end{fig}
\end{minipage}
\bigskip

\bigskip
\begin{minipage}{6in}
\begin{fig}\label{ko^*A4}

{\bf $\At_4^*$}

\begin{center}

\begin{\tz}[scale=.4]
\draw (-.5,0) -- (34.5,0);
\draw (4,0) -- (4,1);
\draw (6,1) -- (8,3) -- (8,1) -- (10,3);
\draw (14,4) -- (16,6) -- (16,5) -- (18,7);
\draw (23,1) -- (24,2);
\draw (24,9) -- (26,11);
\draw [dashed] (26,3) -- (26,11);
\draw [dashed] (30,0) -- (30,12);
\draw (30,0) -- (31,1);
\draw [dashed] (34,0) -- (34,15);
\draw (34,15) -- (32,13);
\node at (0,0) {\lb};
\node at (4,0) {\lb};
\node at (4,1) {\lb};
\node at (6,1) {\lb};
\node at (7,2) {\lb};
\node at (8,3) {\lb};
\node at (8,2) {\lb};
\node at (8,1) {\lb};
\node at (9,2) {\lb};
\node at (10,3) {\lb};
\node at (10,0) {\lb};
\node at (12,4) {\lb};
\node at (14,1) {\lb};
\node at (14,4) {\lb};
\node at (15,5) {\lb};
\node at (16,6) {\lb};
\node at (16,5) {\lb};
\node at (17,6) {\lb};
\node at (18,7) {\lb};
\node at (20,0) {\lb};
\node at (22,1) {\lb};
\node at (23,1) {\lb};
\node at (24,2) {\lb};
\node at (22,8) {\lb};
\node at (26,3) {\lb};
\node at (24,9) {\lb};
\node at (25,10) {\lb};
\node at (26,11) {\lb};
\node at (26,3) {\lb};
\node at (30,0) {\lb};
\node at (31,1) {\lb};
\node at (30,4) {\lb};
\node at (30,12) {\lb};
\node at (34,0) {\lb};
\node at (34,1) {\lb};
\node at (34,3) {\lb};
\node at (34,7) {\lb};
\node at (34,15) {\lb};
\node at (33,14) {\lb};
\node at (32,13) {\lb};
\node at (0,-.6) {$40$};
\node at (4,-.6) {$36$};
\node at (8,-.6) {$32$};
\node at (12,-.6) {$28$};
\node at (16,-.6) {$24$};
\node at (20,-.6) {$20$};
\node at (24,-.6) {$16$};
\node at (28,-.6) {$12$};

\node at (34,-.6) {$6$};

\end{\tz}
\end{center}
\end{fig}
\end{minipage}
\bigskip

\begin{minipage}{6in}
\begin{fig}\label{ko^*A5}

{\bf $\At^*_5$}

\begin{center}

\begin{\tz}[scale=.44]
\draw (-1.5,0) -- (34,0);
\node at (30,21) {\lb};
\node at (0,-.6) {$72$};
\node at (4,-.6) {$68$};
\node at (8,-.6) {$64$};
\node at (12,-.6) {$60$};
\node at (16,-.6) {$56$};
\node at (20,-.6) {$52$};
\node at (24,-.6) {$48$};
\node at (28,-.6) {$44$};
\draw (-1,0) -- (0,1) -- (0,0);
\draw (4,0) -- (4,1);
\draw (6,2) -- (7,3);
\draw (8,3) -- (8,1) -- (10,3);
\draw (14,5) -- (16,7) -- (16,5) -- (18,7);
\draw (22,1) -- (24,3);
\draw [dashed] (26,0) -- (26,4);
\draw (24,10) -- (24,9) -- (26,11);
\draw (30,12) -- (32,14) -- (32,13) -- (34,15);
\node at (-1,0) {\lb};
\node at (0,0) {\lb};
\node at (0,1) {\lb};
\node at (4,0) {\lb};
\node at (4,1) {\lb};
\node at (8,0) {\lb};
\node at (8,1) {\lb};
\node at (8,2) {\lb};
\node at (8,3) {\lb};
\node at (7,3) {\lb};
\node at (6,2) {\lb};
\node at (10,1) {\lb};
\node at (12,4) {\lb};
\node at (14,5) {\lb};
\node at (15,6) {\lb};
\node at (16,7) {\lb};
\node at (16,6) {\lb};
\node at (16,5) {\lb};
\node at (17,6) {\lb};
\node at (18,7) {\lb};
\node at (16,0) {\lb};
\node at (18,0) {\lb};
\node at (20,1) {\lb};
\node at (22,1) {\lb};
\node at (23,2) {\lb};
\node at (24,3) {\lb};
\node at (26,0) {\lb};
\node at (26,4) {\lb};
\node at (30,1) {\lb};
\node at (20,8) {\lb};
\node at (22,9) {\lb};
\node at (24,9) {\lb};
\node at (24,10) {\lb};
\node at (25,10) {\lb};
\node at (26,11) {\lb};
\node at (30,12) {\lb};
\node at (31,13) {\lb};
\node at (32,14) {\lb};
\node at (32,13) {\lb};
\node at (33,14) {\lb};
\node at (34,15) {\lb};
\draw (-.5,17) -- (35,17);
\node at (0,16.4) {$40$};
\node at (4,16.4) {$36$};
\node at (8,16.4) {$32$};
\node at (12,16.4) {$28$};
\node at (16,16.4) {$24$};
\node at (20,16.4) {$20$};
\node at (24,16.4) {$16$};
\node at (28,16.4) {$12$};
\node at (32,16.4) {$8$};
\node at (0,17) {\lb};
\node at (4,17) {\lb};
\node at (4,18) {\lb};
\node at (6,18) {\lb};
\node at (7,19) {\lb};
\node at (8,19) {\lb};
\node at (8,20) {\lb};
\node at (10,17) {\lb};
\node at (10,20) {\lb};
\node at (12,21) {\lb};
\node at (14,18) {\lb};
\node at (14,21) {\lb};
\node at (15,22) {\lb};
\node at (16,23) {\lb};
\node at (18,24) {\lb};
\node at (22,25) {\lb};
\node at (26,28) {\lb};
\node at (30,29) {\lb};
\node at (34,32) {\lb};
\node at (20,17) {\lb};
\node at (22,18) {\lb};
\node at (23,18) {\lb};
\node at (24,19) {\lb};
\node at (26,20) {\lb};
\node at (30,17) {\lb};
\node at (31,18) {\lb};
\node at (34,17) {\lb};
\node at (34,18) {\lb};
\node at (34,20) {\lb};
\node at (34,24) {\lb};
\node at (32,-.6) {$40$};
\node at (9,2) {\lb};
\node at (10,3) {\lb};
\draw (4,17) -- (4,18);
\draw (6,18) -- (8,20) -- (8,19);
\draw (14,21) -- (16,23);
\draw (23,18) -- (24,19);
\draw (30,17) -- (31,18);
\draw (34,17) -- (34,18);
\draw (8,27) -- (10,29);
\draw (16,31) -- (18,33);
\draw (24,35) -- (26,37);
\draw (32,39) -- (34,41);
\node at (17,27) {\rm upper edge should be 7 higher};
\node at (6,26) {\lb};
\node at (8,27) {\lb};
\node at (9,28) {\lb};
\node at (14,30) {\lb};
\node at (10,29) {\lb};
\node at (16,31) {\lb};
\node at (17,32) {\lb};
\node at (18,33) {\lb};
\node at (22,34) {\lb};
\node at (24,35) {\lb};
\node at (25,36) {\lb};
\node at (26,37) {\lb};
\node at (30,38) {\lb};
\node at (32,39) {\lb};
\node at (33,40) {\lb};
\node at (34,41) {\lb};
\node at (30,21) {\lb};
\draw [dashed] (22,1) -- (22,9);
\draw [dashed] (10,20) -- (10,29);
\draw [dashed] (14,21) -- (14,30);
\draw [dashed] (18,24) -- (18,33);
\draw [dashed] (22,25) -- (22,34);
\draw [dashed] (26,20) -- (26,37);
\draw [dashed] (30,17) -- (30,38);
\draw [dashed] (34,18) -- (34,41);

\end{\tz}
\end{center}
\end{fig}
\end{minipage}

We close this section by discussing 
  Theorem \ref{dualitythm}, which establishes duality between $ko_*(K_2)$ and $ko^*(K_2)$.  Here is the proof. Note that \cite[Corollary 9.3]{MR} says that $R=ko$ satisfies the hypothesis of \cite[Theorem 3.1]{DG} with $a=6$,
  while $X=K_2$ has $ko_*(X)$ torsion-free by our calculation  or by \cite{AH}, implying Theorem \ref{dualitythm} by \cite[Theorem 3.1]{DG}. 
  
  Similarly to an observation in Section \ref{kudifflsec}, the actions of $\cdot2$ and $\eta$ in the Pontryagin dual appear backwards from the usual interpretation.
  The duality in Theorem \ref{dualitythm} applies to each $A_k$. For example compare Figures \ref{A5} and \ref{ko^*A5}.
  The actual $A_5$ would have gradings of each increased by 64. The generator $g$ of $(ko^6(\At_5^*))^\vee$ is the class in filtration 31. Then $2^5g$ is the class in filtration 0, and $\eta^2g\ne0$. This corresponds to the class in position $(0,0)$ in Figure \ref{A5}. At the other end of the charts, the class $g'$ in $(66,0)$ in Figure \ref{A5} of order 4 with $\eta g'\ne0$ corresponds to the class in $(72,1)$ in the Pontryagin dual of Figure \ref{ko^*A5}.

  The duality in Theorem \ref{dualitythm} also applies to each $z^t\Bt_{k,\ell}$. Combining Theorems \ref{dualitythm} and \ref{B^*thm}  yields
  \begin{equation}\label{isos}(z^t\Bt_{k,\ell})_*\approx(z^t\Bt_{k,\ell})^{*+6}\approx(z^{4i-\ell-t}\Bt_{k,\ell})_{2^{k+2}+8i-2-*}.\end{equation}
  This gives duality relations among $z^t\Bt_{k,\ell}$ charts. For example, for $1\le t\le3$, let $i=2$ and obtain
  $$(z^t\Bt_{3,4})_*\approx(z^{4-t}\Bt_{3,4})_{46-*},$$
  an isomorphism as $ko_*$-modules, using the dual action of  $ko_*$ on the right hand side. Refer to Figure \ref{ziB34}.
  You can observe that the charts for $z\Bt_{3,4}$ and $z^3\Bt_{3,4}$ are 46-dual, and the chart for $z^2\B_{3,4}$ is 46-self-dual. Note how the exotic $\cdot2$ and $\eta$ extensions between gradings 32 and 34 in $z^2\Bt_{3,4}$ are dual to the nice part of the chart between gradings 12 and 14. Using (\ref{isos}) with $\ell=2\ell'$ and $t=2i-\ell'$, we obtain the following generalization.

  \begin{prop} If $\ell$ is even and $t\equiv\frac\ell2$ {\rm mod} $2$, then $z^t\Bt_{k,\ell}$ is $(2^{k+2}+4t+2\ell-2)$-self-dual.\end{prop}

  \section{Proofs of $ko^*(K_2)$ theorems}\label{cohpfsec}
In this section, we prove Theorems  \ref{B^*thm} and \ref{Ak^*thm}. We begin by recalling from \cite{DW2} results about $\ext_{A(1)}(\zt,-)$ applied to the $A(1)$-modules in the splitting of $H^*(K_2)$.

For $k\ge4$ , let $\Mt_k=\Si^{-2^k}\M_k$ be the $A(1)$-module introduced early in Section \ref{kosummandsec}, and for $r\ge0$ let $M_k^{r*}$ denote the chart obtained by applying $\ext_{A(1)}(\zt,-)$ to $\Si^{-4D}z_J\Mt_k$, where $z_J\Mt_k=z_{j_1}\cdots z_{j_r}\Mt_k$ is a module defined in \cite[Definition 3.3]{DW2} and $D=\prod 2^{j_i}$. In Figure \ref{*0}, we repeat \cite[Figure 4.1]{DW2}, which shows $M_k^{0*}$ for $5\le k\le7$, illustrating what we hope is an obvious pattern.

\bigskip
\begin{minipage}{6in}
\begin{fig}\label{*0}

{\bf $M_k^{0*}$}

\begin{center}

\begin{\tz}[scale=.4]
\draw (-.5,0) -- (6.5,0);
\draw (8.5,0) -- (19.5,0);
\draw (21.5,0) -- (33.5,0);
\node at (0,0) {\lb};
\node at (2,1) {\lb};
\node at (3,2) {\lb};
\node at (4,3) {\lb};
\node at (4,2) {\lb};
\node at (4,1) {\lb};
\node at (5,2) {\lb};
\node at (6,3) {\lb};
\node at (9,0) {\lb};
\node at (13,0) {\lb};
\node at (13,1) {\lb};
\node at (15,2) {\lb};
\node at (16,3) {\lb};
\node at (17,4) {\lb};
\node at (17,3) {\lb};
\node at (17,2) {\lb};
\node at (17,1) {\lb};
\node at (18,2) {\lb};
\node at (19,3) {\lb};
\node at (22,0) {\lb};
\node at (23,1) {\lb};
\node at (23,0) {\lb};
\node at (27,0) {\lb};
\node at (27,1) {\lb};
\node at (27,2) {\lb};
\node at (29,3) {\lb};
\node at (30,4) {\lb};
\node at (31,5) {\lb};
\node at (31,4) {\lb};
\node at (31,3) {\lb};
\node at (31,2) {\lb};
\node at (31,1) {\lb};
\node at (32,2) {\lb};
\node at (33,3) {\lb};
\draw (2,1) -- (4,3) -- (4,1) -- (6,3);
\node at (7,4) {$\iddots$};
\draw (13,0) -- (13,1);
\draw (15,2) -- (17,4) -- (17,1) -- (19,3);
\node at (20,4) {$\iddots$};
\draw (22,0) -- (23,1) -- (23,0);
\draw (27,0) -- (27,2);
\draw (29,3) -- (31,5) -- (31,1) -- (33,3);
\node at (34,4) {$\iddots$};
\node at (0,-.6) {$4$};
\node at (4,-.6) {$0$};
\node at (9,-.6) {$8$};
\node at (17,-.6) {$0$};
\node at (13,-.6) {$4$};
\node at (22,-.6) {$9$};
\node at (27,-.6) {$4$};
\node at (31,-.6) {$0$};
\node at (3,-1.6) {$k=5$};
\node at (14,-1.6) {$k=6$};
\node at (28,-1.6) {$k=7$};

\end{\tz}
\end{center}
\end{fig}
\end{minipage}
\bigskip

In \cite[Proposition 2a]{DW2} it is shown that $M_k^{r*}$ is formed from $M_k^{0*}$ by increasing all filtrations by $r$ and extending to the left using $v_1^4$-periodicity. See Figure \ref{M5} for $M_5^{2*}$. Clearly $M_k^{(r+4)*}=\Si^8M_k^{r*}$.

The tableau for $(z^t\Bt_{k,\ell})^*$ is obtained from that of $z^t\Bt_{k,\ell}$ by applying $*$ to all the summands, reversing the order of the list of summands, and reversing the direction of the arrows. For example, the first (of four) rows of $(\Bt_{5,9})^*$ is obtained from (\ref{B59}) to be
$$\Si^{128}M_7^{1*}\lar \Si^{121}M_4^{5*}\quad \Si^{120}M_4^{3*}\lar \Si^{113}M_5^{4*}\quad \Si^{112}M_5^{2*}\lar \Si^{105}M_4^{4*}\quad \Si^{104}M_4^{2*}\lar \Si^{97}M_6^{3*}.$$

The following result says that if $i+j+k\equiv0$ mod 4, then the charts $M_k^i$ and $M_k^{j*}$ are $(2(i+j+k)-5)$-dual.
We illustrate in Figure \ref{M5}. 

\bigskip
\begin{minipage}{6in}
\begin{fig}\label{M5}

{\bf Dual $M_5$'s.}

\begin{center}

\begin{\tz}[scale=.55]
\draw (1.5,0) -- (13,0);
\draw (15.5,0) -- (27,0);
\node at (3,-.5) {$3$};
\node at (7,-.5) {$7$};
\node at (11,-.5) {$11$};
\node at (18,-.5) {$8$};
\node at (22,-.5) {$4$};
\node at (26,-.5) {$0$};
\draw (2,0) -- (3,1) -- (3,0);
\draw (4,0) -- (5,1);
\draw (9,3) -- (11,5) -- (11,3) -- (13,5);
\draw (17,0) -- (18,1) -- (18,0);
\draw (19,0) -- (20,1);
\draw (24,3) -- (26,5) -- (26,3) -- (28,5);
\node at (2,0) {\lb};
\node at (3,0) {\lb};
\node at (3,1) {\lb};
\node at (4,0) {\lb};
\node at (5,1) {\lb};
\node at (7,2) {\lb};
\node at (9,3) {\lb};
\node at (10,4) {\lb};
\node at (11,5) {\lb};
\node at (11,4) {\lb};
\node at (11,3) {\lb};
\node at (12,4) {\lb};
\node at (13,5) {\lb};
\node at (17,0) {\lb};
\node at (18,0) {\lb};
\node at (18,1) {\lb};
\node at (19,0) {\lb};
\node at (20,1) {\lb};
\node at (22,2) {\lb};
\node at (24,3) {\lb};
\node at (25,4) {\lb};
\node at (26,5) {\lb};
\node at (26,4) {\lb};
\node at (26,3) {\lb};
\node at (27,4) {\lb};
\node at (28,5) {\lb};
\node at (4,3.5) {$M_5^1$};
\node at (19,3.5) {$M_5^{2*}$};

\end{\tz}
\end{center}
\end{fig}
\end{minipage}
\bigskip

\begin{prop}\label{Mdual} Let $S=i+j+k$. If $S\equiv0$ {\rm mod} $4$, then $(M_k^i)_x\approx (M_k^{j*})^{2S-5-x}$.
\end{prop}
\begin{proof} First note that $M_k^{0*}$ leaves filtration 0 in the same way that $M_k^i$ does when $k+i\equiv0$ mod 4, by Figures \ref{same} and \ref{*0}. Since they have the same $k$, hence the same maximal heights, we deduce that if $k+i\equiv0$ mod 4, then $M_k^i$ and $M_k^{0*}$ (drawn cohomologically) are isomorphic up to grading. Next note that $M^3_{5+4t}$, $M^2_{6+4t}$, $M^1_{7+4t}$, and $M^0_{8+4t}$ all leave filtration 0 in grading $11+8t$, while, for any $k$, $M_k^{0*}$ does so in grading 0. With $S=8+4t$ and $x=11+8t$, $2S-5-x=0$, so the gradings in the proposition work if $i\le3$. Increasing $i$ by 4 adds 8 to where $M_k^i$ leaves filtration 0, so the proposition is true if $j=0$.

If the proposition is true for $(i,j,k)$, then it is true for all $(i',j,k)$ with $i'\equiv i$ mod 4, since increasing $i$ by 4 suspends $M_k^i$ by 8. Let $(i,j,k)$ be arbitrary with $i+j+k\equiv0$ mod 4. The proposition is true for $(i',0,k)$ with $i'\equiv i+j$ mod 4 and $i'>j$. Since changing $(i,j)$ from $(i',0)$ to $(i'-j,j)$ increases filtration by $j$ in both 
$M_k^i$ and $M_k^{j*}$, the proposition is true for $(i'-j,j,k)$ and hence also for $(i,j,k)$.
\end{proof}

The proof of Theorem \ref{B^*thm} follows from combining Proposition \ref{Mdual} with the observation preceding it about the tableau for $z^t\Bt_{k,\ell}$. The details are a bit delicate, so we provide two examples.

The first four terms of the tableau for $(z^t\Bt_{5,9})^*$ are
$$\Si^{128}M_7^{(1+t)*}\lar \Si^{121}M_4^{(5+t)*}\quad \Si^{120}M_4^{(3+t)*}\lar \Si^{113}M_5^{(4+t)*}.$$
By Proposition \ref{Mdual}, this is isomorphic to
$$(M^{8-t}_7)_{155-*}\lar (M_4^{7-t})_{148-*}\quad (M_4^{9-t})_{147-*}\lar (M_5^{7-t})_{140-*}.$$
This could also written as
$$(\Si M_7^{8-t})_{156-*}\lar (\Si^8M_4^{7-t})_{156-*}\quad (\Si^9M_4^{9-t})_{156-*}\lar (\Si^{16}M_5^{7-t})_{156-*},$$
which are the first four terms of the tableau for $(z^{3-t}\Bt_{5,9})_{156-*}$ (see (\ref{B59})), consistent with Theorem \ref{B^*thm}.

For a more general verification, we compare dual terms in (\ref{Bkltableau}). The $(2i+1)$st term in the first half of the tableau for $(z^t\Bt_{k,\ell})^*$ is $$\Si^{8(2^{k-1}-i-1)+8}(M_{4+\nu(2^{k-1}-1-i+1)}^{\a(2^{k-1}-1-i+1)-1+t})^*=\Si^{2^{k+2}-8i}(M_{4+\nu(i)}^{k-2-\a(i-1)+t})^*.$$ Since $\a(i)=\a(i-1)+1-\nu(i)$, Proposition \ref{Mdual} implies that this equals
\begin{equation}\label{kl}(M_{4+\nu(i)}^{-k-t+\a(i)+1+4j})_{2^{k+2}-8i+8j+3-*},\end{equation}
where $j$ is chosen so that $-k-t+\a(i)+1+4j\ge0$.
The $(2i+1)$st term in the first half of the tableau for $(z^{4j-\ell-t}\Bt_{k,\ell})_{2^{k+2}+8j+4-*}$ is $$(\Si^{8i+1}M_{4+\nu(i)}^{\ell-k+\a(i)+1+4j-\ell-t})_{2^{k+2}+8j+4-*}.$$
  With a little manipulation, this equals (\ref{kl}).

\begin{proof}  [Proof of Theorem \ref{Ak^*thm}]
It is easy to verify, using Proposition \ref{Mdual}, that the tableau for $\At^*_{2t+\eps}$ is $(2^{2t+\eps+1}+8t+4)$-dual to that of $\E'_{2t+2-\eps,4t+1}$, except at the right end. We illustrate with $\At_5^*$. We write each row of the tableau for $\E'_{5,9}$ directly beneath that of $\At^*_5$, and observe the 84-duality, since $84=65+(2\cdot12-5)$.
\begin{align}\nonumber\Si^{64}M_7^{0*}\lar\Si^{57}M_4^{4*}\quad\Si^{56}M_4^{2*}\lar\Si^{49}M_5^{3*}&\quad\Si^{48}M_5^{1*}\lar\Si^{41}M_4^{3*}\quad\Si^{40}M_4^{1*}\lar\Si^{33}M_6^{2*}\\
\Si M_7^5\ml{5,9}\Si^8M_4^4\qquad \Si^9M_4^6\ml{6,7}\Si^{16}M_5^4&\qquad\Si^{17}M_5^6\ml{6,8}\Si^{24}M_4^5\qquad\Si^{25}M_4^7\ml{7,8}\Si^{32}M_6^4\label{175}\end{align}

\begin{align}\Si^{32}M_6^{0*}\lar \Si^{25}M_4^{3*}\quad\Si^{24}M_4^{1*}\lar \Si^{17}M_5^{2*}&\quad \Si^{16}M_5^{0*}\lar\Si^9M_4^{2*}\quad\Si^8M_4^{0*}\lar V_5^*\nonumber\\
\Si^{33}M_6^6\ml{6,9}\Si^{40}M_4^5\qquad \Si^{41}M_4^7\ml{7,8}\Si^{48}M_5^5&\qquad \Si^{49}M_5^7\ml{7,9}\Si^{56}M_4^6\quad\Si^{57}M_4^8\ml{8,9}\Si^{64}M_{11}^0\label{E59}
\end{align}

The chart $V_k^*$ is formed from $\ext_{A(1)}(\zt,\zt\oplus\Si^{-8}NU)$ with $d^{k+1}$ differentials, similarly to $V_k$ discussed early in Section \ref{A4sec}. The module $NU$ is shown in Figure \ref{NU}, and $\ext_{A(1)}(\zt,\Si^{-8}NU)$ is easily calculated to be the black part of Figure \ref{Vk*}, which includes $\ext_{A(1)}(\zt,\zt)$ in red, and a $d^6$ differential, relevant for $V_5^*$. The right side of Figure \ref{Vk*} shows the result after filtrations are increased to make the differential a $d^1$. Except for the classes in grading 6 and 7, it looks like a Moore spectrum $M(2^{k-1})$. The $\eta$ extension from $-2$ to $-3$ follows from (\ref{bracket}).

\bigskip
\begin{minipage}{6in}
\begin{fig}\label{Vk*}

{\bf Forming $V_5^*$}

\begin{center}

\begin{\tz}[scale=.45]
\draw (-.5,0) -- (12,0);
\draw (15.5,0) -- (27,0);
\draw (0,0) -- (1,1) -- (1,0);
\draw [->] (6,2) -- (6,9);
\draw (6,2) -- (8,4);
\draw [->] (10,5) -- (10,9);
\draw [red] [->] (7,0) -- (7,9);
\draw [red] (7,0) -- (9,2);
\draw [red] [->] (11,3) -- (11,9);
\draw [blue] (7,0) -- (6,6);
\draw [blue] (7,2) -- (6,8);
\draw [blue] (11,3) -- (10,9);
\draw [blue] (11,5) -- (10.33,9);
\draw (16,0) -- (17,1) -- (17,0);
\draw (22,5) -- (22,2) -- (24,4);
\draw (26,5) -- (26,8) -- (24,6);
\node at (12,7) {$\iddots$};
\node at (27,7) {$\iddots$};
\node at (0,-.6) {$7$};
\node at (1,-.6) {$6$};
\node at (6,-.6) {$1$};
\node at (10,-.6) {$-3$};
\node at (16,-.6) {$7$};
\node at (22,-.6) {$1$};
\node at (26,-.6) {$-3$};
\node at (0,0) {\lb};
\node at (1,1) {\lb};
\node at (1,0) {\lb};
\node at (6,2) {\lb};
\node at (6,3) {\lb};
\node at (6,4) {\lb};
\node at (6,5) {\lb};
\node at (7,3) {\lb};
\node at (8,4) {\lb};
\node at (10,5) {\lb};
\node at (10,6) {\lb};
\node at (10,7) {\lb};
\node at (10,8) {\lb};
\node at (16,0) {\lb};
\node at (17,1) {\lb};
\node at (17,0) {\lb};
\node at (22,2) {\lb};
\node at (22,3) {\lb};
\node at (22,4) {\lb};
\node at (22,5) {\lb};
\node at (23,3) {\lb};
\node at (24,4) {\lb};
\node at (26,5) {\lb};
\node at (26,6) {\lb};
\node at (26,7) {\lb};
\node at (26,8) {\lb};
\node at (25,7) {\lb};
\node at (24,6) {\lb};
\node [red] at (8,1) {\lb};
\node [red] at (9,2) {\lb};

\end{\tz}
\end{center}
\end{fig}
\end{minipage}
\bigskip

We illustrate the proof using $k=5$, but the argument clearly generalizes. It seems easiest to use the gradings of $\E'_{5,9}$ and dualize ($84-*$) the gradings of $\At_5^*$. Until the end part of their tableaux ($V_5^*$ or $\Si^{64}M_{11}^0$), the charts and differentials in forming $\At_5^*$ and $\E'_{5,9}$ are the same. In Figure \ref{comp}, we compare $V_5^*$ (with the dual gradings) and $\Si^{64}M^0_{11}$. Note especially that their lower edges agree in grading $\ge83$.

\bigskip
\begin{minipage}{6in}
\begin{fig}\label{comp}

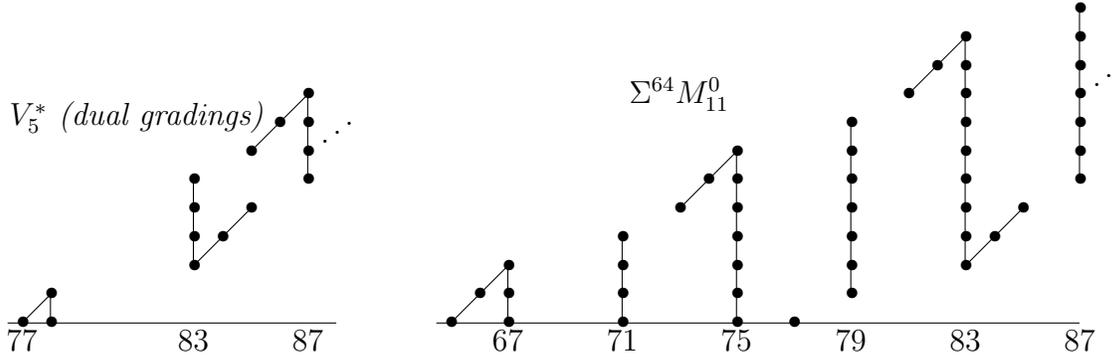
{\bf Comparing $V_5^*$ and $\Si^{64}M_0^{11}$}

\begin{center}

\begin{\tz}[scale=.38]
\draw (-.5,0) -- (11,0);
\draw (14.5,0) -- (37,0);
\draw (0,0) -- (1,1) -- (1,0);
\draw (6,5) -- (6,2) -- (8,4);
\draw (8,6) -- (10,8) -- (10,5);
\node at (11,7) {$\iddots$};
\draw (15,0) -- (17,2) -- (17,0);
\draw (21,0) -- (21,3);
\draw (25,0) -- (25,6) -- (23,4);
\draw (29,1) -- (29,7);
\draw (31,8) -- (33,10) -- (33,2) -- (35,4);
\draw (37,5) -- (37,11);
\node at (38,9) {$\iddots$};
\node at (0,-.6) {$77$};
\node at (6,-.6) {$83$};
\node at (10,-.6) {$87$};
\node at (17,-.6) {$67$};
\node at (21,-.6) {$71$};
\node at (25,-.6) {$75$};
\node at (29,-.6) {$79$};
\node at (33,-.6) {$83$};
\node at (37,-.6) {$87$};
\node at (4,7.2) {$V_5^*$ (dual gradings)};
\node at (23,8) {$\Si^{64}M_{11}^0$};
\node at (0,0) {\lb};
\node at (1,1) {\lb};
\node at (1,0) {\lb};
\node at (6,2) {\lb};
\node at (6,3) {\lb};
\node at (6,4) {\lb};
\node at (6,5) {\lb};
\node at (7,3) {\lb};
\node at (8,4) {\lb};
\node at (10,5) {\lb};
\node at (10,6) {\lb};
\node at (10,7) {\lb};
\node at (10,8) {\lb};
\node at (9,7) {\lb};
\node at (8,6) {\lb};
\node at (15,0) {\lb};
\node at (16,1) {\lb};
\node at (17,2) {\lb};
\node at (17,1) {\lb};
\node at (17,0) {\lb};
\node at (21,0) {\lb};
\node at (21,1) {\lb};
\node at (21,2) {\lb};
\node at (21,3) {\lb};
\node at (23,4) {\lb};
\node at (24,5) {\lb};
\node at (25,6) {\lb};
\node at (25,5) {\lb};
\node at (25,4) {\lb};
\node at (25,3) {\lb};
\node at (25,2) {\lb};
\node at (25,1) {\lb};
\node at (25,0) {\lb};
\node at (27,0) {\lb};
\node at (29,1) {\lb};
\node at (29,2) {\lb};
\node at (29,3) {\lb};
\node at (29,4) {\lb};
\node at (29,5) {\lb};
\node at (29,6) {\lb};
\node at (29,7) {\lb};
\node at (31,8) {\lb};
\node at (32,9) {\lb};
\node at (33,10) {\lb};
\node at (33,9) {\lb};
\node at (33,8) {\lb};
\node at (33,7) {\lb};
\node at (33,6) {\lb};
\node at (33,5) {\lb};
\node at (33,4) {\lb};
\node at (33,3) {\lb};
\node at (33,2) {\lb};
\node at (34,3) {\lb};
\node at (35,4) {\lb};
\node at (37,5) {\lb};
\node at (37,6) {\lb};
\node at (37,7) {\lb};
\node at (37,8) {\lb};
\node at (37,9) {\lb};
\node at (37,10) {\lb};
\node at (37,11) {\lb};

\end{\tz}
\end{center}
\end{fig}
\end{minipage}
\bigskip

 In forming $\E'_{5,9}$, the $\Si^{64}M_{11}^0$ will be modified four times, to $\Phi\Si^{64}M_{10}^0$, then $\Phi^3\Si^{64}M_9^0$, then $\Phi^7\Si^{64}M_8^0$, and finally $\Phi^{15}\Si^{64}M_7^0$, by having its lower edge hit a lightning flash with a $d^1$. These lightning flashes are first $\Si^{57}M_4^8$, then the result of the $\ml{7,9}$ preceding it, then the result of the first two arrows of the second row of  (\ref{E59}), and finally the result of the arrows in the second row of (\ref{175}).
This hitting into lightning flashes will also take place in forming $\At_5^*$, each time reducing the $M(2^4)$ by one 2-power, until after four such $d^1$'s, it has vanished (except for the classes in 77 and 78). 

The difference between forming $\E'_{5,9}$ and $\At_5^*$ is that prior to grading 83 these four stable lightning flashes, which in forming $\E'_{5,9}$ were interacting with parts of the modifications of $\Si^{64}M_{11}^0$, have nothing to interact with in forming $\At^*_5$. So they are present in $\At_5^*$, as extensions of sequences of lightning flashes that were modified in forming $\E'_{5,9}$. Beginning in grading 82, they are terminated.
The differentials and extensions among these lightning flashes follow for many reasons: $v_1^4$-periodicity from earlier ones, comparison with $ku$, or Toda bracket arguments.

\end{proof}

\section{Appendix: Notation and some terminology}
In this appendix, we list our specialized notation and a few items of terminology, in roughly alphabetical order.

\bigskip
\setlength{\tabcolsep}{10pt}
\renewcommand{\arraystretch}{1.3}
\begin{longtable}{cll}
$A_k$&\S1& summand of $ku^*(K_2)$\\
$\A_k$&\S1&$ko_*$ analogue of $A_k$\\
$\At_k$&\S2, Thm.\ref{Akthm}&$\Si^{-2^{k+1}}\A_k$\\
$\as_k$&\S1, \S12&$ku_*$ analogue of $A_k$\\
$\Ast_k$&\S13.3&$\Si^{-2^{k+1}}\as_k$\\
$\A^*_k$&Def.\ref{Ak*def}&$ko^*$ analogue of $A_k$\\
$\At^*_k$&Def.\ref{Ak*def}&$\Si^{-2^{k+1}}\A^*_k$\\
$\Ah_k$&Def.\ref{ABdef}&$E_1$-submodule of $H^*(K_2)$ corresponding to $A_k$\\
$\Ah^o_k$&\S10&$A(1)$-module analogue of $\Ah_k$\\
$A(1)$&\S10&subalgebra of Steenrod algebra generated by $\sq^1$ and $\sq^2$\\
$\a(n)$&\S1&number of 1's in binary expansion of $n$\\
ASS&\S1&2-primary Adams spectral sequence\\
$B_k$&\S1&summand of $ku^*(K_2)$\\
$B_{k,\ell}$&\S4&summand of $ku^*(K_2)$ corresponding to two $B_k$'s and $S_{k,\ell}$\\
$\B_{k,\ell}$&\S1, \S4&$ko_*$ analogue of $B_{k,\ell}$\\
$\Bt_{k,\ell}$&\S4&$\Si^{-2^{\ell+2}}\B_{k,\ell}$\\
$\Bh_{k,\ell}$&Def.\ref{ABdef}&$E_1$-submodule of $H^*(K_2)$ corresponding to $B_{k,\ell}$\\
$\Bh^o_{k,\ell}$&\S10&$A(1)$-submodule of $H^*(K_2)$ corresponding to $\B_{k,\ell}$\\
$b_i(M)$&\S9&$i$th generator or $v$-tower from bottom of an $E_1$-module $M$\\
$c$&\S11&complexification $ko\to ku$\\
$C_k$&Def.\ref{Ckdef}&a certain $E_1$-submodule of $H^*(K_2)$\\
$C_{i,k}$&Def.\ref{hdef}&a chart which appears in $z^i\Bt_{k,\ell}$\\
$d^r$&\S8&differential in ASS for $ku^*(K_2)$, increasing $y$ by $r$\\
$d_r$&\S8&differential in ASS for $ku_*(K_2)$, increasing $y$ by $r$\\
$\E_{e,\ell}$&\S3, Thm.\ref{edges}&edges which form $ko_*(K_2)$ and $ko^*(K_2)$\\
$\E'_{e,\ell}$&Def.\ref{E'def}&pre-edges: $\E_{e,\ell}$ prior to certain differentials\\
$E_1$&\S8&subalgebra of Steenrod algebra generated by $Q_0$ and $Q_1$\\
$E_2$ page&\S7, \S9&the initial stage of a spectral sequence\\
$\eta$&\S2&nonzero element of $ko_1$ (corresponds to Hopf map)\\
$f_b(i)$&Thm.\ref{edges}, Table \ref{tbl}&function telling how many classes hit by differentials\\
$\Phi^i$&Def.\ref{Phidef}&increases filtration of a chart by $i$\\
$\G_k$&\S9&exterior algebra on $\{y_i:i\ge k\}$\\
$h_0$, $h_1$&\S11&elements in Ext corresponding to $2$ and $\eta$\\
$h_k$&Def.\ref{hdef}&a height function\\
$H^*(K_2)$&\S8&cohomology with $\zt$ coefficients\\
$I_k$&Def.\ref{Ldef}&$\E_{1,k}$ in grading $<8$\\
$J(k,i)$&Def.\ref{Ckdef}&a useful function\\
$K_2$&\S1&Eilenberg-MacLane space $K(\zt,2)$\\
$ko_*$&\S1&connective $KO$ homology\\
$ko^*$&\S1&connective $KO$ cohomology\\
$ku_*$&\S1&connective $KU$ homology\\
$ku^*$&\S1&connective $KU$ cohomology\\
$K_{t,\eps}(x,y)$&Def.\ref{Kdef}&charts depicted in Figure \ref{ch2}\\
$L(8a)$&\S13.1&label of $W^{8a-7}\lar W^{8a}$\\
$L_{t,\eps}(x,y)$&Def.\ref{Ldef}&charts depicted in Figure \ref{Lchart}\\
$\L_k$&\S8&exterior algebra on $\{z_i:i\ge k\}$\\
lg$(i)$&Thm.\ref{upperthm}&$[\log_2(i)]$\\
lightning flash&\S1&a frequently-occurring chart with 6 classes\\
$M_k$&\S8&$E_1$-modules from \cite{DW}\\
$M_k$&\S9&$ku^*$ chart for the module $M_k$\\
$\M_k$&\S10&$A(1)$-module corresponding to $M_k$\\
$\Mt_k$&\S10&$\Si^{-2^k}\M_k$\\
$M_k^i$&\S3&$ko_*$ chart for $\Mt_k$ with filtrations decreased by $i$\\
$M_k^i$&\S15&$A(1)$-module with chart $M_k^i$\\
$M_k^i$&\S13.2&$ku$ analogue of the chart $M_k^i$\\
$M_k^{i*}$&\S18&$ko$-cohomology analogue of $M_k^i$\\
$\Mh_4^s$&Def.\ref{Mhdef}&$M_4^s$ with 0 or 1 classes adjoined\\
stably $M_k$&Def.\ref{stable}, \S13&the behavior of any $M_k^i$ in large gradings\\
$M(n)$&Def.\ref{E1'def}, \S12&mod $n$ Moore spectrum or $ko_*(M(n))$\\
MSpin&\S16&Spin cobordism spectrum\\
$N$&Figure \ref{N}&an $E_1$-module which admits structure of $A(1)$-module\\
$NU$&\S10&an $A(1)$-module corresponding to the $E_1$-module $y_1N$\\
$\nu(n)$&\S9&the exponent of 2 in $n$\\
$P$&\S10&the Poincar\'e series of a graded vector space\\
pre-edge&Def.\ref{E'def}&generic term for $\E'_{e,\ell}$\\
$M^\vee$&\S1, \S8, \S17&The Pontryagin dual, Hom$(M,\Z/2^\infty)$, of $M$\\
$q$&\S8&an element of $H^9(K_2)$\\
$Q_0$&\S8&another name for $\sq^1$\\
$Q_1$&\S8&the Milnor primitive $\sq^3+\sq^2\sq^1$\\
$R$&\S6&$ko_*$ in grading 0, 1, and 2\\
$\Si^i$&throughout&$i$-fold suspension, raises grading by $i$\\
$\sq^i$&\S10, \S13.2, \S16&generators of the Steenrod algebra\\
$S_{k,\ell}$&\S1&a factor of summands of $ku^*(K_2)$\\
$\tau_i(M)$&\S9&$i$th generator or $v$-tower from top of an $E_1$-module $M$\\
$U$&\S10&$A(1)$-module corresponding to $y_1$\\
$v$&\S8&generator of $ku^{-2}$ or $ku_2$; Bott periodicity element\\
$v$-tower&\S1, \S8&portion of a $ku^*$ or $ku_*$ chart consisting of elements $v^ix$\\
$v_1^4$&\S3&filtration-4 element of $ko_8$; Adams periodicity element\\
$V_k$&Def.\ref{E1'def}&initial step in forming $\E_{1,k}$\\
$\overline{w}_i$&Thm.\ref{SWthm}, \ref{1.2}&dual Stiefel-Whitney classes\\
$W^n$&\S13.1&$\Si^nM$ term in a tableau\\
$W^n_i$&\S13.1&$W^n$ with subscript of $M$ decreased by $i$\\
$\chi$&\S16&the antiautomorphism of the Steenrod algebra\\
$y_1$&\S8&element of $H^4(K_2)$\\
$y_k$&\S9&$y_1^{2^{k-1}}$\\
$Y_{j,k}$&\S9&term in $C_k$ containing $y_1^j$ as factor\\
$\zt$&\S1&group of order 2\\
$z_j$&\S8, \S9&an element of $H^{2^{j+2}+2}(K_2)$\\
$z_{i,j}$&\S9&$z_i^2z_{i+1}\cdots z_{j-1}$\\
$Z_k^\ell$&\S9&$z_kz_{k+1}\cdots z_{\ell-1}$\\
$z_J\M_k$&\S10&$=z_{j_1}\cdots z_{j_r}\M_k$ is an $A(1)$-module from \cite{DW2}\\
$z^i$&\S4&product of $i$ distinct $\Si^{-2^{j+2}}z_j$'s\\
$z^i\B_{k,\ell}$&\S1, \S4, \S10&chart like $\B_{k,\ell}$ but built from $M_b^{a+i}$\\
$z^i\Bt_{k,\ell}$&\S4, \S10, \S14&$\Si^{-2^{\ell+2}}z^i\B_{k,\ell}$\\
$z^i\B^*_{k,\ell}$&\S17&$ko^*$ analogue of $z^i\B_{k,\ell}$\\
$z^i\Bt_{k,\ell}^*$&\S17&$ko^*$ analogue of $z^i\Bt_{k,\ell}$\\
$z(n)$&Def.\ref{Ckdef}&$\prod z_{j_i}$ where $n=\sum 2^{j_i}$

\end{longtable}
\bigskip

\def\line{\rule{.6in}{.6pt}}

\end{document}